\documentclass[12pt,a4paper]{article}

\usepackage{fontenc}
\usepackage[utf8]{inputenc}
\usepackage[english]{babel}
\usepackage[hidelinks, hypertexnames = false]{hyperref}

\usepackage[svgnames]{xcolor}
\usepackage{abstract}
\usepackage{anyfontsize}
\usepackage{mathrsfs}
\usepackage{csquotes}

\usepackage{amssymb}
\usepackage{environ}
\usepackage{mathtools, nccmath}
\usepackage{amsfonts}
\usepackage{amsthm}
\usepackage[ruled, vlined, german, linesnumbered]{algorithm2e}
\usepackage{bm} 
\usepackage{dsfont} 
\usepackage{gauss}
\usepackage{enumitem} 
\usepackage{chngcntr}
\counterwithout{equation}{section}

\usepackage{libertine,libertinust1math} 

\usepackage{graphicx}
\graphicspath{{./figures/}}
\usepackage{wrapfig} 
\usepackage{caption} 
\usepackage{float} 
\usepackage{tikz}
\usetikzlibrary{angles}
\usepackage{tkz-euclide}
\usepackage{pgfplots}
\usepgfplotslibrary{fillbetween}
\pgfplotsset{compat=1.11}
\tikzset{
  schraffiert/.style={pattern=horizontal lines,pattern color=#1},
  schraffiert/.default=black
}

\usepackage{tabularx}
\newcolumntype{L}[1]{>{\raggedright\arraybackslash}p{#1}} 
\newcolumntype{C}[1]{>{\centering\arraybackslash}p{#1}} 
\newcolumntype{R}[1]{>{\raggedleft\arraybackslash}p{#1}} 
\makeatletter
\newcommand{\thickhline}{%
    \noalign {\ifnum 0=`}\fi \hrule height 3pt
    \futurelet \reserved@a \@xhline
}
\newcolumntype{/}{@{\hskip\tabcolsep\vrule width 1.5pt\hskip\tabcolsep}}
\makeatother
\usepackage{boldline}
\usepackage{colortbl}
\usepackage[includefoot,margin=2.54cm]{geometry}
\usepackage[bottom]{footmisc}

\newcommand{\OU}{{Ornstein\textendash Uhlenbeck }}
\newcommand{\R}{\mathbb{R}}

\newcommand{\e}{\mathrm{e}}
\newcommand{\ii}{\mathrm{i}}
\newcommand{\C}{\mathbb{C}}

\newcommand{\N}{\mathbb{N}}

\newcommand{\E}{\mathbb{E}}
\newcommand{\F}{\mathscr{F}}
\newcommand{\PP}{\mathbb{P}}

\newcommand{\oo}{\mathrm{o}}
\newcommand{\otheta}{{\vartheta}}
\newcommand{\dd}{\:\mathrm{d}}

\newcommand{\Cov}{\mathrm{Cov}}
\newcommand{\kX}{\widehat{X}^\theta}
\newcommand{\wX}{V^\theta}
\newcommand{\oX}{W^\theta}
\newcommand{\langX}{\smash{\widebar X}^\theta}
\newcommand{\langhomX}{\smash{\widebar X}^{\theta,\mathrm{hom}}}
\newcommand{\laengerX}{\underline X^\theta}
\newcommand{\laengerhomX}{\underline X^{\theta,\mathrm{hom}}}
\newcommand{\langoX}{\smash{\widebar X}^\otheta}
\newcommand{\langhomoX}{\smash{\widebar X}^{\otheta,\mathrm{hom}}}

\newcommand{\homkX}{\widehat{X}^{\theta,\mathrm{hom}}}
\newcommand{\homX}{V^{\theta,\mathrm{hom}}}
\newcommand{\homoX}{W^{\theta,\mathrm{hom}}}

\newcommand{\wY}{\smash{\overline{Y}}^{\theta, \mathrm{stat}}}
\newcommand{\wYhat}{\smash{\overline{Y}}^{\theta}}

\newcommand{\chiX}{\smash{\chi^{\theta}}}
\newcommand{\chiXprime}{\smash{\chi^{\theta'}}}

\newcommand{\homY}{\smash{\overline{Y}}^{\theta, \mathrm{hom}}}
\newcommand{\primeX}{\smash{\widetilde X^{\theta}}}
\newcommand{\auxX}{\smash{\widetilde X^{\theta,\mathrm{aux}}}}
\newcommand{\auxaX}{\smash{\widetilde X^{\theta,\mathrm{aux},a}}}
\newcommand{\primemu}{\smash{\widetilde \mu}_{\theta}}
\newcommand{\wmu}{\smash{\widebar{\mu}}_{\theta}}
\newcommand{\omu}{\smash{\underline{\smash{\mu}}}_{\theta}}
\newcommand{\wmuo}{\smash{\widebar{\mu}}_{\otheta}}
\newcommand{\alphaj}{\smash{\alpha_{\lambda, j}^{\otheta}}}
\newcommand{\alphaall}{\smash{\alpha_{\lambda}^\otheta}}
\newcommand{\alphajx}{\smash{\alpha_{\lambda, j}^{\otheta}}}
\newcommand{\alphaallx}{\smash{\alpha_{\lambda}^{\otheta}}}
\newcommand{\tildeP}{\smash{\widebar{P}_{\theta}}}
\newcommand{\tildePtwo}{\smash{\widebar{P}_{\otheta}}}
\newcommand{\hatP}{\smash{\widetilde{P}_\theta}}
\newcommand{\wA}{\smash{\widebar A}^{\theta}}
\newcommand{\wAY}{\smash{\widebar A}^{\theta}_Y}
\newcommand{\wAtwo}{\smash{\underline A^{\theta}}}

\newcommand{\wNtwo}{\smash{\underline N}^{\otheta}}
\newcommand{\wa}{\smash{\widebar a}^{\theta}}

\newcommand{\watwo}{\smash{\underline a}^{\theta}}
\newcommand{\wbtwo}{\smash{\underline b}^{\theta}}
\newcommand{\homA}{\smash{\widebar A}^{\theta}}
\newcommand{\homAtwo}{\smash{\widebar A}^{\theta}_{\otimes 2}}
\newcommand{\norm}[1]{\left\lVert#1\right\rVert} 

\newcommand*{\vsep}{\kern-\tabcolsep\vrule height\arraystretch\ht\strutbox depth\arraystretch\dp\strutbox\kern-\tabcolsep} 

\newcommand{\vertiii}[1]{{\vert\kern-0.25ex\vert\kern-0.25ex\vert #1
    \vert\kern-0.25ex\vert\kern-0.25ex\vert}}
\allowdisplaybreaks

\DeclareMathOperator{\Tr}{Tr}

\usepackage[style = authoryear, citestyle = numeric, maxnames = 20, maxcitenames = 2]{biblatex}

\DeclareFieldFormat*[article, thesis, incollection, inproceedings]{title}{#1}

\makeatletter
\input{numeric.bbx}
\makeatother

\numberwithin{equation}{section}

\usepackage{bm}

\newenvironment{acks}
  {\section*{Acknowledgments}}
  {}

\NewEnviron{fitequation}{%
\begin{equation} \resizebox{0.86\hsize}{!}{$\displaystyle
  \BODY
$}\end{equation}
}

\NewEnviron{fitequationtwo}{%
\begin{equation} \resizebox{0.92\hsize}{!}{$\displaystyle
  \BODY
$}\end{equation}
}

\NewEnviron{fitequationtwo*}{%
\begin{equation*} \resizebox{0.94\hsize}{!}{$\displaystyle
  \BODY
$}\end{equation*}
}

\NewEnviron{fitequation*}{%
\begin{equation*} \resizebox{\hsize}{!}{$\displaystyle
  \BODY
$}\end{equation*}
}

\makeatletter
\newcommand{\customlabel}[2]{%
   \protected@write \@auxout {}{\string \newlabel {#1}{{#2}{\thepage}{#2}{#1}{}} }%
   \hypertarget{#1}{#2}
}
\makeatother

\NewEnviron{fitalign}{
  \sbox2{\parbox{\dimexpr \wd0 + 2\wd1}%
    {\begin{align*}\BODY\end{align*}}}
  \noindent\resizebox{\textwidth}{!}{\usebox2}%
}

\newtheorem{definition}{Definition}[section]
\newtheorem{theorem}[definition]{Theorem}
\newtheorem{proposition}[definition]{Proposition}
\newtheorem{lemma}[definition]{Lemma}
\newtheorem{corollary}[definition]{Corollary}
\newtheorem{assumption}{Assumption}

\newtheorem{assumpA}{Assumption}

\newtheorem{assumpD}{Assumption}

\theoremstyle{definition}
\newtheorem{example}[definition]{Example}
\newtheorem{remark}[definition]{Remark}
\newtheorem{alg}[definition]{Algorithm}

\makeatletter
\newcommand{\partsmash}[2][tb]{%
  \def\mb@t{\ht\z@ #2\ht\z@}\def\mb@b{\dp\z@ #2\dp\z@}%
  \def\mb@tb{\mb@t \mb@b}%
  \edef\finsm@sh{\csname mb@#1\endcsname\box\z@}%
  \ifmmode \@xp\mathpalette\@xp\mathsm@sh
  \else \@xp\makesm@sh
  \fi}
\makeatother

\begin{document}
\title{Parameter Estimation for Partially\\ Observed Affine and Polynomial Processes\footnote{Python code for the simulation studies and examples from Section \ref{sec5: app} is made available as a public GitHub repository accessible via the link \href{http://github.com/Ambossadore/Quasi-Maximum-Likelihood}{http://github.com/Ambossadore/Quasi-Maximum-Likelihood}.}}
\author{Jan Kallsen \quad Ivo Richert\\ Christian-Albrechts-Universität zu Kiel\footnote{Mathematisches Seminar, Christian-Albrechts-Universität zu Kiel, Heinrich-Hecht-Platz 6, 24118 Kiel, Germany, Email: kallsen@math.uni-kiel.de, richert@math.uni-kiel.de}}
\date{}
\maketitle

\begin{abstract}
This paper is devoted to parameter estimation for partially observed polynomial state space models.
This class includes discretely observed affine or more generally polynomial Markov processes.
The polynomial structure allows for the explicit computation of
a Gaussian quasi-likelihood estimator and its asymptotic covariance matrix.
We show consistency and asymptotic normality of the estimating sequence and provide explicitly computable expressions for the corresponding asymptotic covariance matrix.\\

\noindent Key words: polynomial processes, polynomial state space models, parameter estimation, quasi-maximum likelihood estimation, affine processes, filtering\\

\noindent MSC (2020) subject classification: 62M05, 62M20
\end{abstract}

\tableofcontents

\section{Introduction}\label{sec1: Intro}
In light of their computational tractability and their flexibility in modelling empirical properties of common financial time series, the class of \textit{affine processes} has become one of the most prominently used stochastic models in modern financial mathematics. Affine processes are continuous-time Markovian semimartingales whose characteristic function features an exponen\-tial-affine dependence on the initial state vector of the process. First studied in the 1970s in the form of continuous-time limits of Galton--Watson branching processes with immigration in \cite{Kawazu1971}, affine processes have henceforth appeared in many financial applications, including the classic interest rate models of \cite{Vasicek1977} or \cite{Cox1985}, the popular stochastic volatility model of \cite{Heston1993} or the broad class of affine jump-diffusions introduced by \cite{Duffie2000}. In addition, the affine class nests both the class of Lévy processes as well as the class of Lévy-driven \OU processes. Finally, in 2003 affine processes were systematically characterised by \cite{Duffie2003}, combining the two predominantly studied types of processes from \cite{Kawazu1971} and \cite{Duffie2000} in a single broad framework.

The family of affine processes features the convenient property that higher order conditional moments can be expressed by polynomials of the same order in the current state vector. Driven by this observation, the class of affine processes has been subsequently generalised to the broader class of \textit{polynomial processes}, which is defined by this very property and which first appeared by name in \cite{Cuchiero2011} and \cite{Cuchiero2012}, although the application of polynomial processes in financial modelling dates back at least to the early 2000s with the works of \cite{Delbaen2002} and \cite{Zhou2003}. From a practical point of view, polynomial processes become inherently tractable once their infinitesimal generator or their semimartingale characteristics have been specified because the computation of higher-order moments then only requires the simple computation of matrix exponentials. This convenient property also carries over to the computation of joint moments at multiple points in time, see \cite{Benth2021}.

Apart from the affine subclass, polynomial processes cover for example the class of Pearson diffusions, which are easily able to model non-trivial dynamics on compact state spaces, or the class of Generalised \OU processes, which have been applied in financial mathematics e.g.\ in the context of the COGARCH model of \cite{Kluppelberg2004}.
Other examples of the vast applications of polynomial processes in finance include interest rate theory (\cite{Delbaen2002}, \cite{Filipovic2020} or \cite{Zhou2003}), stochastic volatility and option pricing models (\cite{Ackerer2018}), credit risk theory (\cite{Ackerer2020}), life insurance liability (\cite{Biagini2016}), energy commodity market models (\cite{KleisingerYu2020} or \cite{Ware2019}), variance swaps (\cite{Filipovic2016b}) and stochastic portfolio theory (\cite{Cuchiero2019}). Most examples of polynomial processes in the aforementioned applications can be subsumed in the simpler classes of \emph{polynomial diffusions} or \emph{polynomial jump-diffusions} introduced in \cite{Filipovic2016} and \cite{Filipovic2020b}, respectively. For definitions and important properties of affine and polynomial processes we refer to \cite{Eberlein2019}.

From the point of view of these applications, it is inevitable to statistically estimate the underlying parameters of the corresponding parameterised polynomial process. Given the necessity of this task, it is rather surprising that the availability of literature concerning the general estimation of parameterised polynomial processes is sparse. Many of the aforementioned applications establish certain parameter estimation methods that are specifically tailored to the particular employed example of a polynomial process instead of deriving general estimation methods that only use the affine or polynomial structure of the Markov process. The aim of this paper is to close this gap and provide a general framework for the estimation of polynomial --- and hence also of affine --- processes observed at discrete time points. In particular, we study asymptotic properties of our estimators if the length of the observation period tends to infinity. The interval between observations, however, is supposed to be fixed.

The most prominent methods for estimating particular examples of polynomial processes include the maximum likelihood estimation methodology (see for example \cite{Ware2019} for an exact likelihood estimation of the Jacobi process) or a classical least-squares minimisation between simulated model quantities and quantities observed on the market, as it is done for example in \cite{Ackerer2020} or \cite{Filipovic2020}. Our paper is devoted to an estimation under the ``true'' physical measure rather than an equivalent risk-neutral measure, so we will not further discuss any least-squares calibration procedures here. Instead, our main focus is the estimation of parameters using quasi-maximum likelihood (QML) procedures. If the underlying polynomial process is of diffusion-type and observed continuously in time, exact expressions for the likelihood of the process can be obtained by absolutely continuous changes of measure (see for example \cite{Basawa1980}) or by infinitesimal-time expansions of the transition dynamics (see for example \cite{Azencott1981}). If, on the other hand, the process is sampled only at discrete times, exact representations of the likelihood often become unavailable and have to be approximated as for example in \cite{AitSahalia2002} or \cite{Filipovic2013}. Usually, however, not even all components of the underlying process are observable for estimation, which renders even these procedures infeasible. This happens for example in stochastic volatility models or in latent-factor interest rate models. In this case, parameter estimation effectively involves a stochastic filtering problem. 

In this study we tackle estimation in the context of partial observations by first deriving a canonical discrete-time representation of polynomial processes in the form of a potentially higher-dimensional vector-autoregressive model of order one, and then approximating this representation by a Gaussian process with matched first and second moments. This Gaussian process is used as a proxy of the true underlying polynomial process in a classical maximum likelihood procedure. In particular, this allows to express the likelihood function of the observed parts in terms of a Kálmán filter for the latent components. Since this likelihood depends on the underlying process only through its first and second moments, the resulting QML estimator can be interpreted in the widest sense as some kind of generalised method of moments (GMM) approach. We prove that, under some mild conditions on the transition dynamics of the model, our estimator is weakly consistent and asymptotically normal. In addition, we derive an explicit representation of the  asymptotic covariance matrix of the QML estimator which can be used for constructing asymptotic confidence intervals or rejection regions of hypothesis tests. Concerning the parameter estimation of polynomial processes from the applications mentioned above, our methodology is best comparable with that of \cite{Biagini2016}, who are however not concerned with the statistical properties of their obtained QML estimator, and with that of \cite{Cacace2019}, who apply a filtering methodology comparable to ours for estimation, however only for the Heston stochastic volatility model.

As already mentioned, our canonical representation of a polynomial process is that of a vector-autoregressive model, which is itself a special case of a general vector-autoregressive moving average (VARMA) model. In this sense, our results can be seen as a specific application of the general, conceptually rich literature on maximum likelihood estimation of VARMA models. Consistency and asymptotic normality of QML estimators of VARMA models (however with homoskedastic noise) have been shown before for example by \cite{Mainassara2011} and by \cite{Dunsmuir1976} in the case of fully observed processes and by \cite{Schlemm2012} covering the case of partial observations. The key difference between these studies and our results lies in the fact that due to our added Markovian structure and the availability of explicit polynomial moments, we are able to provide a closed-form representation of the asymptotic covariance matrix and to considerably simplify the required assumptions on the underlying process. 

The remainder of this paper is organised as follows. In Section \ref{sec4: Estimation} we introduce the setup for our estimation problem, in particular the notion of a parametric polynomial state space model and the corresponding quasi-likelihood. The main results on consistency and asymptotic normality
are to be found in Section \ref{sec4.2: Proof}. Subsequently, we apply our results to the popular Heston stochastic volatility model from finance and to Lévy-driven \OU models in Section \ref{sec5: app}. The proofs of the main results are to be found in Section \ref{s:proofs}. The \hyperref[appn]{Appendix} provides some technical tools that are needed in different parts of the paper. Python code for the simulation studies included in this paper is available as supplementary material. 

\subsection*{Frequently used notation throughout the paper}
We finish this section by introducing and recalling the basic notation and some theoretical preliminaries. Throughout, we let $\N \coloneq \{0, 1, 2, \dots\}$ denote the set of natural numbers including 0 and $\N^* \coloneq \N \setminus \{0\}$. For any space $E$ we let id denote the identity map on $E$ and $\mathscr{P}(E)$ the power set of $E$. For a multi-index $\lambda = (\lambda_1, \dots, \lambda_d) \in \N^d$ and $x = (x_1, \dots, x_d) \in \R^d$, we write $|\lambda| = \sum_{j=1}^d \lambda_j$ and $\smash{x^\lambda \coloneq \prod_{j=1}^d x_j^{\lambda_j}}$. The notation $\N^d_k$ stands for the set of all multi-indices $\lambda$ with $|\lambda| \leq k$. For $\lambda, \mu \in \N^d$ we write $\mu \leq \lambda$ if $\lambda - \mu \in \N^d$ and in this case we set $\binom{\lambda}{\mu} \coloneq \prod_{j = 1}^d \binom{\lambda_j}{\mu_j}$. The notation $\mathbf{1}_j \in \N^d$ is used for the multi-index containing only zeros except of a one in the $j$-th entry and by convention, we set $\mathbf{0} \coloneq (0, \dots, 0)$. If $E \subset \R^d$, we say that $g: E \to \R$ is a polynomial of order $p$ if it is of the form $g(x) = \sum_{|\lambda| \leq p} \alpha_\lambda x^\lambda$ for some $\alpha_\lambda \in \R$. We denote the set of all such polynomials up to order $p$ by $\mathscr{P}_p(E)$. Additionally, we say that $g: E \to \R^k$ is a $k$-dimensional polynomial of order $p$ if each component $g_i \in \mathscr{P}_p(E)$ for $i \in \{1, \dots, k\}$.

When working on $\R^d$ we assume that it has been endowed with the Borel $\sigma$-algebra $\mathscr{B}(\R^d)$. For any $d \in \N^*$ and $E \in \mathscr{B}(\R^d)$, we write $\lambda_d$ for the Lebesgue measure on $\mathscr{B}(E)$. If $j$ denotes the Hausdorff dimension of $E$, we write $\lambda_E$ for the $j$-dimensional Hausdorff measure on $\mathscr{B}(E)$ and we call $\lambda_E$ the proper Hausdorff measure on $E$. This concept extends the $j$-dimensional volume to arbitrary Borel subsets of $\R^d$ that have topological dimension less than $d$. 

We let $\mathrm{I}_d$ stand for the identity matrix in $\R^{d\times d}$. Whenever the dimension $d$ is unambiguous, we just write $\mathrm{I}$ instead of $\mathrm{I}_d$ and whenever the underlying space $\R^m$ or $\R^{m \times n}$ is unambiguous, we let 0 stand for the respective zero vector or zero matrix. Moreover, $e_j^{(k)}$ denotes the $j$-th unit vector in $\R^k$ for $j,k \in \N^*$. If $A \in \R^{d \times d}$, we let $\mathrm{diag}_k(A)$ denote the $(kd \times kd)$-block matrix that contains the matrix $A$ on the diagonal blocks and elsewhere 0. For any matrix $A \in \R^{d \times d}$, we write $\rho(A)$ for the spectral radius of $A$, i.e.\ $\rho(A) = \max_{i \in \{1, \dots, d\}} |\lambda_i(A)|$ with $\lambda_i(A)$, $i \in \{1, \dots, d\}$, denoting the eigenvalues of $A$. Moreover, we let $\alpha(A)$ stand for the spectral abscissa of $A$, i.e.\ $\alpha(A) = \max_{i \in \{1, \dots, d\}} \mathrm{Re}(\lambda_i(A))$. In multiple occasions, we will use matrix products of the form $\prod_{m=1}^n A_m$ for conformable matrices $A_1, \dots, A_n$. In any of these situations the product is to be read from right to left, that is, $\prod_{m=1}^n A_m \coloneq A_n A_{n-1} \dots A_1$. For some $x \in \R^d$ and $A \in \R^{m \times n}$ and if not specified any further, we let $\norm{x}$ stand for the Euclidean norm of $x$ and $\norm{A}$ stand for the spectral norm of $A$, i.e.\ $\norm{A} = \sup_{\norm{x} = 1} \norm{Ax}$. In particular, $\norm{A}$ is given by the largest singular value of $A$ and $\norm{A} = \rho(A)$ if $A$ is symmetric. \linebreak If $A \in \R^{d \times d}$ is any square matrix, we let $\mathrm{Sym}(A) \coloneq A + A^\top$ denote its symmetrisation.

We say that a sequence $(a_t)_{t \in \N^*}$ in a normed space $E$ converges to $a \in E$ \emph{at a geometric rate} if there are $c \in  \R_+$ and $\gamma \in [0, 1)$ with $\norm{a_t - a} \leq c \gamma^t$ for all $t \in \N^*$. Uniform convergence at a geometric rate of functions $f: E \to F$ with a normed space $F$ is defined analogously.

We let $\mathrm{C}(E, F)$ denote the space of continuous functions $E \to F$ for topological spaces $E$ and $F$. If moreover $E \subset \R^k$ and $F\subset \R^d$, we let $\mathrm{C}^m(E, F)$ denote the space of $m$-times continuously differentiable functions $E \to F$. If $E \subset \R^k$ and $f \in \mathrm{C}^2(E, \R)$, we write $\partial_{x_j} f$ or $\partial_j f$ for $\smash{\frac{\partial f}{\partial x_j}}$ and $\partial_{ij} f$ for $\smash{\frac{\partial^2 f}{\partial x_i x_j}}$. We let $\smash{\nabla_x f \coloneq (\partial_{x_1} f, \dots, \partial_{x_k} f)^\top}$ stand for the gradient of $f$. If $\alpha \in \N^k$ is a multi-index, the operator $\partial_x^\alpha$ is used as a shorthand for $\partial^{\alpha_1}_{x_1} \dots \partial^{\alpha_k}_{x_k}$. We let $\nabla_x^2 f$ stand for the $(k \times k)$ Hessian matrix of $f \in \mathrm{C}^2(E, \R)$. If $f \in \mathrm{C}^1(E, \R^d)$, we write $\nabla_x f: E \to \R^{d \times k}$ to denote the Jacobian matrix of $f$.

Let $E$ be a subset of $\R^d$. We call $E$ a proper state space if $E$ is a closed and connected embedded smooth manifold (possibly with boundary), contains $0 \in\R^d$, and is not contained in any proper linear subspace of $\R^d$. The assumption that $E$ is a smooth manifold ensures that the Hausdorff dimension of $E$ is an integer and coincides with its topological covering dimension, which will be needed later.

Since we make heavy use of the Kronecker product, we recall its basic properties here. For $A \in \R^{m \times n}$ and $B \in \R^{p \times q}$ the Kronecker product $A\otimes B \in \R^{mp \times nq}$ is the block matrix
\[A\otimes B \coloneq \begin{pmatrix}
    A_{11} B & \dots & A_{1n} B \\
    \vdots & \ddots & \vdots \\
    A_{n1} B & \dots & A_{mn} B
\end{pmatrix}.\]
In this sense, the Kronecker product is a non-commutative, associative bilinear map $\otimes: \R^{m \times n} \times \R^{p \times q} \to \R^{mp \times nq}$. We define the Kronecker product for vectors $x \in \R^n$ by interpreting $x$ as an element of $\R^{n \times 1}$. The Kronecker product shares a distributivity property with the matrix product in the sense that $(A \otimes B)(C \otimes D) = (AC) \otimes (BD)$ for any matrices $A$, $B$, $C$, $D$ of conformable size. The Kronecker product gets along well with the vectorisation operator $\mathrm{vec}: \R^{m \times n} \to \R^{mn}$ that stacks the columns of a matrix on top of each other. Indeed, we have $(A \otimes B)\mathrm{vec}(C) = \mathrm{vec}(BCA^\top)$. This implies that $\mathrm{vec}(xy^\top) = y \otimes x$ for $x \in \R^m$ and $y \in \R^n$. By associativity, we can also consider repeated Kronecker multiplication and we will write $A^{\otimes n}$ for the $n$-fold repeated Kronecker product $A \otimes \dots \otimes A$.

The Kronecker product commutes with the operations of transposing and inverting matrices in the sense that $(A \otimes B)^\top = A^\top \otimes B^\top$ and $(A \otimes B)^+ = A^+ \otimes B^+$, where $+$ denotes the Moore--Penrose pseudoinverse. If $A$ and $B$ are square matrices with respective spectrum $\{\lambda_1, \dots, \lambda_m\}$, $\{\mu_1, \dots, \mu_n\} \subset \C$, then $A \otimes B$ has the spectrum $\{\lambda_i \mu_j: \: i \in \{1, \dots, m\}, j \in \{1, \dots, n\}\}$. It follows that $\rho(A \otimes B) = \rho(A) \rho(B)$. Additionally, one has $\vertiii{A \otimes B} = \vertiii{A} \cdot \vertiii{B}$ for any matrix norm $\vertiii{\cdot}$ induced by a vector $p$-norm, see Theorem 8 and the notes on page 413 in \cite{Lancaster1972}. Sometimes, it is convenient to define the Kronecker sum $A \oplus B \in \R^{mp \times np}$ for matrices $A \in \R^{m \times n}$ and $B \in \R^{p \times q}$ such that $np = mq$ by $A \oplus B \coloneq \mathrm{I}_p \otimes A + B \otimes \mathrm{I}_m$. This has the advantage that it behaves well in conjunction with matrix exponentials since $\e^A \otimes \e^B = \e^{A \oplus B}$ for square matrices $A$ and $B$, see \cite{Neudecker1969}. Finally, if $A: E \to \R^{m \times n}$ and $B: E \to \R^{p \times q}$ are continuously differentiable matrix-valued functions on some subset $E$ of $\R$, then one has the convenient Kronecker product rule $\partial_x \big[A(x) \otimes B(x)\big] = \big[\partial_x A(x)\big] \otimes B(x) + A(x) \otimes \big[ \partial_x B(x)\big]$.

Further unexplained notation concerning in particular probability theory, stochastic processes, Markov processes and semimartingales is used as in \cite{Jacod2003, Eberlein2019, Meyn2009, Ethier1986}.

\section{The setup}\label{sec4: Estimation}

As noted in the introduction, the present paper is devoted to estimating an unknown parameter $\theta \in \Theta$ of the law of an affine or polynomial process that is sampled at discrete time points. More specifically, our setup deals with the case of \textit{partial observations} in the sense that only a subset of the coordinates of the multi-dimensional process is actually available for parameter estimation. We start by introducing a flexible discrete-time stochastic framework termed \textit{parametric polynomial state space model}, which naturally generalises the concept of a discretely sampled polynomial process. Since the likelihood of the observed part of the polynomial process is generally unknown in closed form, it is approximated by the likelihood of a Gaussian process with the same first and second moments. Consistency and asymptotic normality of the so-obtained QML estimator is stated in Section \ref{sec4.2: Proof} under suitable ergodicity assumptions on the polynomial process. 

\subsection{Parametric polynomial state space models}\label{sec4.1: StateSpaceModels}
We fix a filtered measurable space $(\Omega, \F,(\F_{t})_{t \in \N})$ along with some $E$-valued adapted process $X = (X(t))_{t \in \N}$ on this space, where $E \subset \R^d$ denotes a proper state space. The following definition is inspired by Definition 2.8 in \cite{KallsenRichert2025}.

\begin{definition}\label{def: polynomial_ssm}
Let $\Theta \subset \R^k$ denote a compact convex set and $(\PP_\theta)_{\theta \in \Theta}$ a family of probability measures on $(\Omega, \mathscr{F})$. We call
$(X,(\PP_\theta)_{\theta \in \Theta})$ a \textbf{parametric polynomial state space model of order 1} if the following holds.
\begin{enumerate}
\item $X$ is a Markov chain with respect to any $\PP_\theta$, $\theta\in\Theta$. We denote its one-step transition probabilities by $P_\theta(x,\cdot)$ for $x \in E$,
i.e.\ $\PP_\theta(X(t)\in B|\F_{t-1})=P_\theta(X(t-1),B)$ holds $\PP_\theta$-almost surely for any $B\in\mathscr B(E)$.
\item $\E_\theta(\lVert X(t)\rVert \mid X(0) = x) < \infty$ for any $t \in \N$, any $x \in E$ and any $\theta\in\Theta$.
\item For any $t \in \N^*$ we have
\begin{equation}\label{Poly}
X(t) = a^\theta + A^\theta X(t-1) + N^\theta(t)
\end{equation}
with some deterministic $a^\theta \in \R^d$ and $A^\theta \in \R^{d\times d}$ and some $\PP_\theta$-martingale difference sequence $N^\theta = (N^\theta(t))_{t \in \N}$,
i.e.\ $\E_\theta(N^\theta(t) |\F_{t-1}) = 0$ for any $t \in \N^*$. We call $A^\theta$ the \textbf{state transition matrix} and $a^\theta$ the \textbf{state transition vector} of $X$. We write $C^\theta(t)\coloneq\Cov_\theta(N^\theta(t))$ for the unconditional covariance matrix under $\PP_\theta$ if second moments exist.
\end{enumerate}
\end{definition}

If the state space $E$ is degenerate, $A^\theta$ above may not be uniquely specified. We will later always choose one possible matrix $A^\theta$ with $\rho(A^\theta) < 1$, see Proposition \ref{prop: ergod}. 

In analogy to polynomial processes, we will now develop the notion of higher-order polynomial state space models. It is the translation of the corresponding Definition 2.9 in \cite{KallsenRichert2025} to our parametric Markovian setup.
\begin{definition}\label{d:orderr}
We call $(X,(\PP_\theta)_{\theta \in \Theta})$ a \textbf{parametric polynomial state space model of order $r\in\N^*$} if
$((X^\lambda)_{\lambda\in\N^d_\ell\setminus\{0\}},(\PP_\theta)_{\theta \in \Theta})$ or, equivalently, $(\mathrm{vec}_{\otimes \ell}(X),(\PP_\theta)_{\theta \in \Theta})$ is a parametric polynomial state space model of order 1 for $\ell=1,\dots,r$, where 
$\mathrm{vec}_{\otimes \ell}(X)\coloneq (X,X\otimes X,\dots,X^{\otimes \ell})$ denotes the $\ell$-fold stacked Kronecker product of $X$ which contains all powers up to order $\ell$ of the elements of $X$ (with some entries occurring multiple times).
\end{definition}

\begin{remark}\label{rem: Mean_Recursion}
As observed in Proposition 2.11 of \cite{KallsenRichert2025}, polynomial state space models satisfy a moment formula similar to the one for polynomial processes in continuous time.
Indeed, if $\smash{a^\theta_{\otimes r}}$ and $\smash{A^\theta_{\otimes r}}$ respectively denote the state transition vector and state transition matrix of $\mathrm{vec}_{\otimes r}(X)$ relative to $\PP_\theta$, then a simple recursion shows that
\begin{equation}\label{eq: Mean Recurrence}
\E_\theta\big[\mathrm{vec}_{\otimes r}(X(t))\big| \F_s\big] = \big(A^\theta_{\otimes r}\big)^{t-s} \mathrm{vec}_{\otimes r}(X(s)) + \bigg(\sum_{u=0}^{t-s-1} \big(A^\theta_{\otimes r}\big)^{u}\bigg) a^\theta_{\otimes r}
\end{equation} 
for any $\theta \in \Theta$ and any $s, t \in \N^*$, $s \leq t$. In particular, since by definition $X$ is also polynomial of all orders $\ell$ less than $r$, it follows that the matrix $A^\theta_{\otimes r}$ has the general block triangular form
    \[
    A^\theta_{\otimes r} = \begin{pmatrix}
        A^\theta & 0 & 0 & \dots & 0 \\
        \circ & \circ & 0 & \dots & 0 \\
        \circ & \circ & \ddots & \ddots & \vdots \\
        \vdots & \vdots & \ddots & \circ & 0 \\
        \circ & \circ& \dots &  \circ & \circ
    \end{pmatrix},
    \]
    where 0 stands for a zero matrix of appropriate size and $\circ$ stands for some matrix of appropriate size. Equation \eqref{eq: Mean Recurrence} can be written in a more concise way as a simple matrix product by $\E_\theta\big[\widetilde{\mathrm{vec}}_{\otimes r}(X(t))\big| \F_s\big] = \big(\widetilde{B}^\theta_{\otimes r}\big)^{t-s}\widetilde{\mathrm{vec}}_{\otimes r}(X(s))$ for $s, t\in \N^*$, $s \leq t$, where we define
    \begin{equation}\label{eq: Mean Recursion Better}
    \widetilde{\mathrm{vec}}_{\otimes r}(X(t)) \coloneq \begin{pmatrix}
        1 \\ \mathrm{vec}_{\otimes r}(X(t))
    \end{pmatrix} \qquad\quad \text{ and } \qquad\quad \widetilde{B}^\theta_{\otimes r} \coloneq \begin{pmatrix}
        1 & 0 \\
        a^\theta_{\otimes r} & A^\theta_{\otimes r}
    \end{pmatrix}.        
    \end{equation} 
    Here, $\widetilde{B}^\theta_{\otimes r}$ is the natural analogue of the matrix exponential in the moment formula for polynomial processes in continuous time and of the matrix $\widetilde B$ from Lemma 2.15 of \cite{KallsenRichert2025}. Similar formulas hold for $(X^\lambda)_{\lambda\in\N^d_\ell\setminus\{0\}}$ in place of $\mathrm{vec}_{\otimes r}(X)$.
\end{remark}

\begin{remark}
If second moments of $N^\theta$ exist, the state space model \eqref{Poly} is a particular case of a vector-autoregressive moving average (VARMA) model. A VARMA($p$, $q$) model is an $\R^d$-valued process $(X(t))_{t \in \N}$ satisfying
\[A_0 X(t) - \sum_{i=1}^p A_i X(t-i) = B_0 \varepsilon(t) - \sum_{i=1}^q B_i \varepsilon(t-i) \]
for some matrices $A_i$, $B_i \in \R^{d \times d}$ and where $(\varepsilon(t))_{t \in \N}$ is an uncorrelated sequence of mean-zero variables. As noted in the introduction, consistency and asymptotic normality for QML estimators of VARMA models (however with homoskedastic noise sequence) have been shown before for example in \cite{Mainassara2011} or \cite{Dunsmuir1976}. The key difference between the polynomial state space model and the VARMA model lies in the fact that due to the Markovian nature and the computability of polynomial moments explicit expressions for the estimator's asymptotic covariance matrix can be obtained in our setup.
\end{remark}

\begin{example}\label{ex: Polynomial}
If $X = (X(t))_{t\in \N}$ is obtained by sampling an $r$-polynomial process $(X(t))_{t \in \R_+}$ at non-negative integer time points, then it is a polynomial state space model by Lemma 2.15 in \cite{KallsenRichert2025}. Conversely, however, there need not be an obvious choice of a continuous-time polynomial process from which $X$ is obtained by discrete sampling. One such example is the GARCH(1, 1) model $(y(t), \sigma^2(t))_{t \in \N}$ of \cite{Bollerslev1986} specified by
\begin{align*}
    y(t) &= \omega_0 + \sigma(t) \varepsilon(t) \\
    \sigma^2(t) &= \omega + \alpha y(t - 1)^2 + \beta \sigma^2(t-1),
\end{align*}
for $t \in \N^*$ for some $\omega_0, \omega > 0$,  $\alpha, \beta \geq 0$ and some adapted i.i.d. sequence $(\varepsilon(t))_{t \in \N}$ with $\varepsilon(t)$ being independent of $\F_{t-1}$ for any $t \in \N^*$. It is not hard to see that $(y(t)^2, \sigma^2(t))_{t \in \N}$ is a polynomial state space model of order $r$ as long as $\varepsilon$ has finite moments of order $2r$ and $\E(\varepsilon(t)^n) = 0$ for odd $n \in \N$. 
\end{example}

As already noted, proper maximum likelihood estimation of the parameter $\theta$ of a polynomial state space model is often infeasible in practice because the likelihood of the model is not known in closed form. We propose here to approximate the unknown likelihood by that of a Gaussian vector-autoregressive process whose first and second moments are matched with those of the given polynomial process. The existence of such a Gaussian equivalent is assured by the following proposition.

\begin{proposition}\label{prop: Gauss_equiv}
Consider a polynomial state space model $(X,(\PP_\theta)_{\theta \in \Theta})$ of order 1 with finite second moments and state transition vector and matrix $a^\theta$ and $A^\theta$. 
Let $Y=(Y(t))_{t \in \N}$ denote an adapted $\R^d$-valued process on some filtered space $(\Gamma, \mathscr{G},(\mathscr G_{t})_{t \in \N})$ and 
$(\mathbb Q_\theta)_{\theta \in \Theta}$ a family of probability measures on that space such that
\begin{equation}\label{Gauss}
Y(t) = a^\theta + A^\theta Y(t-1) + B^\theta(t)W^\theta(t),\quad t \in \N^*,
\end{equation}
where
\begin{enumerate}
\item $Y(0)$ is Gaussian under $\mathbb Q_\theta$ with the same mean and covariance as $X(0)$ under $\PP_\theta$,
\item $B^\theta(t)\in \R^{d\times d}$ is deterministic and satisfies $B^\theta(t)B^\theta(t)^\top = C^\theta(t)$ for any $\theta \in \Theta$,
\item $W^\theta(t)$ is $\R^d$-valued standard Gaussian under $\mathbb Q_\theta$, independent of $\mathscr G_{t-1}$.
\end{enumerate}
Then $\E_\theta(X(t)) = \E_\theta(Y(t))$ and $\E_\theta(X(t)X(s)^\top) = \E_\theta(Y(t)Y(s)^\top)$ for any $s, t \in \N$. We call $(Y,(\mathbb Q_\theta)_{\theta \in \Theta})$ the \textbf{Gaussian equivalent} of the polynomial state space model $(X,(\PP_\theta)_{\theta \in \Theta})$.
\end{proposition}

\begin{remark}
\begin{itemize}
\item[(i)] It is easy to see that $Y$ in Proposition \ref{prop: Gauss_equiv} is a Gaussian process under any $\mathbb Q_\theta$ and hence its law is uniquely determined by the first and second moments.
\item[(ii)] Strictly speaking, Proposition \ref{prop: Gauss_equiv} does not claim that this Gaussian equivalent actually exists. But this is verified easily by choosing $(\Gamma, \mathscr{G},(\mathscr G_{t})_{t \in \N})$ to be the canonical state space, $Y$ the canonical process on this space, and $\mathbb Q_\theta$ as the Gaussian law with the desired first and second moments. 
\item[(iii)] Since $N^\theta(t)$ is in general not independent of $\F_{t-1}$, the conditional second moments of $X$ and $Y$ can differ substantially. In particular, $Y$ is in general not a (time-homogeneous) Markov
chain because its transition function depends on $t$ through $B^\theta(t)$.
\item[(iv)] Contrary to what Proposition \ref{prop: Gauss_equiv} suggests, our approach heavily relies on higher-order polynomial models. Firstly, an order 2 resp.\ 4 is needed for the asymptotic theory contained in the following chapters and in particular for the computation of the asymptotic covariance matrix. But just as importantly, first and second moments of the original process are typically not enough to estimate high-dimensional parameter vectors of a polynomial process. Instead, we will later consider the Gaussian equivalent of the higher-dimensional extension $((X^\lambda)_{\lambda\in L},(\PP_\theta)_{\theta \in \Theta})$ from Definition \ref{d:orderr} for some suitable subset $L \subseteq \N^d_r\setminus\{0\}$, whose first and second moments carry enough information to estimate the unknown model parameters. For concrete examples we refer to Section \ref{sec5: app}.
\end{itemize}
\end{remark}

By slight abuse of notation, we write $(\Omega, \mathscr{F},(\F_{t})_{t \in \N})$, $(\PP_\theta)_{\theta \in \Theta}$ instead of $(\Omega, \mathscr{G},(\mathscr G_{t})_{t \in \N})$ and $(\mathbb Q_\theta)_{\theta \in \Theta}$ in the following even though the spaces and measures will generally differ.

\subsection{Gaussian quasi-likelihood estimation}\label{s:GQLE}
We generally make the following
\begin{assumption}\label{assump: A}
$(X,(\PP_\theta)_{\theta \in \Theta})$ is a parametric polynomial state space model of order 2 and the mappings $\theta\mapsto a^\theta_{\otimes 2},\, A^\theta_{\otimes 2},\, \E_\theta(X(0)),\,\Cov_\theta(X(0))$ are thrice continuously differentiable in $\theta$.
\end{assumption}
For the explicit representation of the asymptotic covariance matrix of the QML estimator in Section \ref{su:explicit} we will need in fact the slightly stronger
\begin{assumpA}\label{assump: AN}
$(X,(\PP_\theta)_{\theta \in \Theta})$ is a parametric polynomial state space model of order 4 and the mappings $\theta\mapsto a^\theta_{\otimes 2},\, A^\theta_{\otimes 2},\, \E_\theta(X(0)),\,\Cov_\theta(X(0))$ are thrice continuously differentiable in $\theta$.
\end{assumpA}
We denote by $\vartheta\in\Theta$ the unknown ``true'' parameter which is to be estimated based on data. The first $m<d$ components of $X$ are supposed to be unobservable. Put differently, the goal is to estimate $\vartheta$ from observing $X_\oo(s)\coloneq(X_{m+1}(s),\dots,X_d(s))$ in the time interval $s=0,1,\dots,t$. We generally let the subscript $\oo$ stand for the observable part of a vector $x \in \R^d$ and set $H \coloneq (\delta_{m+i, j})_{i=1, \dots, d-m; j=1, \dots, d}$, i.e.\ $x_\oo \coloneq Hx = (x_{m+1}, \dots, x_d)$, and for any matrix $\Sigma\in \R^{d\times d}$ we write $\Sigma_{:, \oo} \coloneq \Sigma H^\top =(\Sigma_{ij})_{i=1,\dots,d; j=m+1,\dots,d}$ as well as $\Sigma_{\oo} \coloneq H \Sigma H^\top = (\Sigma_{ij})_{i=m+1,\dots,d; j=m+1,\dots,d}$. Note that $m=0$ corresponds to the case where all components of $X$ are available for estimation.

Denote by $(Y,(\PP_\theta)_{\theta \in \Theta})$ the Gaussian equivalent of our parametric polynomial state space model, see Proposition \ref{prop: Gauss_equiv}.
Since $Y$ is a Gaussian process under any $\PP_\theta$, the induced measures $\PP_\theta^{(Y_\oo(1), \dots, Y_\oo(t))}$ are absolutely continuous with respect to Lebesgue measure $\lambda_{t(d-m)}$ on $\R^{t(d-m)}$ if the matrices $C^\theta(t)$ are positive definite. This will be warranted by Assumption \ref{normassump}(1) below, see Proposition \ref{prop: B-assump}. Hence it makes sense to define the densities $q_t^{\theta} \coloneq \mathrm{d}\PP_\theta^{(Y_\oo(1), \dots, Y_\oo(t))} / \mathrm{d}\lambda_{t(d-m)}$ as well as the conditional Gaussian densities $\smash{q^{\theta}_{t\,|\,t-1}(\cdot |y_1, \dots, y_{t-1})}$ of the conditional distribution of $Y_\oo(t)$ under $\PP_\theta$ given $Y_\oo(1)=y_1, \dots, Y_\oo(t-1)=y_{t-1}$ for $t \in \N^*$. Then we have the decomposition
\[q_t^{\theta}(Y_\oo(1), \dots, Y_\oo(t)) = \prod_{s = 1}^t q_{s \, | \, s- 1}^{\theta}\big(Y_\oo(s) |Y_\oo(1), \dots, Y_\oo(s-1)\big),\]
where $\smash{q_{1\,|\,0}^{\theta}}\coloneq q_1^\theta$ is set to be the unconditional density of $Y_\oo(1)$ under $\PP_\theta$. It is well known that the log-likelihood $\log q_{t \, | \, t-1}^{\theta}$ can be expressed in terms of the Kálmán filter.

Since we do not assume that the state space $E$ is of positive Lebesgue measure on $\R^d$, proper maximum likelihood estimation for $X$ may not even be well-defined because it is not assured that the analogous Lebesgue densities $p_t^{\theta}$ and $\smash{p_{t \, | \, t - 1}^{\theta}}$ for $X$ exist at all. But even if these densities exist, statistical inference cannot be based on the generally unknown $p_t^{\theta}$. The key idea of the QML approach is therefore to replace $p_t^{\theta}$ \mbox{by the Gaussian density $q_t^{\theta}$, see (\ref{eq: log-lik0}, \ref{eq: log-lik}) below.}
\begin{proposition}\label{prop: log-lik}
Consider the Kálmán filter equations for the Gaussian process $Y$ under $\PP_\theta$, but applied to the process $X$ given the observations $X_\oo$. These read as follows:
\begin{align*}
\widehat X^\theta(0, -1) &\coloneq \E_\theta(X(0)),\\
\widehat \Sigma^\theta(0, -1)&\coloneq \Cov_\theta(X(0)),\\
\widehat{X}^\theta(t + 1, t) &\coloneq a^\theta + A^\theta\widehat X^\theta(t, t),   \\
\widehat{X}^\theta(t, t) &\coloneq \widehat{X}^\theta(t, t-1) + \widehat{\Sigma}^\theta_{ :, \oo} (t, t-1) \widehat{\Sigma}^\theta_{\oo}(t, t-1)^{+}\big(X_\oo(t) - \widehat{X}_{\oo}^\theta(t, t-1)\big), \\
\widehat{\Sigma}^\theta(t + 1, t) &\coloneq A^\theta \widehat{\Sigma}^\theta(t, t) (A^\theta)^\top + C^\theta(t + 1), \\
\widehat{\Sigma}^\theta(t, t) &\coloneq \widehat{\Sigma}^\theta(t, t-1) - \widehat{\Sigma}^\theta_{ :, \oo}(t, t-1) \widehat{\Sigma}^\theta_{\oo}(t, t-1)^{+}\widehat{\Sigma}^\theta_{:, \oo}(t, t-1)^\top
\end{align*} 
for $t \in \N$, where the superscript ${}^+$ stands for the Moore--Penrose pseudoinverse. If $\widehat \Sigma_\oo^\theta(t, t-1)$ is invertible for all $t \in \N$, then, setting $\varepsilon^\theta(t) \coloneq X_\oo(t) - \widehat{X}_{\oo}^\theta(t, t-1)$ for $t \in \N^*$, the Gaussian quasi log-likelihood is of the form
\begin{equation}\label{eq: log-lik0}
L^\theta(t)\coloneq \log q_t^{\theta}(X^\vartheta_\oo(1), \dots, X^\vartheta_\oo(t)) = \sum_{s = 1}^t L^\theta(s, s-1)
\end{equation}
for $t \in \N^*$, where (up to an additive constant independent of $\theta$)
\begin{align}\label{eq: log-lik}
L^\theta(t, t-1) = -\frac{1}{2} \Big[ \log \big\lvert\det \widehat{\Sigma}_{\oo}^\theta(t, t-1)\big\rvert + \varepsilon^\theta(t)^\top \widehat{\Sigma}_{\oo}^\theta(t, t-1)^{-1} \varepsilon^\theta(t)\Big].
\end{align}
\end{proposition}

We naturally refer to a maximiser $\widehat\theta(t)$ of \eqref{eq: log-lik0} as a QML estimator, see Definition \ref{def: Gn-Estimator} below.
Instead of considering maximisers of $L^\theta(t)$, one could alternatively define an estimator to be a solution to the score equation $Z^\theta(t) = 0$, where $Z^\theta(t)\coloneq \nabla_\theta L^\theta(t)$ denotes the quasi-score function corresponding to our quasi-likelihood.  This approach is for example pursued in \cite{Jacod2018} in the context of general estimating functions. As is explained in \cite{Heyde1997} or \cite{Soerensen2012}, desirable asymptotic properties of the QML estimator often hold if the  estimating function is unbiased in the sense that $\E_\vartheta(Z^\vartheta(t)) = 0$. More specifically, we may hope for consistency and asymptotic normality with asymptotic covariance of the form $V_\vartheta=W(\vartheta)^{-1}U_\vartheta W(\vartheta)^{-1}$, where $\smash{W(\vartheta)=\lim_{t\to\infty}\frac{1}{t}\nabla_\vartheta Z^\vartheta(t)}$ is a kind of Fisher information matrix for the quasi-likelihood under consideration and $\smash{U_\vartheta=\lim_{t\to\infty}\frac{1}{t}\Cov_\vartheta(Z^\vartheta(t))}$ denotes the asymptotic covariance matrix of the estimating function.

Unbiasedness of the score function typically holds for the true likelihood rather than the quasi-likelihood. But since our Gaussian quasi-score is a quadratic function of the observations and since the first two moments of our model coincide with those of its Gaussian equivalent, we naturally obtain that unbiasedness of the estimating function holds in our quasi-likelihood setup as well.

A number of challenges have to be met. Firstly, we must establish sufficient ergodic properties for the above programme to work in our setup. Secondly, we need to express $W(\vartheta), U_\vartheta$ in a way that is computable in practice. Two different approaches are provided in Sections \ref{sec5.1: CovMat} and \ref{su:explicit} below. The first one is based on observations --- from either data or Monte-Carlo simulations --- of $\smash{\nabla_\theta Z^\theta(t)|_{\theta=\widehat\theta(t)}}$ as well as of $\smash{Z^\theta(t)|_{\theta=\widehat\theta(t)}}$.
The alternative solution from Section \ref{su:explicit} relies on expressing $W(\vartheta), U_\vartheta$ in closed-form as a function of $\vartheta$ and using the fact that $\widehat\theta(t)\to\vartheta$ as $t\to\infty$ in $\PP_\theta$-probability. Due to the complexity of this explicit representation, the implementation of this approach is slightly more involved. But it is rewarded by a much higher computational efficiency because no simulation is needed.

In order to establish the outlined asymptotic properties of the QML estimating sequence, we naturally have to impose a couple of additional assumptions, summarised below.

\begin{assumption}\label{normassump}
We assume the following two conditions:
    \begin{enumerate}
        \item\label{itm: Irreducibility} For any $\theta \in \Theta$, the transition measures $P_\theta(x, \cdot)$ of $X$ are equivalent to the Hausdorff measure $\lambda_E$ for any $x \in E$ and $X$ is weakly Feller, i.e. $P_\theta(x, \cdot) \xrightarrow{x \to x_0} P_\theta(x_0, \cdot)$ weakly for all $x_0 \in E$.
        \item\label{itm: L6-bounded} There is some $\delta > 0$ such that $X$ is bounded in $L^{4 + \delta}(\PP_\theta)$ for any $\theta \in \Theta$.
    \end{enumerate}
\end{assumption}

\begin{remark}
    Assumption \ref{normassump}(\ref{itm: Irreducibility}) implies that the conditional distributions of $X(t)$ given $X(t-1)$ under any $\PP_\theta$ are equivalent to the Hausdorff measure $\lambda_E$ on $E$. By the Chapman--Kolmogorov equation then also the unconditional laws of $X(t)$ are equivalent to $\lambda_E$.
\end{remark}

We note that the conditions in Assumption \ref{normassump} can be slightly weakened if $X$ is obtained by sampling from an affine process with existing finite sixth moments. In this case, the Riccati differential equations for the characteristic function of $X$ can be used to yield particularly simple expressions for the matrices $A^\theta_{\otimes r}$ for orders $r \in \{2, \dots, 6\}$, see Remark \ref{rem: Mean_Recursion}:

\begin{proposition}\label{prop: affine_lp}
    Suppose that $X = (X(t))_{t \in \N}$ is obtained by sampling from an affine process $(X(t))_{t \in \R_+}$ and fix an even $p \in \N^*$. If $X(t) \in L^p(\PP_\theta)$ and $\rho(A^\theta)<1$ for any $t \in \N$ and $\theta \in \Theta$, then $X$ is bounded in $L^p(\PP_\theta)$. Hence, if $p \geq 6$ and $\rho(A^\theta) < 1$, \mbox{condition \ref{itm: L6-bounded} of Assumption \ref{normassump} is fulfilled.}
\end{proposition}

We can now prove an ergodic theorem for the process $X$ which moreover gives a statement about the speed of convergence for expectations of quadratic resp.\ quartic polynomials in $X$ and which justifies the condition $\rho(A^\theta) < 1$ from Proposition \ref{prop: affine_lp}. The following proposition uses the terminology of strong $f$-ergodicity of Markov chains in the following sense:
\begin{definition}\label{defi: f-ergod}
    Let $X$ be an $E$-valued discrete-time process on some space $(\Omega, \F, \PP)$ and $\mu$ a probability measure on $\mathscr{B}(E)$. Moreover, let $f: E \to \R$ be $\mu$-integrable. We call $X$ \textbf{strongly $\bm{f}$-ergodic} (with respect to $\mu$ under $\PP$) if the following conditions hold:
    \begin{enumerate}
        \item $\PP^{X(t)} \to \mu$ weakly as $t \to \infty$.
        \item $\frac{1}{t} \sum_{s=1}^t f\big(X(s)\big) \xrightarrow{t \to \infty} \int f \dd \mu$ $\PP$-almost surely.
    \end{enumerate}
    If convergence in (ii) holds only in probability, then $X$ is called weakly $f$-ergodic with respect to $\mu$ under $\PP$. Whenever there is no ambiguity concerning the measures $\PP$ and $\mu$ under consideration, we just call $X$ strongly (resp. weakly) $f$-ergodic in the above definition.
\end{definition}

\begin{proposition}\label{prop: ergod}
    Suppose that Assumptions \ref{assump: A} (resp. \ref{assump: AN}) and \ref{normassump} hold and let $\theta\in\Theta$. Then there exists a unique stationary law $\mu_\theta$ for $X$ under $\PP_\theta$ and $X$ is strongly $f$-ergodic for any $\mu_\theta$-integrable function $f: E \to \R$. Moreover, $\mu_\theta$ has finite moments of order $4 + \delta$ and $\rho(A^\theta_{\otimes r}) < 1$ for $r = 1, 2$ (resp. $r = 1, \dots, 4$). Finally, $\E_\theta(g(X(t))) \to \int g \dd\mu_\theta$ at a geometric rate for any quadratic (resp. quartic) polynomial $g: \R^d \to \R$.
\end{proposition}

The following proposition is a simple consequence of Proposition \ref{prop: ergod}.

\begin{proposition}\label{prop: B-assump}
    Suppose that Assumptions \ref{assump: A} and \ref{normassump} hold. Then $C^\theta(t)$ is positive definite and thrice continuously differentiable in $\theta$ for any $t \in \N^*$ and $\theta \in \Theta$. Moreover, for any multi-index $\alpha \in \N^k_3$ the matrices $\partial_\theta^\alpha C^\theta(t)$ converge at a geometric rate uniformly in $\theta \in \Theta$ to some limits and $\smash{C^\theta} \coloneq \lim_{t \to \infty} \smash{C^\theta(t)}$ is positive definite.
\end{proposition}

Note that since $C^\theta(t)$ are positive definite by Proposition \ref{prop: B-assump}, the matrices $\smash{\widehat \Sigma^\theta(t, t-1)}$ in Proposition \ref{prop: log-lik} are positive definite for all $t \in \N$, hence the inverse in \eqref{eq: log-lik} is well-defined.

Let $\alpha^\theta = (\mathrm{I}_d - A^\theta)^{-1} a^\theta$ denote the solution to the fixed point equation $\alpha^\theta = a^\theta + A^\theta \alpha^\theta$ and $\Lambda^\theta$ the solution to the Lyapunov equation $\Lambda^\theta = A^\theta \Lambda^\theta A^{\theta^\top} + C^\theta$. Since $\rho(A^\theta) < 1$ by Proposition \ref{prop: ergod} and since $C^\theta$ is positive definite by Proposition \ref{prop: B-assump}, these solutions exist uniquely and $\Lambda^\theta$ is positive definite by Lemma \ref{lem: lin_systems}.2. Moreover, Lemma \ref{lem: lin_systems} also yields that $\alpha^\theta$ and $\Lambda^\theta$ are the stationary mean and covariance matrix of $X$ under its invariant probability measure $\mu_\theta$ from Proposition \ref{prop: ergod}. As a further condition for the main results, we need some form of identifiability for the parameter $\theta \in \Theta$, i.e.\ a sufficient distinguishability of different parameters $\theta \neq \theta'$. To this end we consider the Markov process
\begin{equation}\label{eq: Gauss2}
Y^{\mathrm{stat}}(t) = a^\theta + A^\theta Y^{\mathrm{stat}}(t-1) + B^\theta W^\theta(t)    
\end{equation} 
with independent standard normal $W^\theta(t), t\in\N^*$ under $\PP_\theta$ and initial distribution $\PP_\theta^{Y^{\mathrm{stat}}(0)} = \mathscr{N}(\alpha^\theta, \Lambda^\theta)$. This process differs from the Gaussian equivalent \eqref{Gauss} in Proposition \ref{prop: Gauss_equiv} only by the moments at time 0 and in the sense that the time-dependent covariance matrices $B^\theta(t)B^\theta(t)^\top=C^\theta(t)$ are replaced by the limiting covariance matrix $\smash{C^\theta} =: B^\theta {B^\theta}^\top$. It is easy to see that $Y^{\mathrm{stat}}$  is in fact stationary, i.e.\ the unique invariant probability measure for the Gaussian process $Y^{\mathrm{stat}}$ is given by $\mathscr{N}(\alpha^\theta, \Lambda^\theta)$. Let $\smash{\widehat Y^{\theta,\mathrm{stat}}(t, t-1)}$ and $\smash{\widehat \Sigma^{\theta,\mathrm{stat}}(t, t-1)}$ denote the Kálmán filter with corresponding covariance matrices for $Y^{\mathrm{stat}}$ instead of $Y$, see Proposition \ref{prop: log-lik}. In particular, $\widehat Y^{\theta,\mathrm{stat}}(0, -1) = \alpha^\theta$ and $\widehat \Sigma^{\theta,\mathrm{stat}}(0, -1) = \Lambda^\theta$. We will show in Section \ref{sec4.2.3: Ident} that the following assumption suffices for identifiability of the parameter $\theta$ in a suitable sense:

\begin{assumption}\label{assump: Identifiability}
For every $\theta \neq \theta'$ with $\theta, \theta' \in \Theta$, there exists some integer $u \in \N$ such that $\smash{\widehat \Sigma^{\theta,\mathrm{stat}}_\oo(u, u-1) \neq \widehat \Sigma^{\theta',\mathrm{stat}}_\oo(u, u-1)}$ or $\smash{\widehat Y_\oo^{\theta,\mathrm{stat}}(u, u-1)}$ differs from $\smash{\widehat Y_\oo^{\theta',\mathrm{stat}}(u, u-1)}$ as a deterministic function of the observations $\smash{Y^{\mathrm{stat}}_\oo(0),\dots, Y^{\mathrm{stat}}_\oo(u-1)}$.
\end{assumption}

\begin{remark}
    The above assumption may be somewhat laborious to verify in concrete applications. Most of the time when Assumption \ref{assump: Identifiability} is fulfilled, it is however actually the case that it is fulfilled already with $u = 0$, i.e.\ the parameter $\theta$ is identifiable from the observed part of the stationary mean and covariance matrix. We could therefore replace Assumption \ref{assump: Identifiability} by the following stronger assumption, which is easier to verify in practice:
\end{remark}
\begin{assumpD}\label{assump: Identifiability2}
For every $\theta \neq \theta'$ with $\theta, \theta' \in \Theta$ one has $\alpha^\theta_\oo \neq \alpha^{\theta'}_\oo$ or $\Lambda^\theta_\oo \neq \Lambda^{\theta'}_\oo$.
\end{assumpD}

\section{Main results}\label{sec4.2: Proof}
\subsection{Consistency and asymptotic normality}\label{su: mainmain}
We now turn to the quasi-maximum likelihood (QML) estimation of the unknown parameter $\vartheta$ by focusing on a sequence of maximisers $\smash{\widehat \theta(t)}$ of $L^\theta(t)$ from \eqref{eq: log-lik0}. As in Section \ref{s:GQLE} we fix a parametric state space model $(X,(\PP_\theta)_{\theta \in \Theta})$ with ``true'' parameter $\vartheta\in\Theta$. Moreover, we assume that Assumptions \ref{assump: A}, \ref{normassump} and \ref{assump: Identifiability} hold. Since the parameter space $\Theta$ is assumed to be compact and since the quasi-likelihood $L^\theta(t)$ is a continuous function of $\theta$ by Assumption \ref{assump: A}, the existence of such maximisers is always guaranteed in our setting.
\begin{definition}\label{def: Gn-Estimator}
    Let $t \in \N $. A \textbf{quasi-maximum likelihood estimator} for the polynomial state space model $(X,(\PP_\theta)_{\theta \in \Theta})$ is any $\F_t $-measurable random variable $\smash{\widehat \theta(t)}$ with values in the parameter space $\Theta$ such that $\smash{\widehat \theta(t) \in \mathrm{arg}\max_{\theta \in \Theta} L^\theta(t)}$, where $L^\theta(t)$ is defined in \eqref{eq: log-lik0}.
\end{definition}

As usual in the statistical literature, we call a sequence $(\widehat \theta(t))_{t\in\N }$ \textbf{(weakly) $\vartheta$-consistent} if $\smash{\widehat \theta(t) \xrightarrow{\PP_\vartheta} \vartheta}$. Recall that $\nabla_\theta Z^\theta(t)$ denotes the $(k \times k)$-matrix whose $ij$-th entry is $\partial_{\theta_i} Z^\theta(t)_j$. 

\begin{remark}
    It can be shown that $\F_t $-measurability in Definition \ref{def: Gn-Estimator} of a QML estimator is not a restrictive assumption: In particular, since $\Theta$ is compact and $L^\theta(t)$ is continuous in $\theta$, \cite{Jennrich1969}, Lemma 2, shows existence of an $(\F_t )_{t \in \N }$-adapted sequence $\smash{(\widehat \theta(t))_{t \in \N }}$ of QML estimators as in Definition \ref{def: Gn-Estimator}. 
\end{remark}

\begin{theorem}\label{theo: main}
Any sequence of QML estimators $(\widehat \theta(t))_{t \in \N }$ is $\vartheta$-consistent. 
\end{theorem}
The proof of this and the following theorems is to be found in Section \ref{s:proofs}.

\begin{theorem}\label{theo: main2}
Suppose that $\vartheta \in \mathrm{int}(\Theta)$. Then we have:
\begin{enumerate}
\item The limit $\lim_{t \to \infty} \frac{1}{t} \nabla_\theta Z^\theta(t)\big|_{\theta = \widehat \theta(t)} = W(\vartheta) \in \R^{k \times k}$ exists in $\PP_\vartheta$-probability. 
\item There exists a positive semidefinite matrix $U_\vartheta \in \R^{k \times k}$ such that $\smash{\frac{1}{\sqrt{t}} Z^\vartheta(t)}$ converges in $\PP_\vartheta$-law to some Gaussian random variable with mean 0 and covariance matrix $U_\vartheta$.
\item If the matrix $W(\vartheta) \in \R^{k \times k}$ is invertible, there exists some $\R^k$-valued Gaussian random variable $\smash{Z \sim \mathcal{N}\big(0, V_\vartheta\big)}$ with $V_\vartheta\coloneq W(\vartheta)^{-1}U_\vartheta W(\vartheta)^{-1}$ such that we have
\[\sqrt{t}\big(\widehat \theta(t) - \vartheta\big) \xrightarrow{\PP_\vartheta\text{-}d} Z.\]
\end{enumerate}
\end{theorem}

In order to be applicable in practice, we need more explicit estimates or representations of $V_\vartheta$. As noted in Section \ref{s:GQLE}, two conceptually different ones are provided in the following Sections \ref{sec5.1: CovMat} and \ref{su:explicit}. The former approach only requires Assumption \ref{assump: A}, while the latter relies on the stronger Assumption \ref{assump: AN}.

\subsection{Estimation of the asymptotic covariance matrix}\label{sec5.1: CovMat}
The approach in this section is based on the quasi-score process $Z^\theta(t) = \nabla_\theta L^\theta(t)$ and its derivative $\nabla_\theta Z^\theta(t)$. If we define the conditional score process $\smash{Z^\theta(t, t-1) \coloneq \nabla_\theta L^\theta(t, t-1)}$, one has $Z^\theta(t) = \sum_{s=1}^t Z^\theta(s, s-1)$. Differentiating the Gaussian quasi-log-likelihood \eqref{eq: log-lik} now yields the following expression.

\begin{proposition}\label{prop: Quasi-Score} For $j \in \{1, \dots, k\}$, the $j$-th component of $Z^\theta(t, t-1)$ is given by
    \begin{equation}\label{eq: Score}
        Z^\theta(t, t-1)_j = -\frac{1}{2}\Big[\kappa^{\theta}_j(t) - 2 V^{\theta}_{\oo,j}(t,t-1)^\top\widehat \Sigma^{\theta}_{\oo}(t, t-1)^{-1} \varepsilon^\theta(t) - \varepsilon^\theta(t)^\top \widetilde S^\theta_{\oo,j}(t) \varepsilon^\theta(t) \Big]
    \end{equation}
    with $\kappa^{\theta}_j(t) \coloneq \Tr\big(\widehat \Sigma^\theta_{\oo}(t, t-1)^{-1}S^{\theta}_{\oo,j}(t, t-1)\big)$, initial values $V^{\theta}_j(0, -1) \coloneq \partial_{\theta_j}\E_\theta(X(0))$ as well as $S^{\theta}_j(0, - 1) \coloneq \partial_{\theta_j}\mathrm{Cov}_\theta(X(0))$, where $V^{\theta}_j(t, t-1)$ and $S^{\theta}_j(t, t-1)$ are given by
    \begin{align*}
        V^{\theta}_j(t, t) &\coloneq (\mathrm{I_d} - \widetilde{K}^\theta(t)H) V^{\theta}_j(t, t-1) + N^\theta_j(t) \varepsilon^\theta(t),\\
        V^{\theta}_j(t + 1,t) &\coloneq \partial_{j} a^\theta + (\partial_{j} A^\theta) \widehat X^\theta(t, t) + A^\theta V^{\theta}_j(t, t), \\
        S^{\theta}_j(t, t) &\coloneq (\mathrm{I_d} - \widetilde{K}^\theta(t)H) \big(S^{\theta}_j(t, t-1) - S^{\theta}_{:, \oo, j}(t, t-1)\widetilde K^\theta(t)^\top\big) \\
        S^{\theta}_j(t + 1, t) &\coloneq A^\theta \widehat \Sigma^\theta(t, t) (\partial_{j} A^\theta)^\top + (\partial_{j} A^\theta) \widehat \Sigma^\theta(t, t) {A^\theta}^\top + A^\theta S^{\theta}_j(t, t) {A^\theta}^\top + \partial_{j}C^\theta(t), \\
        \widetilde S^{\theta}_{\oo,j}(t) &\coloneq \widehat \Sigma^\theta_{\oo}(t, t-1)^{-1} S^{\theta}_{\oo,j}(t, t-1) \widehat \Sigma^\theta_{\oo}(t, t-1)^{-1}, \\
        N^{\theta}_j(t) &\coloneq \widehat \Sigma_{:, \oo}^\theta(t,t-1) \widetilde S_{\oo,j}^{\theta}(t) - S_{:, \oo, j}^{\theta}(t,t-1) \widehat \Sigma_{\oo}^\theta(t,t-1)^{-1} ,\\
        \widetilde K^\theta(t) &\coloneq \widehat{\Sigma}_{:, \oo}^\theta(t, t-1)\widehat{\Sigma}_{\oo}^\theta(t, t-1)^{-1}
    \end{align*}
for $t \in \N$. In particular, $S^{\theta}_j(t, t - 1) = \partial_{\theta_j} \widehat \Sigma^\theta(t, t - 1)$ and $S^{\theta}_j(t, t) = \partial_{\theta_j} \widehat \Sigma^\theta(t, t)$ for $t \in \N$. Moreover, we have $\smash{V^{\theta}_j(t,t -1)=\partial_{\theta_j}\widehat X^\theta(t, t-1)}$ and $\smash{V^{\theta}_j(t,t)=\partial_{\theta_j}\widehat X^\theta(t, t)}$ for $t \in \N$.
\end{proposition}

Differentiating the Gaussian quasi-log-likelihood once more, we obtain the following representation for the $\R^{k \times k}$-valued process $\nabla_\theta Z^\theta(t) = \sum_{s=1}^t \nabla_\theta Z^\theta(s, s-1)$.
\begin{proposition}\label{prop: Fisher_Information} For $1 \leq i, j \leq k$, the $ij$-th component of $\nabla_\theta Z^\theta(t, t-1)$ is given by
    \begin{equation}\label{eq: Fisher-Information}
        \nabla_\theta Z^\theta(t, t-1)_{ij} = -\frac{1}{2}\Big[\mu^\theta_{ij}(t) + 2 \nu^\theta_{ij}(t)\varepsilon^\theta(t) + \varepsilon^\theta(t)^\top \psi^\theta_{ij}(t) \varepsilon^\theta(t) + 2 \phi^\theta_{ij}(t)\Big]
    \end{equation}
    with the abbreviations $\mu^\theta_{ij}(t) \coloneq \Tr\big[\widehat \Sigma_{\oo}^{\theta}(t, t-1)^{-1} R_{\oo,ij}^\theta(t, t-1) - \widetilde S_{\oo,i}^\theta(t) S^{\theta}_{\oo,j}(t, t-1)\big]$, $\nu^\theta_{ij}(t) \coloneq V_{\oo,i}^\theta(t, t-1)^\top \widetilde S_{\oo,j}^{\theta}(t) + V_{\oo,j}^{\theta}(t, t-1)^\top \widetilde S_{\oo,i}^\theta(t) - W_{\oo,ij}^\theta(t, t-1)^\top \widehat \Sigma_{\oo}^\theta(t, t-1)^{-1}$,  $\psi^\theta_{ij}(t) \coloneq \smash{\widehat \Sigma_{\oo}^{\theta}(t, t-1)^{-1}}\linebreak\widehat S_{\oo,ij}^\theta(t) - \widehat R_{ij}^\theta(t)$, $\phi^\theta_{ij}(t) = V_{\oo,j}^{\theta}(t, t-1)^\top \widehat \Sigma_{\oo}^{\theta}(t, t-1)^{-1} V_{\oo,i}^\theta(t, t-1)$. Here, $W^\theta_{ij}(0, -1) \coloneq \partial_{\theta_i}\partial_{\theta_j}\E_\theta(X(0))$ and $R^\theta_{ij}(0, -1) \coloneq \partial_{\theta_i} \partial_{\theta_j} \Cov_\theta(X(0))$, and we set\\
    \begin{align*}
        W^\theta_{ij}(t + 1, t) &\coloneq \partial_{ij} a^\theta + (\partial_{ij} A^\theta) \widehat X^\theta(t, t) + (\partial_j A^\theta) V^{\theta}_i(t, t) + (\partial_i A^\theta) V^{\theta}_j(t, t) + A^\theta W^\theta_{ij}(t,t),\\
        W^\theta_{ij}(t, t) &\coloneq (\mathrm{I}_d - \widetilde K^\theta(t)H)W^\theta_{ij}(t,t-1) + M^\theta_{ij}(t) \varepsilon^\theta(t) + \widetilde V^\theta_{ij}(t),\\
        R^\theta_{ij}(t + 1, t) &\coloneq \mathrm{Sym}\big[ \widetilde R^\theta_{ij}(t, t)\big] + A^\theta R^\theta_{ij}(t,t) A^{\theta^\top} + \partial_{ij}C^\theta(t),\\
        R^\theta_{ij}(t, t) &\coloneq R^\theta_{ij}(t, t-1) + \mathrm{Sym}\big[\overline R^\theta_{ij}(t, t)\big]  + \widehat R_{ij}^\theta(t) - \widetilde K^\theta(t) \widehat S_{\oo,ij}^\theta(t) \widehat \Sigma_{:, \oo}^\theta(t,t-1)^\top,\\
        \widetilde R^\theta_{ij}(t, t) &\coloneq (\partial_i A^\theta) \widehat \Sigma^{\theta}(t, t) (\partial_j A^\theta)^\top + \overline{S}^\theta_{ij}(t) + A^\theta \widehat \Sigma^\theta(t,t) (\partial_{ij} A^\theta)^\top ,\\
        \widehat R_{ij}^\theta(t) &\coloneq \widehat \Sigma_{:, \oo}^{\theta}(t, t-1)\widehat \Sigma^\theta_{\oo}(t,t-1)^{-1} R^\theta_{\oo,ij}(t,t-1) \widehat \Sigma^\theta_{\oo}(t,t-1)^{-1}\widehat \Sigma_{:, \oo}^{\theta}(t,t-1)^\top ,\\
        \overline R^\theta_{ij}(t, t) &\coloneq \widehat S_{ij}^\theta(t) \widehat \Sigma_{:, \oo}^\theta(t, t-1)^\top - \widetilde K^\theta(t) R_{:, \oo, ij}^\theta(t, t-1)^\top - \check{S}^\theta_{ij}(t) ,\\
        M^\theta_{ij}(t) &\coloneq (\mathrm{I}_d - \widetilde K^\theta(t) H) \big(R_{:, \oo, ij}^\theta(t, t-1)\widehat \Sigma_{\oo}^\theta(t, t-1)^{-1} - \widehat S_{ij}^\theta(t)\big),\\
        \widetilde V^\theta_{ij}(t) &\coloneq N^{\theta}_j(t)V_{\oo,i}^\theta(t, t-1) + N^\theta_i(t)V_{\oo,j}^{\theta}(t, t-1),\\
        \widehat S_{ij}^\theta(t) &\coloneq S_{:, \oo, j}^{\theta}(t,t-1) \widetilde S_{\oo,i}^\theta(t) + S_{:, \oo, i}^\theta(t,t-1) \widetilde S_{\oo,j}^{\theta}(t) ,\\
        \overline{S}^\theta_{ij}(t) &\coloneq A^\theta S^{\theta}_i(t, t) (\partial_j A^\theta)^\top + A^\theta S^{\theta}_j(t, t) (\partial_i A^\theta)^\top \\
        \check{S}^\theta_{ij}(t) &\coloneq S_{:, \oo, j}^{\theta}(t,t-1) \widehat \Sigma_{\oo}^\theta(t,t-1)^{-1}S_{:, \oo, i}^\theta(t,t-1)^\top
    \end{align*}
    
    \noindent for any $t\in\N$, where we adopt the notation and abbreviations from Proposition \ref{prop: Quasi-Score}. In particular, $R^\theta_{ij}(t + 1, t) = \partial_{\theta_i} \partial_{\theta_j} \widehat \Sigma^\theta(t + 1, t)$ and $R^\theta_{ij}(t, t) = \partial_{\theta_i} \partial_{\theta_j} \widehat \Sigma^\theta(t, t)$ for $t \in \N$. Moreover, we have $W^{\theta}_{ij}(t,t -1)=\partial_{\theta_i}\partial_{\theta_j}\widehat X^\theta(t, t-1)$ and $W^{\theta}_{ij}(t,t)=\partial_{\theta_i}\partial_{\theta_j}\widehat X^\theta(t, t)$ for $t \in \N$.
\end{proposition}

Note that the coefficients in the preceding propositions stabilise in the following sense:
\begin{lemma}\label{limitcov}
    There exists a positive definite matrix $\widehat\Sigma^\theta$ as well as symmetric matrices $S^\theta_j$, $R^\theta_{ij}$ for $\theta \in \Theta$ and $i,j \in \{1, \dots, k\}$ such that, uniformly in $\theta \in \Theta$, we have $\smash{\widehat \Sigma^\theta(t, t-1) \to \widehat\Sigma^\theta}$, $\smash{S^\theta_j(t, t-1) \to S^\theta_j}$, and $\smash{R^\theta_{ij}(t, t-1) \to R^\theta_{ij}}$ at a geometric rate as $t \to \infty$.
\end{lemma}

Since $\smash{W(\vartheta)=\PP_\vartheta\text{-}\lim_{t \to \infty} \frac{1}{t} \nabla_\theta Z^\theta(t)\big|_{\theta = \widehat \theta(t)}}$, we can use $\smash{\widehat W_t[\widehat\theta(t)]\coloneq\frac{1}{t} \nabla_\theta Z^\theta(t)\big|_{\theta = \widehat \theta(t)}}$ as a consistent estimator, where $\smash{\nabla_\theta Z^\theta(t)}$ is obtained using the recursion in Proposition \ref{prop: Fisher_Information}. The situation is less obvious for $U_\vartheta$. Since the increments $Z^\theta(s,s-1)$, $s=1,\dots,t$ of $Z^\theta(t)$ are dependent, it may not be immediately obvious how to express its limiting covariance matrix in terms of a law of large numbers. Due to the Markovian structure of the process $(X(t),\widehat X^\theta(t,t-1),V^\theta(t,t-1))_{t\in \N}$, however, such a representation can be established by an approach going back to \cite{Gordin1978}.

To do so we first need the following matrix product introduced by \cite{Tracy1972}, which basically behaves like a block-wise Kronecker product for partitioned matrices:

\begin{definition}
    Let the matrix $A = (A_{ij})_{ij} \in \R^{m \times n}$ be partitioned into submatrices $A_{ij} \in \R^{m' \times n'}$, and let the matrix $B = (B_{kl})_{kl} \in \R^{p \times q}$ be partitioned into submatrices $B_{kl} \in \R^{p' \times q'}$. Then we define the \textbf{Tracy--Singh product} $A \circ^{m' \times n'}_{p'\times q'} B \coloneq ((A_{ij} \otimes B_{kl})_{ij})_{kl} \in \R^{mp \times nq}$ of $A$ and $B$, i.e. the $kl$-th subblock of the $ij$-th subblock of $A\circ^{m' \times n'}_{p'\times q'} B$ is given by $A_{ij} \otimes B_{kl}$.
\end{definition}

 It is clear from the above definition that the Tracy--Singh product of two partitioned matrices depends crucially on the chosen partition of the matrices into submatrices. If, for example, 

\begin{equation*}
    A = \left(\begin{array}{c | c}
        A_{11} & A_{12} \\ \hline A_{21} & A_{22}
    \end{array}\right) \quad \text{ and } \quad B = \left(\begin{array}{c | c}
        B_{11} & B_{12} \\ \hline B_{21} & B_{22}
    \end{array}\right)
\end{equation*}
are partitioned into $(m\times n)$-blocks $(A_{ij})_{i, j \in \{1, 2\}}$ and $(p\times q)$-blocks $(B_{kl})_{k, l \in \{1, 2\}}$, then

\begin{equation*}
    A\circ^{m \times n}_{p\times q} B = \left(\begin{array}{c|c|c|c}
        A_{11} \otimes B_{11} & A_{11} \otimes B_{12} & A_{12} \otimes B_{11} & A_{12} \otimes B_{12} \\ \hline
        A_{11} \otimes B_{21} & A_{11} \otimes B_{22} & A_{12} \otimes B_{21} & A_{12} \otimes B_{22} \\ \hline
        A_{21} \otimes B_{11} & A_{21} \otimes B_{12} & A_{22} \otimes B_{11} & A_{22} \otimes B_{12} \\ \hline
        A_{21} \otimes B_{21} & A_{21} \otimes B_{22} & A_{22} \otimes B_{21} & A_{22} \otimes B_{22}
    \end{array}\right).
\end{equation*}

 For a given $k$-dimensional polynomial $f: \R^d \to \R^k$ of order $r$, we say that $\alpha_f \in \smash{\R^{k \times \frac{d}{d-1}(d^r - 1)}}$ and $\beta_f \in \R^k$ denote \textit{a pair of coefficients} of $f$ if $f(x) = \alpha_f \mathrm{vec}_{\otimes r}(x) + \beta_f$ for all $x \in \R^d$. Of course, since the stacked Kronecker product $\mathrm{vec}_{\otimes r}(x)$ contains some powers of $x$ multiple times, the choice of coefficients $\alpha_f$ and $\beta_f$ for a given $f$ is not unique, which will however not pose a problem in the sequel. Since $X$ is a polynomial state space model of order 2 with corresponding martingale difference sequence $N^\theta$, there exist some matrices $Q^\theta_{\otimes 2} \in \R^{d^2 \times d^2}$ and $Q^\theta \in \R^{d^2 \times d}$ as well as some $q^\theta \in \R^{d^2}$ such that
\[\E_\theta\big[ N^\theta(t) \otimes N^\theta(t) \: \big|\: X(t-1) = x\big] = Q^\theta_{\otimes 2} (x \otimes x) + Q^\theta x + q^\theta\]
for any $x \in E$ and $t \in \N^*$. Note that by Proposition \ref{prop: affine_lp} and in particular by equation \eqref{eq: otimes-matrix} in its proof, $Q^\theta_{\otimes 2} = 0$ if $X$ is obtained by discrete sampling from an affine process.

For the following statement we introduce some more notation. Set
$\kX(t)\coloneq \widehat{X}^\theta(t, t-1)$, $\wX(t)\coloneq V^\theta(t,t-1)$, $\oX(t)\coloneq W^\theta(t,t-1)$, $\langX(t) \coloneq (X(t), \kX(t), \wX(t))$ as well as
\begin{align}
K^\theta  &\coloneq A^\theta \widehat \Sigma^\theta_{:, \oo} (\widehat\Sigma^\theta_{\oo})^{-1}, \qquad F^\theta  \coloneq A^\theta - K^\theta H ,\nonumber\\
   \smash{\overline Q}^\theta_{\otimes 2} &\coloneq \Big(e_1^{(k+2)}e_1^{(k+2)\top}\Big) \: \circ^{(k+2) \times (k+2)}_{d \times d} \:
        \Big[\Big(e_1^{(k+2)}e_1^{(k+2)\top}\Big) \otimes Q^\theta_{\otimes 2}\Big],  \nonumber\\
    \smash{\overline Q}^\theta &\coloneq \Big(e_1^{(k+2)}e_1^{(k+2)\top}\Big) \: \circ^{(k+2) \times (k+2)}_{d \times d} \:
        \Big[e_1^{(k+2)} \otimes Q^\theta\Big],  \nonumber\\
    \smash{\overline q}^\theta &\coloneq e_1^{(k+2)} \: \circ^{(k+2) \times 1}_{d \times 1} \:
        \Big[e_1^{(k+2)} \otimes q^\theta\Big],\nonumber\\
    \smash{\overline a}^\theta &\coloneq \begin{pmatrix}
        a^\theta \\ a^\theta \\ \nabla_\theta a^\theta
    \end{pmatrix},  \qquad\,\quad\smash{\overline A}^\theta \coloneq \begin{pmatrix}
        A^\theta & 0 & 0 \\
        K^\theta H & F^\theta & 0 \\
        \nabla_{\theta} K^\theta H & \nabla_{\theta} F^\theta & \mathrm{diag}_k(F^\theta)
    \end{pmatrix}, \label{eq: Ahat1_terms}\\
    \Pi_\theta &\coloneq \smash{\overline a}^\theta \otimes \smash{\overline A}^\theta + \smash{\overline A}^\theta \otimes \smash{\overline a}^\theta + \smash{\overline Q}^\theta, \;\;\;\;\: O_\theta \coloneq \smash{\overline A}^\theta \otimes \smash{\overline A}^\theta + \smash{\overline Q}_{\otimes 2}^\theta, \nonumber \\
    \smash{\overline a}^\theta_{\otimes 2} &\coloneq \:\begin{pmatrix}\;
        \smash{\overline a}^\theta \hphantom{,} \\ \;\smash{\overline a}^\theta \otimes \smash{\overline a}^\theta + \smash{\overline q}^\theta\hphantom{,}
    \end{pmatrix}, \qquad \quad\;\;\;\smash{\overline A}^\theta_{\otimes 2} \coloneq \:\begin{pmatrix}
        \;\smash{\overline A}^\theta\; & 0 \hphantom{,} \\
        \;\Pi_\theta\; & O_\theta \hphantom{,}
    \end{pmatrix} \label{eq: Ahatotimes2_terms}.
\end{align}

\begin{proposition}\label{prop: Formula_Cov}
Define the matrices 
\begin{align*}
\Gamma^\theta_1 &\coloneq H^\top (\widehat\Sigma^\theta_{\oo})^{-1} H, \qquad \Gamma_{2,j}^{\theta} \coloneq H^\top (\widehat\Sigma^\theta_{\oo})^{-1} S^\theta_{j, \oo} (\widehat\Sigma^\theta_{\oo})^{-1} H,\\
 \Gamma_\theta^{(j)} &\coloneq \frac{1}{2}\begin{pmatrix}
        \hphantom{-}\Gamma_{2,j}^{\theta} & - \Gamma_{2,j}^{\theta} & \hphantom{-}\big(e_j^{(k)} \otimes \Gamma^\theta_1\big)^\top \\
        - \Gamma_{2,j}^{\theta} & \hphantom{-}\Gamma_{2,j}^{\theta} & - \big(e_j^{(k)} \otimes \Gamma^\theta_1\big)^\top \\
        e_j^{(k)} \otimes \Gamma^\theta_1 & - e_j^{(k)} \otimes \Gamma^\theta_1 & 0
    \end{pmatrix} \in \R^{(k+2)d \times (k+2) d}
\end{align*}
for $j \in \{1, \dots, k\}$, as well as the matrix $\Gamma^\theta \in \R^{k \times (k+2)^2d^2}$ such that its $j$-th row is given by $\mathrm{vec}(\Gamma_\theta^{(j)})$.
Then the matrix $U_\vartheta$ from Theorem \ref{theo: main2} is the limit in $\PP_\vartheta$-probability of the matrices
\begin{equation}\label{eq: Formula_Cov}
    \widehat U_t[\widehat \theta(t)] \coloneq \frac{1}{t}\sum_{s=1}^t \Big[g^\theta\big(\langX(s)\big)g^\theta\big(\langX(s)\big)^\top - h^\theta\big(\langX(s)\big)h^\theta\big(\langX(s)\big)^\top\Big]\bigg|_{\theta = \widehat \theta(t)},
\end{equation}
where $g^\theta: \R^{(k+2)d} \to \R^k$ and $h^\theta: \R^{(k+2)d} \to \R^k$ are $k$-dimensional polynomials with coefficients $\alpha^\theta_g = \big(- \Gamma_\theta (O_\theta - \mathrm{I})^{-1} \Pi_\theta (\overline A^\theta_1 - \mathrm{I})^{-1} \;\: \big| \;\: \Gamma^\theta (O_\theta - \mathrm{I})^{-1} \big)$, $\beta^\theta_g = 0$, $\alpha^\theta_h = \alpha^\theta_g \overline A^\theta_{\otimes 2}$ as well as $\beta^\theta_h = \alpha^\theta_g \overline a^\theta_{\otimes 2}$. 
\end{proposition}

 Since the covariance matrix in Theorem \ref{theo: main2} equals $V_\vartheta = W(\vartheta)^{-1} U_\vartheta W(\vartheta)^{-1}$, it can be consistently estimated by $\smash{\widehat V_t[\widehat \theta(t)] \coloneq \widehat W_t[\widehat \theta(t)]^{-1} \widehat U_t[\widehat \theta(t)] \widehat W_t[\widehat \theta(t)]^{-1} \mathbf{1}_{\{\widehat W_t[\widehat \theta(t)] \text{ invertible}\}}}$ with $\smash{\widehat W_t[\theta]} \coloneq \smash{\frac{1}{t}\nabla_\theta Z^\theta(t)}$ and where $\widehat U_t[\widehat \theta(t)]$ is given in \eqref{eq: Formula_Cov}. If $U_\vartheta$ is invertible, then $\PP(\widehat U_t[\widehat \theta(t)] \text{ invertible}) \to 1$ as $t \to \infty$ and we obtain, for $Y\sim \mathcal{N}(0, \mathrm{I}_k)$,
 \[\sqrt{t}\widehat V_t[\widehat \theta(t)]^{-1/2} \bigl(\widehat \theta(t) - \vartheta\bigr) \mathbf{1}_{\{\widehat V_t[\widehat \theta(t)] \text{ invertible}\}} \xrightarrow{\PP_\vartheta\text{-}d} Y.\]

\begin{figure}[!b]
    \centering
    \includegraphics[width=12cm]{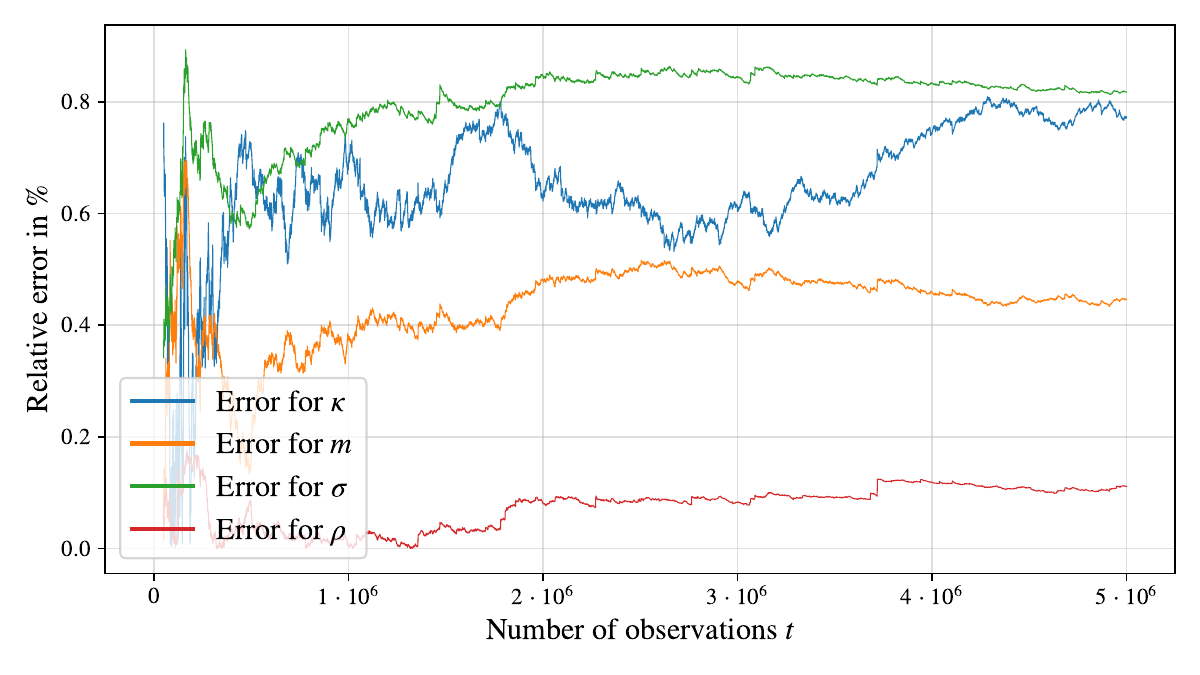}
    \caption{Two sequences of standard deviation estimates for the four parameters $\kappa$, $m$, $\sigma$, $\rho$ of the Heston model (\ref{eq: Heston_v}, \ref{eq: Heston_S}) with unobservable stochastic variance are obtained from the sequence of matrices $\widehat V_t[\widehat \theta(t)]$. For the first sequence of standard deviation estimates, the estimator $\widehat U_t[\widehat \theta(t)]$ from \eqref{eq: Formula_Cov} is computed using the knowledge of the unobservable variance component, while for the second sequence of estimates, the unobservable variance is replaced by its Kálmán filter during computation of $\widehat U_t[\widehat \theta(t)]$. The figure displays the relative error in percent between these two standard deviation estimates for all four parameters.}
    \label{fig: Error_Heston_Covariance}
\end{figure}

\begin{remark}
    The estimator \eqref{eq: Formula_Cov} still cannot be used to estimate $U_\vartheta$ in practice because it depends on the unobservable part of the data $X$. This can be remedied for example by replacing the unobservable components in $X(s)$ in \eqref{eq: Formula_Cov} by those of the filter $\widehat X^\theta(s, s)\big|_{\theta = \widehat\theta(t)}$. Figure \ref{fig: Error_Heston_Covariance} shows that --- at least in the example of the Heston model from Section \ref{sec5.2: Heston} --- the relative error made by this additional approximation is negligible and below one percent for the estimation of asymptotic standard deviations of the QML estimator components. 
    
    As an alternative, one could use Monte-Carlo samples of $X(s)$, $ s=0,1,2,\dots$ rather than real data. This allows to use an arbitrarily large sample size which exceeds the one of the given data set. But since the simulation will be based on the estimated parameter $\smash{\widehat\theta(t)}$ rather than the unknown $\vartheta$, a sufficiently large sample of data is needed in either case.
\end{remark}

 The easy computability of $\widehat U_t$ stems from the polynomial Markov structure which is not present in the estimation of VARMA models with arbitrary uncorrelated noise sequences as in \cite{Mainassara2011} or \cite{Schlemm2012}. The former suggest estimating $\smash{U_\vartheta = \lim_{t \to \infty} \mathrm{Cov}_\vartheta\bigl[ \frac{1}{\sqrt{t}} \sum_{s=1}^t Z^\vartheta(s, s-1)\bigr]}$ by approximating 
 \[\mathrm{Cov}_\vartheta\biggl[ \frac{1}{\sqrt{t}} \sum_{s=1}^t Z^\vartheta(s, s-1)\biggr] =\frac{1}{t}\sum_{s=1}^t\sum_{u=1}^t \mathrm{Cov}_\vartheta\big[Z^\vartheta(s, s-1), Z^\vartheta(u, u-1)\big]\] 
 using a non-parametric kernel estimator, also called heteroskedastic autocorrelation consistent (HAC) estimator or Newey--West estimator, see \cite{Andrews1991} or \cite{Newey1987}. Alternatively, extending a result of \cite{Franq2005}, \cite{Mainassara2011} prove consistency of an estimator of $U_\vartheta$ that builds on a parametric estimate of the spectral density of $Z^\vartheta(t, t-1)$ under the additional assumption of finite $(8+\delta)$th moments of $X(t)$. None of these more advanced techniques are needed in the computation of our covariance estimator $\smash{\widehat U_t[\widehat \theta(t)]}$ from \eqref{eq: Formula_Cov}, which only requires the evaluation of simple matrix manipulations, coming from the fact that the polynomial model assumptions permit the use of the Markov chain central limit result \ref{theo: MCLT} in the proof of Theorem \ref{theo: norm}.

\subsection{Explicit computation of the asymptotic covariance matrix}\label{su:explicit}
In contrast to Section \ref{sec5.1: CovMat}, we now provide explicit representations of $U_\vartheta, W(\vartheta)$ and hence $V_\vartheta$ in Theorem \ref{theo: main2}. To this end, we reconsider some of the processes from the previous section in state-space form. For brevity, set $\smash{\widehat\Sigma^\theta(t)\coloneq \widehat\Sigma^\theta(t, t-1)}$ and $\smash{S^{\theta}_j(t)\coloneq S^{\theta}_j(t,t-1)}$ as well as $\smash{K^\theta (t) \coloneq A^\theta \widetilde{K}^\theta(t)}$ and $\smash{F^\theta (t) \coloneq A^\theta - K^\theta (t)H}$. Observe that $\kX(0) = \E_\theta(X(0))$, $\wX(0)=\nabla_\theta \E_\theta(X(0))$, $\oX(0)=\nabla_\theta^2 \E_\theta(X(0))$, and

\begin{fitalign}
X(t + 1)&=a^\vartheta\hspace*{5.6pt}+\hspace*{18pt}A^\vartheta X(t)+N^\vartheta(t + 1),\\
\kX(t + 1)&=a^{\theta,0}+A^{\theta,0}(t)X(t)+A^{\theta,0,0}(t)\kX(t),\\
\wX_j(t + 1)&=a^{\theta,1}_j+A^{\theta,1}_j(t)X(t)+A^{\theta,1,0}_j(t)\kX(t)+A^{\theta,1,1}(t)\wX_j(t),\\
\oX_{ij}(t + 1)&=a^{\theta,2}_{ij}+A^{\theta,2}_{ij}(t)X(t)+A^{\theta,2,0}_{ij}(t)\kX(t)+A_{i}^{\theta,2,1}(t)\wX_j(t) +A_{j}^{\theta,2,1}(t)\wX_i(t) +A^{\theta,2,2}(t)\oX_{ij}(t)
\end{fitalign}

\noindent for $t \in \N$, where we define the shorthands $a^{\theta,0} \coloneq a^\theta$, $a^{\theta,1}_j \coloneq \partial_{j} a^{\theta}$, $a^{\theta,2}_{ij} \coloneq \partial_{ij} a^{\theta}$ as well as
\begin{align*}
A^{\theta,0}(t) &\coloneq K^\theta(t)H,\\
A^{\theta,1}_j(t) &\coloneq \bigl[\bigl(\partial_j A^\theta\bigr)\widetilde{K}^\theta(t) - N^\theta_j(t)\bigr]H,\\
A^{\theta,2}_{ij}(t) &\coloneq \bigl[\bigl(\partial_{ij}A^\theta\bigr)\widetilde{K}^\theta(t) -\bigl(\partial_{i}A^\theta\bigr)N^\theta_j(t) -(\partial_{j}A^\theta)N^\theta_i(t)+A^\theta M^\theta_{ij}(t)\bigr]H,\\
A^{\theta,0,0}(t) &\coloneq A^{\theta, 1, 1}(t) \coloneq A^{\theta, 2, 2}(t) \coloneq F^\theta(t),\\
A^{\theta,1,0}_j(t) &\coloneq \partial_{j}A^\theta - A^{\theta, 1}_j(t),\\
A^{\theta,2,0}_{ij}(t) &\coloneq \bigl(\partial_{ij}A^\theta\bigr) - A^{\theta, 2}_{ij}(t),\\
A^{\theta,2,1}_{j}(t) &\coloneq \bigl(\partial_{j}A^\theta\bigr) \bigl(\mathrm{I}_d-\widetilde K^\theta(t) H\bigr) + A^\theta N_j^\theta(t) H.
\end{align*}

We let $\widetilde K^\theta$, $K^\theta$, $F^\theta$, $\kappa^\theta_j$, $\mu^\theta_{ij}$,  $\psi_{ij}^\theta$, $\widetilde S^\theta_{\oo,j}$, $M^\theta_{ij}$ and $N^\theta_j$ denote the limits of $\widetilde K^\theta(t)$, $K^\theta(t)$, $F^\theta(t)$, $\kappa^\theta_j(t)$, $\mu^\theta_{ij}(t)$,  $\psi_{ij}^\theta(t)$$\smash{\widetilde S^\theta_{\oo,j}(t)}$, $M^\theta_{ij}(t)$ and $N^\theta_j(t)$, which exist by Lemma \ref{limitcov}. Then we can define time-homogeneous counterparts $\smash{\homkX}$, $\homX$, $\homoX$ of $\smash{\kX}$, $\wX$, $\oX$ with the same initial values. These are obtained from Propositions \ref{prop: log-lik}, \ref{prop: Quasi-Score} and \ref{prop: Fisher_Information} if we use the above limits instead of their time-dependent counterparts and $C^\theta$ \mbox{instead of $C^\theta(t)$, i.e.}

\begin{fitalign}
\homkX(t + 1)&\coloneq a^{\theta,0}+A^{\theta,0}X(t)+A^{\theta,0,0}\homkX(t),\\
\homX_j(t + 1)&\coloneq a^{\theta,1}_j+A^{\theta,1}_jX(t)+A^{\theta,1,0}_j\homkX(t)+A^{\theta,1,1}\homX_j(t),\\
\homoX_{ij}(t + 1)&\coloneq a^{\theta,2}_{ij}+A^{\theta,2}_{ij}X(t)+A^{\theta,2,0}_{ij}\homkX(t)+A_{i}^{\theta,2,1}\homX_j(t) +A_{j}^{\theta,2,1}\homX_i(t) +A^{\theta,2,2}\homoX_{ij}(t)
\end{fitalign}

\noindent for $t \in \N$, where we analogously define the shorthands \pagebreak
\begin{align*}
A^{\theta,0} &\coloneq K^\theta H,\\
A^{\theta,1}_j &\coloneq \bigl[\bigl(\partial_{j}A^\theta\bigr)\widetilde{K}^\theta - N^\theta_j\bigr]H,\\
A^{\theta,2}_{ij} &\coloneq \bigl[\bigl(\partial_{ij}A^\theta\bigr)\widetilde{K}^\theta -\bigl(\partial_{i}A^\theta\bigr)N^\theta_j -(\partial_{j}A^\theta)N^\theta_i+A^\theta M^\theta_{ij}\bigr]H,\\
A^{\theta,0,0} &\coloneq A^{\theta, 1, 1} \coloneq A^{\theta, 2, 2} \coloneq F^\theta,\\
A^{\theta,1,0}_j &\coloneq \partial_{j}A^\theta - A^{\theta, 1}_j,\\
A^{\theta,2,0}_{ij} &\coloneq \partial_{ij}A^\theta - A^{\theta, 2}_{ij},\\
A^{\theta,2,1}_{j} &\coloneq \bigl(\partial_{j}A^\theta\bigr) \bigl(I-\widetilde K^\theta H\bigr) + A^\theta N_j^\theta H.
\end{align*}
Note that $\homkX,\homX,\homoX$ differ from $\kX,\wX,\oX$ only in the sense that the time-depen\-dent coefficients $A^{\theta,1}(t) ,\dots, A^{\theta,2,2}(t)$ are replaced by their limits $A^{\theta,1}$, \dots, $A^{\theta,2,2}$. Observe also that $\langhomX = (X,\homkX,\homX)$ and $\laengerhomX \coloneq (X,\homkX,\homX,\homoX)$ are polynomial state space models of order 2 by Proposition 2.14 in \cite{KallsenRichert2025}.
Proposition \ref{ergod2} and \ref{coro: ergod_oX} show that they are ergodic with respect to unique stationary laws $\overline\mu_\theta$ resp.\ $\underline{\smash{\mu}}_\theta$.

Define $f^\theta:E\times\R^d\times\R^{d\times k}\to\R^{k}$ and $\widetilde f^\theta:E\times\R^d\times\R^{d\times k}\times\R^{d\times k\times k}\to\R^{k\times k}$ by
\begin{align}
f^\theta(\xi,\alpha,\beta)_j&\coloneq -\frac{1}{2}\left(\kappa_j+2\beta_{\oo,j}^\top(\widehat\Sigma^\theta_{\oo})^{-1}\varepsilon-\varepsilon^\top\widetilde{S}^\theta_{\oo,j}\varepsilon\right) \label{eq: f_def} \\
\widetilde f^\theta(\xi,\alpha,\beta,\gamma)_{ij}&\coloneq -\frac{1}{2}\left(\mu_{ij}+2\nu_{ij}\varepsilon+\varepsilon^\top\psi_{ij}\varepsilon+2\beta_{\oo,j}^\top(\widehat\Sigma^\theta_{\oo})^{-1}\beta_{\oo,i}\right), \label{eq: f_tilde_def}
\end{align}
where we define the shorthand $\varepsilon \coloneq \xi_\oo - \alpha_\oo$ and where $\nu \in \R^{d \times k \times k}$ is defined by $\nu_{ij} = \beta_{\oo, i}^\top \widetilde S_{\oo, j}^\theta + \beta_{\oo, j}^\top \widetilde S_{\oo, i}^\theta - \gamma_{\oo, ij}^\top (\widehat\Sigma^\theta_{\oo})^{-1}$. These functions stand for the time-homogeneous counterparts of the quadratic functions defining $Z^\theta(t, t-1)$ and $\nabla_\theta Z^\theta(t, t-1)$ in Propositions \ref{prop: Quasi-Score} and \ref{prop: Fisher_Information}. Moreover, we define the quadratic function $g^\theta:E\times\R^d\times\R^{d\times k}\to\R^{k}$ as the solution to the Poisson equation $f^\theta=\tildeP g^\theta-g^\theta$, where $\tildeP$ denotes the transition operator of $\langhomX = (X,\homkX,\homX)$ under $\PP_\vartheta$. Finally, set $h^\theta\coloneq\tildeP g^\theta,$ i.e. $h^\theta(x) = \tildeP g^\theta(x) = \int g^\theta(y)\tildeP(x,dy)$. These functions actually coincide with the ones from Proposition \ref{prop: Formula_Cov} and can be computed using the Tracy--Singh product. In particular, $g^\theta$ exists in the first place, \mbox{see the proof of Proposition \ref{prop: Formula_Cov}.}
\begin{theorem}\label{theo: main3}
Suppose that $\vartheta \in \mathrm{int}(\Theta)$. The matrices in Theorem \ref{theo: main2} have the form 
\begin{align*}
U_\vartheta&= \int \left(g^\vartheta(x)g^\vartheta(x)^\top
-h^\vartheta(x)h^\vartheta(x)^\top\right)\overline\mu_\vartheta(dx),\\
W(\vartheta)&= \int  \widetilde f^\vartheta(x) d\underline{\smash{\mu}}_\vartheta(dx).
\end{align*}
\end{theorem}
As a consequence, these matrices can be computed as follows if Assumption \ref{assump: AN} holds:

\begin{alg}\label{alg}
\emph{Step 1:}
As $X$ is a polynomial state space model of order 4, 
we have
\begin{equation}\label{e:altpssm}
\E_\theta\bigl(X(t)^\lambda\big|\F_{t-1}\bigr) = \sum_{\mu\in\N^d_4} b_{\lambda,\mu}X(t-1)^\mu,\quad t\in\N^*
\end{equation}
for any $\lambda\in\N^d_4$, with some  $b_{\lambda,\mu}\in\R$ such that
$b_{\lambda,\mu}=0$ if $|\mu|>|\lambda|$, see \cite{KallsenRichert2025}, Lemma 2.10.
If $X = (X(t))_{t\in \N}$ is obtained by sampling a $4$-polynomial process $(X(t))_{t \in \R_+}$ at non-negative integer time points, the $b_{\lambda,\mu}$ are obtained from the corresponding continuous-time coefficients by using the moment formula, see \cite{KallsenRichert2025}, Lemma 2.15. These continuous-time coefficients can in turn be calculated with equation (6.52) in \cite{Eberlein2019} from the local semimartingale characteristics of $X$. This is particularly easy for affine processes, see \cite{Eberlein2019}, Example 6.30.

\emph{Step 2:}
By \cite{KallsenRichert2025}, Proposition 2.14, \eqref{e:altpssm} holds also for $\langhomX = (X,\homkX,\homX)$ and $\laengerhomX = (X,\homkX,\homX,\homoX)$. This result also shows how to obtain the corresponding coefficients $b_{\lambda,\mu}$ for these processes from the ones for $X$.

\emph{Step 3:}
Since $\widetilde f^\vartheta$ is an explicitly given polynomial of order 2, its expectation $W(\vartheta)$ under the stationary law of $\smash{\underline{X}^{\vartheta, \mathrm{hom}} = (X,\widehat{X}^{\vartheta, \mathrm{hom}},V^{\vartheta, \mathrm{hom}},W^{\vartheta, \mathrm{hom}})}$ can now be computed explicitly with the moment formula, see \cite{KallsenRichert2025}, Proposition 2.11 and equation (2.5) therein.

\emph{Step 4:}
The coefficients of the quadratic polynomials $g^\vartheta$ and $h^\vartheta$ are provided explicitly in Proposition \ref{prop: Formula_Cov}. The multinomial theorem then yields the coefficients of the quartic polynomial $x\mapsto g^\vartheta(x)g^\vartheta(x)^\top-h^\vartheta(x)h^\vartheta(x)^\top$ for the computation of $U_\vartheta$, see Theorem \ref{theo: main3}.

\emph{Step 5:}
The expectation $U_\vartheta$ of this quartic polynomial under the unique stationary law of $\smash{\smash{\overline{X}}^{\vartheta, \mathrm{hom}} = (X,\widehat{X}^{\vartheta, \mathrm{hom}},V^{\vartheta, \mathrm{hom}})}$ can now also be obtained explicitly from the moment formula, see \cite{KallsenRichert2025}, Proposition 2.11 and equation (2.5) therein.
\end{alg}

Since the true parameter $\vartheta$ is unknown, we approximate it by its consistent estimator $\widehat\theta(t)$ to apply Algorithm \ref{alg}. A Slutsky-type argument yields the following central limit theorem.
\begin{proposition}\label{p:slutsky}
Let $\vartheta \in \mathrm{int}(\Theta)$ and Assumption \ref{assump: AN} hold. Suppose that $U_\vartheta$ and $W(\vartheta)$ are invertible. Then $V_{\widehat\theta(t)}\xrightarrow{\PP_\vartheta} V_\vartheta$ and
\[\sqrt{t}V_{\widehat \theta(t)}^{-\frac{1}{2}} \bigl(\widehat \theta(t) - \vartheta\bigr) \xrightarrow{\PP_\vartheta\text{-}d} Y\]
as $t\to\infty$ with $Y\sim \mathcal{N}(0, \mathrm{I}_k)$. Here, we set $V_{\widehat \theta(t)}^{-\frac{1}{2}}\coloneq 0 \in \R^{k \times k}$ if $V_{\widehat \theta(t)}$ is not invertible.
\end{proposition}

\subsection{Asymptotic confidence intervals and tests}\label{su:tests}
Once the covariance matrix $V_\vartheta$ has been consistently estimated by $\smash{\widehat V(t)\coloneq V_{\widehat \theta(t)}}$ from Section \ref{su:explicit} or $\smash{\widehat V(t)\coloneq \widehat V_t[\widehat \theta(t)]}$ from Section \ref{sec5.1: CovMat}, it can be used to calculate asymptotic confidence intervals for the components of $\vartheta$. In particular, for $j \in \{1, \dots, k\}$, the symmetric interval
\[\Big[\widehat \theta(t)_j - \sqrt{\textstyle \frac{\widehat V(t)_{jj}}{t}} \Phi^{-1}\big(1 - \textstyle \frac{\alpha}{2}\big), \;\widehat \theta(t)_j + \sqrt{\frac{\widehat V(t)_{jj}}{t}} \Phi^{-1}(1 - \frac{\alpha}{2})\Big]\]
covers the true parameter value $\vartheta_j$ with a probability converging to $1 - \alpha$, where $\Phi^{-1}(\alpha)$ denotes the standard normal $\alpha$-quantile. One-sided intervals and multivariate elliptic confidence regions can be constructed in a similar manner. Alternatively, it might be of interest to simultaneously test $m \leq k$ non-linear constraints on $\vartheta$ in the form of a hypothesis
\[H_0: R(\vartheta) = r,\]
where $r \in \R^m$ and $R: \R^k \to \R^m$ is continuously differentiable with Jacobian $\nabla_\theta R: \R^k \to \R^{m \times k}$ such that $\mathrm{rank}(\nabla_\theta R(\theta)) = m$ for all $\theta \in \Theta$. Most of the testing principles used in applied mathematics to test such hypotheses are based on the Wald test, the Lagrange multiplier (LM) test, also called score test or Rao-score test, and the likelihood-ratio (LR) test, also called Wilks test, see \cite{Engle1984}. The latter two require the use of a constrained estimator which we conveniently define to be any $\F_t $-measurable $\Theta$-valued random variable $\smash{\widehat \theta^c(t)}$ such that $\smash{\widehat \theta^c(t) \in \mathrm{arg}\max_{\theta \in \Theta_c} L^\theta(t)}$, where $\Theta_c \coloneq \{\theta \in \Theta: \: R(\theta) = r\}$.\footnote{Since $R$ is a continuous function and since $\Theta$ is compact, the constrained set $\Theta_c$ is compact as well, so the existence of measurable constrained QML estimators once more follows from \cite{Jennrich1969}, Lemma 2.} We further make the assumption that, under $H_0$, the sequence $\smash{\big(\widehat{\theta}^c(t)\big)_{t \in \N }}$ is weakly $\vartheta$-consistent and moreover that $\smash{\big(\sqrt{t}\bigl(\widehat \theta(t) - \vartheta\bigr), \sqrt{t}(\widehat \theta^c(t) - \vartheta)\big)_{t \in \N }}$ is jointly asymptotically normal. In this case, since $\vartheta \in \mathrm{int}(\Theta)$, the Lagrange multiplier theorem (see e.g.\ \cite{Fitzpatrick2006}, Theorem 17.17) and consistency imply that, with probability converging to 1, $\widehat \theta^c(t)$ solves the Lagrangian equations $\smash{\frac{1}{t}Z^\theta(t)} = \smash{\nabla_\theta R(\theta)^\top \widehat \lambda(t)}$ and $R(\theta) = r$ together with some $\R^m$-valued Lagrange multiplier variable $\smash{\widehat \lambda(t)}$.

\begin{proposition}\label{prop: significance_tests}
    Assume that $\vartheta \in \mathrm{int}(\Theta)$ and that both $U_\vartheta$ and $W(\vartheta)$ are invertible. Then the following conclusion holds under $H_0$, where the occurring statistic $\xi_\mathrm{Wald}(t)$ is defined to take the value 0 on $\{\widehat V(t) \text{ is singular}\}$, the probability of which converges to 0.
    \begin{itemize}
        \item[1.] The Wald statistic $\xi_{\mathrm{Wald}}(t) \coloneq t \big[R(\theta) - r\big]^\top \big[\nabla_\theta R(\theta) \widehat V(t) \nabla_\theta R(\theta)^\top\big]^{-1} \big[R(\theta) - r\big]\big|_{\theta = \widehat \theta(t)}$
        satisfies $\smash{\xi_{\mathrm{Wald}}(t) \xrightarrow{\PP_\vartheta\text{-}d} \Xi}$, where $\Xi$ follows a $\chi^2$-distribution with $m$ degrees of freedom.
    \end{itemize}
    Set $\widehat V^c(t) \coloneq \widehat V_t[\widehat \theta^c(t)]$, $\widehat W^c(t) \coloneq \widehat W_t[\widehat \theta^c(t)]$ or $\widehat V^c(t) \coloneq V_{\widehat \theta^c(t)}$, $\widehat W^c(t) \coloneq W_{\widehat \theta^c(t)}$, depending on whether the approach in Section \ref{sec5.1: CovMat} or Section \ref{su:explicit} is chosen. Suppose further that $R \in \mathrm{C}^2(\Theta, \R^m)$ and that the objects in Assumption \ref{assump: A} are $\mathrm{C^4}$-functions. Let $\Sigma_{\widehat V} \coloneq \nabla_\theta R(\widehat \theta^c(t)) \widehat V^c(t) \nabla_\theta R(\widehat \theta^c(t))^\top$, $\Sigma_{V} \coloneq \nabla_\theta R(\vartheta) V_\vartheta \nabla_\theta R(\vartheta)^\top$, and further define $\Sigma_{W} \coloneq \nabla_\theta R(\vartheta) W(\vartheta)^{-1} \nabla_\theta R(\vartheta)^\top$. Then the following two statements hold under $H_0$, where we set $\xi_\mathrm{LM}(t) \coloneq 0$ on the set $\{\widehat V^c(t) \text{ is singular}\}$:
    \begin{itemize}
        \item[2.] The LM statistic $\xi_\mathrm{LM}(t) \coloneq \frac{1}{t} Z^\theta(t) \widehat W^c(t)^{-1} \nabla_\theta R(\theta) \Sigma_{\widehat V}^{-1} \nabla_\theta R(\theta)^\top \widehat W^c(t)^{-1}Z^\theta(t) \big|_{\theta = \widehat \theta^c(t)}$ satisfies $\smash{\xi_\mathrm{LM}(t) \xrightarrow{\PP_\vartheta\text{-}d} \Xi}$, where $\Xi$ follows a $\chi^2$-distribution with $m$ degrees of freedom.
        
        \item[3.] Similarly, the LR statistic $\xi_\mathrm{LR}(t) \coloneq -2 \log(\varphi_t)$ satisfies $\xi_\mathrm{LR}(t) \xrightarrow{\PP_\vartheta\text{-}d} \sum_{j=1}^m \lambda_j Z_j^2$, where $\varphi_t \coloneq \smash{\exp\big(L^{\widehat \theta^c(t)} - L^{\widehat \theta(t)}\big)}$ is the quasi-likelihood ratio, $\smash{Z_j \stackrel{\mathrm{iid}}{\sim} \mathcal{N}(0, 1)}$ and $\lambda_j$ are the eigenvalues of $[-W(\vartheta)]^{-1/2} \nabla_\theta R(\vartheta)^\top \Sigma_W^{-1} \Sigma_V \Sigma_W^{-1} \nabla_\theta R(\vartheta) [-W(\vartheta)]^{-1/2}$.
    \end{itemize}
\end{proposition}

\begin{remark}
    Above we assumed that the sequence of constrained QML estimators $\widehat \theta^c(t)$ is consistent and, jointly with $\smash{\widehat \theta(t)}$, asymptotically normal under $H_0$. This can be naturally motivated by the following construction, which assumes that the hypothesis can be parameterised by some lower-dimensional space: Suppose that there exists a convex and compact set $\widetilde \Theta_c \in \R^{k'}$ with non-empty interior, where $0 \leq k' \leq k$, and some bijection $h \in \mathrm{C}^3(\widetilde \Theta_c, \Theta)$ such that $\mathrm{rank}(\nabla_\eta h(\eta)) = k'$ for any $\eta \in \widetilde \Theta_c$. Then it is easy to check that all conditions of the main theorems \ref{theo: main} and \ref{theo: main2} are fulfilled for the estimation model with lower-dimensional parameter space $\widetilde \Theta_c$ and so, under $H_0$, $\smash{\widehat \eta(t) \in \mathrm{arg}\max_{\eta \in \widetilde \Theta_c} L^{h(\eta)}(t)}$ is consistent and asymptotically normal for the unique value $\eta^*$ such that $h(\eta^*) = \vartheta$. Consequently, also $\smash{\widehat \theta^c(t) = h(\widehat \eta(t))}$ is $\vartheta$-consistent and asymptotically normal, and the joint asymptotic normality assumption follows by considering estimation of $(\vartheta, \eta^*)$ on the parameter space $\Theta \times \widetilde \Theta_c$.
\end{remark}

 Compared to the likelihood-ratio test, the Lagrange multiplier test has the advantage that only the constrained estimator $\smash{\widehat \theta^c(t)}$ has to be evaluated instead of both the constrained and unconstrained estimator. Additionally and in contrast to the well-known simple results obtained for maximum likelihood estimation using independent and identically distributed random variables, the asymptotic distribution of the likelihood-ratio test is more complicated than that of the Wald and Lagrange multiplier test. It can be computed using various numerical algorithms, see for example \cite{Davies1980, Farebrother1984, Imhof1961, Liu2009}.

\section{Examples}\label{sec5: app}
In this section we illustrate the QML approach by applying it to two affine models from finance. In Section \ref{sec5.2: Heston} we study the popular Heston model, where the variance process is unobservable. In Section \ref{sec5.3: OU} we consider a Lévy-driven two-factor short rate interest model by viewing it as an instance of a general multivariate Lévy-driven Ornstein-Uhlenbeck process.

\subsection{The Heston stochastic volatility model}\label{sec5.2: Heston}
In the stochastic volatility model of \cite{Heston1993}, the asset $S = (S(t))_{t \in \R_+}$ incorporates a stochastic variance component $v = (v(t))_{t \in \R_+}$ driven by a square-root diffusion. More specifically, it is assumed that $(v, S)$ follows the dynamics
\begin{align}
    \mathrm{d}v(t) &= \kappa(m - v(t)) \mathrm{d}t + \sigma \sqrt{v(t)} \mathrm{d}W^{(2)}(t), \label{eq: Heston_v} \\
        \mathrm{d}S(t) &= \big(\mu + \delta v(t)\big) S(t) \mathrm{d}t + \sqrt{v(t)} S(t) \mathrm{d}W^{(1)}(t) \label{eq: Heston_S}
\end{align}
under the physical probability measure $\PP_\theta$, where $W^{(1)}$ and $W^{(2)}$ are two Brownian motions with correlation $\rho \in [-1, 1]$. This way, the variance process exhibits a mean reversion behaviour at the rate $\kappa \in (0,\infty)$ towards the long run variance $m \in (0,\infty)$. The parameter $\sigma \in (0,\infty)$ acts as a so-called volatility of volatility, while $\mu \in \R$ describes the drift rate of the asset $S$ and $\delta \in \R$ is a parameter controlling the volatility response to the drift. Existence and uniqueness of a non-negative strong solution to the stochastic differential equation \eqref{eq: Heston_v} for any $v(0) \geq 0$ follows for example as in \cite{Ikeda1989}, Example IV.8.2. Since $S$ can be expressed as the exponential of a stochastic integral of $v$, \eqref{eq: Heston_S} also admits a unique strong solution. Instead of the equation \eqref{eq: Heston_S} for the asset $S$ itself, we can also consider the log-spot process $Y=(Y(t))_{t \in \R_+}$ with $Y(t) \coloneq \log S(t)$ for $t \in \R_+$, described by the equation
\begin{align}\label{eq: Heston_Y}
    \mathrm{d}Y(t) = \Big( \mu + \Big[\delta - \frac{1}{2}\Big]v(t)\Big) \mathrm{d}t + \sqrt{v(t)} \mathrm{d}W^{(1)}(t).
\end{align}
Then $(v, Y)$ is a bivariate affine process described by the local semimartingale characteristics 
\[(b^{\mathrm{id}}(t), c(t), K(t)) = \bigg[ \begin{pmatrix}
    \kappa m \\ \mu
\end{pmatrix} + \begin{pmatrix}
    -\kappa \\ \delta - \frac{1}{2}
\end{pmatrix} v(t) ,\; \begin{pmatrix}
    \sigma^2 & \sigma \rho \\ \sigma\rho & 1
\end{pmatrix}v(t), \; 0\bigg].\]
It is not hard to see that the corresponding generalised Riccati equations are explicitly solvable in the case of the Heston model $(v, Y)$, see \cite{Eberlein2019}, Section 8.2.4. 

Since both the joint likelihood of $Y(t)$ and $v(t)$ as well as the marginal likelihood of $Y(t)$ are unknown in explicit form in the Heston model, estimation of the model parameters $\mu$, $\delta$, $\kappa$, $m$, $\sigma$, $\rho$ is a non-trivial task. Since the Heston model belongs to the class of affine stochastic volatility models, its characteristic function can be derived from the affine model structure, and the corresponding density can be obtained by classic Fourier inversion techniques from the characteristic function (see e.g.\ \cite{Bates2006}), which is however a time-consuming numerical task. Alternatively, estimation in the Heston model can solely be based on empirical estimates of the characteristic function as for example in \cite{Jiang2002} or \cite{Singleton2001}, or Fourier inversion-based methods can be applied to calibrate the model parameters directly to observed option prices. The latter approach however determines the model parameters under an equivalent martingale measure, a so-called risk-neutral measure, see for example \cite{Eberlein2019}, Section 11.2.3, while estimation routines based on observed returns like maximum likelihood or least squares minimisation procedures fit the model parameters under the true physical measure driving the returns.

Estimation in the latter case essentially boils down to a filtering problem and often consists in computing an approximate density of discretely observable quantities. Due to their simplicity, Gaussian density approximations enjoy particular popularity in the field of estimating multivariate diffusions or affine models like the Heston model, see for example \cite{Duffie2002}, \cite{Fisher1996}, \cite{Hurn2013} or \cite{Wang2018}. Alternatively, method of moments estimators for the latent volatility process can be employed as in \cite{Bollerslev2002}. Of all comparative studies found in the literature, the approach followed by \cite{Cacace2019}, which includes a polynomial filtering method similar to ours, is most comparable to our estimation framework for the Heston model.

We now aim to express the Heston model as a polynomial state space model in the sense of Definition \ref{def: polynomial_ssm}. For computational simplicity, we consider fixed drift parameter $\mu = 0$ and volatility response parameter $\delta = \frac{1}{2}$ so that the log-spot process $Y$ from \eqref{eq: Heston_Y} has the form
\begin{align}\label{eq: Heston_Ytilde}
        \mathrm{d} Y(t) = \sqrt{v(t)} \mathrm{d}W^{(1)}(t),
\end{align}
which is a martingale. Then $\Delta Y = \big(Y(t) - Y(t - 1)\big)_{t \in \N}$ is the discretely observed log-returns process in this Heston model.\footnote{We assume $Y(-1)$ to be given for this purpose.} Since $v$ only affects the returns $\Delta Y$ through their conditional variance, it makes sense to base the parameter estimation in the Heston model not only on the observed returns $\Delta Y$, but also on the squared returns $(\Delta Y)^2 \coloneq (\Delta Y(t)^2)_{t \in \N}$. Then $X \coloneq (v, \Delta Y, (\Delta Y)^2)$ with $v = (v(t))_{t \in \N}$ and $\theta = (\kappa, m, \sigma, \rho)$ is a polynomial state space model with proper state space $E = \R_+ \times \{(x, y) \in \R^2: y = x^2\}$. The compact parameter space $\Theta$ is chosen as $\Theta = I_\kappa \times I_m \times I_\sigma \times I_\rho \coloneq [10^{-4}, 10] \times [10^{-8}, 1] \times [10^{-4}, 1] \times [-1, 1]$. For the following simulation studies, we assume $X(0)$ to be externally known, which corresponds to a Dirac distribution as the initial law for $X(0)$. Of course, our setup from Section \ref{sec4: Estimation} also admits other arbitrary choices of initial distributions.

Since the distribution of $v(s + t)$ given $v(s)$ with $s, t \in \R_+$ is that of $c_t$ times a non-central chi-square distribution with $\frac{4\kappa m}{\sigma^2}$ degrees of freedom and non-centrality parameter $\smash{\frac{v(s)}{c_t} \mathrm{e}^{-\kappa t}}$, where $\smash{c_t = \frac{1}{4 \kappa} \sigma^2 (1 - \mathrm{e}^{-\kappa t})}$ (see \cite{Alfonsi2010}), the state transition vector and matrix of $X$ are
\begin{align*}
    a^\theta = m \begin{pmatrix}
        1 - \e^{-\kappa} \\ 0 \\ 1 - \frac{1}{\kappa}(1 - \e^{-\kappa})
    \end{pmatrix} \qquad \qquad A^\theta = \begin{pmatrix}
        \e^{-\kappa} & 0 & 0 \\
        0 & 0 & 0 \\
        \frac{1}{\kappa}(1 - \e^{-\kappa}) & 0 & 0
    \end{pmatrix}.
\end{align*}

We now turn to a brief justification of the assumptions introduced in Section \ref{sec4.1: StateSpaceModels} in the case of the Heston model. Since the affine characteristics of $(X, Y)$ are smooth functions of $\theta$, the same holds for the state transition vector and matrix $a^\theta_{\otimes r}$ and $A^\theta_{\otimes r}$ by the moment formula \cite{Eberlein2019}, Theorem 6.26, for any $r \in \N$. Since $\Delta Y(1)$ possesses a positive density with respect to Lebesgue measure given\footnote{See for example the considerations in \cite{DelBanoRollin2010} or use methods from Malliavin calculus, for example Corollary 3.3 in \cite{Nourdin2009} together with the expressions from Section 3.3 in \cite{Alos2021}.} $\Delta Y(0)$ and since the distribution of $v(1)$ given $v(0)$ under any $\PP_\theta$ is explicitly known and equivalent to Lebesgue measure on $\R_+$ if $\kappa$, $m$, $\sigma$ are positive, $X(1)$ has a density with respect to two-dimensional Hausdorff measure on $E$ given $X(0)$. Moreover, since $(v(t), Y(t))_{t \in \R_+}$ is affine, it possesses the Feller property (see \cite{Duffie2003}, Theorem 2.7), which then also holds for $X$. Hence, Assumptions \ref{assump: AN} and \ref{normassump}(\ref{itm: Irreducibility}) are fulfilled provided $\E_\theta(v(0))$ and $\mathrm{Var}_\theta(v(0))$ are $\mathrm{C}^3$-functions of $\theta$. Concerning Assumption \ref{normassump}(\ref{itm: L6-bounded}), note that from the knowledge of the distributions of $v(t)$ it follows that $v$ has bounded moments of all orders, and so, by the Burkholder--Davis--Gundy inequality (see e.g. \cite{Revuz1999}, Corollary 4.2, \mbox{or Theorem \ref{theo: BDG}), there is $C_p > 0$ such that}
\begin{align*}
\E_\theta(|\Delta Y(t)|^p) &= \E_\theta \Big( \Big\lvert \int_{t-1}^{t} \sqrt{v(s)} \dd W^{(1)}(s)\Big\rvert^p \Big) \leq C_p \E_\theta \Big( \Big\lvert \int_{t-1}^{t} v(s) \dd s \Big\rvert^{\frac{p}{2}} \Big) \\
&\leq C_p \int_{t-1}^{t} \E_\theta\big(v(s)^{\frac{p}{2}}\big) \dd s
\end{align*}
for any $p \geq 2$. Hence $X$ is bounded in any $L^p(\PP_\theta)$ and Assumption \ref{normassump}(\ref{itm: L6-bounded}) is fulfilled. Likewise, Assumption \ref{assump: Identifiability2} holds true for the Heston model because the asymptotic covariance matrix $\Lambda^\theta_{\oo}$ of $X_\oo$ under $\PP_\theta$ is identifiable, where $\Lambda^\theta$ solves $\Lambda^\theta = A^\theta \Lambda^\theta A^{\theta^\top} + C^\theta$. Explicit formulas for the matrix $C^\theta$ in the Heston model can be found in Example 4.8 of the companion paper \cite{KallsenRichert2025}.

\begin{remark}
    Instead of considering squared returns $(\Delta Y)^2$ in the specification of the Heston model as a polynomial state space model, one could also consider the quantities $\Delta Y^2$, i.e.\ differences of the squared log-spot, or both. Additionally, higher order versions of both quantities could be added as additional components to $X$, incorporating higher powers of returns in the estimation framework. We did not test the effect of these model modifications concerning the efficiency of the resulting quasi-maximum likelihood estimator, say by examining the resulting differences in size of the asymptotic variances of the estimator components. However, the simultaneous use of first and second powers $\Delta Y$ and $(\Delta Y)^2$ is crucial for identifiability of the parameters $\sigma$ and $\rho$: if $X = (v, \Delta Y)$, then neither $\sigma$ nor $\rho$ influence the asymptotic distribution of the Gaussian equivalent, which are hence non-identifiable and cannot be estimated. If, on the other hand, estimation is based only on the squared returns, i.e.\ $X = (v, (\Delta Y)^2)$, then $\sigma$ is estimable but only $\rho^2$ is identifiable, i.e.\ the strength of correlation between the asset and the volatility becomes estimable but not its sign. If the parameter $\delta$ is set to a different value than $\frac{1}{2}$, then $(v, \Delta Y, (\Delta Y)^2)$ ceases to be a polynomial state space model and $(v, v^2, \Delta Y, (\Delta Y)^2)$ needs to be used instead. In these state space specifications, all components containing $v$ are assumed to be unobservable.
\end{remark}

\begin{figure}[!b]
    \centering
    \includegraphics[width=15cm]{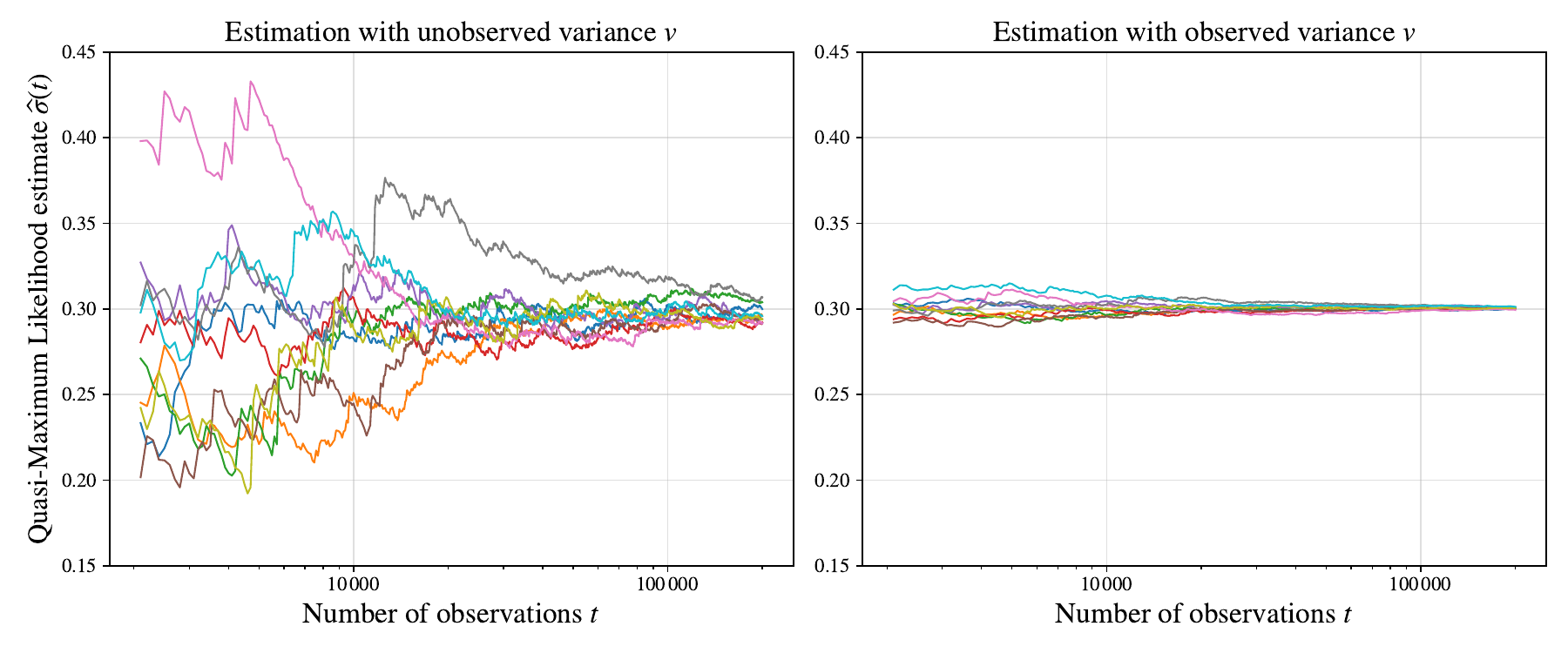}
    \caption{Ten sequences of quasi-maximum likelihood estimators for the parameter $\sigma^* = 0.3$. On the left $v$ is assumed to be unobservable, while on the right $v$ is assumed to be observable.}
    \label{fig: QML_Sequences}
\end{figure}

\begin{figure}[!b]
    \centering
    \includegraphics[width=13.cm]{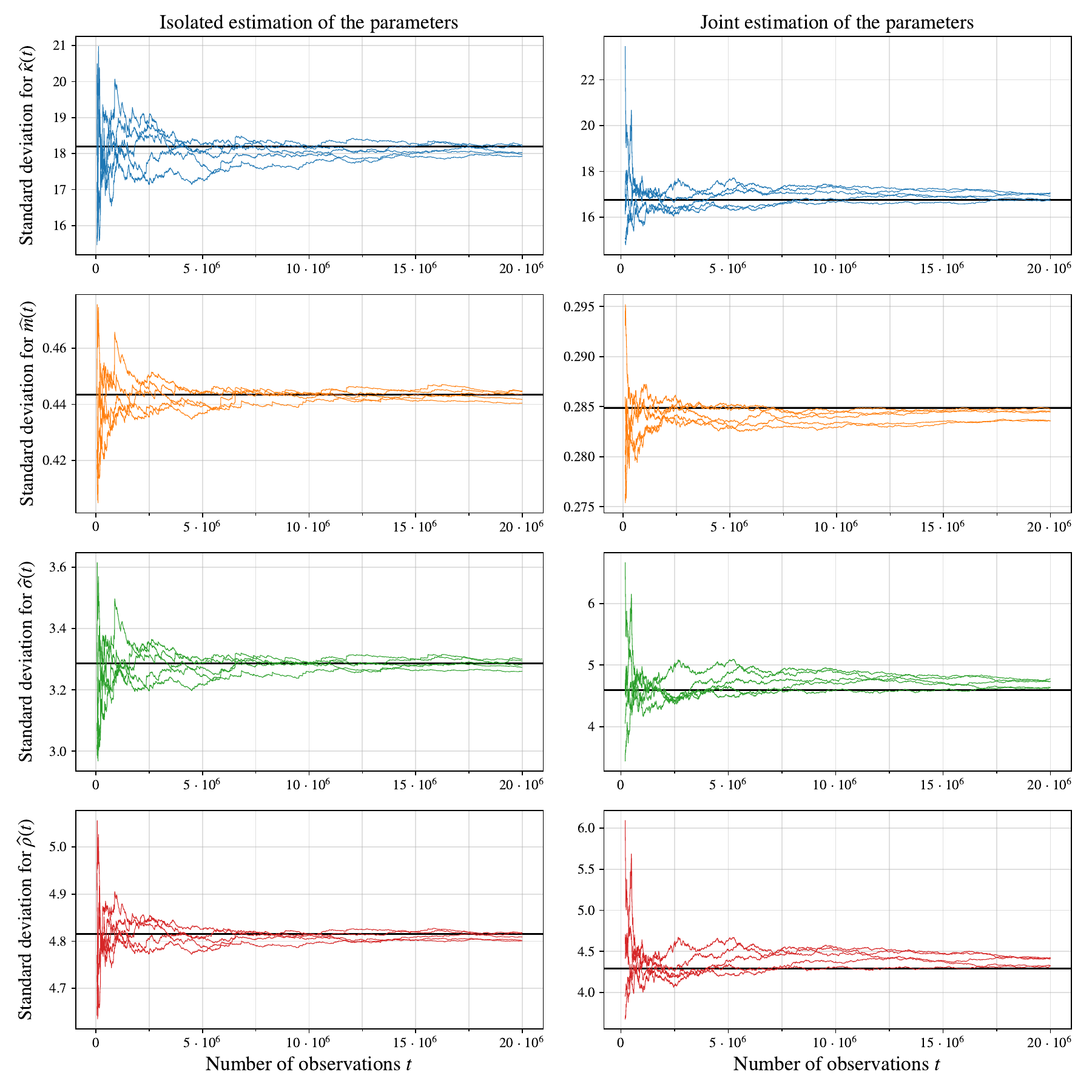}
    \caption{Five sequences of estimator standard deviations, obtained from the covariance estimates $\smash{\widehat V_t[\widehat \theta(t)]}$ for $t$ between 1 and $20 \cdot 10^6$. The left column displays the standard deviation in the case of an isolated estimation of the parameters $\kappa^*$, $m^*$, $\sigma^*$ and $\rho^*$, assuming all other parameters to be known. The right column displays the standard deviation in the case of a joint estimation of the parameter $\vartheta = (\kappa^*, m^*, \sigma^*, \rho^*)$. Black lines show values obtained from the explicit calculations detailed in Section \ref{su:explicit}.}
    \label{fig: Covariance_Heston_firstobserved1}
\end{figure}

Since the Assumptions \ref{assump: AN}, \ref{normassump} and \ref{assump: Identifiability2} are fulfilled in the case of the Heston model, Theorems \ref{theo: main} and \ref{theo: main2} yield consistency of any quasi-maximum likelihood estimator sequence and also asymptotic normality, provided that the asymptotic Fisher information matrix $W(\vartheta)$ is invertible. This condition can be justified easily using the explicit calculations from Algorithm \ref{alg}. To visualise the content of Theorem \ref{theo: main}, we let $\vartheta = (\kappa^*, m^*, \sigma^*, \rho^*) \coloneq (1,\, 0.4^2,\, 0.3,\, -0.5)$ and $X(0) = (v(0),\Delta Y(0),\Delta Y(0)^2) = (0.3^2,\, 0,\, 0)$ and simulate $N = 10\, 000$ independent trajectories $(X^\vartheta(t))_{t \in \{0, 1, \dots, T\}}$ of size $T = 200\, 000$ using the programming language Python. The employed discretisation scheme for the equations \eqref{eq: Heston_v} and \eqref{eq: Heston_Ytilde} is the standard Euler--Maruyama method with step size $\Delta t = \frac{1}{250}$, with the addition that at each step of the discretisation the function $x^+ = \max\{x, 0\}$ is used to ensure non-negativity of the variance (this is termed the absorption scheme in \cite{Lord2010}, who compare different simulation schemes for the Heston model). We first focus solely on the estimation of the parameter $\sigma$, assuming all other parameters to be known. Figure \ref{fig: QML_Sequences} shows ten independent sequences $(\widehat \sigma(t))_{t \in \{1, \dots, T\}}$ of the resulting quasi-maximum likelihood estimator, first if the variance component $v$ is assumed to be unobservable and secondly if the variance component $v$ can be observed.

Here it is visible that the absence of additional observable components significantly increases the estimator standard deviation by approximately a factor of 10. We suffer from another increase in variance if the whole parameter $\vartheta$ rather than only $\sigma$ is unknown. The tiles of Figure \ref{fig: Covariance_Heston_firstobserved1} each contain five standard deviation sequences for the estimator components, obtained from $\widehat V_t[\widehat \theta(t)]$ introduced in Section \ref{sec5.1: CovMat}, first for an isolated estimation of the parameter components and secondly for the joint estimation. The horizontal black lines show the explicit computation of the standard deviation components using Algorithm \ref{alg} from Section \ref{su:explicit}. In the case of the joint estimation, the asymptotic standard deviation for $\widehat \sigma(t)$ is larger by approximately a factor of 1.5 compared to the isolated estimation. Interestingly however, a reversed effect can be observed for the other parameters $\kappa$, $m$, $\rho$, the estimation of which becomes more accurate in the joint estimation case. In particular, the estimation accuracy achieved by the quasi-maximum likelihood estimators seems to be comparable to the accuracy reported in \cite{Cacace2019}, Table 5.

\begin{figure}[!b]
    \centering
    \includegraphics[width=11.9cm]{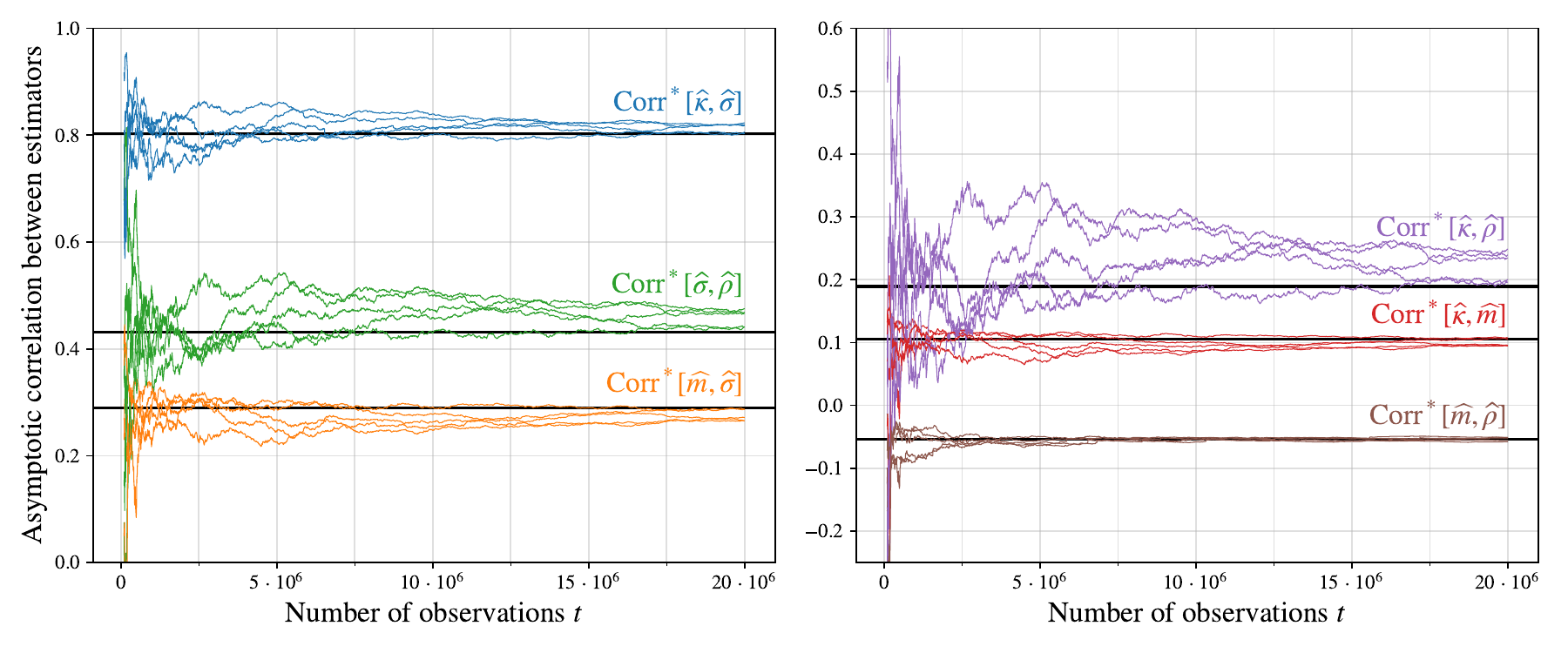}
    \caption{For each pair of estimator components, the figure displays five independent sequences of estimator correlations, respectively obtained from the covariance estimates $\smash{\widehat V_t[\widehat \theta(t)]}$ for $t$ varying between 1 and $20 \cdot 10^6$, and a little less rigorously called $\mathrm{Corr}^*[\cdot, \, \cdot]$ to denote the asymptotic correlation between the estimator components in the square brackets. Black lines show values obtained from the explicit calculations detailed in Section \ref{su:explicit}.}
    \label{fig: Correlation_Heston_firstobserved1}
\end{figure}

For the case of joint estimation, Figure \ref{fig: Correlation_Heston_firstobserved1} contains the sequences of correlations between the estimator components, obtained from $\smash{\widehat V_t[\widehat \theta(t)]}$, again accompanied by the explicit calculations from Section \ref{su:explicit} using black horizontal lines. Another visible fact from Figures \ref{fig: Covariance_Heston_firstobserved1} or \ref{fig: Correlation_Heston_firstobserved1} is the potentially low speed of convergence of the asymptotic covariance estimator. In Figure \ref{fig: Covariance_Heston_firstobserved1} about $t = 10^7$ observations are necessary for the covariance estimator to reach a stable state with acceptable accuracy compared to the explicit calculations. This is a potential drawback of the easily computable estimators $\widehat V_t[\theta]$ that we introduced in Section \ref{sec5.1: CovMat}, which have a rather high variance due to the occurrence of empirical fourth moments of $X$ in the expression of the estimator $\widehat U_t[\theta]$. In practice, it is rather unlikely that such a large history of data is available for estimation. In this case, once the parameter $\vartheta$ of the model has been consistently estimated from data by $\widehat \theta(t)$, a bootstrap procedure can be used. To this end, one simulates a large history of synthetic data from the presumed model using the fitted parameter and uses the simulated history to estimate the asymptotic covariance matrix of the quasi-maximum likelihood estimator as in Figure \ref{fig: Covariance_Heston_firstobserved1}. This discussion once more emphasises the advantages of our explicit calculations from Algorithm \ref{alg} in Section \ref{su:explicit}.

\begin{figure}[!t]
    \centering
    \includegraphics[width=12cm]{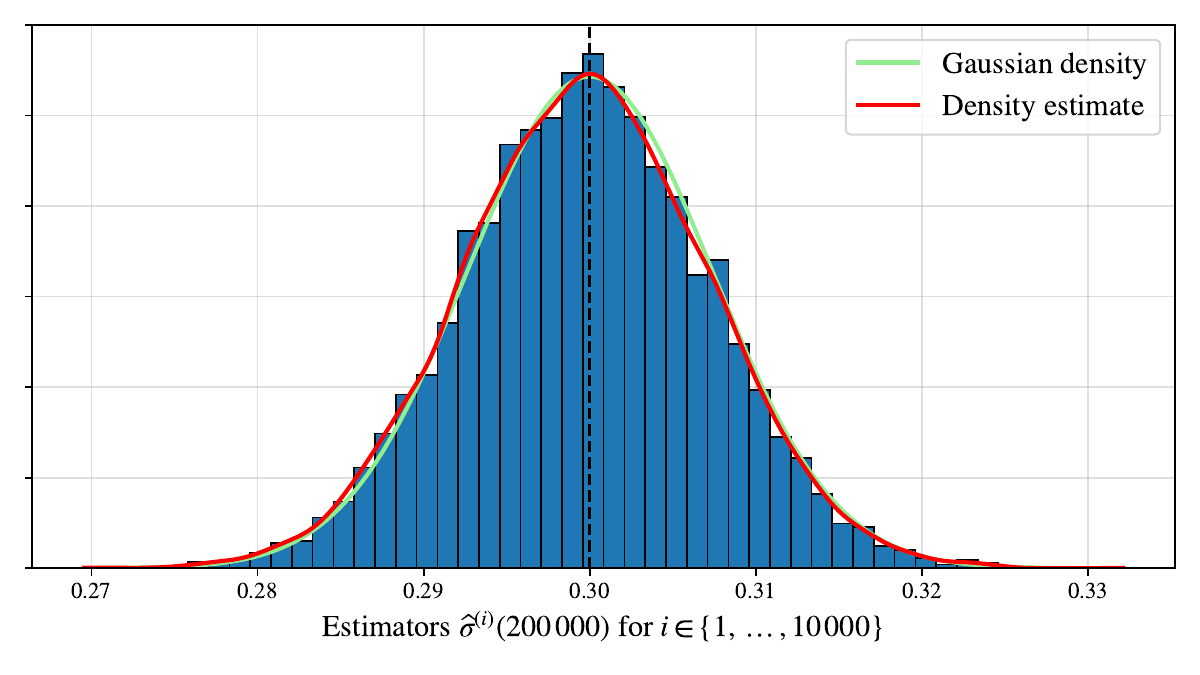}
    \caption{Histogram of $N = 10\,000$ independent quasi-maximum likelihood estimators $\widehat \sigma^{(i)}(200\,000)$. On top we draw the density of the normal distribution with mean $\sigma^* = 0.3$ and variance given by the explicit calculation of $V_\vartheta$ in green as well as a kernel density estimate of the estimators $\widehat \sigma^{(i)}(T)$ in red. The black dotted line indicates the true value $\sigma^* = 0.3$.}
    \label{fig: QMLDistribution_Heston_firstobserved1}
\end{figure}

In the case of an isolated estimation of the parameter $\sigma$, the explicit calculation of the asymptotic covariance matrix reveals an asymptotic estimator standard deviation of approximately $\mathrm{Std}^*[\widehat \sigma] \colonapprox 3.2871$ for the asymptotic standard deviation of the estimator $\widehat \sigma(t)$. We can validate this calculation by visualising the empirical distribution of the $N=10\,000$ independent estimators $\widehat\sigma^{(i)}(T)$ for $T = 200\,000$ and $i \in \{1, \dots, N\}$. This is shown in Figure \ref{fig: QMLDistribution_Heston_firstobserved1}, which displays a histogram of the $10 \, 000$ estimators along with the Gaussian density with mean $\sigma^*=0.3$ and standard deviation $\frac{\mathrm{Std}^*[\widehat \sigma]}{\sqrt{200\,000}}$ and next to a density estimate of the empirical distribution using a Gaussian kernel density estimator with bandwidth selection according to \cite{Scott1992}. A D'Agostino--Pearson test (see \cite{DAgostino1973}) does not reject the normality of the sample at any reasonable level with a $p$-value of approximately 0.6479. The mean of the empirical distribution is approximately $0.2997$ with a standard deviation of $\smash{\frac{3.2974}{\sqrt{200\,000}}}$, validating the standard deviation calculation: a two-sided chi-square test does not reject the hypothesis $\mathrm{Var}(\widehat \sigma(200\,000)) = \frac{3.2871^2}{200\,000}$ at any reasonable level with a $p$-value\footnote{Since the asymptotic distribution of the test is sufficiently symmetric, we use the definition $p = 2 \min\{F(S), \linebreak 1-F(S)\}$ for the two-sided $p$-value of a test with asymptotic distribution function $F$ and test statistic $S$.} of approximately $0.6512$.

A further validation of the covariance calculations from Section \ref{sec5.1: CovMat} and Section \ref{su:explicit} can be obtained by evaluating the empirical size of the Wald, Lagrange multiplier, and likelihood-ratio tests, that is, the frequency of committing a Type I error using $N = 10\,000$ independent repetitions of the tests. Even though these tests become rather simplistic in the case of an isolated estimation of a single parameter, we include the test results for $H_0: \sigma^* = 0.3$ for the sake of exposition. Table \ref{tab: tests} includes the empirical size of the tests, i.e.\ the frequency under $H_0$ of the test statistics exceeding the $(1 - \alpha)$-quantile of the respective limiting distribution from Proposition \ref{prop: significance_tests}. The number of observations included for the tests is again $T = 200\,000$. By the central limit theorem, these sizes are approximately distributed as $\mathcal{N}\big(\alpha, \frac{\alpha(1- \alpha)}{N}\big)$ so that an asymptotic 95\% confidence interval (CI) and a two-sided $p$-value for the sizes are given in Table \ref{tab: tests} as well. The latter is again calculated as $2\min\{F(S), 1 - F(S)\}$, where $S$ is the empirical size and $F$ is the distribution function of the $\mathcal{N}\big(\alpha, \frac{\alpha(1- \alpha)}{N}\big)$ distribution. Evidently, all tests at the levels $\alpha = 5\%$ and $\alpha = 10\%$ yield satisfying rejection frequencies. However, the rejection frequencies at $\alpha = 1\%$ exceed the confidence range for the size, indicating a slightly downward-biased variance estimate $\widehat V(t) = V_{\widehat \theta(t)}$. These findings are roughly in line with the test results obtained by \cite{Mainassara2014} for the estimation of weak VARMA models.

\begin{remark}
    It is also possible to estimate the parameters of a Heston model on a small time scale. To this end, one may define the polynomial state space model $X = (\widetilde v, \Delta Y, (\Delta Y)^2)$ with $\widetilde v(k) \coloneq v(k\Delta t)$ and $\Delta Y(k) \coloneq Y(k \Delta t) - Y(k \Delta t - \Delta t)$ for some small fixed time increment $\Delta t$. For $\Delta t = \frac{1}{24 000}$, which roughly corresponds to a five-minute interval of observations if time is measured in years, one obtains the asymptotic standard deviation $\mathrm{Std}^*[\widehat \sigma] \approx 160.0723$ for the example of an isolated estimation of the volatility of volatility $\sigma$. The asymptotic accuracy can potentially be improved by incorporating higher powers of the variance and the returns to the state space model in order to leverage the information contained in higher moments of $v$ and $Y$. For example, if one uses $X = (\widetilde v, \widetilde v^2, \Delta Y, (\Delta Y)^2, (\Delta Y)^4)$ with the first two components treated as unobservable, one obtains $\mathrm{Std}^*[\widehat \sigma] \approx 81.8182$. In terms of calendar time, this means that using this larger state space model in conjunction with a smaller time-discretisation requires a time series of data that is approximately 39 times shorter, albeit with a higher data frequency, to achieve the same accuracy of estimation compared to using the time-discretisation $\Delta t = 1$.
\end{remark}

\begin{table}[!b]
    \centering
    \caption{Frequency of rejection under the null hypothesis $H_0: \sigma^* = 0.3$ for the Wald, the LM and the LR test from Proposition \ref{prop: significance_tests}. Asymptotically, these sizes of the hypothesis tests fall theoretically within the given confidence intervals in the fifth column of the table with a probability of 95\%. $p$-values of the sizes are given below.} 
    \small
    \begin{tabular}{|C{1.4cm}|C{1.8cm}|C{1.8cm}|C{1.8cm}|C{3cm}|} \hline
    Level & Wald test & LM test & LR test & $95\%$ CI for the size \\ \hline
    \vspace*{-0.1cm}$\alpha = 1\%$ & 1.27\% \linebreak {\scriptsize ($p = 0.0067$)} & 1.30\% \linebreak {\scriptsize ($p = 0.0026$)} & 1.30\% \linebreak {\scriptsize ($p = 0.0026$)} & \vspace*{-0.1cm}[0.805\%, 1.195\%] \\ \hline
    \vspace*{-0.1cm}$\alpha = 5\%$ & 5.32\% \linebreak {\scriptsize ($p = 0.1420$)} & 5.39\% \linebreak {\scriptsize ($p = 0.0735$)} & 5.39\% \linebreak {\scriptsize ($p = 0.0735$)} & \vspace*{-0.1cm}[4.573\%, 5.427\%] \\ \hline
    \vspace*{-0.1cm}$\alpha = 10\%$ & 10.44\% \linebreak {\scriptsize ($p = 0.1425$)} & 10.54\% \linebreak {\scriptsize ($p = 0.0719$)} & 10.53\% \linebreak {\scriptsize ($p = 0.0773$)} & \vspace*{-0.1cm}[9.412\%, 10.588\%] \\ \hline
    \end{tabular}
    \label{tab: tests}
\end{table}

\begin{remark}
    One might pose the question how the asymptotic standard deviations for the four parameters $\kappa$, $m$, $\sigma$ and $\rho$ reported in Figure \ref{fig: Covariance_Heston_firstobserved1} compare to asymptotic standard deviations of other estimators in the literature. The smallest possible estimator variance is given by the inverse Fisher information, which is achieved by the true maximum likelihood estimator (based on partial observations). Even though the latter is too expensive to compute in practice, one can obtain particle approximations of the observed Fisher information, see e.g. \cite{Nemeth2016}. Using Algorithm 2 of \cite{Nemeth2016} with $\lambda = 0.999$, $N = 10\,000$ particles and $200\,000$ observations of $X(t) = (v, \Delta Y, (\Delta Y)^2)$ with $\Delta t = 1$, one obtains asymptotic standard deviations of the maximum likelihood estimators for the isolated estimation of the parameters $\kappa$, $m$, $\sigma$ and $\rho$. These amount to approximately $3.4676$, $0.1910$, $0.5344$ and $1.7780$, respectively, which is smaller than the asymptotic standard deviation of the QML estimators by a factor approximately between 2 and 6. However, we reckon that different choices of the state space components of $X$ including higher powers of $\Delta Y$ would \mbox{lead to improved QML estimator standard deviations in practice.}
\end{remark}

\subsection{Multivariate Lévy-driven \OU processes}\label{sec5.3: OU}

In this section we study the estimation of \OU processes driven by a multivariate Lévy process. The univariate Gaussian \OU process was initially proposed by \cite{Uhlenbeck1930} in a physical environment and motivated in this context by its mean-reverting properties. Since then, multiple generalisations to the classic \OU process have appeared in the literature including those driven by a general Lévy process, as initially proposed by \cite{Sato1983} and \cite{Sato1985}. As in \cite{Masuda2004} we define a $d$-dimensional analogue of a univariate Lévy-driven \OU process as a strong solution to the stochastic differential equation
\begin{equation}\label{eq: OU_equation}
    \mathrm{d}X(t) = -QX(t) \dd t + \mathrm{d}L(t)
\end{equation}
on some probability space $(\Omega, \F, \PP)$, where $Q \in \R^{d\times d}$ and $L = (L(t))_{t \in \R_+}$ is a $d$-dimensional Lévy process, usually called the \textit{background driving Lévy process}. 
This family of models has a wide range of applications. A particularly prominent one in financial economics lies in the stochastic volatility specification of the well-known \cite{BarndorffNielsen2001} model, but other motivations for studying processes of the form \eqref{eq: OU_equation} include the modelling of energy prices and derivatives (see for example \cite{Benth2006}) or of the so-called storage equation from storage theory (see \cite{Cinlar1971}). Multivariate Lévy-driven \OU processes have for example been applied to the phenomenon of co-integration in econometrics by \cite{Fasen2013}. We can state the following standard result from \cite[Section 6.3]{Applebaum2009}:

\begin{proposition}\label{prop: OU_strongsolution}
    Given $X(0)$, the unique strong solution to the equation \eqref{eq: OU_equation} is of the form $X = (X(t))_{t\in\R_+}$ with $X(t) = \e^{-Qt} X(0) + \int_0^t \e^{-Q(t-s)} \dd L(s)$ for $t \in \R_+$.
\end{proposition}

For an \OU process $X$, it can be shown that the law of $X(t)$ is infinitely divisible and that the conditional characteristic function of $X(t)$ given $\F_s$ is given by
\begin{equation}\label{eq: OU_char}
    \E\big[\e^{\ii u^\top X(t)} \:\big|\: \F_s\big] = \exp\bigg[\ii u^\top \e^{-Q(t-s)} X(s) + \int_s^t \psi^L\big(\e^{-Q^\top(t-r)} u\big) \dd r \bigg]
\end{equation}
(see e.g.\ \cite{Sato1984}, Theorem 3.1), where $\psi^L$ denotes the characteristic exponent of $L$ and $\F$ is the natural filtration of $X$. Thus, $X$ is an affine process and we can apply the theory developed in Sections \ref{sec4: Estimation} and \ref{sec4.2: Proof} to the estimation of an \OU process $X$ defined in terms of some $Q^\theta \in \R^{d \times d}$ and some $d$-dimensional Lévy process $L$ whose law under $\PP_\theta$ depends on $\theta$. Concerning the latter assumption we have the following simple condition:

\begin{proposition}\label{prop: OU_lp}
    Let $X$ be an \OU process with background driving Lévy process $L$. Let $p \in [1,\infty)$ and suppose $\E\big(\lVert L(1) \rVert^p\big) < \infty$ (or equivalently $\int_{\lVert x \rVert \geq 1} \lVert x \rVert^p \nu^L(\mathrm{d}x) < \infty$, where $\nu^L$ is the Lévy measure of $L$). Then $\smash{\E\big(\lVert X(t) \rVert^p\big) < \infty}$ for any $t \in \R_+$. If moreover $\alpha(Q) > 0$, i.e.\ $Q$ has only eigenvalues with positive real parts, then $X$ is bounded in $L^p$.
\end{proposition}

The preceding proposition shows that the \OU polynomial state space model $(\widetilde X(t))_{t \in \N}$ with $\widetilde X(k) = X(k \Delta t)$ defined by some matrix $Q^\theta \in \R^{d\times d}$, some $d$-dimensional Lévy process $L$, and some step size $\Delta t$ is a parametric polynomial state space model of order $r$ if $\E_\theta(\lVert L(1) \rVert^r) < \infty$ for any $\theta \in \Theta$. By slight abuse of notation, we will write $X$ for both the continuous-time \OU process $(X(t))_{t \in \R_+}$ and the state space model $\widetilde X$ with step size $\Delta t$, even though the time scales generally differ. It is then possible to express the state transition matrix and state transition vector $A^\theta$ and $a^\theta$ as well as the noise covariance matrices $C^\theta(t)$ from Proposition \ref{prop: Gauss_equiv} in closed form if $Q^\theta$ is assumed to be non-singular. To do so one only needs expressions for the expectation and the covariance matrix of $L(1)$, which are easily obtained for most common Lévy processes:

\begin{proposition}\label{prop: OU_expressions}
    Assume that $\E_\theta(\lVert L(1)\rVert^2) < \infty$ for any $\theta \in \Theta$ and set $a^L_\theta \coloneq \E_\theta(L(1))$ as well as $c^L_\theta \coloneq \mathrm{Cov}_\theta(L(1))$. Assume that the matrix $Q^\theta$ is non-singular for any $\theta \in \Theta$ and define the long-run mean $\mu^\theta_\infty \coloneq (Q^\theta)^{-1} a^L_\theta$. Then the state transition vector and state transition matrix of the \OU polynomial state space model $X$ are given by $\smash{a^\theta = (\mathrm{I}_d - \e^{-Q^\theta \Delta t}) \mu^\theta_\infty}$ and $\smash{A^\theta = \e^{-Q^\theta \Delta t}}$. Moreover, the noise covariance matrices $C^\theta$ from Proposition \ref{prop: Gauss_equiv} are independent of $t$. With $\Xi_\theta \coloneq \mathrm{I}_{d^2} - \e^{-(Q^\theta \oplus Q^\theta) \Delta t}$ they are given by
    \begin{equation}\label{eq: OU_BB}
        \mathrm{vec}(C^\theta(t)) = (Q^\theta \oplus Q^\theta)^{-1} \Xi_\theta\big([\mu^\theta_\infty \oplus \mu^\theta_\infty] a^L_\theta + \mathrm{vec}(c^L_\theta)\big) -\Xi_\theta \big(\mu^\theta_\infty \otimes \mu^\theta_\infty\big).
    \end{equation}
\end{proposition}

Since the expressions for the mean and covariance matrix of the Lévy process $L$ together with $Q^\theta$ fully determine the quasi-maximum likelihood estimator for the model, this once more puts an emphasis on the interpretation of the Gaussian quasi-maximum likelihood estimator as a sort of generalised method of moments estimator, as discussed in Section \ref{sec1: Intro}.

\begin{remark}
\OU processes are natural generalisations of discrete-time autoregressive processes of order one and hence belong to the class of continuous-time autoregressive moving average (CARMA) processes as introduced in \cite{Doob1944}. In particular, it has been shown in \cite{Schlemm2012b}, Corollary 3.4, that multivariate CARMA processes and Lévy-driven \OU processes are equivalent. Estimation of these processes is a well-developed field in the literature, including \cite{Taufer2009} using empirical characteristic functions and \cite{Valdivieso2009} as well as \cite{Lu2022} using maximum likelihood methods in the univariate and in the multivariate setting, respectively. If the \OU process has finite second moments and is observed at discrete points in time, consistency and asymptotic normality have been shown for a standard least-squares estimator in \cite{Fasen2013}, Proposition 3.2, for a method of moments estimator in \cite{Spiliopoulos2009}, Theorem 4.1, and for a quasi-maximum likelihood estimator in \cite{Schlemm2012} under slightly different conditions than ours. Of these authors, only \cite{Schlemm2012} deal with partially observed \OU processes. On the other hand, neither of these studies seems to use the affine polynomial moment structure of the \OU process, and so the presumed conditions for consistency and asymptotic normality in \cite{Schlemm2012} are based on strong mixing properties and all in all less elementary than those introduced in Section \ref{sec4.1: StateSpaceModels}.
\end{remark}

We will now briefly treat the fulfilment of the Assumptions \ref{assump: AN}, \ref{normassump}, \ref{assump: Identifiability2} from Section \ref{s:GQLE} for a parametric \OU polynomial state space model $X$ parameterised by some $\theta \in \Theta$, where $\Theta$ is a convex and compact subset of $\R^k$. By Proposition \ref{prop: OU_expressions}, Assumption \ref{assump: AN} is fulfilled whenever $Q^\theta$, $a_\theta^L$, and $c_\theta^L$ are $\mathrm{C}^3$-functions of $\theta$. The identification of a suitable state space $E$ for $X$ such that the transition measures of $X$ are equivalent to $\lambda_E$ is difficult to determine in generality and depends on the chosen background-driving Lévy process $L$. Conditions ensuring the existence of a Lebesgue density of $X(t)$ given $X(t-\Delta t)$ can e.g.\ be found in \cite{Priola2007} or \cite{Simon2011}. For example, Theorem 1.1 in the former establishes that the law of $X(t)$ for any $t > 0$ is absolutely continuous with respect to $\lambda_{d}$ if the Lévy measure $\nu^\theta$ of $L$ under $\PP_\theta$ is infinite and has a Lebesgue density on some open neighbourhood of zero. However, these conditions do not warrant positivity of the density on all of $\R^d$ or a suitable subset of $\R^d$. Under a condition which ensures that the Lévy measure $\nu^\theta$ has sufficient mass around zero (see \cite{Priola2007}, Hypothesis 1.2), the characteristic function \eqref{eq: OU_char} is integrable, and so, by the Fourier inversion formula \mbox{(cf. \cite{Sato1999}, Proposition 2.5(xii)), the conditional density of $X(t)$ is}
\begin{equation*}
    f^\theta_{X(t) \mid X(t - \Delta t) = x}(y) = \frac{1}{(2 \pi)^d} \int_{\R^d} \exp\Big[\ii u^\top \big(\e^{-Q^{\theta} \Delta t}x - y\big) + \int_{0}^{\Delta t} \psi^{L}_\theta \big(\e^{-Q^{\theta^\top} s} u\big) \dd s  \Big] \dd u
\end{equation*}
for $y \in \R^d$, whence positivity conditions can be derived to fully establish Assumption \ref{normassump}(\ref{itm: Irreducibility}) if a specific Lévy process $L$ has been chosen. If $a_\theta^L = 0$, i.e.\ $L$ is a centred Lévy process under $\PP_\theta$, then \eqref{eq: OU_BB} simplifies to $\mathrm{vec}(C^\theta(t)) = (Q^\theta \oplus Q^\theta)^{-1} (\mathrm{I}_d - \e^{-(Q^\theta \oplus Q^\theta) \Delta t}) \mathrm{vec}(c_\theta^L)$, which can also be expressed as $\smash{C^\theta(t) = \int_0^{\Delta t} \e^{-Q^\theta s} c_\theta^L \e^{-Q^{\theta^\top} s} \dd s}$, see the proof of Proposition \ref{prop: OU_expressions}. Together with Proposition \ref{prop: B-assump} this shows that Assumption \ref{normassump}(\ref{itm: Irreducibility}) can fail if $c^L_\theta$ is not positive definite for some $\theta \in \Theta$. The Feller property for $X$ follows from the fact that $X$ is obtained by sampling from an affine process or by an application of the dominated convergence theorem to the representation of $X(t)$ from Proposition \ref{prop: OU_strongsolution}. Concerning Assumption \ref{normassump}(\ref{itm: L6-bounded}), $L^{4 + \delta}$-boundedness follows from Proposition \ref{prop: OU_lp} if $\alpha(Q^\theta) > 0$ and $\E_\theta(\lVert L(1)\rVert^{4 + \delta}) < \infty$ for some $\delta > 0$ and all $\theta \in \Theta$.

\begin{remark}
    Even though we focus on Lévy-driven \OU processes here, it should be rather straightforward to extend the results of this section to \textit{Generalised \OU processes}, i.e.\ càdlàg $\R^d$-valued processes solving $\mathrm{d}X(t) = \mathrm{d}V(t)X(t-) + \mathrm{d}L(t)$ for some $\R^{d\times d} \times \R^d$-valued Lévy process $(V, L)$ (see \cite{Behme2012}), which  are in general non-affine polynomial. In this case explicit expressions for the state transition vector and matrix can be derived using similar methods as in the proof of Proposition \ref{prop: OU_expressions} or as in \cite{Eberlein2019}, Example 6.31. Concerning Assumption \ref{normassump}, the Feller property follows again from an explicit expression for the process as in Proposition \ref{prop: OU_strongsolution} (see \cite{Behme2012}, Theorem 3.4) and conditions for boundedness in $L^{4 + \delta}$ can be derived by methods similar to \cite{Behme2011}, Theorem 3.1, or \cite{Lindner2005}, Proposition 4.1. 
\end{remark}

As a specific application of the above theory inspired by financial mathematics, we consider a Lévy-driven two-factor interest rate model for the short rate, where the short rate can be thought of as some sort of instantaneous forward interest rate with infinitesimal horizon. For details about interest rate models in general and the example considered here see Chapter 14 and Example 14.10 in \cite{Eberlein2019}, respectively. In our case we assume that the short rate $r = (r(t))_{t \in \R_+}$ follows an \OU process that incorporates a mean-reverting behaviour towards some average $m = (m(t))_{t \in \R_+}$ that is itself random. Specifically,
\begin{align}
    \mathrm{d}m(t) &= -\lambda m(t) \dd t + \mathrm{d}L_1(t), \label{eq: ou_interest1} \\
    \mathrm{d}r(t) &= \kappa\big(m(t) - r(t)\big) \dd t + \mathrm{d}L_2(t) \label{eq: ou_interest2}
\end{align}
under $\PP_\theta$, where $\lambda, \kappa \in (0,\infty)$ and $L = (L_1, L_2)^\top$ is a bivariate Lévy process and where $\theta = (\lambda, \kappa, \theta_0) \in \Theta \subseteq \R^k$ with some additional parameters $\theta_0 \in \R^{k-2}$. It is convenient for our purposes to assume that the short rate $r$ can be observed from bond prices in the market and that the stochastic mean $m$ is unobservable. This model $X = (m, r)$ is of the form \eqref{eq: OU_equation} with
\begin{equation*}
    Q^\theta = \begin{pmatrix}
        \lambda & 0 \\
        -\kappa & \kappa
    \end{pmatrix}.
\end{equation*}

As a parametric family of Lévy processes $L$  we consider bivariate normal-inverse Gaussian L\'evy processes with parameters $\mu,\beta \in \R^2$, $\alpha, \delta \in (0,\infty)$, where $\Delta \in \R^{2\times 2}$ is positive definite with $\det(\Delta) = 1$ and $\beta^\top \Delta \beta < \alpha^2$, see for example \cite{Eberlein2019}, Section 2.4.8. The mean and covariance matrix of $L(1)$ are given by $\smash{a^L_\theta = \mu + \frac{\delta}{\sqrt{\psi}} \Delta \beta}$ and $\smash{c^L_\theta = \frac{\delta}{\sqrt{\psi}} \big(\Delta + \frac{1}{\psi} \Delta \beta \beta^\top \Delta \big)}$, where $\psi \coloneq \alpha^2 - \beta^\top \Delta \beta$.
For ease of exposition we simplify our calculations by focusing on centred and symmetric normal-inverse Gaussian Lévy processes with uncorrelated components as the background driving Lévy process, i.e.\ we focus on $\mu=0$, $\beta = 0$, $\Delta = \mathrm{I}_2$, for which $a^L_\theta = 0$, $c^L_\theta = \frac{\delta}{\alpha} \mathrm{I}_2$. In this case estimation of the four-dimensional parameter $\theta = (\lambda, \kappa, \alpha, \delta)$ using the Gaussian quasi-maximum likelihood approach is not possible because the parameters $\alpha$ and $\delta$ are not identifiable from the second moments of $L^\theta$ but only their quotient $\frac{\delta}{\alpha}$, so Assumption \ref{assump: Identifiability2} is not satisfied for a separate estimation of $\alpha$ and $\delta$. Hence, we confine ourselves to the estimation of the three-dimensional parameter $\theta = (\lambda, \kappa, \delta)$ with fixed $\alpha = 1$.

\begin{remark}
As a limitation introduced by Assumption \ref{assump: Identifiability}, it is only possible to estimate those parameters of a state space model by QML estimation that are identifiable from the first and second moments of the state space model. In order to make inference on the parameters $\alpha$ and $\delta$ separately, one could allow for second-order polynomials of $r(t)$, $t\in\N$ rather than only $r(t)$. This can be achieved by replacing $X=(m,r)$ with the higher-dimensional polynomial state-space model $X=(r, r^2, m, mr, m^2)$, whose components $m, mr, m^2$ are unobservable. In the same vein, one could also estimate a non-zero parameter $\beta$ of the underlying Lévy process. In that case all parameters of the state space model $X=(m, r)$ would be identifiable.
\end{remark}

\begin{figure}[!b]
    \centering
    \includegraphics[width=13.1cm]{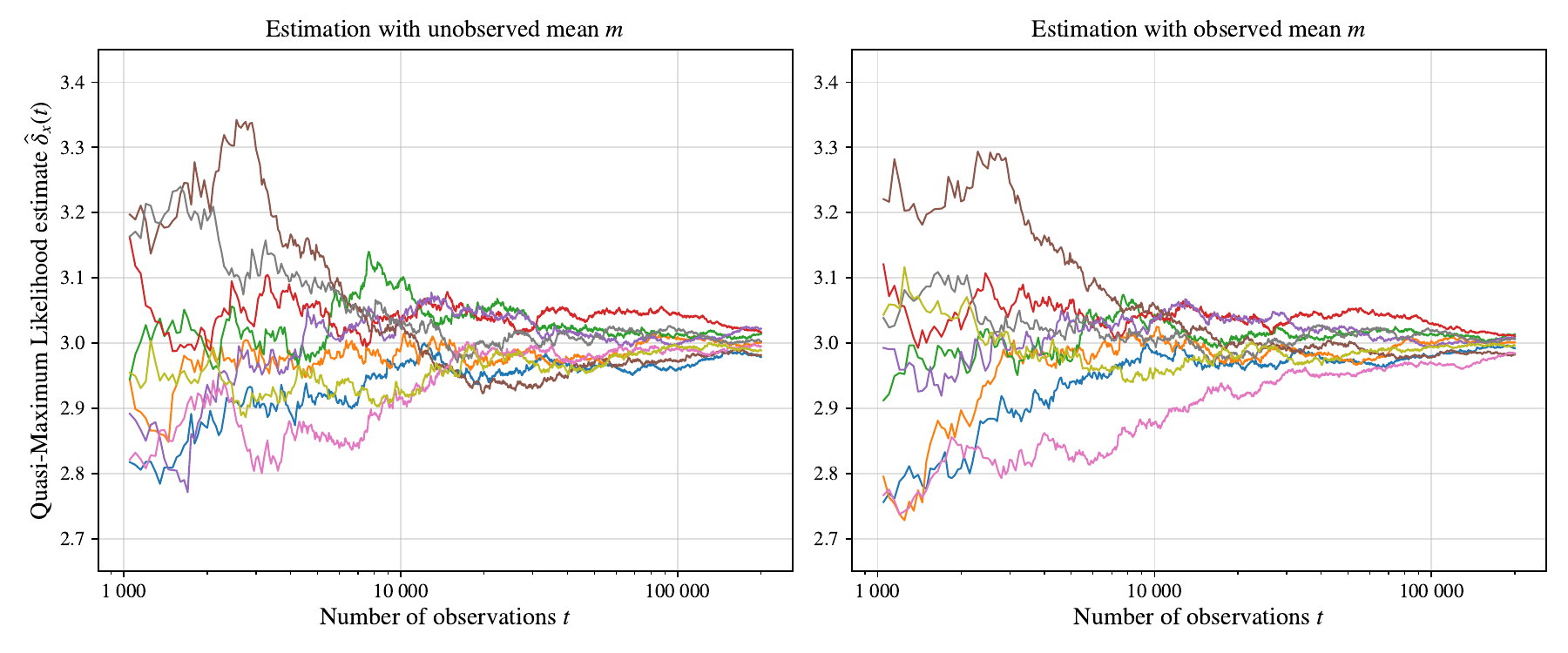}
    \caption{Ten sequences of quasi-maximum likelihood estimators for the parameter $\delta^* = 3$. On the left $m$ is assumed to be unobservable, while it is observed on the right.}
    \label{fig: OU_Sequences}
\end{figure}

Since the Assumptions \ref{assump: AN}, \ref{normassump}, \ref{assump: Identifiability2} are fulfilled for the \OU model described above, Theorems \ref{theo: main} and \ref{theo: main2} yield again consistency and asymptotic normality of any sequence of quasi-maximum likelihood estimators, provided that the Fisher information matrix $W(\vartheta)$ is invertible. We again visualise this result by setting $\vartheta = (\lambda^*, \kappa^*, \delta^*) \coloneq (1,\, 0.5,\, 3)$ and by simulating $N = 10\, 000$ independent sample trajectories $(X^\vartheta(t))_{t \in \{0, 1, \dots, T\}}$ of size $T = 200\, 000$ with step size $\Delta t = 1$, where the stochastic integral from Proposition \ref{prop: OU_strongsolution} is discretised by an Euler method with mesh size $\frac{1}{5000}$. Similar to the Heston model from Section \ref{sec5.2: Heston}, the initial distribution for $X(0)$ is set to a Dirac distribution in the point $x = (0.5, 1)$ under any $\PP_\theta$. The parameter space $\Theta$ is chosen as $\Theta = [10^{-4}, 10] \times [10^{-4}, 10] \times [10^{-4}, 1000]$. As in the case of the Heston model, we compare our results for (i) a joint estimation of $\theta$ and (ii) an isolated estimation of the separate parameter component $\delta$ with all other components assumed to be known. Figure \ref{fig: OU_Sequences} shows ten independent sequences $\smash{(\widehat \delta(t))_{t \in \{1, \dots, T\}}}$ of the quasi-maximum likelihood estimator, first if the latent mean $m$ is assumed to be unobservable and secondly if $m$ is assumed to be observable.

\begin{figure}[!b]
    \centering
    \includegraphics[width=14.5cm]{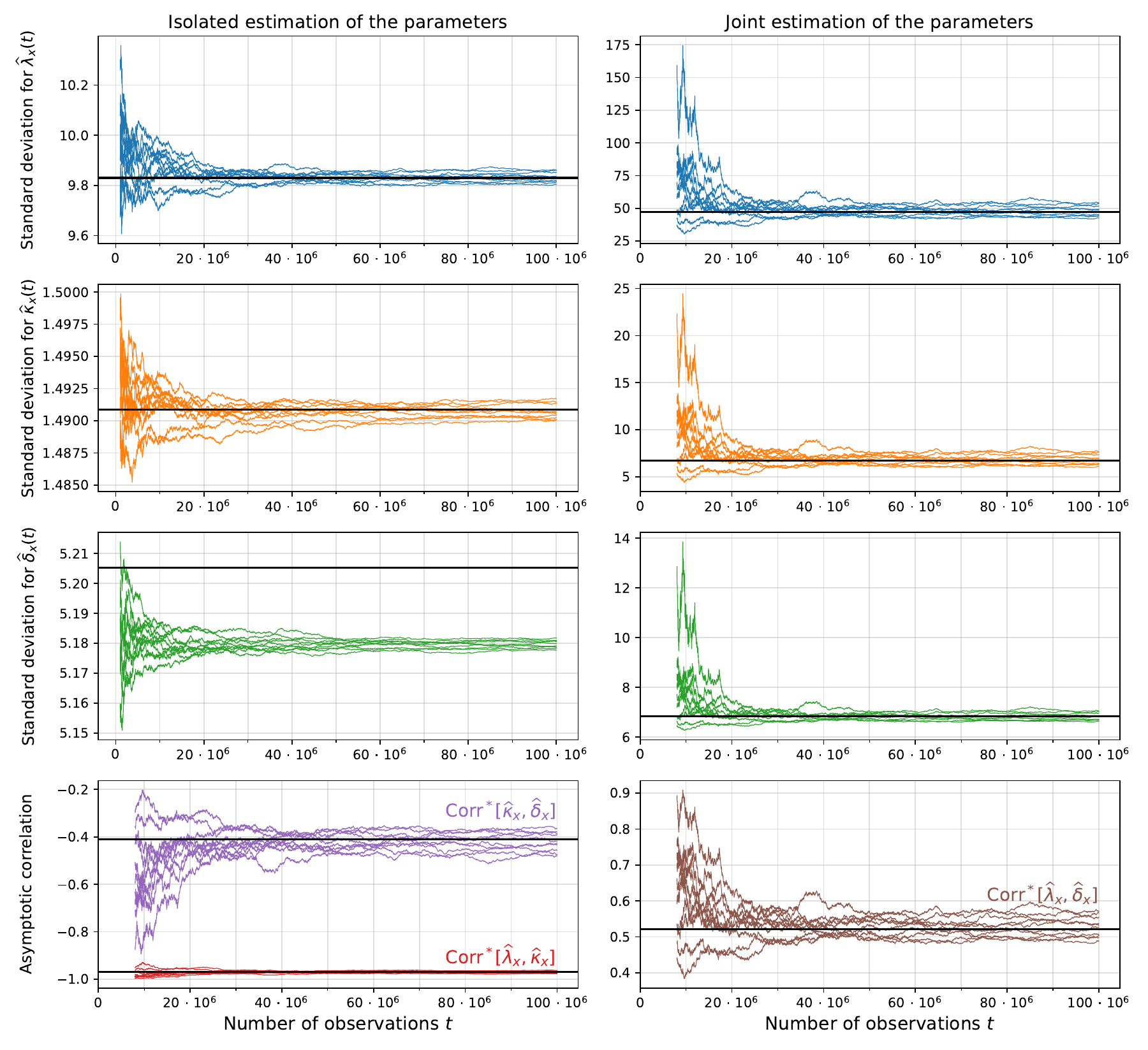}
    \caption{The upper three rows contain ten sequences of estimator standard deviations, respectively obtained from the estimates $\widehat V_t[\widehat \theta(t)]$ for $t$ between 1 and $100 \cdot 10^6$. The left column displays the standard deviation for an isolated estimation of the parameter components, while the right shows the same if the whole parameter $\vartheta$ is estimated. The last row contains ten sequences of asymptotic estimator correlations obtained from $\widehat V_t[\widehat \theta(t)]$. Black lines show values obtained from the explicit calculations detailed in Section \ref{su:explicit}.}
    \label{fig: Covariance_OU_firstobserved1}
\end{figure}

Again as in the case of the Heston model, introducing additional observable components for estimation reduces the estimator's standard deviation. In the case of estimating the parameter $\delta$ of the \OU model described above, the reduction in variance is however much smaller compared to the estimation of the volatility of volatility parameter $\sigma$ in the Heston model (see Figure \ref{fig: QML_Sequences}). This could potentially originate from the fact that the dispersion parameter $\delta$ influences both components of our model equally, while in the case of the Heston model, the parameter $\sigma$ affects the observed returns only indirectly. Similar to Figures \ref{fig: Covariance_Heston_firstobserved1} and \ref{fig: Correlation_Heston_firstobserved1}, the tiles of  Figure \ref{fig: Covariance_OU_firstobserved1} each contain ten standard deviation (or correlation) sequences for the estimator's components next to the exact calculation of the corresponding asymptotic covariance matrix components, which is depicted by a black line.

Recall the discussion in Section \ref{sec5.2: Heston} of the potential drawback of the covariance estimator $\widehat V_t[\widehat \theta(t)]$. It becomes even more pronounced here because the large amount of around $50 \cdot 10^6$ observations is necessary for the estimator to reach a sufficiently stable state, as it is visible in Figure \ref{fig: Covariance_OU_firstobserved1}. Again, this demands a rather high calculation effort in practice, speaking once more for the explicit covariance matrix calculations from Section \ref{su:explicit}. Apart from that, all covariance estimates visible in Figure \ref{fig: Covariance_OU_firstobserved1} are well in line with the theoretical computations from Section \ref{su:explicit}, except for the case of an isolated estimation of $\delta$, in which case the standard deviation estimate has a downward bias of around 5\% of the true asymptotic standard deviation.

\begin{figure}[!b]
    \centering
    \includegraphics[width=12cm]{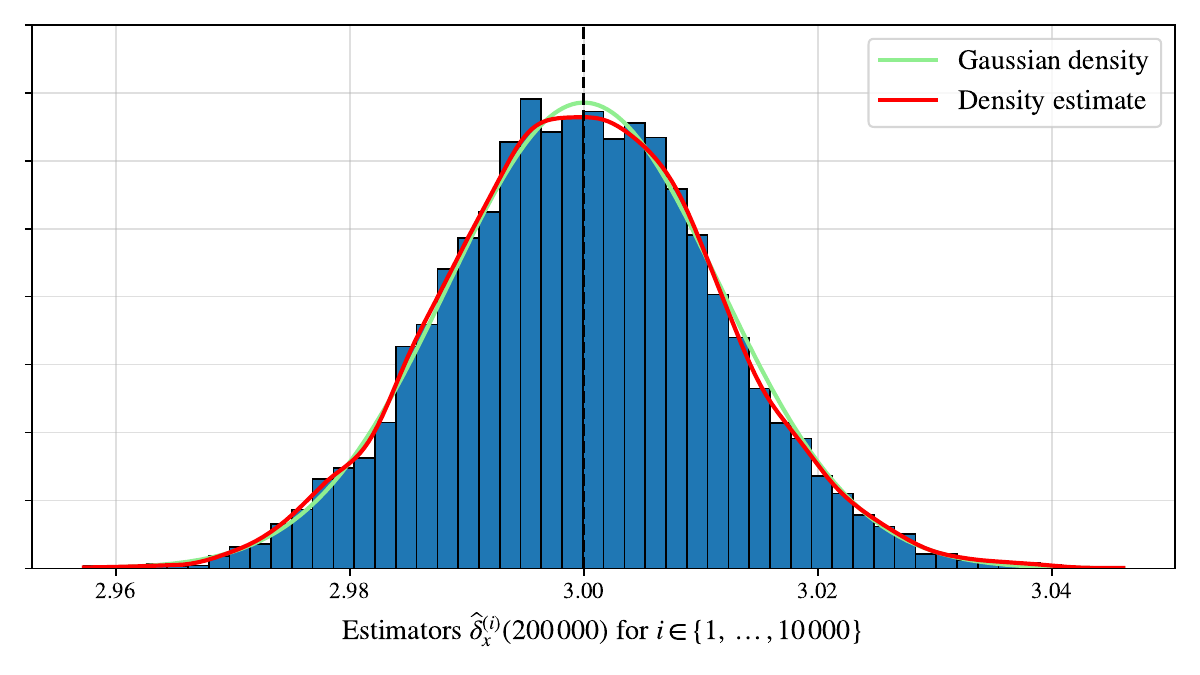}
    \caption{Histogram of $N = 10\,000$ independent quasi-maximum likelihood estimators $\widehat \delta^{(i)}(T)$, $i \in \{1, \dots, N\}$, for $\delta$ with $T = 200\,000$. On top we draw the Gaussian density corresponding to the normal distribution with mean $\delta^* = 3$ and variance given by the explicit calculation of $V_\vartheta$ in green, and we draw a kernel density estimate of the estimators $\smash{\widehat \delta^{(i)}(T)}$ in red. The black dotted line indicates the true parameter value $\delta^* = 3$.}
    \label{fig: QMLDistribution_OU_firstobserved1}
\end{figure} 

In Figure \ref{fig: QMLDistribution_OU_firstobserved1} we again display a histogram of $N = 10 \, 000$ independent quasi-maximum likelihood estimators $\widehat \delta^{(i)}(T)$ for $T = 200\, 000$ and $i \in \{1, \dots, N\}$ in the case of an isolated estimation of $\delta$ with the latent mean $m$ treated as unobservable. Next to a kernel density estimate, we again show the Gaussian density corresponding to a mean of $\delta^* = 3$ and a standard deviation of $\smash{\frac{\mathrm{Std}^*[\widehat \delta]}{\sqrt{200\,000}}}$, where the asymptotic estimator standard deviation $\smash{\mathrm{Std}^*[\widehat \delta] :\approx 5.2054}$ has been calculated explicitly using the methods from Section \ref{su:explicit}. As in the case of the Heston model, the computed standard deviation of $\smash{\frac{5.2054}{\sqrt{200\,000}}}$ matches well with the standard deviation of the empirical distribution of the $N$ estimators, which is approximately $\frac{5.1953}{\sqrt{200\,000}}$. A two-sided chi-square test does not reject the hypothesis $\smash{\mathrm{Var}(\widehat \delta(200\,000)) = \frac{5.2054^2}{200\,000}}$ at any reasonable level with a $p$-value of approximately $0.7890$, which is again calculated as $2\min\{F(S), 1 - F(S)\}$, where $S$ and $F$ are the corresponding test statistic and the chi-square cumulative distribution function, respectively.

\begin{remark}
    As in the case of the Heston model, the asymptotic accuracy of the estimation depends on the chosen step size $\Delta t$. Choosing $\Delta t = \frac{1}{24000}$ corresponds again approximately to a five-minute increment if time is measured in years. In this case the asymptotic standard deviation for estimating $\delta$ becomes $\smash{\mathrm{Std}^*[\widehat \delta]} \approx 464.7759$. As discussed earlier in the case of the Heston model, a higher accuracy of estimation can be achieved by incorporating higher powers of the components into the state space model.
\end{remark}

\section{Proofs}\label{s:proofs}

\subsection{Proofs for Section \ref{sec4: Estimation}}

\begin{proof}[Proof of Proposition \ref{prop: Gauss_equiv}]
This follows from Proposition 2.18 in \cite{KallsenRichert2025}.
\end{proof}

\begin{proof}[Proof of Proposition \ref{prop: log-lik}]
For the Kálmán filter recursions see for example \cite{KallsenRichert2025}, Proposition 3.1. Equation \eqref{eq: log-lik} follows then from the fact that $Y(t)$ is conditionally Gaussian given $Y_\oo(1), \dots, \linebreak Y_\oo(t-1)$ under $\mathbb{Q}_\theta$, with conditional mean given by the Kálmán filter mean and conditional covariance matrix given by the Kálmán filter covariance matrix. 
\end{proof}

\begin{proof}[Proof of Proposition \ref{prop: affine_lp}]
    As noted e.g. in \cite{Eberlein2019}, Example 6.30, the affine process $(X(t))_{t \in \R_+}$ is a $p$-polynomial process. By Example \ref{ex: Polynomial}, $(X,(\PP_\theta)_{\theta \in \Theta})$ is a polynomial state space model of order $p$. Hence, $(\mathrm{vec}_{\otimes p}(X),(\PP_\theta)_{\theta \in \Theta})$ is a polynomial state space model of order 1 and can be represented as $\mathrm{vec}_{\otimes p}(X(t)) = a^\theta_{\otimes p} + A^\theta_{\otimes p} \mathrm{vec}_{\otimes p}(X(t-1))  + N^\theta_{\otimes p}(t)$ with some vector $a^\theta_{\otimes p}$, some matrix $A^\theta_{\otimes p}$, and a martingale difference $N^\theta_{\otimes p}$.
    By equation \eqref{eq: Mean Recurrence} and Corollary \ref{linalg}, it suffices to show that $\rho(A^\theta_{\otimes p}) < 1$ for boundedness of $p$-th moments of $X$, which yields that $X$ is bounded in $L^p(\PP_\theta)$ because $p$ is even. Without loss of generality we assume $p=2$. The general case follows along the same lines. By \cite{Eberlein2019}, Theorem 6.6,
    \begin{equation}\label{eq: varphi}
        \varphi^\theta_t(u) \coloneq \E_\theta\big[\e^{\ii u^\top X(t + 1)} \big| \F_{t}\big] = \e^{\Psi^\theta_0(u) + \sum_{j=1}^d\Psi^\theta_j(u) X_j(t)} =: \e^{\psi^\theta(u, t)} , \quad u \in \R^d
    \end{equation} 
    for $t \in \N$ and some complex-valued functions $\Psi^\theta_j$, $j \in \{0, \dots, d\}$. Since $\E_\theta[X(t + 1) |\F_{t}] = -\ii \nabla_u \varphi^\theta_t(0)$, the state transition matrix $A^\theta$ of $X$ is given by $A^\theta_{ij} = -\ii \,\partial_{u_i} \Psi^\theta_j(0)$ for $i, j \in \{1, \dots, d\}$, which follows by differentiation in \eqref{eq: varphi} and because $\varphi^\theta_t(0) = 1$. Moreover, 
    \begin{align*}
        \E_\theta\big[X_i(t+1)X_j(t+1) |\F_{t}] &= - \partial_{u_i}\partial_{u_j} \varphi^\theta_t(0) = - \partial_{u_i} \Big[ \big(\partial_{u_j} \psi^\theta(u, t)\big)\varphi^\theta_t(u)\Big]\Big\rvert_{u = 0} \\
        &= -\partial_{u_i}\partial_{u_j} \psi^\theta(0, t) - \big(\partial_{u_i} \psi^\theta(0, t)\big) \big(\partial_{u_j} \psi^\theta(0, t)\big)
    \end{align*}
    for $i, j \in \{1, \dots, d\}$ which, up to an affine-linear function $f(X(t))$, is quadratic in $X(t)$ with $$\E_\theta\big[X_i(t+1)X_j(t+1) |\F_{t}] = \sum_{(n, m) \in \{1, \dots, d\}^2} -\partial_{u_i}\Psi^\theta_n(0)\partial_{u_j}\Psi^\theta_m(0) X_n(t)X_m(t) + f(X(t)).$$
    The quadratic coefficients of the above function are precisely the entries of the matrix $\smash{A^\theta \otimes A^\theta}$. It follows that the state transition matrix $A^\theta_{\otimes 2}$ has the block triangular form
    \begin{equation}\label{eq: otimes-matrix}
    A^\theta_{\otimes 2} = \begin{pmatrix}
        A^\theta & 0 \\
        \circ & A^\theta \otimes A^\theta
    \end{pmatrix}.        
    \end{equation} 
    Since $\rho(A^\theta)<1$, we have $\rho(A^\theta \otimes A^\theta) < 1$ and hence $\rho(A^\theta_{\otimes 2}) < 1$, which proves that $X$ is bounded in $L^p(\PP_\theta)$ for $p=2$, as desired. The case $p > 2$ can be handled analogously.
\end{proof}

\begin{proof}[Proof of Proposition \ref{prop: ergod}]
    Existence and uniqueness of a stationary law for $X$ and ergodicity of $X$ follow immediately by Theorem \ref{theo: f-ergod}  because Assumption \ref{normassump}.\ref{itm: Irreducibility} implies $\psi$-irreducibility and aperiodicity for $\psi=\lambda_E$. The final part of Theorem \ref{theo: f-ergod} together with Assumption \ref{normassump}.\ref{itm: L6-bounded} yields finiteness of moments of order $4+\delta$ or less under $\mu_\theta$. Moreover, Theorem \ref{theo: f-ergod} states that for any polynomial $g: \R^d \to \R$ of order 4 or less $\E_{\theta}\big[g(X(t)) \mid X(0) = x\big] \to \int g \dd\mu_\theta$ for $\lambda_E$-almost all $x \in E$. It remains to verify that $\E_{\theta}\big[g(X(t))\big] \to \int g \dd\mu_\theta$ at a geometric rate if $g: \R^d \to \R$ is a quadratic resp.\ quartic polynomial. From the above and \eqref{eq: Mean Recurrence}, we deduce that
    \begin{equation}\label{e:aus2.2}
    \E_{\theta}(X(t) \mid X(0) = x) = (A^\theta)^t x + \Bigl(\sum_{s=0}^{t-1} (A^\theta)^{s}\Bigr) a^\theta
    \end{equation}
    converges to some constant independent of $x$ as $t \to \infty$ for all $x \in D$, where $D$ is a dense subset of $E$. In particular, $(A^\theta)^t (x-y)  = \E_{\theta}(X(t) \mid X(0)= x) - \E_{\theta}(X(t) \mid X(0) = y) \to 0$ for all $x, y \in D$. By Lemma \ref{lem: power_convergence} we can find $\widetilde A^\theta \in \R^{d\times d}$ with $\rho(\widetilde A^\theta) < 1$ such that $\widetilde A^\theta (x-y) = A^\theta(x-y)$ for all $x, y \in D$. Since $0 \in E$ and $D$ is dense in $E$, we can first send $y$ to 0 to deduce $\widetilde A^\theta x = A^\theta x$ on the dense subset $D \subset E$, which implies $\widetilde A^\theta x = A^\theta x$ for all $x \in E$. So without loss of generality, we can assume that $A^\theta = \widetilde A^\theta$ and hence $\rho(A^\theta) < 1$. This yields 
     \begin{fitequationtwo}\label{eq: EXt_decomp}
    \E_{\theta}(X(t))=\int\E_{\theta}(X(t) \mid X(0) = x)\PP_\theta^{X(0)}(\mathrm{d} x)=(A^\theta)^t\E_{\theta}(X(0))+ \Bigl(\sum_{s=0}^{t-1} (A^\theta)^{s}\Bigr) a^\theta \to   \Bigl(\sum_{s=0}^{\infty} (A^\theta)^{s}\Bigr) a^\theta,
    \end{fitequationtwo}
    which proves the claim for first order polynomials. The parallel statement for quadratic resp.\ quartic polynomials follows from the same argument applied to $A^\theta_{\otimes 2}$ resp.\ $A^\theta_{\otimes 4}$.
\end{proof}

\begin{proof}[Proof of Proposition \ref{prop: B-assump}]
    First note that $\E_\theta(X(t))$ and $\E_\theta(\mathrm{vec}_{\otimes 2}(X(t)))$ converge uniformly in $\theta$ at a geometric rate to their respective limits from Proposition \ref{prop: ergod}: Indeed, since $a^\theta$ and $A^\theta$ are continuous in $\theta$ by Assumption \ref{assump: A}, Lemma \ref{lem: linalg_unif} shows that there are $c \in \R_+$, $\gamma \in [0, 1)$ such that $\sup_{\theta\in\Theta} \lVert (A^\theta)^t\rVert \leq c \gamma^t$. Combined with the fact that $\E_\theta(X(0))$ is bounded in $\theta$, we get that $\E_\theta(X(t))$ converges uniformly at a geometric rate. The same holds for $\E_\theta(\mathrm{vec}_{\otimes 2}(X(t)))$ by applying the previous argument to $a^\theta_{\otimes 2}$ and $A^\theta_{\otimes 2}$. Now, since $C^\theta(t)$ contains by definition only first and second moments of $X(t-1)$, it also converges uniformly in $\theta$ at a geometric rate to some limiting matrix $C^\theta$ by Proposition \ref{prop: ergod}. To show that the same holds also for the derivatives $\partial_\theta^\alpha C^\theta(t)$ for any multi-index $\alpha \in \N_3^k$, it suffices to show that the derivatives of $\E_\theta(X(t))$ and $\E_\theta(\mathrm{vec}_{\otimes 2}(X(t)))$ converge uniformly in $\theta$ at a geometric rate. For $\E_\theta(X(t))$, this follows quickly by differentiating \eqref{eq: EXt_decomp} and applying Corollary \ref{coro: linalg_diff}. The same holds for $\E_\theta(\mathrm{vec}_{\otimes 2}(X(t)))$ if one exchanges $a^\theta$ and $A^\theta$ by $a^\theta_{\otimes 2}$ and $A^\theta_{\otimes 2}$ in the previous argument.
    
    It remains to show that $C^\theta(t)$ and $C^\theta = \lim_{t \to \infty} C^\theta(t)$ are positive definite. Since $P_\theta(x, \cdot)$ is equivalent to $\lambda_E$ for any $x \in E$ and $E$ is not contained in any proper subspace of $\R^d$, $\smash{\mathrm{Cov}_\theta(X(t) \mid X(t-1) = x) = \mathrm{Cov}_\theta(N(t) \mid X(t-1)=x)}$ is positive definite for any $x \in E$. Now, because $\E_\theta(N(t) \mid X(t-1)) = 0$ almost surely and $$\smash{\mathrm{Cov}_\theta(N(t)) = \E_\theta\big[\mathrm{Cov}_\theta(N(t) \mid X(t-1))\big] + \mathrm{Cov}_\theta\big[\E_\theta(N(t) \mid X(t-1))\big]},$$ we have $\mathrm{Cov}_\theta(N(t)) = \int \mathrm{Cov}_\theta(N(t) \mid X(t-1) = x) \PP_\theta^{X(t-1)}(\mathrm{d}x)$, so $\mathrm{Cov}_\theta(N(t))$ is also positive definite. This result holds irrespective of the initial distribution of $X(0)$, so by choosing the invariant law $X(0) \sim \mu_\theta$, we get that $C^\theta = \lim_{t \to \infty} C^\theta(t)$ is also positive definite.
\end{proof}

\subsection{Proofs for Section \ref{su: mainmain}}

\begin{proof}[Proof of Proposition \ref{prop: Quasi-Score}]
    For $t \in \N$ the matrices $\widehat \Sigma^\theta(t + 1, t)$  are positive definite because the matrices $C^\theta(t)$ are positive definite (see Proposition \ref{prop: B-assump}) and hence all pseudoinverses occurring in Proposition \ref{prop: log-lik} are proper inverses. The result now follows by differentiating the log-likelihood \eqref{eq: log-lik} with respect to $\theta_j$, using $\partial_{\theta_j} M_\theta^{-1} = -M_\theta^{-1} (\partial_{\theta_j} M_\theta)M_\theta^{-1}$  and $\partial_{\theta_j} \log \det M_\theta = \Tr\big[M_\theta^{-1} \partial_{\theta_j} M_\theta\big]$ for a matrix-valued differentiable function $\theta \mapsto M_\theta$ with $M_\theta$ being positive definite for all $\theta$, see \cite{Petersen2012} and its references.
\end{proof}

\begin{corollary}\label{coro: EZ=0}
    $\E_\vartheta(Z^\vartheta(t)) = 0$ for any $t \in \N$.
\end{corollary}
\begin{proof}
    Since $Z^\vartheta(t)$ is quadratic in $X_{\oo}$ and the second moments of $X_{\oo}$ and $Y_{\oo}$ under $\PP_\vartheta$ (resp. $\mathbb{Q}_\theta$) match, we have $\smash{\E_\vartheta(Z^\vartheta(t)) = \E_\vartheta\big[\nabla_\theta \log q_t^{\vartheta}(Y_\oo(1), \dots, Y_\oo(t))\big] = \int \nabla_\theta q_t^{\vartheta}(y) \dd y}$. This vanishes if we can argue that $\smash{\nabla_\theta \int q_t^{\vartheta}(y) \dd y = \int \nabla_\theta q_t^{\vartheta}(y) \dd y}$, which holds by dominated convergence if we find Lebesgue-integrable functions $h_j$ with $\smash{h_j(y) \geq |\partial_{\theta_j} q_t^\theta(y)|}$ for all $\theta \in \Theta$. Now we have that $\smash{\partial_{\theta_j} q_t^\theta(y)} = \big(\partial_{\theta_j} \log q_t^\theta(y)\big) q_t^\theta(y)$. The first factor on the right-hand side is quadratic in $y_1, \dots, y_t$ with continuous coefficients in $\theta$ while the second factor on the right hand side is a Gaussian density with continuous mean and covariance in $\theta$. Since $\Theta$ is compact, all these coefficients and Gaussian densities are uniformly bounded in $\theta$.
\end{proof}
The proof of Theorems \ref{theo: main} and \ref{theo: main2} relies on two classic results from estimation theory:

\begin{proposition}\label{prop: Jacod-Consistency}
    Assume that there is some continuous function $Q: \Theta \to \R$ such that \linebreak $\smash{\sup_{\theta \in \Theta} \lVert \frac{1}{t} L^\theta(t) - Q(\theta)\rVert \xrightarrow{\PP_\vartheta} 0}$ and such that $Q$ has a unique maximum at $\vartheta$. Then every sequence $\smash{(\widehat \theta(t))_{t \in \N}}$ of quasi-maximum likelihood estimators is $\vartheta$-consistent.
\end{proposition}
\begin{proof}
    Since $\Theta$ is separable and both $L^\theta(t)$ and $Q(\theta)$ are continuous functions of $\theta$, it follows from Lemma \ref{lem: measurability} that the uncountable supremum above is measurable and hence a well-defined random variable.
    Since convergence in probability is equivalent to the fact that each subsequence has a further subsequence converging almost surely, we see that any subsequence of $\smash{(\widehat \theta(t))_{t \in \N}}$ has a further subsequence which converges $\PP_\vartheta$-almost surely to $\vartheta$ by Lemma \ref{lem: argmax}. Hence $\smash{\widehat \theta(t) \xrightarrow{\PP_\vartheta} \vartheta}$ for any such sequence $\smash{(\widehat \theta(t))_{t \in \N}}$ of estimators.
\end{proof}

Under the conditions above, the existence of a $\vartheta$-consistent sequence of quasi-maximum likelihood estimators becomes very natural: since $\frac{1}{t}L^\theta(t)$ converges to some continuous function $Q$ in probability uniformly on $\Theta$, and since $Q$ is uniquely maximised at $\vartheta$, the quasi-likelihood process $L^\theta(t)$ must also have a maximum near $\vartheta$ for sufficiently large $t \in \N$.

\begin{proposition}\label{prop: Jacod-Normality}
    Suppose that $\vartheta \in \mathrm{int}(\Theta)$ and let the conditions of Proposition \ref{prop: Jacod-Consistency} hold with twice continuously differentiable $Q$. Set $W \coloneq \nabla_\theta^2 Q$. Assume that $\smash{\frac{1}{\sqrt{t}} Z^\vartheta(t) \xrightarrow{\PP_\vartheta\text{-}d} Z}$ for some $\R^k$-valued random variable $Z$, that $W(\vartheta)$ is invertible, and that $\sup_{\theta \in \Theta} \lVert \frac{1}{t} \nabla_\theta Z^\theta(t) - W(\theta) \rVert \smash{\xrightarrow{\PP_\vartheta}} 0$. Then $\smash{\sqrt{t}(\widehat \theta(t) - \vartheta)\xrightarrow{\PP_\vartheta\text{-}d} -W(\vartheta)^{-1} Z}$ for any $\vartheta$-consistent sequence $\smash{(\widehat \theta(t))_{t \in \N}}$ of QML estimators.
\end{proposition}
\begin{proof}
    This follows as in \cite{Jacod2018}, Theorem 2.11. In particular, since $\vartheta \in \mathrm{int}(\Theta)$ and since $\smash{\widehat \theta(t) \xrightarrow{\PP_\vartheta} \vartheta}$, we find that $\smash{\PP_\vartheta(B_t) \coloneq \PP_\vartheta(\widehat \theta(t) \in \mathrm{int}(\Theta)) \to 1}$ and so $\smash{\PP_\vartheta(Z^{\theta}(t)\big|_{\theta = \widehat \theta(t)} = 0) \to 1}$. Moreover, by the same argument as in the proof of Proposition \ref{prop: Jacod-Consistency}, the uncountable supremum occurring in this proposition is measurable and hence a well-defined random variable. By the multivariate mean value theorem it follows that
    \[Z^{\widehat \theta(t)}(t) - Z^\vartheta(t) = \underbrace{\Big(\int_0^1  \nabla_\theta Z^{\lambda \vartheta + (1 - \lambda) \widehat \theta(t)}(t) \dd \lambda \Big)}_{\eqqcolon \widetilde Z^\vartheta(t)} (\widehat \theta(t) - \vartheta).\]
    Since the integrand above is continuous in $\lambda$, $\widetilde Z^\vartheta(t)$ is measurable. Moreover, we can deduce
    $$\lVert \frac{1}{t}\widetilde Z^\vartheta(t) - W(\vartheta) \rVert  \leq \sup_{\theta \in \Theta} \lVert \frac{1}{t} \nabla_\theta Z^\theta(t) - W(\theta) \rVert + \sup_{\lVert \theta - \vartheta\rVert \leq \lVert \widehat \theta(t) - \vartheta\rVert} \big\lbrace \lVert W(\theta) - W(\vartheta)\rVert \big\rbrace,$$
 which implies that $\smash{\frac{1}{t}\widetilde Z^\vartheta(t) \xrightarrow{\PP_\vartheta} W(\vartheta)}$
 because $W$ is continuous and $\smash{\widehat \theta(t)}$ is consistent. Defining the set $\smash{C_t \coloneq \{\widetilde Z^\vartheta(t) \text{ is invertible}\} \cap B_t}$, we have $\PP_\vartheta(C_t) \to 1$. On the set $C_t$ it holds that
    $\smash{\sqrt{t}(\widehat \theta(t) - \vartheta)} = \smash{- t\widetilde Z^\vartheta(t)^{-1} \frac{1}{\sqrt{t}} Z^\vartheta(t)}$ and the right-hand side converges in $\PP_\vartheta$-law to $-W(\vartheta)^{-1}Z$ by Slutsky's theorem. Since $\PP_\vartheta(C_t) \to 1$ as $t \to \infty$, this concludes.
\end{proof}

The main goal of Sections \ref{sec4.2.1: CLT}--\ref{sec4.2.3: Ident} is to establish the various conditions from Propositions \ref{prop: Jacod-Consistency} and \ref{prop: Jacod-Normality} in order to prove Theorems \ref{theo: main} and \ref{theo: main2}. The following Section \ref{sec4.2.1: CLT} is concerned with the condition $\smash{\frac{1}{\sqrt{t}} Z^\vartheta(t) \xrightarrow{\PP_\vartheta\text{-}d} Z}$ from Proposition \ref{prop: Jacod-Normality}. In Section \ref{sec4.2.2: ULLN} we deal with the uniform convergence $\smash{\sup_{\theta \in \Theta} \lVert \frac{1}{t}\nabla_\theta Z^\theta(t) - W(\theta)\rVert \xrightarrow{\PP_\vartheta} 0}$ from Proposition \ref{prop: Jacod-Normality} before Section \ref{sec4.2.3: Ident} fills in missing details and finishes the proof of Theorems \ref{theo: main} and \ref{theo: main2}.

\subsubsection{A functional central limit theorem for the quasi-score process}\label{sec4.2.1: CLT}

The general strategy to prove the condition from Proposition \ref{prop: Jacod-Normality} is rather straightforward: Since $\smash{\frac{1}{\sqrt{t}} Z^\vartheta(t) = \frac{1}{\sqrt{t}}\sum_{s=1}^t Z^\vartheta(s, s-1)}$, the goal is to represent $Z^\vartheta(t, t-1)$ as a function of a Markov chain and use a suitable central limit theorem for stationary Markov chains. The proof of this result will require several ergodicity properties of $Z^\vartheta$. Thus it will be necessary to use notions from the ergodic theory of discrete-time Markov chains. These details are confined to Appendix \ref{appA: Markov}, where we adapt our exposition to the classic book of \cite{Meyn2009}. The difficulty of our approach lies in the fact that $Z^\theta(t, t-1)$ is in general a function of all observations $X_\oo(1), \dots, X_\oo(t-1)$ because the observed part $X_\oo$ is typically not a Markov process. As in Section \ref{sec5.1: CovMat}, we therefore introduce the $E \times \R^{d(k + 1)}$-valued processes $\smash{\langX(t) \coloneq \big(X(t), \widehat X^\theta(t, t-1), V^\theta(t, t-1)\big)}$, which contain the polynomial state space model $X$ augmented by the filtered component and by the first partial derivatives of the filtered component. The dynamics of $\smash{\langX}$ under $P_\vartheta$ can be written as
\begin{fitequation}\label{eq: Xtilde_dynamics}
    \langX(t + 1) = \begin{pmatrix}
        a^\vartheta \\ a^\theta \\ \partial_1 a^\theta \\ \partial_2 a^\theta \\ \vdots \\ \partial_k a^\theta
    \end{pmatrix} + \begin{pmatrix}
        A^\vartheta & 0 & 0 & 0 & \dots & 0 \\
        K^\theta(t) H & F^\theta(t) & 0 & 0 & \dots & 0 \\
        \partial_{1} K^\theta(t) H & \partial_{1} F^\theta(t) & F^\theta(t) & 0 & \dots & 0\\
        \partial_{2}K^\theta(t)H & \partial_{2} F^\theta(t) & 0 & F^\theta(t) & \dots & 0 \\
        \vdots & \vdots & \vdots & \vdots & \ddots & \vdots \\
        \partial_{k} K^\theta(t)H & \partial_{k}F^\theta(t) & 0 & 0 &\dots & F^\theta(t)
    \end{pmatrix} \langX(t) + \begin{pmatrix}
        N^\vartheta(t + 1) \\ 0 \\ 0 \\ 0 \\ \vdots \\ 0
    \end{pmatrix}
\end{fitequation}

\noindent for $t \in \N^*$, where $K^\theta(t) = A^\theta \widehat \Sigma^\theta_{:, \oo}(t, t-1) \widehat\Sigma^\theta_{\oo}(t, t-1)^{-1}$, where $F^\theta(t) \coloneq A^\theta - K^\theta(t)H$, and $H \coloneq (\delta_{m+i, j})_{i\in\{1, \dots, d-m\}, j \in \{1, \dots, d\}}$, i.e., $H \in \R^{(d-m) \times d}$ is the matrix such that $x_\oo = Hx$ for $x \in \R^d$. Moreover, $\smash{\langX(0) = (X(0), \E_\theta(X(0)), \nabla_\theta \E_\theta(X(0)))}$ holds almost surely under any measure $\PP_\theta$. The advantage of studying the process $\smash{\langX}$ in place of $X$ is that the conditional quasi-score $Z^\theta(t, t-1)$ is a function just of $\smash{\langX(t)}$ and not of any other past values as it is the case for $X$.

Note that $\smash{\langX}$ is not a time-homogeneous Markov chain because its transition dynamics depend on time $t$ through $C^\theta(t)$. A direct application of a suitable Markov chain central limit theorem to the quasi-score is hence out of reach. However, Lemma \ref{limitcov} in Section \ref{sec5.1: CovMat} establishes that the Kálmán filter covariance matrices $\smash{\widehat \Sigma^\theta(t, t-1)}$ and their derivatives $\smash{S^\theta_j(t, t-1)}$ converge as $t \to \infty$ at a geometric rate. We show in Proposition \ref{prop: asymptotic_equiv} that this implies that $\smash{\langX}$ can be asymptotically well-approximated by a time-homogeneous Markov chain because its transition dynamics become asymptotically independent of time. 

\begin{remark}
    In comparison to comparable literature for estimation of (partially-observed) VARMA models as \cite{Mainassara2011} or \cite{Schlemm2012}, our approach of using Markovianity of an augmented version of the original process seems to be novel in the context of quasi-maximum likelihood estimation. The usual path to proving asymptotic normality of the quasi-score process consists in using certain strong mixing properties for the observed part of the process as well as results from spectral time series analysis, see for example the summability condition on the mixing coefficients cited in \cite{Schlemm2012}. Our conditions from Section \ref{s:GQLE} seem sufficient to prove such mixing properties by first establishing geometric ergodicity of $X$ in a similar manner as in Proposition \ref{prop: geom_ergod_y} and then using the well-known fact that geometric ergodicity implies an exponential mixing behaviour, see for example \cite{Liebscher2005}, Proposition 4. Nonetheless, we follow the different route described above because it enables us to derive closed-form expressions for asymptotic estimator covariances without using spectral analysis.
\end{remark}

\begin{proof}[Proof of Lemma \ref{limitcov}]
    For notational convenience we fix some $j \in \{1, \dots, k\}$ and write $\widehat \Sigma^\theta(t)$ for $\smash{\widehat \Sigma^\theta(t,t-1)}$ and  $S^\theta(t)$ for $S^\theta_j(t, t-1)$. From the Kálmán filter recursions it is apparent that $\smash{\widehat \Sigma^\theta(t)}$ obeys the discrete-time algebraic Riccati difference equation
    \begin{equation}\label{Riccati}
        \widehat \Sigma^\theta(t + 1) = A^\theta \Big[\widehat \Sigma^\theta(t) - \widehat \Sigma^\theta_{:, \oo}(t) \widehat \Sigma^\theta_\oo(t)^{-1}\widehat \Sigma^\theta_{:, \oo}(t)^\top\Big]A^{\theta^\top} + C^\theta(t+1)
    \end{equation}     
     for $t \in \N^*$. Indeed, since $C^\theta(t)$ is positive definite by Proposition \ref{prop: B-assump}, $\widehat \Sigma^\theta(t)$ and hence also $\widehat \Sigma^\theta_\oo(t)$ are positive definite for $t \in \N^*$, whence the inverse in \eqref{Riccati} is well-defined. Since $\rho(A^\theta) < 1$ by Proposition \ref{prop: ergod}, \cite{Anderson1979}, pp. 78--82, show that there exists a unique positive definite fixed point $\widehat \Sigma^\theta$ for \eqref{Riccati} with $\smash{C^\theta }$ in place of $\smash{C^\theta(t)}$. We proceed to show that $\widehat \Sigma^\theta(t)$ in \eqref{Riccati} converges to $\widehat \Sigma^\theta$ at a geometric rate uniformly in $\theta \in \Theta$. First, we argue that the sequence ($\widehat \Sigma^\theta(t))_{t \in \N^*}$ is bounded uniformly in $\theta$. Consider the auxiliary filter $\widehat X^{\theta,\mathrm{aux}}$ with
    \[\widehat X^{\theta,\mathrm{aux}}(t, t-1) \coloneq a^\theta + A^\theta \widehat X^{\theta,\mathrm{aux}}(t-1, t-2).\]
    Since the Kálmán filter $\widehat X^\theta(t, t-1)$ is the optimal filter in the equivalent Gaussian state space model, \cite{KallsenRichert2025}, Proposition 3.2, implies that
    \begin{align}\label{e:ast}
    \Tr \,\widehat \Sigma^\theta(t) &= \E_\theta(\lVert X(t) - \widehat X^\theta(t, t-1)\rVert^2) \leq \E_\theta(\lVert X(t) - \widehat X^{\theta, \mathrm{aux}}(t, t-1)\rVert^2) \\
    &= \Tr\,\widehat\Sigma^{\theta,\mathrm{aux}}(t), \nonumber
    \end{align}
    where $\widehat \Sigma^{\theta,\mathrm{aux}}(t)$ denotes the auxiliary filter error covariance matrix. But $\widehat \Sigma^{\theta,\mathrm{aux}}(t)$ obeys the discrete Lyapunov equation
    \begin{equation}\label{eq: auxfilter}
    \widehat \Sigma^{\theta,\mathrm{aux}}(t) = A^\theta \widehat \Sigma^{\theta,\mathrm{aux}}(t-1) {A^\theta}^\top + C^\theta(t),       
    \end{equation} 
    with $\rho(A^\theta) < 1$. Since $C^\theta(t)$ converges to $\smash{C^\theta}$ uniformly, it is bounded uniformly in $\theta$, so Lemma \ref{lem: lin_systems}(1) shows that the sequence $(\widehat \Sigma^{\theta,\mathrm{aux}}(t))_{t \in \N^*}$ is bounded uniformly in $\theta$. A fortiori, the sequence $(\widehat \Sigma^\theta(t))_{t\in \N^*}$ is also bounded uniformly in $\theta$ by \eqref{e:ast}. Now, a simple calculation shows that the Kálmán filter covariance matrix can be written in the Lyapunov form
    \begin{equation}\label{eq: Lyapunov}
        \widehat \Sigma^\theta(t + 1) = (A^\theta - K^\theta(t) H) \widehat \Sigma^\theta(t) (A^\theta - K^\theta(t) H)^\top + C^\theta(t + 1)
    \end{equation} 
    for $t \in \N^*$, where $K^\theta(t) = A^\theta \widehat \Sigma^\theta_{:, \oo}(t) \widehat \Sigma^\theta_{\oo}(t)^{-1}$. Fix $s \in \{1, \dots, t-1\}$ and define the matrix $\smash{\Psi^\theta_{t, s}} \coloneq (A^\theta - K^\theta(t-1)H) \cdots (A^\theta - K^\theta(s) H)$ . Since $C^\theta(u)$ is positive definite for any $u \in \N^*$, it follows that
    \[\widehat \Sigma^\theta(t) = \Psi^\theta_{t, s} \widehat \Sigma^\theta(s) \Psi^{\theta^\top}_{t, s} + \text{positive definite terms} \geq \Psi^\theta_{t, s} \widehat \Sigma^\theta(s) \Psi^{\theta^\top}_{t, s},\]
    and hence the matrices on the right-hand side are bounded uniformly in $\theta$. By the Courant-Fischer theorem, the $j$-th diagonal element of the matrix on the right hand side is bounded from below by $\lambda_{\min}(\widehat \Sigma^\theta(s)) \lVert(\Psi^\theta_{t, s})_j\rVert^2$, where $(\Psi^\theta_{t, s})_j$ denotes the $j$-th row of $\Psi^\theta_{t, s}$ and $\lambda_{\min}(\widehat \Sigma^\theta(s))$ the smallest eigenvalue of $\smash{\widehat \Sigma^\theta(s)}$. Since $\Theta$ is compact and $\smash{\lambda_{\min}(\widehat \Sigma^\theta(s))}$ is continuous in $\theta$ and positive for all $\theta \in \Theta$, it is uniformly bounded from below by a positive constant $\lambda$. Hence $\lambda \lVert(\Psi^\theta_{t, s})_j\rVert^2$ is uniformly bounded in $\theta$ and $t$. In particular, the matrices $\Psi^\theta_{t, s}$ are bounded uniformly in $\theta$, $t$, and $s < t$. Following a lengthy calculation, one can show
    \[\widehat \Sigma^\theta(t) - \widehat \Sigma^\theta = (A^\theta - K^\theta H) (\widehat \Sigma^\theta(t-1) - \widehat \Sigma^\theta) (A^\theta - K^\theta(t-1) H)^\top + D^\theta(t),\]
    where $K^\theta \coloneq A^\theta \widehat \Sigma^\theta_{ :, \oo} (\widehat \Sigma^\theta_{\oo})^{-1}$ is the limiting Kálmán gain and $D^\theta(t) \coloneq C^\theta(t) - \smash{C^\theta }$, see \cite{Anderson1979}, Problem 4.5. By iterating this equation one obtains
    \begin{fitequation}\label{eq: Essential_Decomp}
    \widehat\Sigma^\theta(t) - \widehat \Sigma^\theta = (A^\theta - K^\theta H)^{t-1} (\widehat \Sigma^\theta(1) - \widehat \Sigma^\theta) \Psi_{t, 1}^{\theta^\top} + \sum_{s=0}^{t-2} (A^\theta - K^\theta H)^{s} D^\theta(t-s) \Psi^{\theta^\top}_{t, t-s}.
    \end{fitequation}
    Let $F^\theta \coloneq A^\theta - K^\theta H$. We now argue that $\rho(F^\theta) \leq \alpha$ for any $\theta \in \Theta$ and some $\alpha  \in [0,1)$. Analogously to equation \eqref{eq: Lyapunov}, $\smash{\widehat \Sigma^\theta}$ is a fixed point of the matrix-valued Lyapunov equation
    \begin{equation}\label{eq: Lyapunov_stable}
        \widehat \Sigma^\theta = F^\theta \widehat \Sigma^\theta F^{\theta^\top} + C^\theta.
    \end{equation} 
    Suppose now that for any $\varepsilon > 0$ there exists some $\theta \in \Theta$ such that $\smash{F^{\theta^\top}}$ has some eigenvalue $\lambda_\theta$ with $|\lambda_\theta|^2 > 1 - \varepsilon$, and let $v_\theta \neq 0$ denote a corresponding complex eigenvector. It follows
    \begin{equation}\label{eq: eigenvaluebound}
    (1 - |\lambda_\theta|^2) \overline{v_\theta} \widehat \Sigma^\theta v_\theta = \overline{v_\theta}C^\theta  v_\theta,
    \end{equation} 
    where the bar denotes the complex conjugate transpose. If $|\lambda_\theta| \geq 1$, the left-hand side is non-positive while the right-hand side is positive by the positive definiteness of $\smash{C^\theta }$, which yields a contradiction. Hence $|\lambda_\theta| < 1$. Since $\smash{C^\theta}$ is positive definite for any $\theta$ by Proposition \ref{prop: B-assump}, the right-hand side of \eqref{eq: eigenvaluebound} is bounded from below by $\lambda_{\min}(C^\theta) |v_\theta|^2$. As $C^\theta(t)$ is a sequence of continuous functions converging uniformly to $\smash{C^\theta}$, $\smash{C^\theta}$ is continuous in $\theta$ as well. It follows that also $\lambda_{\min}(C^\theta)$ is continuous in $\theta$ and hence uniformly bounded from below by a positive constant $\eta$ because $\Theta$ is compact. In a similar manner, the left-hand side of \eqref{eq: eigenvaluebound} is bounded from above by $\varepsilon \rho(\widehat \Sigma^\theta) |v_\theta|^2$. By cancelling $|v_\theta|^2$ we obtain all in all that
    \[\varepsilon \rho(\widehat \Sigma^\theta) \geq \eta.\]
    But this is impossible because $\varepsilon$ was arbitrary and $\widehat \Sigma^\theta$ is bounded in $\theta$ as the pointwise limit of the uniformly bounded sequence $\widehat \Sigma^\theta(t)$. So $\rho(F^\theta) \leq \alpha < 1$ for all $\theta \in \Theta$ and some $\alpha < 1$. Finally, $\widehat \Sigma^\theta \geq C^\theta$ by \eqref{eq: Lyapunov_stable} in the Loewner order. Since the eigenvalues of $C^\theta$ are uniformly bounded from below on $\Theta$ by some positive $\eta$ as before, the eigenvalues and hence the determinant of $\widehat \Sigma^\theta$ is uniformly bounded from below on $\Theta$ by some positive constant. By Lemma \ref{lem: strange_bound} this shows that $(\widehat \Sigma^\theta_{\oo})^{-1}$ and hence also $F^\theta$ is bounded in $\theta$. Lemma \ref{lem: linalg_unif_bounded} thus yields that there are $c \in \R_+$ and $\gamma \in [0, 1)$ such that $\lVert (F^\theta)^t \rVert \leq c \gamma^t$ for all $\theta \in \Theta$. For sufficiently large $c$ and some $\widehat \gamma \in (\gamma,1)$ we also have $\lVert D^\theta(t) \rVert \leq c \widehat \gamma^t$ for all $\theta$ and some $p \geq 0$, see Proposition \ref{prop: B-assump}. By \eqref{eq: Essential_Decomp} we obtain for all $\theta \in \Theta$
    \begin{align*}
        \lVert \widehat\Sigma^\theta(t) - \widehat \Sigma^\theta \rVert \leq c \gamma^{t-1} \lVert \widehat \Sigma^\theta(1) - \widehat \Sigma^\theta\rVert \lVert \Psi^{\theta^\top}_{t, 1}\rVert + c^2 \widehat \gamma^t \sum_{s=0}^{t-2} \Big(\frac{\gamma}{\widehat \gamma}\Big)^s \lVert \Psi^{\theta^\top}_{t, t-s}\rVert.
    \end{align*}
    Since the norms $\lVert \Psi^{\theta^\top}_{t, t-s}\rVert$ and $\lVert \widehat \Sigma^\theta(1) - \widehat \Sigma^\theta\rVert$ are bounded uniformly in $\theta$, this proves $\smash{\widehat \Sigma^\theta(t) \to \widehat \Sigma^\theta}$ uniformly at a geometric rate. As products of uniformly convergent and uniformly bounded sequences are uniformly convergent, Lemma \ref{lem: geometric} 1--3 show that the same holds also for the sequence of matrices $A^\theta - K^\theta(t) H$ as well as $\widehat \Sigma^\theta(t, t)$. In particular, the matrix $F^\theta = A^\theta - K^\theta H$ is continuous on $\Theta$ because it is the uniform limit of continuous functions. 
    Let $F^\theta(t) \coloneq A^\theta - K^\theta(t) H$. To show uniform geometric convergence of the matrices $S^\theta(t)$, note that algebraic manipulations of the recursions from Proposition \ref{prop: Quasi-Score} yield
    \begin{equation*}
        S^\theta(t+1) = F^\theta(t) S^\theta(t) F^\theta(t)^\top + A^\theta \widehat \Sigma^\theta(t, t) (\partial_{j} A^\theta)^\top + (\partial_{j} A^\theta) \widehat \Sigma^\theta(t,t)A^{\theta^\top} + \partial_{\theta_j}\big(C^\theta(t+1)\big).
    \end{equation*}
    Since $A^\theta$ and its partial derivatives are bounded as continuous functions on a compact domain, the  last three summands on the right converge uniformly at a geometric rate. Thus, Lemma \ref{lem: lin_systems}.3 implies that $S^\theta(t)$ converges uniformly at a geometric rate, as desired. The matrices $R^\theta_{ij}(t, t-1)$ can be treated completely analogously.
\end{proof}

As in Section \ref{su:explicit}, Lemma \ref{limitcov} suggests that we can approximate the process $\langX$ by the time-homogeneous Markov chain $\smash{\langhomX = (X, \homkX, \homX)}$ that is defined by transition dynamics identical to $\smash{\langX}$ but with the limiting matrices $\widehat \Sigma^\theta$ and $S^{\theta}_j$ in place of $\smash{\widehat \Sigma^\theta(t, t-1)}$ and $\smash{S^\theta_j(t,t-1)}$, respectively, for $j \in \{1, \dots, k\}$. Then we obtain the following

\begin{corollary}\label{convcoro}
    Let $\smash{\wa}$ denote the state transition vector and $\smash{\wA(t)}$ the (time-dependent) state transition matrix of $\langX$ under $\PP_\vartheta$. Additionally, let $\smash{\overline N}^{\vartheta}$ denote the $d(k+2)$-dimensional martingale difference sequence for $\langX$ from \eqref{eq: Xtilde_dynamics}. Then there exists a matrix $\homA \in \R^{d(k+2) \times d(k+2)}$ with $\smash{\rho(\homA) < 1}$ such that $\smash{\wA(t) \to \homA}$ uniformly in $\theta$ at a geometric rate. Moreover, the homogeneous Markov chain $\langhomX$ with $\smash{\langhomX(0) \coloneq \langX(0)}$ and $\smash{\langhomX(t)} \coloneq \widehat a^\theta_1 + \homA \langhomX(t-1) + \smash{\overline N}^\vartheta(t)$  is a polynomial state space model of order 2 resp., under the stronger Assumption \ref{assump: AN}, of order 4.
\end{corollary}
\begin{proof}
    The uniform convergence $\smash{\wA(t) \to \homA}$ at a geometric rate follows from Lemma \ref{limitcov} because all submatrices present in $\smash{\wA(t+1)}$ are products of constant matrices and the matrices $\widehat \Sigma^\theta(t, t-1)$, $\widehat \Sigma^\theta(t, t)$ as well as $S^\theta_j(t, t-1)$ and $S^\theta_j(t, t)$ for $j \in \{1, \dots, k\}$. Moreover, $\smash{\langhomX}$ has bounded $(4+\delta)$th moments by Lemma \ref{lem: Hamilton} because $X$ does. Since $\smash{\overline N}^\vartheta$ contains only $N^\vartheta$ and zero components, it follows then that conditional expectations of quadratic (resp.\ quartic) polynomials in the components of $\smash{\langhomX(t)}$ given $\F_s$ are quadratic (resp.\ quartic) in the components of $\smash{\langhomX(s)}$, that is $\smash{\langhomX}$ is a polynomial state space model of order 2 (resp.\ 4) and bounded in $L^{4 + \delta}$. Finally, we have $\rho(\homA) < 1$ because $\homA$ is a block upper-triangular matrix whose matrices on the diagonal are $A^\vartheta$ and $F^\theta$, where $F^\theta = A^\theta - K^\theta H$ with $\smash{K^\theta = A^\theta \widehat \Sigma^\theta_{ :, \oo} (\widehat \Sigma^\theta_{\oo})^{-1}}$. Here, $\smash{\rho(A^\vartheta) < 1}$ by Proposition \ref{prop: ergod} and $\rho(F^\theta) < 1$ by the proof of Lemma \ref{limitcov}.
\end{proof}

In order to obtain stochastic stability properties for the inhomogeneous process $\langX$, we are now ready to establish a weak form of ergodicity for the homogeneous counterpart $\smash{\langhomX}$.

\begin{proposition}\label{ergod2}
    For any $\theta\in \Theta$ there is a unique invariant probability measure $\wmu$ for $\langhomX$ such that $\smash{\langhomX}$ is weakly $f$-ergodic with respect to $\wmu$ under $\PP_\vartheta$ for any polynomial $f: E \times \R^{d(k+1)} \to \R$ of degree 4 or less. Moreover, $\smash{\E_\vartheta(f(\langhomX(t))) \to \int f \mathrm{d}\wmu}$ at a geometric rate if $f$ is quadratic (quartic if Assumption \ref{assump: AN} holds). Finally, the process $\smash{\langhomX}$ is bounded in $L^{4 + \delta}(\PP_\vartheta)$.
\end{proposition}
\begin{proof}
    The state transition matrix $\homA$ of $\langhomX$ under $\PP_\vartheta$ has the block triangular form
    \begin{equation}\label{eq: transitiontilde}
    \homA = \begin{pmatrix}
        A^\vartheta & 0 & 0 \\
        G^\theta & F^\theta & 0 \\
        \circ & \circ & \mathrm{diag}_k(F^\theta)
    \end{pmatrix} = \begin{pmatrix}
        A^\vartheta & 0 & 0 & \dots & 0 \\
        G^\theta & F^\theta & 0 & \dots & 0 \\
        \circ & \circ & F^\theta & \dots & 0\\
        \vdots & \vdots & \vdots  & \ddots & \vdots \\
        \circ & \circ & \circ  &\dots & F^\theta
    \end{pmatrix},
    \end{equation}
    where $\circ$ is an irrelevant matrix of suitable size, $G^\theta = K^\theta H$ and $F^\theta = A^\theta - K^\theta H$ as before. Since the proof of this result is surprisingly technical, we start by showing the claim for the process $\primeX \coloneq (X(t), \homkX(t, t-1))_{t \in \N}$. The same strategy of proof can then be extended to the whole $d(k+2)$-dimensional process $\smash{\langhomX}$ in the end. First, $\primeX$ is weakly Feller because the transition measures $\smash{\widetilde{P}}_\theta(\cdot)$ of $\primeX$ are of the form $\smash{\widetilde{P}}_\theta(\cdot) = P_\vartheta(x_{1:d}, \cdot) \otimes \varepsilon(G^\theta x_{1:d} + F^\theta x_{d+1:2d})$, where $\varepsilon(y)$ is the Dirac measure in $y$, as usual.
    In order to prove irreducibility and aperiodicity, we consider arbitrary deterministic starting points $X(0) = x \in E$ and $\smash{\widehat X^{\theta,\mathrm{hom}}(0, -1) = \widehat x \in \R^d}$, which is signified here by writing $\PP_{\vartheta,(x,\widehat x)}$ instead of $\PP_\vartheta$ as usual.
    Iterating the Kálmán filter recursion under the measure $\PP_{\vartheta,(x,\widehat x)}$ 
    yields the closed-form expression
    \begin{fitequation}\label{eq: filterdynamics}
    \homkX(t, t-1) = (F^\theta)^t \widehat X^{\theta,\mathrm{hom}}(0, -1) + \sum_{s=0}^{t-1} (F^\theta)^s G^\theta X(t-s-1) + \Big(\sum_{s=0}^{t-1} (F^\theta)^s\Big) a^\theta.        
    \end{fitequation} 
    The process $\primeX$ itself is in general not irreducible, so a classic ergodicity proof for Markov chains cannot be applied to $\primeX$ directly. However, we are able to approximate $\primeX$ by irreducible processes in a suitable manner. Let $a^\theta_\infty \coloneq \big(\sum_{s=0}^{\infty} (F^\theta)^s\big) a^\theta$, which exists and is finite because $\rho(F^\theta) < 1$, see Lemma \ref{limitcov}. We now approximate $\smash{\homkX}$ by the auxiliary processes
    \begin{align*}
        \widehat X^{\theta,\mathrm{aux},a}(t) &\coloneq (F^\theta)^t \widehat x + \sum_{s=0}^{t-1} (F^\theta)^s G^\theta X(t-s-1) + a^\theta_\infty, \\
        \widehat X^{\theta,\mathrm{aux}}(t) &\coloneq (F^\theta)^t \widehat x + \sum_{s=0}^{t-1} (F^\theta)^s G^\theta X(t-s-1),
    \end{align*}
    and let $\auxaX \coloneq (X, \widehat X^{\theta,\mathrm{aux},a})$ and $\auxX \coloneq (X, \widehat X^{\theta,\mathrm{aux}})$. We proceed to prove irreducibility of $\auxX$ and $\auxaX$ if the state space is suitably restricted to some subset of $E \times \R^d$.
    For any matrix $M \in \R^{d \times d}$ we define $M(E)$ to be the image of $E$ under the linear map defined by $M$, and we let $\oplus$ denote the usual Minkowski sum of subsets of $\R^d$. Define
    \begin{equation}\label{eq: setsum}
     \mathscr{R}^\theta \coloneq \bigoplus_{s=0}^\infty (F^\theta)^{s}G^\theta(E) \in \mathscr B(\R^d),
    \end{equation} 
    and let $\mathscr{R}^\theta_a \coloneq \mathscr{R}^\theta + a^\theta_\infty$. Since $E$ is a connected smooth manifold containing 0, one can show that the topological covering dimension of the connected set $\mathscr R^\theta \ni 0$ coincides with its Hausdorff dimension, and the partial sums in \eqref{eq: setsum} form an increasing sequence of sets. Suppose that the state space for $\auxX$ is restricted from $E\times \R^d$ to $E \times \mathscr R^\theta$ so that $\smash{\widehat X^{\theta,\mathrm{aux}}(0) \in \mathscr{R}^\theta}$, where we equip $E \times \mathscr{R}^\theta$ with the induced topology. This restriction makes sense because as long as $\auxX(0) \in E \times \mathscr R^\theta$, we also have $\auxX(t) \in E \times \mathscr R^\theta$ and so, in the language of \cite{Meyn2009}, $E \times \mathscr{R}^\theta$ is an absorbing set for the Markov chain $\auxX$. Let $j \leq d$ be the topological covering dimension of $\mathscr R^\theta$ and let $l \in \N$ denote the smallest integer such that
    \begin{equation}\label{eq: J}
    \mathscr J^\theta\coloneq G^\theta(E) \oplus F^\theta G^\theta(E) \oplus \dots \oplus (F^\theta)^lG^\theta(E) \subseteq \mathscr R^\theta
    \end{equation}
    has topological covering dimension $j$. In particular,  any topological covering dimension considered here coincides with the Hausdorff dimension of the corresponding set under consideration. Let $\lambda_{\mathscr{R}^\theta}$ denote the $j$-dimensional Hausdorff measure on $\mathscr{R}^\theta$, as usual. Then we can find some $\varepsilon > 0$ and some open $H \subseteq \mathrm{int}(\mathscr J^\theta)$ such that $\lambda_{\mathscr{R}^\theta}(H) > 0$ and such that any $x \in \mathscr{R}^\theta$ with $\norm{x - H} < \varepsilon$ is still an element of $\mathscr J^\theta$. Here, the interior is to be understood in the induced topology on $\mathscr{R}^\theta$. As $\rho(F^\theta) < 1$, we can find $t > l + 1$ such that $\norm{(F^\theta)^{t-1}G^\theta x} < \frac{\varepsilon}{2}$ as well as $\norm{(F^\theta)^{t}\widehat x} < \frac{\varepsilon}{2}$ for any $(x,\widehat x) \in E \times \mathscr R^\theta$ by Corollary \ref{linalg}. Since $\smash{(X(1), \dots, X(t-1))}$ has a positive density with respect to the proper Hausdorff measure \mbox{$\lambda_{E^{t-1}}$ on $E^{t-1}$ by Assumption B(1), we get}
    \begin{align*}
    \PP_{\vartheta,(x,\widehat x)}(\widehat X^{\theta,\mathrm{aux}}(t) \in H) &= \PP_{\vartheta,(x,\widehat x)}((F^\theta)^{t-1}G^\theta x + (F^\theta)^t \widehat x + \zeta \in H) \\ &= \PP_{\vartheta,(x,\widehat x)}(\zeta \in H - (F^\theta)^{t-1}G^\theta x - (F^\theta)^t \widehat x),        
    \end{align*}
    \noindent where $\zeta$ is a random variable with a positive density with respect to $j$-dimensional proper Hausdorff measure on $G^\theta(E) \oplus F^\theta G^\theta(E) \oplus \dots \oplus (F^\theta)^{t-2}G^\theta(E)$. Since $\mathscr{J}^\theta$ has the same dimension as $\mathscr{R}^\theta$, since $(F^\theta)^{t-1}G^\theta x + (F^\theta)^t \widehat x \in \mathscr{R}^\theta$, and since $\big\lVert (F^\theta)^{t-1}G^\theta x + (F^\theta)^t \widehat x \big\rVert<\varepsilon$, it follows that $H - (F^\theta)^{t-1}G^\theta x - (F^\theta)^t\widehat x \subseteq \mathscr{J}^\theta \subseteq G^\theta(E) \oplus F^\theta G^\theta(E) \oplus (F^\theta)^{t-2}G(E)$, and so we obtain that $\smash{\PP_{\vartheta,(x,\widehat x)}(\zeta \in H - (F^\theta)^{t-1}G^\theta x - (F^\theta)^t \widehat x) > 0}$. As we can repeat this construction for any subset of $H$ with positive $\lambda_{\mathscr{R}^\theta}$-measure, this proves that $\auxX$ restricted to $E \times \mathscr{R}^\theta$ is irreducible with respect to $\varphi \coloneq \lambda_E \otimes \lambda_H$, where $\lambda_H$ is the restriction of $\lambda_{\mathscr{R}^\theta}$ to $H$. By $\psi$ we denote a corresponding maximal irreducibility measure for $\auxX$. It follows that the shifted process $\auxaX = \auxX + a^\theta_\infty$ is irreducible as well if the state space is restricted to $E \times \mathscr{R}^\theta_a$.

    In a next step we show that $\auxX$ restricted to $E \times \mathscr{R}^\theta$ is aperiodic. Suppose there exists a $k$-cycle for $\auxX$ with $k \geq 2$, i.e.\ disjoint $D_1, \dots, D_k \in \mathscr{B}(E \times \mathscr{R}^\theta)$ with $\smash{\PP_{\vartheta,(x,\widehat x)}(\auxX(1) \in D_{i+1}) = 1}$ for all $(x,\widehat x) \in D_i$ and $i = 0, \dots, k-1$ (mod $k$) and $\psi\big[(\bigcup_{i=1}^k D_i)^c\big] = 0$. Since $X$ is aperiodic, the (up to $\psi$-null sets) only possible choices for $D_i$ are of the form $D_i = E \times D'_i$ for some $D'_i \subseteq \mathscr R^\theta$, where $D_2' = F^\theta(D_1') \oplus G^\theta(E)$, $\dots$, $\smash{D_k' = (F^\theta)^{k-1}(D_1') \oplus \bigoplus_{i=0}^{k-2} (F^\theta)^{i}G^\theta(E)}$. Then
    \[D_1' = (F^\theta)^{nk}(D_1') \oplus \bigoplus_{i=0}^{nk-1} (F^\theta)^{i}G^\theta(E) \quad \text{ and } \quad D_2' = (F^\theta)^{nk + 1}(D_1') \oplus \bigoplus_{i=0}^{nk} (F^\theta)^{i}G^\theta(E)\]
    for any $n \in \N^*$. Fix any $z \in D_1'$. It follows that $(F^\theta)^{nk}z + \mathscr{J}^\theta \subseteq D_1'$ and $(F^\theta)^{nk + 1}z + \mathscr{J}^\theta \subseteq D_2'$, where $\mathscr{J}^\theta$ is given in \eqref{eq: J}. Choose $n$ large enough so that $\lVert (F^\theta)^{nk}z\rVert < \varepsilon$ and $\lVert (F^\theta)^{nk + 1}z\rVert < \varepsilon$, which is possible because $\rho(F^\theta) < 1$. Since the set $H \subseteq \mathscr{J}^\theta$ is chosen in such a way that any $x \in \mathscr{R}^\theta$ with $\lVert x - H\rVert < \varepsilon$ is still an element of $\mathscr{J}^\theta$, we obtain $H \subseteq D_1'$ and $H \subseteq D_2'$. But $\psi(H) > 0$, which is a contradiction to the disjointness of the sets $D_1'$ and $D_2'$ up to $\psi$-null sets. It follows that the process $\auxX$ is aperiodic and hence also the shifted process $\auxaX$.
    
    The next step of the proof consists in showing that $\auxaX$ has uniformly bounded moments of order $4+\delta$. Since $\lVert a + b\rVert^p \leq 2^{p-1}(\norm{a}^p + \norm{b}^p)$ for $p \in [1,\infty)$, one obtains that
    \[\E_{\vartheta}\big[\lVert \auxaX(t)\rVert^p\big]\leq 2^{p - 1} \E_{\vartheta}\big[\lVert X(t)\rVert^p\big] + 2^{p - 1} \E_{\vartheta}\big[\lVert \widehat X^{\theta,\mathrm{aux},a}(t)\rVert^p\big]\]
    for any $p \in [1,\infty)$. Since $X$ has uniformly bounded moments of order $4+\delta$, the claim follows by establishing that $\widehat X^{\theta,\mathrm{aux},a}$ has uniformly bounded moments of order $4+\delta$. This however holds by Lemma \ref{lem: Hamilton}. By the same argument, the processes \mbox{$\smash{\primeX}$ and $\smash{\langhomX}$ are also bounded in $L^{4 + \delta}$.} 
    
    Since $\mathrm{supp}(\psi) \supseteq H$ has non-empty interior, all conditions from Theorem \ref{theo: f-ergod} are fulfilled for $\auxaX$ restricted to the space $E \times \mathscr{R}^\theta_a$. It follows that $\auxaX$ restricted to $E \times \mathscr{R}^\theta_a$ is strongly ergodic under $\PP_{\vartheta}$, i.e. there exists a unique invariant probability measure $\primemu$ on $E \times \mathscr{R}^\theta_a$ for $\auxaX$ and $\smash{\PP_{\vartheta}^{\widetilde X^{\theta,\mathrm{aux},a}(t)} \xrightarrow{w} \primemu}$ as well as $\frac{1}{t}\sum_{s=1}^t f(\auxaX(s)) \to \int f \dd \primemu$ $\PP_{\vartheta}$-almost surely as long as $\widehat x \in \mathscr{R}^\theta$, and that moreover $\E_{\vartheta,(x,\widehat x)}\big[f(\auxaX(t))\big] \to \int f \dd \primemu$ for $\lambda_E$-almost all $x \in E$ and all $\primemu$-integrable functions $f$. Since $\auxaX$ has uniformly bounded moments of order $4+\delta$, the final part of Theorem \ref{theo: f-ergod} yields that $\primemu$ has finite moments of order $4 + \delta$. The transition dynamics of $\smash{\widehat X_\theta^{\text{aux,$a$}}}$ are given by 
    \begin{align*}
    \smash{\widehat X_\theta^{\text{aux,$a$}}}(t) &= G^\theta X(t-1) + F^\theta(\smash{\widehat X_\theta^{\text{aux,$a$}}}(t-1) - a^\theta_\infty) + a^\theta_\infty \\
    &= G^\theta X(t-1) + F^\theta\smash{\widehat X_\theta^{\text{aux,$a$}}}(t-1) + (\mathrm{I}_d - F^\theta)a^\theta_\infty \\
    &= G^\theta X(t-1) + F^\theta\smash{\widehat X_\theta^{\text{aux,$a$}}}(t-1) + a^\theta
    \end{align*} 
    because $a^\theta_\infty = (\mathrm{I}_d - F^\theta)^{-1} a^\theta$. Hence $\auxaX$ has the same transition dynamics as $\primeX$ and so $\primemu$ is also a unique invariant probability measure for $\primeX$. Finally, we enlarge the state space from $E \times \mathscr{R}_a^\theta$ back to $E \times \R^d$ and consider again the original $E \times \R^d$-valued process $\primeX$. In line with the above, we consider initial distributions $\primeX(0) \sim \nu_\vartheta \otimes \varepsilon(\E_\theta(X(0)))$ for arbitrary distributions $\nu_\vartheta$ on $E$ and $\auxaX(0) \sim \nu_\vartheta \otimes \varepsilon(\widehat x + a^\theta_\infty)$ for $\widehat x \in \mathscr{R}^\theta$, signified by the notation $\PP_{\vartheta, (\nu_\vartheta, \widehat x)}$. Then $\lVert \primeX(t) - \auxaX(t) \rVert$ forms a deterministic sequence converging to 0 under any $\PP_{\vartheta, (\nu_\vartheta, \widehat x)}$. Since $\primeX$ is bounded in $L^{4 + \delta}$, it follows that $f(\primeX(t)) - f(\auxaX(t)) \to 0$ in $L^1$ under any $\PP_{\vartheta, (\nu_\vartheta, \widehat x)}$ and for any polynomial $f$ of order 4 or less by Lemma \ref{lem: Conv_of_Polynomials}. This finally yields for any such polynomial $f$ that $\E_{\vartheta, (x, \widehat x)}\big[f(\primeX(t))\big] \to \int f \dd \primemu$ for $\lambda_E$-almost any $x \in E$ and any $\widehat x \in \R^d$ as well as $\smash{\frac{1}{t}\sum_{s=1}^t f(\primeX(s)) \xrightarrow{L^1} \int f \dd \primemu}$ under $\PP_{\vartheta}$. In particular, we have shown that $\primeX$ is weakly $f$-ergodic with respect to $\primemu$ for any such function $f$. The final part of the theorem for $\primeX$, i.e.\ the fact that convergence of expectations of quadratic (quartic) polynomials holds even for any starting distribution and at a geometric rate, follows then by the same argument as in the proof of Proposition \ref{prop: ergod} because $\primeX$ is a polynomial state space model of order 2 (resp. 4) as $X$ is a state space model of order 2 (resp. 4).
    
    The treatment of the $d(k+2)$-dimensional process $\langhomX$ is analogous. This can be seen by setting $Y^{\theta,\mathrm{hom}} \coloneq (\homkX,\homX)^\top$. From \eqref{eq: transitiontilde}, the dynamics of $\langhomX$ are
    \[\begin{pmatrix}
        X(t) \\ Y^{\theta,\mathrm{hom}}(t)
    \end{pmatrix} = \wa + \begin{pmatrix}
        A^\vartheta & 0 \\ \overline{G}^\theta & \overline F^\theta
    \end{pmatrix}\begin{pmatrix}
        X(t-1) \\ Y^{\theta,\mathrm{hom}}(t-1)
    \end{pmatrix} + \smash{\overline N}^\vartheta(t) \]
    with some matrices $\overline G^\theta \in \R^{d(k+1) \times d}$ and $\overline F^\theta \in \R^{d(k+1) \times d(k+1)}$ satisfying $\rho(\overline F^\theta) < 1$ because $\rho(F^\theta) < 1$ and $\smash{\overline F^\theta}$ is a block-lower triangular matrix whose submatrices on the main diagonal are $F^\theta$, see \eqref{eq: transitiontilde}. Hence $\smash{\langhomX}$ has exactly the same structure as $\smash{\primeX}$ when $\smash{\homkX}$ is replaced by $\smash{Y^{\theta,\mathrm{hom}}}$, and the same proof as above can be applied to $\smash{\langhomX}$.
\end{proof}

The preceding results from Lemma \ref{limitcov}, Corollary \ref{convcoro}, and Proposition \ref{ergod2} are intuitive and constitute to some extent well-known facts in the filtering and engineering literature subsumed under the name of a \textit{steady-state Kálmán filter}, a term which dates over 50 years back to \cite{ODonnel1966}.

\begin{remark}\label{rem: Radius_for_homX}
    By exactly the same argument as in the proof of Proposition \ref{prop: ergod}, we can choose the state transition matrix $\homAtwo$ in such a way that $\rho(\homAtwo) < 1$, where $\homAtwo$ denotes the state transition matrix for $\mathrm{vec}_{\otimes 2}(\langhomX)$ under $\PP_\vartheta$, see Remark \ref{rem: Mean_Recursion}.
\end{remark}

\begin{proposition}\label{prop: asymptotic_equiv}
    $\langX(t) - \langhomX(t) \xrightarrow{t \to \infty} 0$ in $L^{4+\delta}(\PP_\vartheta)$ at a geometric rate.
\end{proposition}
\begin{proof}
    As before let $\wa$ denote the state transition vector of $\langX$ and $\langhomX$, let $\wA(t)$ and $\homA$ denote the state transition matrices of $\langX$ and $\langhomX$, respectively, and let $\smash{\overline N}^\vartheta$ denote the martingale difference sequence in the state space representation of $\langX$ and $\langhomX$ (all under $\PP_\vartheta$). Then $\langX(t) = \wa + \wA(t) \langX(t-1)+ \smash{\overline N}^\vartheta(t)$ and $\langhomX(t) = \wa + \homA \langhomX(t-1) + \smash{\overline N}^\vartheta(t)$ as well as $\langX(0) = \langhomX(0)$. Moreover $\wA(t) \to \homA$ as $t \to \infty$ at a geometric rate and $\rho(\homA) < 1$ by Corollary \ref{convcoro}. By iterating and subtracting these equations, one obtains
    \begin{fitalign}
        \langX(t) - \langhomX(t) = \underbrace{\sum_{s=1}^t \Big[\prod_{r=s+1}^t \wA(r) - (\homA)^{t-s}\Big]\wa}_{(1)} &+ \underbrace{\sum_{s=1}^t \Big[\prod_{r=s+1}^t \wA(r) - (\homA)^{t-s}\Big]  \smash{\overline N}^\vartheta(s)}_{(2)} \\ &+ \underbrace{\Big(\prod_{s=1}^t \wA(s) - (\homA)^t\Big) \langX(0)}_{(3)}.
    \end{fitalign}
    
    \noindent By the same argument as in the proof of Lemma \ref{lem: lin_systems}.3, the deterministic term (1) converges to 0 uniformly in $\theta$ at a geometric rate and (3) converges to 0 in $L^{4 + \delta}$ at a geometric rate. So it only remains to focus on the term (2). By equation \eqref{eq: bound_for_product} there are a constant $c \in \R_+$ and $\gamma \in [0, 1)$ such that $\big\lVert \prod_{r=s+1}^t \wA(r) - (\homA)^{t-s}\big\rVert \leq c \gamma^t.$ Hence
    \begin{align*}
        \lVert (2) \rVert_{L^{4 + \delta}} \leq c \gamma^t \sum_{s=1}^t \lVert \smash{\overline N}^\vartheta(s)\rVert_{L^{4 + \delta}} \leq c \gamma^t
    \end{align*}
    (with varying constants $c,\gamma$ from one expression to the next) because $\smash{\overline N}^\vartheta$ is bounded in $L^{4+\delta}$ as it consists only of $N^\vartheta$ and zero components, and because $N^\vartheta$ is bounded in $L^{4 + \delta}$ as $X$ is bounded in $L^{4 + \delta}$. It follows that $(2) \xrightarrow{t \to \infty} 0$ in $L^{4+\delta}$ at a geometric rate and the same then also holds for $\langX(t) - \langhomX(t)$.
\end{proof}

Using the preceding Proposition \ref{prop: asymptotic_equiv}, we can transfer the ergodicity results in Proposition \ref{ergod2} from the homogeneous Markov process $\langhomX$ to the inhomogeneous process $\langX$:

\begin{corollary}\label{coro: ergod_wX}
    For any $\theta\in \Theta$ the process $\langX$ is bounded in $L^{4+\delta}(\PP_\vartheta)$ and weakly $f$-ergodic with respect to $\wmu$ for any polynomial $f: E \times \R^{d(k+1)} \to \R$ of degree 4 or less. In particular, $\PP_\vartheta^{\langX(t)} \xrightarrow{w} \wmu$ and $\frac{1}{t}\sum_{s=1}^t f(\langX(s)) \to \int f \dd \wmu$ in $L^1(\PP_\vartheta)$. Moreover $\E_\vartheta(f(\langX(t))) \to \int f \mathrm{d}\wmu$ at a geometric rate if $f$ is a quadratic polynomial (resp. quartic polynomial if Assumption \ref{assump: AN} holds).
    \end{corollary}
\begin{proof}
    Boundedness in $L^{4+\delta}$ follows from the preceding Proposition \ref{prop: asymptotic_equiv} because $\langhomX$ is boun\-ded in $L^{4+\delta}$ by Proposition \ref{ergod2}. Since the law of $\langhomX(t)$ converges weakly to $\wmu$, the law of $\langX$ does so too by Slutsky's theorem. Finally, by Lemma \ref{lem: Conv_of_Polynomials}, $f(\langX(t)) - f(\langhomX(t)) \to 0$ in $L^1$, so the convergence of arithmetic means and expectations follows from the respective convergence of $\langhomX$ from Proposition \ref{ergod2}. Likewise, since $\langX(t) - \langhomX(t) \to 0$ in $L^{4+\delta}$ converges even at a geometric rate, we obtain from Lemma \ref{lem: Conv_of_Polynomials} that expectations of quartic functions of the process $\smash{\langX(t)}$ converge at a geometric rate because the same holds for the process $\langhomX$ by Proposition \ref{ergod2}.
\end{proof}

In the following we let $\tildeP$ denote the one-step transition measures for the Markov chain $\langhomX$ under $\PP_\vartheta$ so that $\smash{\tildeP g(x) \coloneq \E_\vartheta\big[g(\langhomX(t)) \mid \langhomX(t-1) = x\big]}$. The last ingredient on the way to a functional central limit theorem for the quasi-score process is the fact that every quadratic polynomial $f: \R^{d(k+2)} \to \R$ gives rise to an explicitly solvable Poisson-type equation involving the operator $\tildeP$. This is the content of the following lemma:

\begin{lemma}\label{range}
    Every quadratic $f: \R^{d(k+2)} \to \R$ can be represented as $f(x) = \tildeP g(x) - g(x) + r$ with some quadratic polynomial $g: \R^{d(k+2)} \to \R$ and some $r \in \R$.
\end{lemma}
\begin{proof}
    Since $f$ is a quadratic polynomial, it can be written in the form $f(x) = \alpha_f^\top \mathrm{vec}_{\otimes 2}(x) + \beta_f$, where $\alpha_f \in \R^{d(k+2) + d^2(k+2)^2}$ is the vector of coefficients of the linear and quadratic terms of $f$ and $\beta_f \in \R$. Suppose now that $g$ is quadratic without a constant, i.e.\ of the form $g(x) = \alpha_g^\top \mathrm{vec}_{\otimes 2}(x)$. Since $\langhomX$ is a state space model of order 2, we have $\tildeP \mathrm{vec}_{\otimes 2}(x) = \smash{\overline a}^\theta_{\otimes 2} + \homAtwo \mathrm{vec}_{\otimes 2}(x)$, where $\smash{\overline a}^\theta_{\otimes 2}$ denotes the state transition vector for $\mathrm{vec}_{\otimes 2}(\langhomX)$. Hence
    \begin{align*}
        \tildeP g(x) - g(x) &= \alpha_g^\top \smash{\overline a}^\theta_{\otimes 2} + \alpha_g^\top \homAtwo \mathrm{vec}_{\otimes 2}(x) - \alpha_g^\top \mathrm{vec}_{\otimes 2}(x) \\
        &= \alpha_g^\top \smash{\overline a}^\theta_{\otimes 2} + \alpha_g^\top \Big[\homAtwo - \mathrm{I}\Big] \mathrm{vec}_{\otimes 2}(x),
    \end{align*}
    
    \noindent where $\mathrm{I}$ denotes the $d(k+2) + d^2(k+2)^2$-dimensional identity matrix. Now set $r \coloneq \beta_f -  \alpha_g^\top \smash{\overline a}^\theta_{\otimes 2}$. Then the Poisson equation $f(x) = \tildeP g(x) - g(x) + r$ is fulfilled whenever $\alpha_g^\top (\homAtwo - \mathrm{I}) = \alpha_f^\top$. This equation has the unique solution $\alpha_g^\top = \alpha_f^\top (\homAtwo - \mathrm{I})^{-1}$, where the inverse is well-defined because $\rho(\homAtwo) < 1$ by Remark \ref{rem: Radius_for_homX}. This finishes the proof.
\end{proof}

The preceding lemma is the main ingredient of the following proof, which establishes the asymptotic normality condition for the quasi-score process from Proposition \ref{prop: Jacod-Normality}. The strategy of the proof consists of recognising that $Z^\theta(t, t-1)$ is a quadratic polynomial in $\langX(t)$ with coefficients depending on $t$, which however converge as $t \to \infty$ to some limiting coefficients. Hence we can approximate this polynomial by a time-invariant polynomial using the limiting coefficients, and we can moreover suitably approximate $\langX$ by $\langhomX$.

\begin{theorem}\label{theo: norm}
There is $U_\vartheta \geq 0 \in \R^{k \times k}$ with $\smash{\frac{1}{\sqrt{t}} Z^\vartheta(t) \xrightarrow{\PP_\vartheta\text{-}d} Z \sim N(0, U_\vartheta)}$.
\end{theorem}
\begin{proof}
We have $\smash{\frac{1}{\sqrt{t}} Z^\vartheta(t) = \frac{1}{\sqrt{t}} \sum_{s=1}^t Z^\vartheta(s, s-1)}$, where each component $\smash{Z^\vartheta(t, t-1)_j}$ for $j = 1, \dots, k$ is a quadratic polynomial with time-varying coefficients in $\langoX(t)$ with 
\[Z^\vartheta(t, t-1)_j = \sum_{|\lambda|\leq 2} \alphajx(t) \langoX(t)^\lambda \eqqcolon f^\vartheta_{j}(t, \langoX(t)),\] 
the sum referring to $d(k+2)$-dimensional multi-indices $\lambda$. 
By the explicit representation in equation \eqref{eq: Score} the coefficients $\alphajx(t)$ depend only on the matrices $\widehat \Sigma^\vartheta(t, t-1)$ and $S^{\vartheta, j}(t, t-1)$, which by Lemma \ref{limitcov} converge as $t \to \infty$ at a geometric rate. It follows quickly that there exist limiting coefficients $\alphaj$ such that $\alphajx(t) \to \alphaj$ at a geometric rate. Let $f^\vartheta_j(y) \coloneq \sum_{|\lambda|\leq 2} \alphaj y^\lambda$ denote the quadratic polynomial with these limiting coefficients, see equation \eqref{eq: f_def}. We use the following decomposition of $\frac{1}{\sqrt{t}}Z^\vartheta(t)$:
    \begin{fitequation*}
        \frac{1}{\sqrt{t}}Z^\vartheta(t) = \underbrace{\frac{1}{\sqrt{t}}\sum_{s=1}^t f^\vartheta( \langhomoX(s))}_{(1)} + \underbrace{\frac{1}{\sqrt{t}}\sum_{s=1}^t \Big(\sum_{|\lambda| \leq 2} \alphaall \big[\langoX(s)^\lambda - \langhomoX(s)^\lambda\big]\Big)}_{(2)} + \underbrace{\frac{1}{\sqrt{t}}\sum_{s=1}^t \Big(\sum_{|\lambda| \leq 2} \big[\alphaallx(s) - \alphaall\big] \langoX(s)^\lambda\Big)}_{(3)},
    \end{fitequation*}
    where $f^\vartheta = (f^\vartheta_1, \dots, f^\vartheta_k)^\top$ and $\smash{\alphaall = (\alpha^{\vartheta}_{\lambda, 1}, \dots, \alpha^\vartheta_{\lambda, k})^\top}$. We show that (1) converges in distribution to a Gaussian random variable $Z$, while (2) and (3) converge in probability to 0. The result follows then from Slutsky's theorem. To treat (1), note that by Corollary \ref{coro: EZ=0} we have $\E_\vartheta\big[Z^\vartheta(t)\big] = \sum_{s=1}^t \E_\vartheta\big[Z^\vartheta(s, s-1)\big] = 0$ for any $t \in \N^*$, so also $\E_\vartheta\big[Z^\vartheta(t, t-1)\big] = 0$ for any $t \in \N^*$. By Corollary \ref{coro: ergod_wX}, this implies that 
    \begin{align*}
        0 = \E_\vartheta\big[Z^\vartheta(t, t-1)\big] = \sum_{|\lambda| \leq 2} \alphaallx(t) \E_\vartheta\big[\langoX(t)^\lambda\big] \xrightarrow{t \to \infty} \sum_{|\lambda \leq 2} \alphaall \int x^\lambda \wmuo(\mathrm{d}x) = \int f^\vartheta \dd \wmuo,
    \end{align*}
    so $\int f^\vartheta \dd \wmuo = 0$. By applying Lemma \ref{range} componentwise, there exists a $k$-dimensional quadratic polynomial $g^{\vartheta}: \R^{d(k+2)} \to \R^k$ such that $f^\vartheta = \tildePtwo g^\vartheta - g^\vartheta + r^\vartheta$ for some constant $r^\vartheta \in \R^k$. Since $\wmuo$ is an invariant probability measure for $\langhomoX$ and $\tildePtwo$ is the transition operator of $\langhomoX$, it follows from the definition of an invariant probability measure that $\smash{\int (\tildePtwo g^\vartheta - g^\vartheta) \dd \wmuo = 0}$. Hence we conveniently obtain $\smash{r^\vartheta = 0}$. Now define the triangular array $\smash{(U_t^{\vartheta}(s))_{t \in \N^*, s \in \{1, \dots, t\}}}$ of random variables by
    \[U_t^{\vartheta}(s) \coloneq \frac{1}{\sqrt{t}}\Big[\tildePtwo g^\vartheta(\langhomoX(s-1)) - g^\vartheta(\langhomoX(s))\Big].\]
    Then $\smash{\E_\vartheta(U_t^{\vartheta}(s) \mid \F_{s-1}) = 0}$ for each $t \in \N^*$ and $s \in \{1, \dots, t\}$, so the process $M^\vartheta$ with
    \begin{align*}
    M^\vartheta(t) \coloneq \sum_{1 \leq s \leq t} U_t^\vartheta(s) = \frac{1}{\sqrt{t}}\Big[g^\vartheta(\langhomoX(0)) - g^\vartheta(\langhomoX(t)) + \sum_{0 \leq s \leq t - 1} f^\vartheta(\langhomoX(s))\Big]
    \end{align*}
    is a $k$-dimensional martingale with respect to $(\F_t)_{t \in \N^*}$. Since $g^\vartheta(\langhomoX(t))$ is bounded in $L^1$, it follows that $\smash{\frac{1}{\sqrt{t}}}\big[g^\vartheta(\langhomoX(0)) - g^\vartheta(\langhomoX(t))\big] \to 0$ in $L^1$. Hence it suffices to prove asymptotic normality of $M^\vartheta$ in order to show asymptotic normality of (1). To do so, we establish the conditions of  the martingale central limit theorem \ref{theo: MCLT},  which tells us that the asymptotic covariance matrix of (1) is given by $U_\vartheta \coloneq \lim_{t\to \infty} \sum_{s=1}^t U^\vartheta_t(s) U^{\vartheta^\top}_t(s)$, where the limit is in $\PP_\vartheta$-probability if it exists. The sum $\sum_{s=1}^t U^\vartheta_t(s) U^{\vartheta^\top}_t(s)$ is given by
    \begin{equation*}
        \frac{1}{t} \sum_{s=1}^t \Big[\tildePtwo g^\vartheta(\langhomoX(s-1)) - g^\vartheta(\langhomoX(s))\Big]\Big[\tildePtwo g^\vartheta(\langhomoX(s-1)) - g^\vartheta(\langhomoX(s))\Big]^\top,
    \end{equation*}
    where each summand is an $\R^{k \times k}$-valued random variable whose entries are quartic polynomials in $\langhomoX (s-1)$ and $\langhomoX (s)$. Since $\langhomoX $ is weakly $f$-ergodic for any quartic polynomial $f$ by Proposition \ref{ergod2}, it is clear that the same holds also true for $(\langhomoX (t), \langhomoX (t-1))_{t \in \N^*}$. Thus the limit $U_\vartheta$ in $\PP_\vartheta$-probability is indeed well-defined and exists. It remains to prove the Lindeberg condition 2. from the martingale central limit Theorem \ref{theo: MCLT} to show asymptotic normality of (1). Since $\big[\tildePtwo g^\vartheta(\langhomoX (s-1)) - g^\vartheta(\langhomoX (s))\big]$ is a $k$-dimensional quadratic polynomial in $\langhomoX (s-1)$ and $\langhomoX (s)$, it is bounded in $L^{2 + \frac{\delta}{2}}$ because $\langhomoX $ is bounded in $L^{4 + \delta}$ by Lemma \ref{lem: Boundedness_of_Polynomials}. Hence there exists some $B \geq 0$ such that $\E_\vartheta\big[\lVert U_t^\vartheta(s) \rVert^{2 + \eta} \big] \leq Bt^{-(1+\frac{\eta}{2})}$ for all $s \in \{1, \dots, t\}$, where $\eta \coloneq \frac{\delta}{2}$. Moreover, $\lVert U_t^\vartheta(s)\rVert^2 \mathbf{1}_{\{\lVert U_t^\vartheta(s) \rVert > \varepsilon\}} \leq \varepsilon^{-\eta} \lVert U_t^\vartheta(s) \Vert^{2 + \eta}$, so
    \begin{equation*}
        \sum_{s=1}^t \E_\vartheta \Big[\lVert U_t^\vartheta(s)\rVert^2 \mathbf{1}_{\{\lVert U_t^\vartheta(s) \rVert > \varepsilon\}}\Big] \leq \varepsilon^{-\eta}\sum_{s=1}^t \E_\vartheta\Big[\lVert U_t^\vartheta(s) \Vert^{2 + \eta}\Big] \leq B \varepsilon^{-\eta} t^{-\frac{\eta}{2}} \xrightarrow{t \to \infty} 0.
    \end{equation*}
    This shows that the Lindeberg condition from Theorem \ref{theo: MCLT} holds with $L^1$ convergence, which establishes that $(1) \xrightarrow{\PP_\vartheta-d} Z \sim N(0, U_\vartheta)$. To show that (2) converges in probability to 0, note that $\sum_{|\lambda| \leq 2} \alphaall \big[\langoX(t)^\lambda - \langhomoX (t)^\lambda\big] = f^\vartheta(\langoX(t)) - f^\vartheta(\langhomoX (t))$ converges in $L^1$ to 0 at a geometric rate by Proposition \ref{prop: asymptotic_equiv} and Lemma \ref{lem: Conv_of_Polynomials}, so $\sum_{s=1}^t \E_\vartheta\big[ \lVert f^\vartheta(\langoX(s)) - f^\vartheta(\langhomoX (s)) \rVert \big]$ converges as $t \to \infty$ and thus $(2) \to 0$ in $L^1$. Finally, we prove that (3) converges to 0 in $L^1$. Since all absolute first and second powers of $\langoX$ are bounded in expectation by some constant $C \in \R_+$ (see Corollary \ref{coro: ergod_wX}), one has
    \begin{align*}
        \lVert (3) \rVert_{L^1} \leq \frac{C}{\sqrt{t}}\sum_{s=1}^t \sum_{|\lambda| \leq 2} \big\lVert \alphaallx(s) - \alphaall \big\rVert.
    \end{align*}
    But $\alphaallx(t) \to \alphaall$ at a geometric rate, so the sum $\sum_{s=1}^t \sum_{|\lambda| \leq 2} \big\lVert \alphaallx(s) - \alphaall \big\rVert$ converges and therefore (3) converges to 0 as $t \to \infty$ in $L^1$, finishing the proof.
\end{proof}

\begin{corollary}\label{coro: limiting_U}
    Let $g^\vartheta$ be the solution to the Poisson equation $\smash{f^\vartheta = \tildePtwo g^\vartheta - g^\vartheta}$, where $f^\vartheta$ is as in equation \eqref{eq: f_def}. Then $U_\vartheta \in \R^{k \times k}$ from Theorem \ref{theo: norm} is given by \begin{align}\label{eq: limiting_U}
    U_\vartheta &= \PP_\vartheta\text{-}\lim_{t\to\infty} \frac{1}{t}\sum_{s=1}^t \Big[g^\vartheta(\langoX(s))g^\vartheta(\langoX(s))^\top - \tildePtwo g^\vartheta(\langoX(s)) \tildePtwo g^\vartheta(\langoX(s))^\top\Big].
    \end{align}
    Additionally, if Assumption \ref{assump: AN} holds, the matrix $U_\vartheta$ can be expressed as
    \begin{equation}\label{eq: limiting_U_exp}
    U_\vartheta = \lim_{t \to \infty}\E_\vartheta\big[g^\vartheta(\langoX(t))g^\vartheta(\langoX(t))^\top - \tildePtwo g^\vartheta(\langoX(t)) \tildePtwo g^\vartheta(\langoX(t))^\top\big].
    \end{equation}
\end{corollary}
\begin{proof}
    By the reasoning from the proof of Theorem \ref{theo: norm}, the given Poisson equation has a solution $g^\vartheta$ within the $k$-dimensional quadratic polynomials and $U_\vartheta$ is given by
    \[U_\vartheta = \E_{\wmuo} \Big(\Big[\tildePtwo g^\vartheta(\langhomoX (0)) - g^\vartheta(\langhomoX (1))\Big]\Big[\tildePtwo g^\vartheta(\langhomoX (0)) - g^\vartheta(\langhomoX (1))\Big]^\top\Big),\]
    where $\E_{\wmuo}$ denotes the expectation under the law on $(\Omega, \mathscr{F})$ for which $\langhomoX (0) \sim \wmuo$. But $$\E_{\wmuo}\big[\tildePtwo g^\vartheta(\langhomoX (0))g^\vartheta(\langhomoX (1))^\top \big] = \E_{\wmuo}\big[\tildePtwo g^\vartheta(\langhomoX (0))\tildePtwo g^\vartheta(\langhomoX (0))^\top\big],$$ so
    \begin{align}\label{e:Uexplizit}
    U_\vartheta = \E_{\wmuo}\Big[ g^\vartheta&(\langhomoX (1))g^\vartheta(\langhomoX (1))^\top \Big] - \E_{\wmuo}\Big[\tildePtwo g^\vartheta(\langhomoX (0))\tildePtwo g^\vartheta(\langhomoX (0))^\top \Big]  \nonumber \\
        &= \int g^\vartheta(x) g^\vartheta(x)^\top \wmuo(\mathrm{d}x) - \int \tildePtwo g^\vartheta(x)\tildePtwo g^\vartheta(x)^\top \wmuo(\mathrm{d}x).
    \end{align} 
    The latter term is the ($\PP_\vartheta$-)limit of the claimed expressions by Corollary \ref{coro: ergod_wX}.
\end{proof}
Corollary \ref{coro: limiting_U} tells us that once an expression for $g^\vartheta$ and $\tildePtwo g^\vartheta$ is found, the asymptotic covariance matrix $U_\vartheta$ from Theorem \ref{theo: norm} can be easily estimated by \eqref{eq: limiting_U} or, if one works with a polynomial state space model of order 4, explicitly calculated using \eqref{eq: limiting_U_exp}. This is done in Sections \ref{sec5.1: CovMat} and \ref{su:explicit}. Finding an analytical expression for $g^\vartheta$ and $h^\vartheta \coloneq \tildePtwo g^\vartheta$, which is the content of Proposition \ref{prop: Formula_Cov}, reduces to knowing the quadratic polynomial coefficients of the conditional quasi-score and to a simple matrix inversion, see the construction in the proof of Lemma \ref{range}. This is in sharp contrast to the complicated methods needed for estimation of $U_\vartheta$ in non-polynomial VARMA models, see for example \cite{Mainassara2011}, Section 4, \cite{Franq2005}, Section 5.2, or \cite{Schlemm2012}, p. 2197. There, an estimator of $U_\vartheta$ is derived from the spectral density of $Z^\vartheta(t, t-1)$ under the assumption that $X$ has $(8 + \delta)$th moments and that the conditional score possesses an infinite-order autoregressive representation, which is then estimated by a least-squares autoregression. In addition, these methods are only applicable for estimating $U_\theta$ at $\theta = \vartheta$ which can turn out to be insufficient, for example in the context of hypothesis testing using the Lagrange multiplier test, see Section \ref{su:tests}. This is remedied by a different estimator in \cite{Mainassara2014} and \cite{Mainassara2012}.

Theorem \ref{theo: norm} establishes the asymptotic normality condition needed in Proposition \ref{prop: Jacod-Normality}, which is taken from \cite{Jacod2018}. For our consistency proof of the quasi-maximum likelihood estimator for the polynomial state space model, we now focus our attention on the uniform convergence condition from Proposition \ref{prop: Jacod-Normality}. This is the goal of the next section.

\subsubsection{A law of large numbers for the observed Fisher information}\label{sec4.2.2: ULLN}

We now turn to establishing the uniform law of large numbers-type result for $\frac{1}{t} \nabla_\theta Z^\theta(t)$ from Proposition \ref{prop: Jacod-Normality}, where $\nabla_\theta Z^\theta(t, t-1)$ is the $(k\times k)$-matrix whose $ij$-th entry is $\partial_{\theta_i} Z^\theta(t)_j$, i.e.\ the Hessian matrix of the quasi-log-likelihood $L^\theta(t) = \log q_t^\theta(X_\oo(1), \dots, X_\oo(t))$ with respect to $\theta$. Due to the obvious similarities to the observed Fisher information matrix from classic statistical theory, we will generally refer to $\frac{1}{t} \nabla_\theta Z^\theta(t)$ as the observed Fisher information. 

Up to a handful of arguments, the convergence proof for the observed Fisher information can be handled by exactly the same techniques as established in the preceding section. After representing $\nabla_\theta Z^\theta(t,t-1)$ explicitly, we proceed by augmenting the process $\langX(t)$ with second derivatives of the Kálmán filter $\widehat X^\theta(t, t-1)$ with respect to $\theta$ and then establishing basic ergodicity properties of this augmented process. Proposition \ref{prop: Fisher_Information} provides such a closed form for the observed Fisher information.

\begin{proof}[Proof of Proposition \ref{prop: Fisher_Information}:]
    The result follows by differentiating the $j$-th component of the quasi-score \eqref{eq: Score} with respect to $\theta_i$, similar to how it was done in Proposition \ref{prop: Quasi-Score}, and using the identities from \cite{Petersen2012}.
\end{proof}

The specific lengthy form of the observed Fisher information in \eqref{eq: Fisher-Information} is not of crucial importance. The main observation here is that its structure is similar to that of the quasi-score from equation \eqref{eq: Score}. In particular, $\nabla_\theta Z^{\theta}(t,t-1)$ is a $(k \times k)$-dimensional quadratic polynomial in $\laengerX(t)$ with time-varying coefficients, where $\laengerX$ is the $E \times \R^{d(k^2+k+1)}$-valued process $\laengerX(t) \coloneq (X(t),\kX(t),\wX(t), W^\theta(t))^\top$, which has the same form as $\langX$ in the preceding section but augmented by all second derivatives of the Kálmán filtered component. Then the dynamics of $\laengerX$ under $\PP_\vartheta$ can be written for $t \in \N^*$ as
\begin{fitequation}\label{eq: Xbar_dynamics}
    \laengerX(t + 1) = \begin{pmatrix}
        \wa \\ \partial_{11} a^\theta \\ \partial_{12} a^\theta \\ \vdots \\ \partial_{kk} a^\theta
    \end{pmatrix} + \begin{pmatrix}
        \wA(t + 1) & 0 & 0 & \dots & 0 \\
        \circ & F^\theta(t) & 0 & \dots & 0 \\
        \circ & 0 & F^\theta(t) & \dots & 0 \\
        \vdots & \vdots & \vdots & \ddots & \vdots \\
        \circ & 0 & 0 &\dots & F^\theta(t)
    \end{pmatrix} \laengerX(t) + \begin{pmatrix}
        \smash{\overline N}^\vartheta(t + 1) \\ 0 \\0 \\ \vdots \\ 0
    \end{pmatrix},
\end{fitequation}
where $\wa$ and $\wA(t)$ are the state transition vector and the time-dependent state transition matrix of $\langX$, where $\smash{\overline N}^\vartheta$ denotes the martingale difference sequence for $\langX$ from \eqref{eq: Xtilde_dynamics}, and where $\circ$ denotes some irrelevant submatrix of suitable dimensionality. For the matrices $R_{ij}^\theta(t, t-1)$ from Proposition \ref{prop: Fisher_Information}, it is not hard to see that by Proposition \ref{prop: B-assump}
\begin{equation}\label{eq: R_dynamics}
    R^\theta_{ij}(t+1,t) = F^\theta(t) R^\theta_{ij}(t,t-1) F^\theta(t)^\top + \text{ terms converging uniformly in $\theta$}.
\end{equation}
Hence, by the same reasoning as in Lemma \ref{limitcov} for the matrices $S^\theta_j(t, t-1)$, there exist matrices $R^{\theta}_{ij}$ for $i, j \in \{1, \dots, k\}$ such that, uniformly in $\theta\in\Theta$, $R^\theta_{ij}(t,t-1) \to R^{\theta}_{ij}$ at a geometric rate. Simple differentiation as in Propositions \ref{prop: Quasi-Score} or \ref{prop: Fisher_Information} and similar algebraic manipulations as in equation \eqref{eq: R_dynamics} show for any $i, j, l \in \{1, \dots, k\}$ we have
\begin{equation}\label{eq: partial_R_dynamics}
\partial_{\theta_l} R^\theta_{ij}(t+1,t) = F^\theta(t) \partial_{\theta_l} R^\theta_{ij}(t,t-1)F^\theta(t)^\top + \text{ terms uniformly bounded in $\theta$}
\end{equation}
because all third derivatives of $C^\theta(t)$ are uniformly bounded in $\theta$ by Assumption \ref{assump: A}. Lemma \ref{lem: lin_systems}.1 then yields that all derivatives $\partial_{\theta_l} R^\theta_{ij}(t,t -1)$ are uniformly bounded in $\theta$ and $t \in \N^*$.

In consequence, we can approximate $\laengerX$ by the time-homogeneous Markov chain $\laengerhomX$ defined by transition dynamics identical to $\laengerX$ but with the limiting matrices $\widehat \Sigma^\theta$, $S^{\theta}_j$ and $R^{\theta}_{ij}$ in place of $\smash{\widehat \Sigma^\theta(t,t-1)}$, $\smash{S^\theta_j(t,t-1)}$ and $\smash{R^\theta_{ij}(t,t-1)}$, respectively, for $i, j \in \{1, \dots, k\}$. Since the structure of the state transition matrix of $\laengerhomX$ is completely analogous to the state transition matrix of $\langhomX$, we can show existence of a stationary distribution $\omu$ for $\laengerhomX$ as well as weak ergodicity and boundedness in $L^{4 + \delta}$ of $\laengerhomX$ by the exact same argument as given in Proposition \ref{ergod2} for $\langhomX$. Since $\laengerhomX$ is a polynomial state space model of order 2 just as $\langhomX$, we have that $\laengerX(t) - \laengerhomX(t) \to 0$ in $L^{4 + \delta}$ under $\PP_\vartheta$ at a geometric rate as in Proposition \ref{prop: asymptotic_equiv}. We obtain the following analogue of Corollary \ref{coro: ergod_wX}:

\begin{corollary}\label{coro: ergod_oX}
     For any $\theta\in \Theta$, the process $\laengerX$ is bounded in $L^{4+\delta}(\PP_\vartheta)$ and weakly $f$-ergodic with respect to $\omu$ for any polynomial $f: E \times \R^{d(k^2+k+1)} \to \R$ of degree 4 or less. In particular, $\smash{\PP_\vartheta^{\laengerX(t)} \xrightarrow{w} \omu}$ as well as $\frac{1}{t}\sum_{s=1}^t f(\laengerX(s)) \to \int f \dd \omu$ in $L^1(\PP_{\vartheta})$ and $\theta \mapsto \int f \dd \omu$ is continuous. Moreover, $\E_\vartheta\big[f(\laengerX(t))\big]\to \int f \dd \omu$ at a geometric rate if $f$ is quadratic (resp. quartic if Assumption \ref{assump: AN} holds).
\end{corollary}
\begin{proof}
    The only new statement that has not been covered before and needs to be verified is the continuity of the function $\theta \mapsto \int f \dd \omu$ for any polynomial $f$ of degree 4. Equivalently, it needs to be proven that the function $\theta \mapsto \int f \dd \omu$ is continuous, where $f: x \mapsto x^{\otimes p}$ for any $p \in \{1, \dots, 4\}$. For such $f$ we first claim that the sequence of expectations $\smash{\E_\vartheta\big[f(\laengerX(t)) \mid X(0) = x\big]}$ is continuous in the parameter $\theta$. In order to see this, let $\watwo$ and $\wAtwo(t)$ denote the state transition vector and time-dependent state transition matrix of $\laengerX$. Moreover, let $\wNtwo$ denote the martingale difference sequence for $\laengerX$ from \eqref{eq: Xbar_dynamics}. If we substitute $\wbtwo(t) \coloneq \watwo + \wNtwo(t)$, iterating the dynamic equation for the process $\laengerX$ yields
    \begin{align*}
        \laengerX(t) = \sum_{s=1}^t \prod_{u = s}^{t-1} \Big(\wAtwo(u)\Big) \wbtwo(s) + \Big(\prod_{s=0}^{t-1} \wAtwo(u)\Big) \laengerX(0).
    \end{align*}
    For simplicity, we will assume that $p = 2$ and that $\laengerX(0) = 0$. The general case of larger $p$ and $\laengerX(0) \neq 0$ can be obtained completely analogously but with higher notational effort. Let $\mu_{s, w}(\theta) \coloneq \E_\vartheta\big[\wbtwo(s) \otimes \wbtwo(w) \mid X(0) = x\big]$, which is continuously differentiable in $\theta$ by thrice continuous differentiability assumed in Assumption \ref{assump: A}. Moreover, since $\wNtwo$ is bounded in $L^{4+\delta}$, $\mu_{s, w}(\theta)$ and its derivatives are uniformly bounded in $s, w$ and $\theta$. Hence
    \begin{align}
        \E_\vartheta\Big[\laengerX(t)^{\otimes 2} \:\big|\: X(0) = x\Big] = \sum_{s=1}^t \sum_{w=1}^t \Big[ \prod_{u=s}^{t-1} \wAtwo(u) \otimes \prod_{u=w}^{t-1}\wAtwo(u)\Big]\mu_{s, w}(\theta).\label{eq: oxkronecker2}
    \end{align}
     Moreover, the proof of Lemma \ref{limitcov} shows that $\wAtwo(t) \to \wAtwo$ uniformly on $\Theta$ for some matrix $\wAtwo$ continuous in $\theta$ with $\rho(\wAtwo) < 1$ because $\wAtwo$ is a block lower-triangular matrix whose submatrices on the main diagonal have spectral radius less than 1. Additionally, the matrices $\partial_{\theta_l} \wAtwo(t)$ are uniformly bounded in $\theta$ and $t$ for any $l \in \{1, \dots, k\}$ because they depend on $t$ only through the matrices $\widehat \Sigma^\theta(t, t-1)$ and its inverse, $S^\theta_j(t,t-1)$, $R^\theta_{ij}(t, t-1)$, and $\partial_{\theta_l}R^\theta_{ij}(t, t-1)$, which are uniformly bounded in $\theta$ and $t$ for any $i, j \in \{1, \dots, k\}$. 
     
     All in all, Lemma \ref{lem: linalg_unif} and Corollary \ref{coro: linalg_diff} can now be applied to show that there exist some constants $c \in \R_+$ and $\gamma \in [0, 1)$ such that for any $s \in \N$ and any $m \in \N^*$ we have $\sup_{\theta \in \Theta}\lVert \wAtwo(s+m) \dots \wAtwo(s+1)\rVert \leq c  \gamma^m$ and $\sup_{\theta \in \Theta}\lVert \partial_{\theta_l}\big[\wAtwo(s+m) \dots \wAtwo(s+1)\big]\rVert \leq c  m\gamma^m$. Differentiating $\eqref{eq: oxkronecker2}$ with respect to $\theta_l$, applying $\norm{A \otimes B} = \norm{A} \cdot \norm{B}$ as well as the Kronecker product rule $\partial(A \otimes B) = (\partial A) \otimes B + A \otimes (\partial B)$, it is then not hard to see that also $\partial_{\theta_l} \E_\vartheta(\laengerX(t)^{\otimes 2} \mid X(0) = x)$ is bounded uniformly in $t$ and $\theta$ for any $l \in \{1, \dots, k\}$. Since all partial derivatives $\partial_{\theta_l}\E_\vartheta\big[f(\laengerX(t)) \mid X(0) = x\big]$ are uniformly bounded in $t$ and $\theta$, the sequence $(\E_\vartheta\big[f(\laengerX(t)) \mid X(0) = x\big])_{t \in \N^*}$ is equicontinuous by the multivariate mean value theorem, see for example Theorem \ref{theo: matrix_mvt}, and it converges pointwise for $\lambda_E$-almost any $x \in E$. Since equicontinuity and pointwise convergence imply uniform convergence on compacta, $\E_\vartheta\big[f(\laengerX(t)) \mid X(0) = x\big] \to \int f \dd \omu$ uniformly on the compact space $\Theta$, and hence $\int f \dd \omu$ is also continuous in $\theta$.
\end{proof}

The continuity of the limiting constant $\int f \dd \omu$ in $\theta$ will mostly be needed for technical reasons to ensure measurability of any uncountable suprema occurring in the sequel. In order to prove the desired uniform law of large numbers for $\nabla_\theta Z^{\theta}$, we need two more auxiliary results concerning the process $\laengerX$. The first one given in the following lemma strengthens the $L^p$-boundedness of $\laengerX$ to an $L^p$-boundedness uniform in $\theta$. The second result then establishes uniformity in the above Corollary \ref{coro: ergod_oX} for the process $\laengerX$.

\begin{lemma}\label{lem: uniform_lp_boundedness}
    We have $\sup_{t \in \N} \E_\vartheta\big[\sup_{\theta \in \Theta} \lVert\laengerX(t)\rVert^{p}\big] < \infty$ for any $1 \leq p \leq 4 + \delta$.
\end{lemma}
\begin{proof}
    Since $\smash{\E_\vartheta\big[\sup_{\theta \in \Theta} \lVert\laengerX(t)\rVert^{p}\big]}^{\frac{1}{p}} \leq \smash{\E_\vartheta\big[\sup_{\theta \in \Theta} \lVert\laengerX(t)\rVert^{q}\big]}^{\frac{1}{q}}$ for $1 \leq p \leq q$ by Hölder's inequality, it suffices to prove the claim for $p = 4 + \delta$. By \eqref{eq: Xtilde_dynamics} and \eqref{eq: Xbar_dynamics}
    \begin{equation}\label{eq: Xbar_dynamics2}
        \laengerX(t) \eqqcolon \begin{pmatrix}
            X(t) \\ \chiX(t)
        \end{pmatrix} = \begin{pmatrix}
            a^\vartheta \\ r^{\theta} 
        \end{pmatrix} + \begin{pmatrix}
            A^\vartheta & 0 \\
            \underline{G}^\theta(t-1) & \underline{F}^\theta(t-1)
        \end{pmatrix} \laengerX(t-1) + \begin{pmatrix}
            N^\vartheta(t) \\0
        \end{pmatrix},
    \end{equation}
    where $r^\theta \in \R^{d(k^2+k+1)}$, $\underline{G}^\theta(t) \in \R^{d(k^2+k+1) \times d}$ and $\underline{F}^\theta(t) \in \R^{d(k^2+k+1) \times d(k^2+k+1)}$ are continuous in $\theta$. Moreover, the proof of Lemma \ref{limitcov} shows that $\underline{F}^\theta(t) \to \underline{F}^\theta$ uniformly on $\Theta$ for some matrix $\underline{F}^\theta$ that is continuous in $\theta$ with $\rho(\underline{F}^\theta) < 1$ because $\underline{F}^\theta$ is a block lower-triangular matrix whose submatrices on the main diagonal are just repeated copies of the matrix $F^\theta$ and $\rho(F^\theta)< 1$ by the proof of Lemma \ref{limitcov}. Lemma \ref{lem: linalg_unif} shows moreover that there exist constants $c \in \R_+$ and $\gamma \in [0, 1)$ such that $\sup_{\theta \in \Theta}\lVert \underline{F}^\theta(s+m) \underline{F}^\theta(s+m - 1) \dots \underline{F}^\theta(s+1)\rVert \leq c \gamma^m$ holds for any $s \in \N$ and $m \in \N^*$. By iterating equation \eqref{eq: Xbar_dynamics2} we obtain the following decomposition for $\chiX(t)$:
    \begin{fitequation}\label{eq: chi_decomp}
        \chiX(t) = \underbrace{\sum_{s=1}^t \Big(\prod_{u=s}^{t-1} \underline{F}^\theta(u)\Big)r^\theta}_{(1)} + \underbrace{\sum_{s=0}^{t-1}\Big(\prod_{u=s+1}^{t-1} \underline{F}^\theta(u)\Big) \underline{G}^\theta(s)  X(s)}_{(2)} + \underbrace{\Bigg[\prod_{u=0}^{t-1} \underline{F}^\theta(u)\Bigg] \chi^\theta(0)}_{(3)}.
    \end{fitequation}
    For a family $Y = (Y(\theta))_{\theta \in \Theta}$ of random variables indexed by $\theta \in \Theta$, we now write $\vertiii{Y}_{L^p}$ for $\smash{\E_\vartheta\big[\sup_{\theta \in \Theta} \lVert Y(\theta) \rVert^{p}\big]^{1/p}}$. By the Minkowski inequality $\vertiii{Y_1 + Y_2}_{L^p} \leq \vertiii{Y_1}_{L^p} + \vertiii{Y_2}_{L^p}$ proven in Lemma \ref{lem: Hamilton} for two such families $Y_1$ and $Y_2$, it follows that $$\vertiii{\laengerX(t)}_{L^p} \leq \vertiii{X(t)}_{L^p} + \vertiii{\chiX(t)}_{L^p} = \E_\vartheta\big[\lVert X(t)\rVert^p\big]^{\frac{1}{p}} + \E_\vartheta\big[\sup_{\theta \in \Theta} \lVert \chiX(t) \rVert^{p}\big]^{\frac{1}{p}}.$$ Since the first summand on the right is bounded in $t$, it suffices to prove boundedness of the second summand on the right. Since the matrices $\smash{\underline{G}^\theta(t)}$ are continuous in $\theta$ and depend on $t$ only through the matrices $\smash{\widehat \Sigma^\theta(t, t-1)}$ and its inverse, $S^\theta_j(t,t-1)$ and $R^\theta_{ij}(t, t-1)$ which are uniformly bounded in $\theta$ and $t$ for any $i, j \in \{1, \dots, k\}$, we can find some $M \in \R_+$ such that $\lVert \underline{G}^\theta(t)\rVert \leq M$ and $\lVert r^\theta \rVert \leq M$ for any $t \in \N^*$ and $\theta \in \Theta$. For the deterministic summand (1) in \eqref{eq: chi_decomp} we see that $\norm{(1)} \leq c M \sum_{s=1}^t \gamma^{t-s}$, which is bounded in $t$. Likewise, also the deterministic summand (3) in \eqref{eq: chi_decomp} is bounded in $t$. To handle (2), observe that $$\sum_{s=1}^t \sup_{\theta \in \Theta} \Big\lVert \Big(\prod_{u=s}^{t-1} \underline{F}^\theta(u)\Big) \underline{G}^\theta(s-1)\Big\rVert \leq cM \sum_{s=1}^t \gamma^{t-s},$$ which is likewise bounded in $t$. Lemma \ref{lem: Hamilton} then yields that $\vertiii{(2)}_{L^p}$ is bounded in $t$.
\end{proof}

\begin{lemma}\label{lem: PötscherPrucha}
    For any $m$-dimensional polynomial $f: \R^{d(k^2+k+2)} \to \R^m$ of order 4 or less we have that $\PP_\vartheta\big(\sup_{\theta \in \Theta}\big\lVert\frac{1}{t}\sum_{s=1}^t f(\laengerX(s)) - \int f \dd\omu \big\rVert > \varepsilon\big) \to 0$ holds for  any $\varepsilon > 0.$
\end{lemma}
\begin{proof}
    Since $f(x)$ is of the form $f(x) = A_4x^{\otimes 4} + \dots + A_1 x^{\otimes 1} + a_0$ for some matrices $A_j \in \R^{m \times d^j(k^2+k+2)^j}$, $j \in \{1, \dots, 4\}$, and some $a_0 \in \R^m$, it suffices to prove the claim for $f(x) = x^{\otimes j}$, $j \in \{1, \dots, 4\}$. We focus on the case $j=4$ as $1 \leq j \leq 3$ follow similarly. By thrice continuous differentiability assumed in Assumption \ref{assump: A}, $r^\theta$, $\underline{G}^\theta(t)$ and $\underline{F}^\theta(t)$ in \eqref{eq: Xbar_dynamics2} are continuously differentiable in $\theta$. As in the proof of Lemma \ref{lem: uniform_lp_boundedness}, fix $c \in \R_+$ and $\gamma \in [0, 1)$ such that $\sup_{\theta \in \Theta}\lVert \underline{F}^\theta(s+m) \underline{F}^\theta(s+m - 1) \dots \underline{F}^\theta(s+1)\rVert \leq c \gamma^m$ for any $s \in \N$ and $m \in \N^*$. Simple differentiation as in Proposition \ref{prop: Quasi-Score} or \ref{prop: Fisher_Information} yields that $\partial_{\theta_l} \underline{F}^\theta(t)$ and $\partial_{\theta_l} \underline{G}^\theta(t)$ depend on $t$ only through $\widehat \Sigma^\theta(t, t-1)$ and its inverse, $S^\theta_j(t, t-1)$, $R^\theta_{ij}(t,t-1)$ as well as $\partial_{\theta_l} R^\theta_{ij}(t,t-1)$ for $i, j, l \in\{1, \dots, k\}$. Since all partial derivatives $\partial_{\theta_l} R^\theta_{ij}(t,t -1)$ are uniformly bounded in $\theta$ and $t$ by \eqref{eq: partial_R_dynamics} and \ref{lem: lin_systems}.1, the sequences $\partial_{\theta_l} \underline{F}^\theta(t)$ and $\partial_{\theta_l} \underline{G}^\theta(t)$ are also bounded uniformly in $t$ and $\theta$. Fix now $\theta, \theta' \in \Theta$ as well as $x \in E$. Letting $\underline{F}^\theta_{s, t} \coloneq \prod_{u=s}^{t-1} \underline{F}^\theta(u)$, \eqref{eq: chi_decomp} \mbox{yields that $\chiX(t) - \chiXprime(t)$ is of the form}
    \begin{equation*}
        \sum_{s=1}^t \Bigg[ \underline{F}^\theta_{s, t}r^\theta - \underline{F}^{\theta'}_{s, t} r^{\theta'} + \bigg[\underline{F}^\theta_{s, t} \underline{G}^{\theta}(s-1) - \underline{F}^{\theta'}_{s, t} \underline{G}^{\theta'}(s-1)\bigg] X(s-1)\Bigg] + \Bigg[\underline{F}^\theta_{0, t} - \underline{F}^{\theta'}_{0, t}\Bigg] \chi^\theta(0).
    \end{equation*}
    \noindent Integration by parts yields that the derivative $\partial_{\theta_l}\big(\prod_{u=s}^{t} \underline{F}^\theta(u)\big)$ for $l \in \{1, \dots, k\}$ is given by 
    \begin{equation*}
        \partial_{\theta_l}\big(\prod_{u=s}^{t} \underline{F}^\theta(u)\big)=\big[\partial_{l} \underline{F}^{\theta}(t)\big] \underline{F}^{\theta}_{s, t} + \underline{F}^{\theta}(t) \big[\partial_{l} \underline{F}^{\theta}(t-1)\big]\underline{F}^{\theta}_{s, t-1} + \dots + \underline{F}^{\theta}_{s+1, t} \big[\partial_{\theta_l} \underline{F}^{\theta}(s)\big].
    \end{equation*}
    If $M \geq 1$ is fixed such that $\lVert \partial_{\theta_l} \underline{F}^{\theta}(t) \rVert \leq M$, $\lVert \underline{G}^{\theta}(t) \rVert \leq M$, $\lVert \partial_{\theta_l} \underline{G}^{\theta}(t) \rVert \leq M$, $\lVert r^\theta \rVert \leq M$, and $\lVert \partial_{\theta_l} r^\theta \rVert \leq M$ for all $\theta \in \Theta$ and $t \in \N^*$ and since $c \geq 1$ can be assumed without loss of generality, we get $\lVert \partial_{\theta_l}\big(\prod_{u=s}^{t} \underline{F}^{\theta}(u)\big) \rVert \leq c^2 M (t - s + 1)\gamma^{t-s}$. \mbox{Hence we have for $s \leq t$ and $\theta \in \Theta$ that}
    \begin{equation*}
        \bigg\lVert \partial_{l}\bigg[\underline{F}^{\theta}_{s, t}r^\theta \bigg]\bigg\rVert \leq M \bigg\lVert \partial_{\theta_l}\Big(\prod_{u=s}^{t-1} \underline{F}^{\theta}(u)\Big) \bigg\rVert + M \bigg\lVert \Big(\prod_{u=s}^{t-1} \underline{F}^{\theta}(u)\Big) \bigg\rVert \leq c^2 M^2 (t-s+1) \gamma^{t-s - 1}
    \end{equation*}
    and $\big\lVert \partial_{\theta_l} \big[\big(\prod_{u=s}^{t-1} \underline{F}^{\theta}(u)\big) \underline{G}^{\theta}(s-1) \big]\big\rVert \leq c^2 M^2 (t-s+1) \gamma^{t-s-1}$. Then, Theorem \ref{theo: matrix_mvt} yields
    \begin{equation*}
        \bigg\lVert \underline{F}^{\theta}_{s, t}r^\theta - \underline{F}^{\theta'}_{s, t} r^{\theta'}\bigg\lVert \leq \sqrt{k} c^2M^2 (t-s+1) \gamma^{t-s-1} \lVert \theta - \theta'\rVert
    \end{equation*}
    and $\big\lVert\underline{F}^{\theta}_{s, t} \underline{G}^{\theta}(s-1) - \underline{F}^{\theta'}_{s, t} \underline{G}^{\theta'}(s-1)\big\rVert \leq \sqrt{kd} c^2M^2 (t-s+1) \gamma^{t-s-1} \lVert \theta - \theta'\rVert$. Since $\chi^\theta(0)$ is deterministic, we can moreover assume that $\lVert \chiX(0)\rVert \leq M$ so that we have $$\bigg\lVert \bigg[\prod_{u=0}^{t-1} \underline{F}^{\theta}(u) - \prod_{u=0}^{t-1} \underline{F}^{\theta'}(u)\bigg] \chiX(0) \bigg\rVert \leq \sqrt{k}c^2M^2 t \gamma^{t-1} \lVert \theta - \theta'\rVert.$$
    Since $\lVert \laengerX(t) - \underline{X}^{\theta'}(t)\rVert = \lVert \chiX(t) - \chiXprime(t)\rVert$ by definition, we altogether obtain
    \begin{align}
        \lVert \laengerX(t) - \underline{X}^{\theta'}(t)\rVert &\leq C \lVert \theta-\theta'\rVert \sum_{s=1}^t (t-s+1) \gamma^{t-s-1} \Big(1 + \big\lVert X(s-1)\big\rVert\Big) + C \lVert \theta-\theta'\rVert  t \gamma^{t-1} \nonumber\\
        &= C \lVert \theta-\theta'\rVert  \bigg[t \gamma^{t-1} + \sum_{s=1}^t s \gamma^{s-2} \Big(1 + \big\lVert X(t-s)\big\rVert\Big)\bigg], \label{eq: Pötscher_bound}
    \end{align}
    for $t \in \N^*$, where $C \coloneq \sqrt{kd}c^2M^2$. Our strategy to prove the claimed uniform convergence consists now of using this upper bound to show the uniform stochastic equicontinuity condition from Theorem \ref{theo: Pötscher_Prucha}. We first show that $k(\varepsilon) \coloneq \limsup_{t \to \infty} \PP_\vartheta\big( \sup_{\lVert \theta - \theta'\rVert < \alpha} \lVert \laengerX(t) - \underline{X}^{\theta'}(t) \rVert > \varepsilon\big) \to 0$ as $\alpha \to 0$ for any $\varepsilon > 0$. Inserting the inequality \eqref{eq: Pötscher_bound} yields
    \begin{align*}
        k(\varepsilon) &\leq \limsup_{t \to \infty} \PP_\vartheta\Big(C \alpha \Big[t \gamma^{t-1} + \sum_{s=1}^t s \gamma^{s-2} \Big(1 + \big\lVert X(t-s)\big\rVert\Big)\Big] > \varepsilon\Big) \\
        &\leq C \frac{\alpha}{\varepsilon} \limsup_{t \to \infty} \Big[t \gamma^{t-1} + \sum_{s=1}^t s \gamma^{s-2} \Big(1 + \E_\vartheta\Big[\big\lVert X(t-s)\big\rVert\Big]\Big)\Big],
    \end{align*}
    where the superior limit on the right-hand side is finite because $X$ is bounded in $L^1$. Hence $k(\varepsilon) \to 0$ as $\alpha \to 0$, as claimed. Let $\eta, \nu > 0$ be arbitrary. Since $\E_\vartheta\big[\sup_{\theta \in \Theta} \lVert \laengerX(t) \rVert \big]$ is bounded by Lemma \ref{lem: uniform_lp_boundedness}, we have $\limsup_{t \to \infty} \PP_\vartheta\big(\sup_{\theta \in \Theta}\lVert \laengerX(t) \rVert > K \big) < \frac{\nu}{2}$ for some $K \in \R_+$. Now, the function $f: x \mapsto x^{\otimes 4}$ is uniformly continuous on compacta, so we can find $\varepsilon > 0$ such that $\lVert f(x) - f(y) \rVert < \eta$ when $\lVert x - y \rVert \leq \varepsilon$ and $\lVert x \rVert \leq K$. Let $\smash{\Xi_\alpha(t) \coloneq \sup_{\lVert \theta - \theta'\rVert < \alpha} \lVert f(\laengerX(t)) - f(\underline{X}^{\theta'}(t)) \rVert}$ and $k(\eta, f) \coloneq \limsup_{t \to \infty} \PP_\vartheta\big( \Xi_\alpha(t) > \eta\big)$. Then
    \begin{equation*}
        k(\eta, f) \leq \limsup_{t \to \infty} \PP_\vartheta\Big( \sup_{\lVert \theta - \theta'\rVert < \alpha} \lVert \laengerX(t) - \underline{X}^{\theta'}(t) \rVert > \varepsilon\Big) + \limsup_{t \to \infty} \PP_\vartheta\big(\sup_{\theta \in \Theta}\lVert \laengerX(t) \rVert > K \big).
    \end{equation*}
    The first summand on the right is $k(\varepsilon)$. Hence we can find $\alpha > 0$ small enough such that $k(\varepsilon) < \frac{\nu}{2}$, so $k(\eta, f) < \nu$. In other words, $k(\eta, f) \to 0$ as $\alpha \to 0$. For $\varepsilon > 0$, we can now write 
    \begin{align*}
        \limsup_{t \to \infty} \E_\vartheta\big[\Xi_\alpha(t)\big] &\leq \varepsilon + \limsup_{t \to \infty} \E_\vartheta(\Xi_\alpha(t) \mathbf{1}_{\{\Xi_\alpha(t) > \varepsilon\}}) \\
        &\leq \varepsilon + \limsup_{t\to\infty} \E_\vartheta\big[ \Xi_\alpha(t)^{\frac{4 + \delta}{4}} \big]^{\frac{4}{\delta + 4}} \PP_\vartheta( \Xi_\alpha(t) > \varepsilon)^{\frac{\delta}{\delta + 4}} \\
        &\leq \varepsilon + \sup_{t\in \N^*} \E_\vartheta\big[ \Xi_\alpha(t)^{\frac{4 + \delta}{4}} \big]^{\frac{4}{\delta + 4}} \limsup_{t \to \infty} \PP_\vartheta(\Xi_\alpha(t) > \varepsilon)^{\frac{\delta}{\delta + 4}},
    \end{align*}
    \noindent where we used Hölder's inequality. The second factor tends to 0 as $\alpha \to 0$ because $k(\varepsilon, f) \to 0$ as $\alpha \to 0$. \mbox{But as $(a + b)^p \leq 2^{p-1} (a^p + b^p)$ for $a, b > 0$ and any $p \geq 1$ we have}
    \begin{align*}
        \E_\vartheta\big[ \Xi_\alpha(t)^{\frac{4 + \delta}{4}} \big] &= \E_\vartheta\Big( \sup_{\lVert \theta - \theta'\rVert < \alpha} \lVert f(\laengerX(t)) - f(\underline{X}^{\theta'}(t)) \rVert^{\frac{4+\delta}{4}}\Big) \leq 2^{\frac{4+\delta}{4}} \E_\vartheta\Big(\sup_{\theta \in \Theta} \lVert f(\laengerX(t)) \rVert^{\frac{4 + \delta}{4}} \Big) \\
        &= 2^{\frac{4+\delta}{4}} \E_\vartheta\Big(\sup_{\theta \in \Theta} \lVert \laengerX(t)^{\otimes 4} \rVert^{\frac{4 + \delta}{4}} \Big) = 2^{\frac{4+\delta}{4}} \E_\vartheta\Big(\sup_{\theta \in \Theta} \lVert \laengerX(t) \rVert^{4 + \delta} \Big),
    \end{align*}
    where we used that $\lVert A \otimes B \rVert = \norm{A} \cdot\norm{B}$ for any matrices $A$ and $B$. The expected value on the right hand side is bounded in $t$ by Lemma \ref{lem: uniform_lp_boundedness}, so $\sup_{t \in \N^*}\E_\vartheta\big[ \Xi_\alpha(t)^{\frac{4 + \delta}{4}} \big] < \infty$. Combined with the above we obtain $\limsup_{\alpha \to 0} \limsup_{t \to \infty}\E_\vartheta\big[\Xi_\alpha(t)\big] \leq \varepsilon$. Since $\varepsilon > 0$ was arbitrary, it follows that $\limsup_{t \to \infty}\E_\vartheta\big[\Xi_\alpha(t)\big] \to 0$ as $\alpha \to 0$. Since $\limsup_{t \to \infty} \frac{1}{t} \sum_{s=1}^t a_s \leq \limsup_{t \to \infty} a_t$ for any sequence $(a_t)_{t \in \N^*}$ of real numbers (see \cite{Goldberg1976}, Exercise 2.9.7), we also have $\limsup_{t \to \infty} \frac{1}{t}\sum_{s=1}^t \E_\vartheta\big[\Xi_\alpha(s)\big] \to 0$ as $\alpha \to 0$. Finally,
    \begin{fitequation*}
        \PP_\vartheta\Big( \sup_{\lVert \theta- \theta'\rVert < \alpha} \Big\lVert \frac{1}{t}\sum_{s=1}^t f(\laengerX(s)) - \frac{1}{t}\sum_{s=1}^t f(\underline{X}^{\theta'}(s)) \Big\rVert > \varepsilon\Big) \leq \PP_\vartheta \Big( \frac{1}{t}\sum_{s=1}^t \Xi_\alpha(s) > \varepsilon\Big) \leq \frac{1}{t\varepsilon}\sum_{s=1}^t \E_\vartheta\big[\Xi_\alpha(s)\big],
    \end{fitequation*}
    so $\limsup_{t \to \infty} \PP_\vartheta\big( \sup_{\lVert \theta- \theta'\rVert < \alpha} \big\lVert \frac{1}{t}\sum_{s=1}^t f(\laengerX(s)) - \frac{1}{t}\sum_{s=1}^t f(\underline{X}^{\theta'}(s)) \big\rVert > \varepsilon\big) \to 0$ whenever $\smash{\alpha \to 0}$. This is the uniform stochastic equicontinuity condition from Theorem \ref{theo: Pötscher_Prucha}. Since we have $\frac{1}{t}\sum_{s=1}^t f(\laengerX(s)) \to \int f \dd \omu$ in $\PP_\vartheta$-probability for any fixed $\theta$ and since $\frac{1}{t} \sum_{s=1}^t f(\laengerX(s))$ and $\int f \dd \omu$ are continuous in $\theta$ by Corollary \ref{coro: ergod_oX}, Theorem \ref{theo: Pötscher_Prucha} can be applied to yield \linebreak $\smash{\sup_{\theta \in \Theta} \lVert \frac{1}{t} \sum_{s=1}^t f(\laengerX(s)) - \int f \dd \omu \rVert \to 0}$ in $\PP_\vartheta$-probability.
\end{proof}

We are now ready to prove the main result of this section, assumed in Proposition \ref{prop: Jacod-Consistency}:

\begin{proposition}\label{prop: uniform_Z}
    There exists a continuous matrix-valued function $W: \Theta \to \R^{k \times k}$ such that $\sup_{\theta \in \Theta} \lVert \frac{1}{t}\nabla_\theta Z^\theta(t) - W(\theta)\rVert \to 0$ in $\PP_\vartheta$-probability.
\end{proposition}
\begin{proof}
    By continuity of $W$ all suprema over $\Theta$ are measurable by Lemma \ref{lem: measurability}. Since $\nabla_\theta Z^\theta(t,t-1)$ is a $(k \times k)$-dimensional quadratic polynomial of $\laengerX(t)$, it is of the form
    \[\nabla_\theta Z^\theta(t, t-1) = \sum_{|\lambda| \leq 2} \beta_{\lambda}^\theta(t) \laengerX(t)^\lambda \coloneq \widetilde f^\theta(t, \laengerX(t)),\]
    where the sum extends over the $d(k^2+k+2)$-dimensional multi-indices $\lambda$ and where $\beta_\lambda^\theta(t) \in \R^{k \times k}$. As in the proof of Theorem \ref{theo: norm}, all coefficients $\beta_\lambda^\theta(t)$ depend on $t$ only through products of the Kálmán filter covariance matrices $\widehat \Sigma^\theta(t,t-1)$ and its inverse, the matrices $S^\theta_j(t,t-1)$ and $R^\theta_{ij}(t,t-1)$, see the explicit representation in equation \eqref{eq: Fisher-Information}. Since these matrices converge as $t \to \infty$ uniformly on $\Theta$, there exist limiting coefficients $\beta_\lambda^\theta$, continuous in $\theta$, such that $\beta_\lambda^\theta(t) \to \beta_\lambda^\theta$ uniformly in $\theta$. Let $\widetilde f^\theta(y) \coloneq \sum_{|\lambda| \leq 2} \beta_\lambda^\theta y^\lambda$ denote the $(k \times k)$-dimensional quadratic polynomial with these limiting coefficients, see also equation \eqref{eq: f_tilde_def}. Define the matrix-valued function $W(\theta) \coloneq \int \widetilde f^\theta \dd \omu$, which is continuous because the limiting coefficients $\beta_\lambda^\theta$ are continuous in $\theta$ and by Corollary \ref{coro: ergod_oX}. We can then decompose
    \begin{fitequation*}
        \sup_{\theta\in\Theta}\Big\lVert \frac{1}{t}\nabla_\theta Z^\theta(t) - W(\theta)\Big\rVert \leq \underbrace{\sup_{\theta \in \Theta} \Big\lVert\frac{1}{t} \sum_{s=1}^t \widetilde f^\theta(\laengerX(s)) - W(\theta)\Big\rVert}_{(1)} + \underbrace{\sup_{\theta \in \Theta}\Big\lVert \frac{1}{t} \sum_{s=1}^t \sum_{|\lambda| \leq 2} \big[\beta_\lambda^\theta(s) - \beta_\lambda^\theta \big] \laengerX(s)^{\lambda}\Big\rVert}_{(2)}.
    \end{fitequation*}
    The summand (1) converges to 0 in probability by Lemma \ref{lem: PötscherPrucha}. For (2) we find that
    \begin{align*}
        \E_\vartheta\big[(2)\big] \leq \frac{1}{t}\sum_{s=1}^t \sum_{|\lambda|\leq 2} \sup_{\theta \in \Theta} \big\lVert \beta_\lambda^\theta(s) - \beta_\lambda^\theta \big\rVert \E_\vartheta \Big[\sup_{\theta \in \Theta} \big\lVert \laengerX(s)^\lambda\big\rVert\Big].
    \end{align*}
    Since $\sup_{t \in \N} \E_\vartheta\big[\sup_{\theta \in \Theta} \lVert\laengerX(t)\rVert^{2}\big] < \infty$ by Lemma \ref{lem: uniform_lp_boundedness}, it follows by the same calculations as in Lemma \ref{lem: Boundedness_of_Polynomials} that also $\sup_{t \in \N} \E_\vartheta\big[\sup_{\theta \in \Theta} \lVert\laengerX(t)^\lambda\rVert\big] < \infty$ for any $\lambda$ with $|\lambda|\leq 2$. It follows that we can find a constant $C \in \R_+$ such that $\E_\vartheta\big[(2)\big] \leq \frac{C}{t}\sum_{s=1}^t \sum_{|\lambda|\leq 2}\sup_{\theta \in \Theta} \big\lVert \beta_\lambda^\theta(s) - \beta_\lambda^\theta \big\rVert$ and the arithmetic mean on the right-hand side converges to 0 because $\beta_\lambda^\theta(t) \to \beta_\lambda^\theta$ uniformly in $\theta$ and because convergence of a sequence implies convergence in Cesàro mean. So we have shown that $(2) \to 0$ in $L^1(\PP_\vartheta)$.
\end{proof}

\subsubsection{Proof of consistency and asymptotic normality}\label{sec4.2.3: Ident}

The final ingredient for the proof of our main Theorems \ref{theo: main} and \ref{theo: main2} is the fact that $\frac{1}{t}L^\theta(t)$ converges in $\PP_\vartheta$-probability uniformly in $\theta$ to some twice continuously differentiable $Q$ that is uniquely maximised at the true parameter value $\vartheta$. Recall from the proof of Proposition \ref{prop: uniform_Z} that the Fisher information matrix $W(\theta)$ is given by $\int \widetilde f^\theta \dd \omu$, where $\widetilde f^\theta$ denotes the quadratic polynomial describing the limiting dependence of $\nabla_\theta Z^\theta(t, t-1)$ on $\laengerX(t)$ given in \eqref{eq: f_tilde_def}. We start by showing that the same limit $W(\theta)$ \mbox{can be achieved by considering the likelihood of}
\begin{equation*}
Y^{\mathrm{stat}}(t) = a^\theta + A^\theta Y^{\mathrm{stat}}(t-1) + B^\theta W^\theta(t),    
\end{equation*} 
which is a stationary Markov chain and was introduced in equation \eqref{eq: Gauss2}. To do so, define the Gaussian densities $\smash{\widetilde q_t^\theta} \coloneq \smash{\mathrm{d}\PP_\vartheta^{(Y^{\mathrm{stat}}_\oo(1), \dots, Y^{\mathrm{stat}}_\oo(t)} / \mathrm{d} \lambda_{t (d - m)}}$ as well as the conditional Gaussian densities $\smash{\widetilde q_{t \mid t-1}^\theta(\cdot \mid y_1, \dots, y_{t-1})}$ of the conditional distribution of $\smash{Y^{\mathrm{stat}}_\oo(t)}$ under $\PP_\theta$ given $Y^{\mathrm{stat}}_\oo(1) = y_1, \dots, Y^{\mathrm{stat}}_\oo(t-1) = y_{t-1}$. The closed form of $\smash{\log \widetilde q_{t \mid t-1}^\theta}$ is similar to the closed form of $\smash{\log q_{t\mid t-1}^\theta}$ given in \eqref{eq: log-lik} with the only exception that the time-dependent matrices $C^\theta(t+1)$ in Proposition \ref{prop: log-lik} are replaced by the limiting matrices $C^\theta$ and that the initial mean and covariance matrix differ. We write $\smash{L^{\theta}_Y(t) \coloneq \log q_t^\theta(Y_\oo(1), \dots, Y_\oo(t))}$ as well as $\smash{\widetilde L^{\theta}_Y(t) \coloneq \log \widetilde q_t^\theta(Y_\oo^{\mathrm{stat}}(1), \dots, Y_\oo^{\mathrm{stat}}(t))}$. Moreover, we define the scores $\smash{Z^{\theta}_Y(t) \coloneq \nabla_\theta L^{\theta}_Y(t)}$ and $\smash{\widetilde Z^{\theta}_Y(t) \coloneq \nabla_\theta \widetilde L^{\theta}_Y(t)}$ as well as the observed Fisher information $\smash{\nabla_\theta Z^{\theta}_Y(t)}$ and $\smash{\nabla_\theta \widetilde Z^{\theta}_Y(t)}$. The closed-form expressions for the score process $\widetilde Z^{\theta}_Y(t)$ and the observed Fisher information $\nabla_\theta \widetilde Z^{\theta}_Y(t)$ are again analogous to Propositions \ref{prop: Quasi-Score} and \ref{prop: Fisher_Information} with $C^\theta$ in place of $C^\theta(t).$

We know that $\smash{L^\theta(t, t-1) = \log q_{t \mid t-1}^\theta(X_\oo(t) \mid X_\oo(1), \dots, X_\oo(t - 1))}$ is a quadratic polynomial in $\langX(t)$ with time-dependent coefficients. If we replace $X$ by $Y$ in this log-likelihood, we can define an $\R^{d(k+2)}$-valued process $\wYhat$ in exactly the same manner as we defined $\langX$ but with $Y$ in place of $X$ so that $\smash{\log q_{t \mid t-1}^\theta(Y_\oo(t) \mid Y_\oo(1), \dots, Y_\oo(t - 1))}$ is a quadratic polynomial in $\wYhat$. In other words, we have $\gamma_\lambda^\theta \in \R$ and multi-indices $\lambda$ such that
\begin{equation}\label{eq: Y_loglik}
    L^{\theta}_Y(t) = \log q_t^\theta(Y_\oo(1), \dots, Y_\oo(t)) = \sum_{s=1}^t \sum_{|\lambda| \leq 2} \gamma_\lambda^\theta(s) \wYhat(s)^\lambda.
\end{equation}
It is easy to see that $Y^{\mathrm{stat}}$ is a parametric polynomial state space model of order 2 that fulfils Assumptions A, B, C. In particular, we can repeat the construction from Section \ref{sec4.2.1: CLT} for $Y^{\mathrm{stat}}$ instead of $X$ to obtain the augmented process $\wY$ and its homogeneous counterpart $\homY$ in place of $\langX$ and $\langhomX$, respectively. We then find that:

\begin{lemma}\label{lem: Y_hat_Y}
    $\wYhat(t) - \wY(t) \xrightarrow{L^p} 0$ under any $\PP_\vartheta$, for any $\theta \in \Theta$ and any $p \geq 1$.
\end{lemma}
\begin{proof}
    By definition the process $\wYhat(t)$ has the same state transition vector $\wa$ and time-de\-pen\-dent state transition matrix $\wA(t)$ as $\langX(t)$ because the Kálmán filter covariance matrices for $Y^\theta$ coincide with those for $X^\theta$. Moreover, the state transition vector of $\wY$ is also given by $\wa$. Let $\wAY(t)$ denote the time-dependent state transition matrix for $\wY(t)$ under $\PP_\vartheta$. It is easy to see that the Kálmán filter covariance matrices for $Y^{\mathrm{stat}}$ converge to the same limit as the Kálmán filter covariance matrices for $Y$ because the limit of the algebraic Riccati difference recursion \eqref{Riccati} remains unchanged if $B^\theta(t+ 1)$ is replaced by $B^\theta$ in \eqref{Riccati}. It follows that $\wA(t)$ and $\wAY(t)$ converge to the same limit $\homA$ at a geometric rate, which is the state transition matrix of $\langhomX$ and $\homY$. By the same calculations that were used to show that $\smash{\langX(t) - \langhomX(t) \to 0}$ in $L^{4+\delta}$ in Proposition \ref{prop: asymptotic_equiv}, it then follows that $\smash{\wYhat(t) - \homY(t) \xrightarrow{L^p} 0}$ and $\smash{\wY(t) - \homY(t) \xrightarrow{L^p} 0}$ for any $p \geq 1$.
\end{proof}

\begin{proposition}\label{prop: equal_limits}
    For any $\theta \in \Theta$ we have:
    \begin{itemize}
        \item[1.] $\lim\limits_{t \to \infty} \frac{1}{t} \E_\vartheta\big[ L^\theta(t)\big] = \lim\limits_{t \to \infty} \frac{1}{t} \E_\vartheta\big[ \widetilde L^{\theta}_Y(t)\big] \eqqcolon Q(\theta)$,
        \item[2.] $\lim\limits_{t \to \infty} \frac{1}{t}\E_\vartheta\big[Z^{\theta}(t)\big] = \lim\limits_{t \to \infty} \frac{1}{t}\E_\vartheta\big[ \widetilde Z^{\theta}_Y(t)\big] \eqqcolon G(\theta)$,
        \item[3.] $\lim\limits_{t \to \infty} \frac{1}{t} \E_\vartheta\big[\nabla_\theta Z^{\theta}(t)\big] = \lim\limits_{t \to \infty} \frac{1}{t} \E_\vartheta\big[\nabla_\theta\widetilde Z^{\theta}_Y(t)\big] = W(\theta)$.
    \end{itemize}
\end{proposition}
\begin{proof}
    We prove the claim for the limits in 1. because 2. and 3. follow completely analogously. By construction we have $\E_\vartheta\big[ \log q_t^\theta(X_\oo(1), \dots, X_\oo^(t))\big] = \E_\vartheta\big[ \log q_t^\theta(Y_\oo(1), \dots, Y_\oo(t))\big]$ for $\smash{t \in \N^*}$ because $\log q_t^\theta(X_\oo(1), \dots, X_\oo(t))$ is linear in the components of $(1, X_\oo(u)) (1, X_\oo(s))^\top$ with $1 \leq u,s \leq t$, which have the same expectation as $(1, Y_\oo(u)) (1, Y_\oo(s))^\top$ by Proposition \ref{prop: Gauss_equiv}. In line with \eqref{eq: Y_loglik} and with the notation introduced further above, we can write
    \begin{align*}
        L^{\theta}_Y(t) \coloneq \log q_t^\theta(Y_\oo(1), \dots, Y_\oo(t)) &= \sum_{s=1}^t \sum_{|\lambda| \leq 2} \gamma_\lambda^\theta(s) \wYhat(s)^\lambda, \\
        \widetilde L^{\theta}_Y(t) \coloneq \log \widetilde q_t^\theta(Y^{\mathrm{stat}}_\oo(1), \dots, Y^{\mathrm{stat}}_\oo(t)) &= \sum_{s=1}^t \sum_{|\lambda| \leq 2} \widetilde \gamma_\lambda^\theta(s) \wY(s)^\lambda,
    \end{align*}
    where the coefficients $\gamma_\lambda^\theta(t)$ and $\widetilde \gamma_\lambda^\theta(t)$ converge to the same limiting coefficients as $t \to \infty$ because the Kálmán filter covariance matrices for $Y$ and the Kálmán filter covariance matrices for $Y^{\mathrm{stat}}$ converge to the same limits. Note that $\wYhat(t)^\lambda - \wY(t)^\lambda \to 0$ in $L^1$ for any multi-index $\lambda$ with $|\lambda| \leq 2$ by Lemma \ref{lem: Y_hat_Y} and Lemma \ref{lem: Conv_of_Polynomials}. Altogether, we obtain that 
    \begin{equation*}
        \frac{1}{t}L^{\theta}_Y(t) - \frac{1}{t}\widetilde L^{\theta}_Y(t) = \frac{1}{t} \sum_{s=1}^t \sum_{|\lambda| \leq 2} \Big[ \gamma_\lambda^\theta(s) \big(\wYhat(s)^\lambda - \wY(s)^\lambda\big) + \big( \gamma_\lambda^\theta(s) - \widetilde \gamma_\lambda^\theta(s)\big) \wY(s)^\lambda\Big]
    \end{equation*}
    converges to 0 in $L^1(\PP_\vartheta)$. Since we know that the limit $\lim_{t\to \infty} \frac{1}{t}\E_\vartheta(\widetilde L^{\theta}_Y(t))$ exists and is independent of $x$ by ergodicity of $\wY$ (see Corollary \ref{coro: ergod_wX}), the identity 1. follows.
\end{proof}

\begin{remark}\label{rem: unif_conv_to_Q}
    By ergodicity of $\langX$ and $\wY$ from Corollary \ref{coro: ergod_oX}, it is clear that the deterministic limits $Q(\theta) \in \R$ and $G(\theta) \in \R^k$ are also limits in $\PP_\vartheta$-probability of $\frac{1}{t}L^\theta(t)$ and $\frac{1}{t} Z^\theta(t)$, respectively. In view of the uniform convergence to $W(\theta)$ proven in Proposition \ref{prop: uniform_Z}, it follows from Lemma \ref{lem: tao} that $W(\theta) = \nabla_\theta G(\theta)$ as well as $G(\theta) = \nabla_\theta Q(\theta)$, and that $\smash{\sup_{\theta \in \Theta} \lVert \frac{1}{t} Z^\theta(t) - G(\theta)\rVert \xrightarrow{\PP_\vartheta} 0}$ as well as $\smash{\sup_{\theta \in \Theta} \lVert \frac{1}{t} L^\theta(t) - Q(\theta)\rVert \xrightarrow{\PP_\vartheta} 0}$. This establishes the uniform convergence property from Proposition \ref{prop: Jacod-Consistency}. It only remains to be verified that $Q(\vartheta) > Q(\theta)$ for any $\theta \neq \vartheta$.
\end{remark}

Our proof of the fact that $Q(\vartheta) > Q(\theta)$ for all $\theta \neq \vartheta$ is based on the results from \cite{Douc2011} and builds on an information-theoretic approach. To this end we make use of the Kullback--Leibler divergence or relative entropy between probability measures, defined by $\mathrm{KL}(\mu \, \lVert \, \nu) \coloneq \int \log\big(\frac{\mathrm{d}\mu}{\mathrm{d}\nu} \big) \dd \mu$ if $\mu$ is absolutely continuous with respect to $\nu$ and $\mathrm{KL}(\mu \, \lVert \, \nu) = \infty$ otherwise. The following notion of \textit{exponential separation} was introduced in \cite{Douc2011}:

\begin{definition}\label{defi: Exponential Separability}
    For $t \in \N^*$ let $\mu_t$ and $\nu_t$ denote probability measures on measurable spaces $(E_t, \mathscr{E}_t)$. The sequence $(\nu_t)_{t \in \N^*}$ is called \textbf{exponentially separated} from $(\mu_t)_{t \in \N^*}$, denoted by $(\nu_t)_{t \in \N^*} \dashv (\mu_t)_{t \in \N^*}$ if there exist sets $A_t \in \mathscr{E}_t$ such that we have $\liminf_{t \to \infty} \mu_t(A_t) > 0$ and  $\limsup_{t\to\infty}\frac{1}{t}\log \nu_t(A_t) < 0$. If $\mu$ and $\nu$ are probability measures on $(E^{\N^*}, \mathscr{E}^{\otimes \N^*})$ for some measurable space $(E, \mathscr{E})$, we write $\nu \dashv \mu$ if $(\nu_t)_{t\in \N^*} \dashv (\mu_t)_{t \in \N^*}$, where $\mu_t \coloneq \mu \circ \pi_t^{-1}$, $\nu_t \coloneq \nu \circ \pi_t^{-1}$ and $\pi_t: E^{\N^*} \to E^t$ denotes the canonical projection on the first $t$ components.
\end{definition}

The connection between the preceding definition of exponential separation and the relative entropy is given in the following lemma, \mbox{which is shown in \cite{Douc2011}, Lemma 10.}

\begin{lemma}\label{lem: Kullback_Leibler}
    If $(\nu_t)_{t \in \N^*} \dashv (\mu_t)_{t \in \N^*}$, then $\liminf_{t \to \infty} \frac{1}{t} \mathrm{KL}(\mu_t \, \lVert \, \nu_t) > 0$. $\hfill \qed$
\end{lemma}

For verifying that $Q$ has a unique maximum at the true parameter $\vartheta$, it suffices to show $\PP_\theta^{Y^{\mathrm{stat}}_\oo} \dashv \PP_\vartheta^{Y^{\mathrm{stat}}_\oo}$ for any $\theta \neq \vartheta$. Indeed, since then also $\PP_\theta^{(Y^{\mathrm{stat}}_\oo(1), \dots)} \dashv \PP_\vartheta^{(Y^{\mathrm{stat}}_\oo(1), \dots)}$, this implies by definition of the Kullback--Leibler divergence that $$\liminf_{t \to \infty} \frac{1}{t} \E_\vartheta \bigg[\log \frac{\widetilde q_t^{\vartheta}(Y^{\mathrm{stat}}_\oo(1), \dots, Y^{\mathrm{stat}}_\oo(t))}{\widetilde q_t^\theta(Y^{\mathrm{stat}}_\oo(1), \dots, Y^{\mathrm{stat}}_\oo(t))}\bigg] = \liminf_{t \to \infty} \frac{1}{t}\E_\vartheta\big[\widetilde L^{\vartheta}_Y(t) - \widetilde L^{\theta}_Y(t)\big] > 0,$$
i.e.\ $Q(\vartheta) > Q(\theta)$ for any $\theta \neq \vartheta$ by Proposition \ref{prop: equal_limits}. The key to this exponential separation lies in the geometric ergodicity of the process $Y^{\mathrm{stat}}$, a concept introduced in Definition \ref{defi: geom_ergod}:

\begin{proposition}\label{prop: geom_ergod_y}
    Under all $\PP_\theta$, $Y^{\mathrm{stat}}$ is $V_\theta$-geometrically ergodic for some $V_\theta: \R^d \to [1,\infty)$.
\end{proposition}
\begin{proof}
    Since $Y^{\mathrm{stat}}$ is a polynomial state space model meeting Assumption B, we obtain from the proof of Proposition \ref{prop: ergod} that $Y^{\mathrm{stat}}$ is a $\psi$-irreducible, aperiodic weak Feller chain, $\psi$ denoting Lebesgue measure on $\R^d$, with unique stationary law $\mu_{\theta, Y} \coloneq \mathcal{N}(\alpha^\theta, \Lambda^\theta)$ on $\mathscr{B}(\R^d)$, see the discussion before Assumption C in Section \ref{s:GQLE}. By Theorem \ref{theo: geom-ergod} it remains to construct a function $V_\theta: \R^d \to [1,\infty)$ such that $\hatP V_\theta - V_\theta \leq \beta V_\theta + b \mathbf{1}_C$ for some $\beta < 0$, some $b \in \R$ and some compact $C\subseteq \R^d$, where $\hatP$ denotes the transition operator for the Markov chain $Y^{\mathrm{stat}}$ under $\PP_\theta$. Here, the constants $\beta$, $b$ and the set $C$ are allowed to depend on $\theta \in \Theta$ as well, which is suppressed in the notation. Let $N \in \R^{d\times d}$ be any positive definite matrix and let $M^\theta$ be the unique positive definite solution to the Lyapunov equation $M^\theta = A^{\theta^\top} M^\theta A^\theta + N$, which exists by Lemma \ref{lem: lin_systems}.2. Let $V_\theta(x) = 1 + x^\top M^{\theta}x$. Then
    \begin{align*}
        \hatP V_\theta(x) - V_\theta(x) &= a^{\theta^\top} M^\theta a^\theta + x^\top \Big(A^{\theta^\top}M^\theta A^\theta - M^\theta\Big) x + \Tr\big(B^{\theta^\top}M^\theta B^\theta\big) + 2 a^{\theta^\top} A^\theta x \\
        &= a^{\theta^\top} M^\theta a^\theta + \Tr\big(B^{\theta^\top}M^\theta B^\theta\big) - x^\top N x + 2 a^{\theta^\top} A^\theta x \\
        &\leq a^{\theta^\top} M^\theta a^\theta + \Tr\big(B^{\theta^\top}M^\theta B^\theta\big) - \frac{\lambda_{\min}(N)}{\lambda_{\max}(M^\theta)} x^\top M^\theta x + 2 a^{\theta^\top} A^\theta x \\
        &= c_\theta + \beta_\theta V_\theta(x) + d_\theta^\top x = \beta_\theta V_\theta(x) \Big( 1 + \frac{c_\theta + d_\theta^\top x}{\beta_\theta V_\theta(x)}\Big),
    \end{align*}
    where $\beta_\theta \coloneq - \frac{\lambda_{\min}(N)}{\lambda_{\max}(M^\theta)} < 0$, $d_\theta \coloneq 2 A^{\theta^\top}a^\theta$, and $c_\theta \coloneq a^{\theta^\top} M^\theta a^\theta + \Tr\big(B^{\theta^\top}M^\theta B^\theta\big) - \beta_\theta$. Since $\smash{\frac{c_\theta + d_\theta^\top x}{\beta_\theta V_\theta(x)} \to 0}$ as $\norm{x} \to \infty$, we can find $K_\theta \geq 0$ such that $\smash{1 + \frac{c_\theta + d_\theta^\top x}{\beta_\theta V_\theta(x)} > \frac{1}{2}}$ for $\norm{x} \geq K_\theta$. Define $C_\theta \coloneq \big\lbrace x \in \R^d: \: \norm{x} \leq K_\theta\big\rbrace$. Then $\hatP V_\theta(x) - V_\theta(x) \leq \frac{\beta_\theta}{2} V_\theta(x)$ for $x \in C_\theta^c$ and $\hatP V_\theta(x) - V_\theta(x) \leq \frac{\beta_\theta}{2} V_\theta(x) + b_\theta$ for $x \in C_\theta$, where $b_\theta \coloneq \sup_{x \in C_\theta} \big[ c_\theta + d_\theta^\top x\big]$.
\end{proof}

Using the geometric ergodicity of $Y^{\mathrm{stat}}$ from above, the key step to establish the exponential separation $\PP_\theta^{Y^{\mathrm{stat}}_\oo} \dashv \PP_\vartheta^{Y^{\mathrm{stat}}_\oo}$ consists in proving a certain large deviations principle for geometrically ergodic Markov chains. This is done by embedding the Markov chain into some wide sense regenerative process. This splitting construction goes back to \cite{Douc2011}, Theorem 17 and culminates in an {Azuma--Hoeffding-type inequality for $Y^{\mathrm{stat}}_\oo$, given below. For brevity, we write $Y^{\mathrm{stat}}_\oo(s, \dots, t)$ for $Y^{\mathrm{stat}}_\oo(s), \dots, Y^{\mathrm{stat}}_\oo(t)$ and $s \leq t$.

\begin{proposition}\label{theo: Azuma-Hoeffding}
    Fix $u \in \N^*$ and $\theta \in \Theta$ and let $f: \R^{(d - m)(u+1)} \to \R$ be bounded and measurable. Then there is $K > 0$ depending on $u$, $f$, $\theta$ such that we have
     \begin{equation}\label{eq: final_azuma_hoeffding}
         \PP_\theta\bigg( \bigg\lvert \frac{1}{t}\sum_{s=1}^t f\big[Y^{\mathrm{stat}}_\oo(s, \dots, s + u)\big] - \E_{\theta}\Big[ f\big( Y^{\mathrm{stat}}_\oo(0, \dots, u)\big)\Big] \bigg\rvert  \geq \lambda\bigg) \leq K \e^{-\frac{\lambda t}{K}\big(\lambda \wedge 1\big)}
     \end{equation}
    for all $\lambda > 0$.
\end{proposition}
\begin{proof}
    By Theorem 17 in \cite{Douc2011}, we can find $K > 0$ for any $V$-geometrically ergodic $\R^d$-valued Markov chain $X$ with invariant measure $\mu$ and initial distribution $\nu$ such that we have
    \begin{equation}\label{eq: theo17}
        \PP\bigg(\bigg\lvert \frac{1}{t}\sum_{s=1}^t f(X(s)) - \int f \dd \mu \bigg\rvert \geq \lambda \bigg) \leq K \Big(\int V \dd \nu\Big) \e^{-\frac{\lambda t}{K \lVert f\rVert_\infty} \big(\frac{\lambda}{\lVert f\rVert_\infty} \wedge 1\big)}.
    \end{equation}
    for all bounded, measurable functions $f: \R^d \to \R$, for all $x \in \R^d$ and for all $\lambda > 0$. Assume $\E_\theta\Big[ f\big( Y^{\mathrm{stat}}_\oo(0, \dots, u)\big)\Big] = 0$ without loss of generality and consider the decomposition
    \begin{fitequation}\label{eq: decomp_AzH}
        \sum_{s=1}^t f\big[Y^{\mathrm{stat}}_\oo(s, \dots, s + u)\big] = \sum_{r=0}^u \bigg(\sum_{s=1}^t \xi^\theta_{s, r} \bigg) + \sum_{s=1}^t \E_\theta\Big(f\big[Y^{\mathrm{stat}}_\oo(s, \dots,s + u)\big] \: \Big| \: \F_{s-1}\Big),
    \end{fitequation}
    where we define
    \begin{equation*}
        \xi^\theta_{s, r} \coloneq \E_\theta\Big( f\big[Y^{\mathrm{stat}}_\oo(s, \dots, s + u)\big] \mid \F_{s+r}\Big) - \E_\theta\Big( f\big[Y^{\mathrm{stat}}_\oo(s, \dots, s + u)\big] \mid \F_{s+r-1}\Big).
    \end{equation*}
    for $r = 0, \dots, u$ and $s = 1, \dots, t$. Since $(\xi^\theta_{s, r})_{s \in \{1, \dots, t\}}$ is a $\PP_\theta$-martingale difference sequence with respect to the filtration $\big(\F_{s + r}\big)_{s \in \{1, \dots, t\}}$ for fixed $r$ and since $\lVert \xi^\theta_{s, r}\rVert_\infty \leq 2 \lVert f \rVert_\infty$ for any $s \in \{1, \dots, t\}$ and $r \in \{0, \dots, u\}$, the classic Azuma--Hoeffding inequality (see for example \cite{Williams1991}, E.14.2(b)) implies
    \begin{equation}\label{eq: Azuma_Hoeff}
        \PP_\theta\bigg( \frac{1}{t}\bigg\lvert \sum_{s=1}^t \xi^\theta_{s, r} \bigg\rvert \geq \lambda \bigg) \leq 2\exp\Big( - \frac{\lambda^2 t}{8 \lVert f \rVert_\infty^2}\Big)
    \end{equation}
    for any $r \in \{0, \dots, u\}$. Furthermore, $\smash{\E_\theta\big(f\big[Y^{\mathrm{stat}}_\oo(s), \dots, Y^{\mathrm{stat}}_\oo(s + u)\big] \: \big| \: \F(s-1)\big)}$ can be written as a measurable function $\smash{g_\theta\big(Y^{\mathrm{stat}}(s-1)\big)}$ by the Markov property, where $\smash{\int g_\theta \dd \mu_{\theta, Y} = 0}$ because we assumed that $\smash{\E_\theta\big[ f\big( Y^{\mathrm{stat}}_\oo(0), \dots, Y^{\mathrm{stat}}_\oo(u)\big)\big] = 0}$. In particular, $\lVert g_\theta \rVert_\infty \leq \lVert f \rVert_\infty$. Since $\smash{Y^{\mathrm{stat}}}$ is $V_\theta$-geometrically ergodic for some $V_\theta: \R^d \to [1,\infty)$ by Proposition \ref{prop: geom_ergod_y}, it follows from \eqref{eq: theo17} that
    \begin{equation}\label{eq: consequence_theo17}
        \PP_\theta\bigg( \frac{1}{t}\bigg\lvert \sum_{s=1}^t g_\theta(Y^{\mathrm{stat}}(s-1)) \bigg\rvert \geq \lambda\bigg) \leq K \Big(\int V_\theta \dd \mu_{\theta, Y}\Big) \e^{-\frac{\lambda t}{K \lVert f\rVert_\infty} \big(\frac{\lambda}{\lVert f\rVert_\infty} \wedge 1\big)}.
    \end{equation}
    In particular, the integral on the right-hand side is finite because $V_\theta$ is a quadratic polynomial by the construction in the proof of Proposition \ref{prop: geom_ergod_y}. Combining the inequalities \eqref{eq: Azuma_Hoeff} and \eqref{eq: consequence_theo17} and inserting these into \eqref{eq: decomp_AzH}, one obtains
    \begin{align*}
        \PP_\theta\bigg(\bigg\lvert \frac{1}{t}&\sum_{s=1}^t f\big[Y^{\mathrm{stat}}_\oo(s, \dots, s + u)\big] \bigg\rvert \geq \lambda\bigg) \\
        &\leq \sum_{r=0}^u \PP_\theta\bigg( \bigg\lvert \frac{1}{t} \sum_{s=1}^t \xi^\theta_{s, r} \bigg\rvert \geq \frac{\lambda}{2(u + 1)}\bigg) + \PP_\theta\bigg(\bigg\lvert \frac{1}{t} \sum_{s=1}^t g_\theta(Y^{\mathrm{stat}}(s-1)) \bigg\rvert \geq \frac{\lambda}{2}\bigg) \\
        &\leq 2(u+1) \e^{-\frac{\lambda^2 t}{32(u+1)^2 \lVert f \rVert^2_\infty}} + K \Big(\int V_\theta \dd \mu_{\theta, Y} \Big)\e^{- \frac{\lambda t}{2 K \norm{f}_\infty}\big( \frac{\lambda}{2\norm{f}_\infty} \wedge 1\big)}.
    \end{align*}
\noindent Rearranging and combining the constants in the sum above yields the claim.
\end{proof}

\begin{remark}
    In principle the proof of the Azuma--Hoeffding-type inequality \eqref{eq: final_azuma_hoeffding} could be shortened by deriving it directly from the simpler concentration inequality \eqref{eq: theo17}, relying on geometric ergodicity of the Markov chain $\big(Y^{\mathrm{stat}}(t), \dots, Y^{\mathrm{stat}}(t + u)\big)_{t \in \N}$. This would however require a more technical construction of a suitable function $V_\theta$ for $V_\theta$-geometric ergodicity of this process in the spirit of the proof of Proposition \ref{prop: geom_ergod_y}. Instead we chose the different strategy of proof presented above, which only requires $V_\theta$-geometric ergodicity of $Y^{\mathrm{stat}}$.
\end{remark}

At its core the preceding Proposition \ref{theo: Azuma-Hoeffding} establishes that the convergence of arithmetic means in $Y^{\mathrm{stat}}$ to corresponding expectations under the stationary distribution is geometrically fast. Finally, the following Lemma \ref{lem: Identifiability} shows that under the identifiability Assumption C, the parameters $\theta$ and $\vartheta$ with $\theta \neq \vartheta$ can be distinguished by observing $Y^{\mathrm{stat}}_\oo$:

\begin{lemma}\label{lem: Identifiability}
    For any $\theta \neq \theta'$ there is $u \in \N$ such that $\PP_{\theta}^{Y^{\mathrm{stat}}_\oo(0, \dots, u)} \neq \PP_{\theta'}^{Y^{\mathrm{stat}}_\oo(0, \dots, u)}$.
\end{lemma}
\begin{proof}
    By the considerations right before Assumption C in Section \ref{s:GQLE}, $\mu_{\theta, Y}$ is a Gaussian distribution and hence admits a Lebesgue density $f_\theta^0$. The Lebesgue density of the joint Gaussian law $\smash{\PP_{\theta}^{(Y^{\mathrm{stat}}_\oo(0), \dots, Y^{\mathrm{stat}}_\oo(u))}}$ is the product of $f_\theta^0(y_0)$ and the conditional densities of $Y^{\mathrm{stat}}_\oo(s)$ given $Y^{\mathrm{stat}}_\oo(0)=y_0, \dots, Y^{\mathrm{stat}}_\oo(s-1)=y_{s-1}$ for $s \in \{1, \dots, u\}$. Hence it suffices to show that for any $\theta \neq \theta'$ there exists some $u \in \N^*$ such that the conditional density $\smash{f^{u \mid u-1}_\theta(y_u \mid y_0, \dots, y_{u-1})}$ of $Y^{\mathrm{stat}}_\oo(u)$ given $Y^{\mathrm{stat}}_\oo(0)=y_0, \dots, Y^{\mathrm{stat}}_\oo(u-1)=y_{u-1}$ differs from the corresponding conditional density $\smash{f_{\theta'}^{u \mid u-1}(y_u \mid y_0, \dots, y_{u-1})}$. Here, $\smash{f_\theta^{0 \mid -1}}$ is set to the unconditional density $f_\theta^0$. Since $Y^{\mathrm{stat}}$ is a Gaussian process, \cite{KallsenRichert2025}, Proposition 3.2, yields that the conditional density $\smash{f_\theta^{s\mid s-1}(y_s \mid y_0, \dots, y_{s-1})}$ is the density of the normal distribution with mean $\smash{\widehat Y_\oo^{\theta, \mathrm{stat}}}$ evaluated at $Y^{\mathrm{stat}}_\oo(0)=y_0, \dots, Y^{\mathrm{stat}}_\oo(s-1)=y_{s-1}$ and covariance matrix $\smash{\widehat \Sigma^{\theta, \mathrm{stat}}_\oo(s, s-1)}$ given by the Kálmán filter recursions. But for any $\theta \neq \theta'$, there exists $u \in \N$ such that $\smash{\widehat Y^{\theta, \mathrm{stat}}_\oo(u, u-1) \neq \widehat Y^{\theta', \mathrm{stat}}_\oo(u, u-1)}$ or $\widehat \Sigma^{\theta, \mathrm{stat}}_\oo(u,u-1) \neq \widehat \Sigma^{\theta', \mathrm{stat}}_\oo(u,u-1)$ by Assumption C.
\end{proof}

The preceding lemma now finally yields the desired exponential separation property:

\begin{corollary}\label{coro: exp_sep}
     For any $\theta \neq \vartheta$ the exponential separation $\PP_\theta^{Y^{\mathrm{stat}}_\oo} \dashv \PP_\vartheta^{Y^{\mathrm{stat}}_\oo}$ holds.
\end{corollary}
\begin{proof}
     Since $\theta \neq \vartheta$, we have $\PP_{\theta}^{(Y^{\mathrm{stat}}_\oo(0, \dots, u))} \neq \PP_{\vartheta}^{(Y^{\mathrm{stat}}_\oo(0, \dots, u))}$ for some $u \in \N$ by the preceding Lemma \ref{lem: Identifiability}. Accordingly, since a measure is determined by its action on bounded measurable functions, there exists $u \in \N$ and a bounded measurable $h: \R^{(d-m)(u + 1)} \to \R$ such that $\E_{\theta}\big[ h(Y^{\mathrm{stat}}_\oo(0, \dots, u))\big] = 0$ and $\E_{\vartheta}\big[ h(Y^{\mathrm{stat}}_\oo(0, \dots, u))\big] = 1$. For $t > u$ define 
     $$A_t \coloneq \Big\lbrace (y_0, \dots, y_t) \in \R^{(d-m)(t+1)}: \;\Big\lvert \frac{1}{t - u}\sum_{s=1}^{t-u} h\big(y_s, \dots, y_{s+u}\big) \Big\rvert \geq \frac{1}{2}\Big\rbrace \in \mathscr{B}(\R^{(d-m)(t + 1)}).$$
     By ergodicity of $Y^{\mathrm{stat}}$ it follows that $\PP_\vartheta^{(Y^{\mathrm{stat}}_\oo(0), \dots, Y^{\mathrm{stat}}_\oo(t))}(A_t) \to 1$ as $t \to \infty$. On the other hand, Theorem \ref{theo: Azuma-Hoeffding} shows that there exists some constant $K > 0$ such that $\PP_\theta^{(Y^{\mathrm{stat}}_\oo(0, \dots, t))}(A_t) \leq K \exp(\frac{u-t}{4K})$, which implies that
     \[\limsup_{t\to \infty} \frac{1}{t} \log \PP_\theta^{(Y^{\mathrm{stat}}_\oo(0,\dots,t))}(A_t) \leq -\frac{1}{4K} < 0.\]
     By definition, this establishes the exponential separation property $\PP_\theta^{Y^{\mathrm{stat}}_\oo} \dashv \PP_\vartheta^{Y^{\mathrm{stat}}_\oo}$.
\end{proof}

As noted earlier, the exponential separation property established in the preceding Corollary \ref{coro: exp_sep} yields that $Q(\theta) < Q(\vartheta)$ for any $\theta \neq \vartheta$. This is the final ingredient for the proof of the main theorem which we will cite as a separate corollary here for the sake of completeness:

\begin{corollary}\label{coro: uniquely_maximised}
    For any $\theta \in \Theta$ with $\theta \neq \vartheta$ it holds that $Q(\theta) < Q(\vartheta)$.
\end{corollary}
\begin{proof}
    Since $\PP_\theta^{Y^{\mathrm{stat}}_\oo} \dashv \PP_\vartheta^{Y^{\mathrm{stat}}_\oo}$ by Corollary \ref{coro: exp_sep}, we have that $(\nu^\theta_t)_{t \in \N} \dashv (\nu^\vartheta_t)_{t \in \N}$, where we set $\smash{\nu^\theta_t \coloneq \PP_\theta^{(Y^{\mathrm{stat}}_\oo(0), \dots, Y^{\mathrm{stat}}_\oo(t))}}$. Hence we obtain $\liminf_{t \to \infty} \frac{1}{t} \mathrm{KL}(\nu_t^{\vartheta}  \, \lVert \, \nu^\theta_t) > 0$ by Lemma \ref{lem: Kullback_Leibler}. From the definition of the Kullback--Leibler divergence it follows that $$\liminf_{t \to \infty} \frac{1}{t} \E_\vartheta \bigg[\log \frac{\widetilde q_t^{\vartheta}(Y^{\mathrm{stat}}_\oo(1), \dots, Y^{\mathrm{stat}}_\oo(t))}{\widetilde q_t^\theta(Y^{\mathrm{stat}}_\oo(1), \dots, Y^{\mathrm{stat}}_\oo(t))}\bigg] = \liminf_{t \to \infty} \frac{1}{t}\E_\vartheta\big[\widetilde L^{\vartheta}_Y(t) - \widetilde L^{\theta}_Y(t)\big] > 0,$$ which is equivalent to $Q(\vartheta) > Q(\theta)$ by Proposition \ref{prop: equal_limits}. This finishes the proof.
\end{proof}

We are now finally ready to assemble all pieces collected throughout the preceding sections to prove our main theorems. To this end, we establish the conditions from Propositions \ref{prop: Jacod-Consistency} and \ref{prop: Jacod-Normality}:

\begin{proof}[Proof of Theorems \ref{theo: main} and \ref{theo: main2}] We have $\sup_{\theta \in\Theta} \lVert \frac{1}{t}L^\theta(t) - Q(\theta) \rVert \xrightarrow{\PP_\vartheta} 0$ by Remark \ref{rem: unif_conv_to_Q} for some continuous $Q$. Since $Q$ has a unique maximum at $\vartheta$ by Corollary \ref{coro: uniquely_maximised}, Proposition \ref{prop: Jacod-Consistency} establishes consistency of any sequence of quasi-maximum likelihood estimators. Moreover we have that $\smash{\frac{1}{\sqrt{t}}Z^\vartheta(t) \xrightarrow{\PP_\vartheta\text{-}d} Z \sim N(0, U_\vartheta)}$ by Theorem \ref{theo: norm} and $\sup_{\theta \in \Theta} \lVert \frac{1}{t}\nabla_\theta Z^\theta(t) - W(\theta)\rVert \smash{\xrightarrow{\PP_\vartheta} 0}$ by Proposition \ref{prop: uniform_Z}, where $W = \nabla_\theta^2 Q$ by Remark \ref{rem: unif_conv_to_Q}. Hence the conditions from Proposition \ref{prop: Jacod-Normality} are met if $W(\vartheta)$ is invertible. It follows that we indeed have convergence \begin{equation}\label{eq: Asymp_norm_in_Proof}
\smash{\sqrt{t}\bigl( \widehat \theta(t) - \vartheta\bigr) \xrightarrow{\PP_\vartheta\text{-}d}} -W(\vartheta)^{-1} Z \sim N(0, W(\vartheta)^{-1}U_\vartheta W(\vartheta)^{-1}).    
\end{equation}
It remains to prove that $\smash{\frac{1}{t} \nabla_\theta Z^\theta(t)\big|_{\theta = \widehat \theta(t)} \xrightarrow{\PP_\vartheta} W(\vartheta)}$. This follows by decomposing
\begin{align*}
    \big\lVert \textstyle\frac{1}{t} \displaystyle\nabla_\theta Z^\theta(t) \big|_{\theta = \widehat \theta(t)} - W(\vartheta) \big\rVert \leq \sup_{\theta \in \Theta} \big\lVert \textstyle\frac{1}{t}\displaystyle \nabla_\theta Z^\theta(t) - W(\theta) \big\rVert + \big\lVert W(\widehat \theta(t)) - W(\vartheta)\big\rVert, 
\end{align*}
where the first summand on the right tends to 0 in probability by Proposition \ref{prop: uniform_Z} and the second one by consistency and the continuous mapping theorem. 
\end{proof}

\begin{remark}
    Even though, from a theoretical point of view, the proof of Theorems \ref{theo: main} and \ref{theo: main2} is complete, the condition of invertibility of $W(\vartheta)$ imposed in the third statement is hard to verify because the true parameter value $\vartheta$ is of course unknown in practice. However, once $\vartheta$ has been consistently estimated by $\widehat \theta(t)$, $W(\widehat \theta(t))$ can be computed explicitly by the calculations in Section \ref{su:explicit} (if Assumption \ref{assump: AN} holds) in order to verify the invertibility assumption from Theorem \ref{theo: main2}. 
    
    A large part of the statistical literature on asymptotic normality of (quasi-)maximum likelihood estimators in hidden Markov models or general ARMA-type models a priori assumes an invertibility condition for some kind of expected Fisher information matrix such as $W(\vartheta)$. This is for example the case in \cite{Alj2017}, \cite{Bickel1998} or \cite{Melard2022}. By contrast, \cite{Schlemm2012} develop a sufficient condition for invertibility of $W(\vartheta)$, which, however, seems at least as hard to verify as invertibility of $W(\vartheta)$ itself. Indeed, it involves verifying a rank condition for a complicated Jacobian matrix evaluated at the true unknown $\vartheta$.
\end{remark}

\subsection{Proofs for Sections \ref{sec5.1: CovMat}, \ref{su:explicit} and \ref{su:tests}}

\begin{proof}[Proof of Proposition \ref{prop: Formula_Cov}]
    We first argue that
    $\E_\vartheta\big[\smash{\overline N}^\vartheta(t) \otimes \smash{\overline N}^\vartheta(t)|\langhomX(t-1)=x\big]=\smash{\overline Q}^\theta_{\otimes 2}\smash{(x\otimes x)}+\smash{\overline Q}^\theta x+\smash{\overline q}^\theta$ holds for any $x\in E\times\R^{(k+1)d}$, $t\in\N$, where $\smash{\overline N}^\vartheta$ is defined in Corollary \ref{convcoro}. This identity follows from a simple calculation and the fact that $\smash{\overline N}^\vartheta = (N^{\vartheta^\top}, 0, \dots, 0)^\top \in \R^{(k+2)d}$, yielding $\smash{\overline N}^\vartheta \otimes \smash{\overline N}^\vartheta = e_1^{(k+2)} \circ_{d\times 1}^{(k+2)\times 1} \big[ e_1^{(k+2)} \otimes N^\vartheta \otimes N^\vartheta\big]$.

    The $k$-dimensional quadratic polynomial $g^\theta$ solves the Poisson equation $f^\theta = \tildeP g^\theta - g^\theta$ (see the proof of Theorem \ref{theo: norm}), where $\tildeP$ denotes the transition operator of $\langhomX$ under $\PP_\vartheta$ and where $f^\theta$ given in equation \eqref{eq: f_def} is the $k$-dimensional quadratic polynomial that captures the dependence of $Z^\theta(t,t-1)$ on $\langX(t)$ but with limiting coefficients instead of time-dependent ones. Note that $\smash{\overline A}^\theta$ and $\smash{\overline a}^\theta$ are the state transition matrix resp.\ vector for $\langhomX$ and $\smash{\overline A}^\theta_{\otimes 2}$, $\smash{\overline a}^\theta_{\otimes 2}$ are the ones for $\mathrm{vec}_{\otimes 2}(\langhomX)$, in line with the notation in Corollary \ref{convcoro}, Remark \ref{rem: Radius_for_homX}, and Lemma \ref{range}. Let $\smash{\alpha^\theta_f}$ and $\smash{\beta^\theta_f}$ denote a pair of coefficients for $f^\theta$. From equation \eqref{eq: Score} it is apparent that $\smash{f^\theta(x)_j = \smash{x^\top \Gamma_\theta^{(j)} x} + (\beta_f^\theta)_j}$, i.e. $\smash{f^\theta(x)_j = \mathrm{vec}(\Gamma_\theta^{(j)}) (x \otimes x) + (\beta_f^\theta)_j}$. Hence $\smash{f^\theta(x) = \Gamma^\theta (x \otimes x) + \beta_f}$, so $\smash{\alpha^\theta_f = (0 \mid \Gamma^\theta)} \in \R^{k \times (k+2) d + (k+2)^2d^2}$. The proof of Lemma \ref{range} yields $\smash{\alpha^\theta_g = \alpha^\theta_f (\smash{\overline A}^\theta_{\otimes 2} - \mathrm{I})^{-1}}$. By writing
\begin{fitequation*}
    \langhomX(t+1)^{\otimes 2} = \smash{\overline a}^\theta \otimes \smash{\overline a}^\theta + (\smash{\overline A}^\theta \otimes \smash{\overline A}^\theta)\langhomX(t)^{\otimes 2} + \smash{\overline a}^\theta \otimes (\smash{\overline A}^\theta \langhomX(t))  + (\smash{\overline A}^\theta\langhomX(t)) \otimes \smash{\overline a}^\theta + \smash{\overline N}^\theta(t)^{\otimes 2}
\end{fitequation*}

\noindent and inserting in the definition of the matrices $\smash{\overline Q}^\theta_{\otimes 2}$, $\smash{\overline Q}^\theta$ and the vector $\smash{\overline q}^\theta$, we obtain the expressions \eqref{eq: Ahatotimes2_terms}. By Remark \ref{rem: Radius_for_homX} it follows once more that $\rho(O_\theta) < 1$. Note that
\begin{equation*}
\begin{pmatrix}
    A & 0 \\
    C & D    
\end{pmatrix}^{-1} = \begin{pmatrix}
    A^{-1} & 0 \\
    -D^{-1}CA^{-1} & D^{-1}
\end{pmatrix}
\end{equation*}
for any block matrix such that $A$ and $D$ are invertible, see \cite{Bernstein2011}, Proposition 3.9.7. Applying this formula to $(\smash{\overline A}^\theta_{\otimes 2} - \mathrm{I})^{-1}$, the given expression for $\alpha_g^\theta = (0 \mid \Gamma^\theta) (\smash{\overline A}^\theta_{\otimes 2} - \mathrm{I})^{-1}$ follows. Now set $h^\theta = \tildeP g^\theta$. Since $\tildeP \mathrm{vec}_{\otimes 2}(x) = \smash{\overline A}^\theta_{\otimes 2} \mathrm{vec}_{\otimes 2}(x) + \smash{\overline a}^\theta_{\otimes 2}$, the given expressions for $\alpha^\theta_h$ and $\beta^\theta_h$ also follow. Corollary \ref{coro: limiting_U} yields that $U_\vartheta$ is the limit in $\PP_\vartheta$-probability of $\widehat U_t[\vartheta]$, where
\begin{equation}\label{eq: utheta}
    \widehat U_t[\vartheta] = \frac{1}{t}\sum_{s=1}^t \Big[g^\vartheta\big(\langoX(s)\big)g^\vartheta\big(\langoX(s)\big)^\top - h^\vartheta\big(\langoX(s)\big)h^\vartheta\big(\langoX(s)\big)^\top\Big].
\end{equation}
To prove that the matrix $U_\vartheta$ is also the limit in $\PP_\vartheta$-probability of the expression \eqref{eq: Formula_Cov}, it remains to show that $U_\vartheta$ is also the limit of $\widehat U_t[\widehat \theta(t)]$ in place of $\widehat U_t[\vartheta]$. We can decompose
\begin{equation}\label{eq: utheta_decomp}
    \big\lVert \widehat U_t[\widehat \theta(t)] - \widehat U_t[\vartheta] \big\rVert \leq \sup_{\theta \in \Theta} \big\lVert \widehat U_t[\theta] - U_\theta \big\rVert + \big\lVert U_{\widehat \theta(t)} - U_\vartheta \rVert.
\end{equation}
Since each summand in the definition \eqref{eq: utheta} of $\widehat U_t[\theta]$ is a $(k\times k)$-dimensional polynomial of order 4 in $\langX(s)$, Lemma \ref{lem: PötscherPrucha} shows that the first summand on the right of \eqref{eq: utheta_decomp} converges to 0 in $\PP_\vartheta$-probability. And since $U_\theta$ is continuous in $\theta$ by the final statement in Corollary \ref{coro: ergod_oX}, the second summand converges to 0 in $\PP$-probability as $(\widehat \theta(t))_{t \in \N }$ is $\vartheta$-consistent.
\end{proof}

\begin{proof}[Proof of Theorem \ref{theo: main3}]
The first statement is \eqref{e:Uexplizit}.
The representation of $W(\vartheta)$ follows from how it is introduced in the proof of Proposition \ref{prop: uniform_Z}.
\end{proof}

\begin{proof}[Proof of Proposition \ref{p:slutsky}]
$W(\theta),U(\theta)$ and $V_\theta$ are continuous by Proposition \ref{prop: uniform_Z} and by the continuity statement in Corollary \ref{coro: ergod_oX}. Since $(\widehat \theta(t))_{t \in \N }$ is $\vartheta$-consistent, the first claim follows from the continuous mapping theorem. The second follows from Slutsky's theorem.
\end{proof}

\begin{proof}[Proof of Proposition \ref{prop: significance_tests}]
    For the Wald test statistic, note that the asymptotic normality property $\sqrt{t}\bigl(\widehat \theta(t) - \vartheta\bigr) \smash{\xrightarrow{\PP_\vartheta\text{-}d}} N(0, V_\vartheta)$ implies together with the Delta method that $$\smash{\sqrt{t} \big[R(\widehat \theta(t)) - R(\vartheta)\big] \stackrel{H_0}{=} \sqrt{t} \big[R(\widehat \theta(t)) - r\big] \smash{\xrightarrow{\PP_\vartheta\text{-}d}} N(0, \nabla_\theta R(\vartheta) V_\vartheta \nabla_\theta R(\vartheta)^\top)}.$$ The conclusion about the limiting distribution of the Wald test statistic hence follows from the continuous mapping theorem and Slutsky's theorem. We now turn to the Lagrange multiplier and likelihood-ratio test. We will denote $Z_t(\theta) \coloneq Z^\theta(t)$ and we will work on the set $\{\widehat \theta(t) \in \mathrm{int}(\Theta)\} \cap \{\widehat \theta^c(t) \in \mathrm{int}(\Theta)\} \cap \{\widehat V^c(t) \text{ invertible}\}$ whose probability converges to 1 under $H_0$. As already noted, the constrained estimator $\widehat \theta^c(t)$ solves the Lagrangian equations $\frac{1}{t} Z_t(\widehat \theta^c(t)) = \nabla_\theta R(\widehat \theta^c(t))^\top \widehat \lambda(t)$ and $R(\widehat \theta^c(t)) = r$ together with some $\R^m$-valued random variables $\widehat \lambda(t)$. (The proof of the Lagrange multiplier theorem constructs $\widehat \lambda(t)$ explicitly as a continuous function of $\widehat \theta^c(t)$ using the implicit function theorem, so $\widehat \lambda(t)$ is measurable and adapted.) Using Taylor's theorem we get
    \begin{align}\label{eq: LM_stat}
        0 = \frac{1}{\sqrt{t}} Z_t(\widehat \theta(t)) = \frac{1}{\sqrt{t}}Z_t(\widehat \theta^c(t)) + \frac{1}{t} \nabla_\theta Z_t(\widehat \theta^c(t)) \sqrt{t}\big[\widehat \theta(t) - \widehat \theta^c(t)\big] + r(t),
    \end{align}
    where the $j$-th element of $r(t)$ is given by $\frac{1}{2}\big[\widehat \theta(t) - \widehat \theta^c(t)\big]^\top \frac{1}{t}\nabla_\theta^2 Z_t(\xi_t^j)_j \sqrt{t}\big[\widehat \theta(t) - \widehat \theta^c(t)\big]$ for $j \in \{1, \dots, k\}$ and where the (random) $\xi_t^j$ lies on the line segment between $\widehat\theta(t)$ and $\widehat\theta^c(t)$. Since Assumption \ref{assump: A} holds with four times continuous differentiability, a similar argument as in Section \ref{sec4.2.2: ULLN} and Lemma \ref{lem: PötscherPrucha} yields that $\smash{\sup_{\theta \in \Theta} \big\lVert \frac{1}{t}\nabla_\theta^2 Z_t(\theta)_j - H_j(\theta)\big\rVert \xrightarrow{\PP_\vartheta} 0}$ for some continuous $H_j: \Theta \to \R^{k \times k}$, and so $\frac{1}{t}\nabla_\theta^2 Z_t(\xi_t^j)_j \xrightarrow{\PP_\vartheta} H_j(\vartheta)$ because both $\widehat \theta(t)$ and $\widehat \theta^c(t)$ are consistent under $H_0$. By joint asymptotic normality, $\sqrt{t} \big[\widehat \theta(t) - \widehat \theta^c(t)\big]$ converges in distribution, so $\smash{r(t) \xrightarrow{\PP_\vartheta} 0}$. Then we get $\sqrt{t}\big[R(\widehat \theta(t)) - r\big] = \sqrt{t}\big[R(\widehat \theta(t)) - R(\widehat \theta^c(t))\big] = \sqrt{t} \nabla_\theta R(\widehat\theta^c(t))\big[\widehat \theta(t) - \widehat \theta^c(t)\big] + \widetilde r(t)$ with some $\smash{\widetilde r(t) \xrightarrow{\PP_\vartheta} 0}$ by applying Taylor's theorem. Inserting \eqref{eq: LM_stat}, we obtain
    \begin{align*}
        \sqrt{t}\big[R(\widehat \theta(t)) - r\big] &=  \nabla_\theta R(\widehat\theta^c(t)) \widehat W^c(t)^{-1} \Big[-\frac{1}{\sqrt{t}} Z_t(\widehat \theta^c(t)) - r(t)\Big] + \widetilde r(t) \\
        &= -\nabla_\theta R(\widehat\theta^c(t)) \widehat W^c(t)^{-1} \sqrt{t} \nabla_\theta R(\widehat\theta^c(t))^\top \widehat \lambda(t) + \mathrm{o}_{\PP_\vartheta}(1),
    \end{align*}
    where we used that $\frac{1}{\sqrt{t}} Z_t(\widehat \theta^c(t)) = \sqrt{t} \nabla_\theta R(\widehat \theta^c(t))^\top \widehat \lambda(t)$ from the Lagrangian equations and where $\mathrm{o}_{\PP_\vartheta}(1)$ denotes a shorthand for any random variable converging to 0 in $\PP_\vartheta$-probability. Thus $\sqrt{t}\widehat \lambda(t) \smash{\xrightarrow{\PP_\vartheta\text{-}d}} N(0, \Sigma_W^{-1} \Sigma_V \Sigma_W^{-1})$, and accordingly it follows that $t \widehat \lambda(t)^\top \Sigma_{\widehat W} \Sigma_{\widehat V}^{-1} \Sigma_{\widehat W}\widehat \lambda(t)$ converges in distribution to a $\chi^2$-distributed random variable with $m$ degrees of freedom, where we define $\smash{\Sigma_{\widehat W} = \nabla_\theta R(\widehat \theta^c(t)) \widehat W^c(t)^{-1} \nabla_\theta R(\widehat \theta^c(t))^\top}$. Inserting the Lagrangian equations, the latter is equal to the Lagrange multiplier test statistic $\xi_{\mathrm{LM}}$, proving the second claim. For the likelihood ratio test, first, note that $(-W(\vartheta))^{-\frac{1}{2}}$ is well defined if $W(\vartheta)$ is invertible since then $-W(\vartheta)$ is positive definite since $W(\vartheta) = \nabla^2_\theta Q(\vartheta)$. Then, let $L_t(\theta)$ stand for $L^\theta(t)$ and apply Taylor's theorem as before to obtain
    \begin{align*}
        L_t(\widehat \theta^c(t)) = L_t(\widehat \theta(t)) + \frac{1}{2} \big[\widehat \theta(t) - \widehat \theta^c(t)\big]^\top \nabla_\theta Z_t(\widehat \theta(t)) \big[\widehat \theta(t) - \widehat \theta^c(t)\big] + \mathrm{o}_{\PP_\vartheta}(1),
    \end{align*}
    so $\xi_\mathrm{LR}(t) = 2\big[ L_t(\widehat \theta(t)) - L_t(\widehat \theta^c(t))\big] = -\sqrt{t} \big[\widehat \theta(t) - \widehat \theta^c(t)\big]^\top \frac{1}{t}\nabla_\theta Z_t(\widehat \theta(t)) \sqrt{t} \big[\widehat \theta(t) - \widehat \theta^c(t)\big] + \mathrm{o}_{\PP_\vartheta}(1)$. By \eqref{eq: LM_stat} we have $\sqrt{t} \big[\widehat \theta(t) - \widehat \theta^c(t)\big] = - \widehat W^c(t)^{-1} \nabla_\theta R(\widehat \theta^c(t))^\top \sqrt{t}\widehat \lambda(t) + \mathrm{o}_{\PP_\vartheta}(1)$. Using $\sqrt{t}\widehat \lambda(t) \smash{\xrightarrow{\PP_\vartheta\text{-}d}} \lambda$ with $\lambda \sim N(0, \Sigma_W^{-1} \Sigma_V \Sigma_W^{-1})$, this yields that $\xi_\mathrm{LR}(t)$ converges in $\PP_\vartheta$-distribution to $\lambda^\top \nabla_\theta R(\vartheta) [-W(\vartheta)]^{-1} \nabla_\theta R(\vartheta)^\top \lambda$. The result then follows by standard arguments concerning quadratic forms in normal vectors, see for example \cite{VanDerVaart1998}, Lemma 17.1.
\end{proof}

\subsection{Proofs for Section \ref{sec5: app}}

\begin{proof}[Proof of Proposition \ref{prop: OU_lp}]
     Fix $t \in \R_+$ and $x \in \R^d$. By Theorem \ref{theo: Sp-Hp-domination} there exists some $c_p > 0$ such that $\big\lVert \int_0^t \e^{Qs} \dd L(s) \big\rVert_{L^p} \leq \big\lVert \int \e^{Q\cdot\mathrm{id}} \mathbf{1}_{[0, t]} \dd L \big\rVert_{S^p} \leq c_p \big\lVert \int \e^{Q\cdot\mathrm{id}} \mathbf{1}_{[0, t]} \dd L \big\rVert_{H^p}$. Let $L = M^L + A^L$ be the special semimartingale decomposition of $L$ into the local martingale part $M^L$ and the predictable part $A^L$ of finite variation with $M^L(0)=A^L(0)=0$. Note that $A^L(t) = a^L t$ almost surely for $a^L = \E[L(1)]$, see \cite{Eberlein2019}, Theorem 2.21. Using Émery's inequality from Theorem \ref{theo: Emery} as well as the Burkholder--Davis--Gundy inequality from Theorem \ref{theo: BDG}, we obtain
    \begin{align*}
        \Big\lVert \int_0^t \e^{Qs} \dd L(s) \Big\rVert_{L^p}^p &\leq c \:\Big\lVert \int \e^{Q\cdot\mathrm{id}} \mathbf{1}_{[0, t]} \dd L \Big\rVert_{H^p}^p \leq c \: \big\lVert \e^{Q\cdot\mathrm{id}} \mathbf{1}_{[0, t]}\big\rVert^p_{S^\infty} \big\lVert L \cdot\mathbf{1}_{[0, t]} \big\rVert^p_{H^p} \\
        &\leq c \: \Big(\sup_{s\leq t} \big\lVert \e^{Qs}\big\rVert^p \Big)\Big\lVert \lVert [M^L, M^L](t)\rVert^{\frac{1}{2}} + \lVert A^L\rVert(t) \Big\rVert_{L^p}^p \\
        &\leq c\: \Big[ \E\big( \lVert [M^L, M^L](t)\rVert^{\frac{p}{2}} \big) + \E\big(\lVert A^L\rVert(t)^p \big) \Big] \\
        &\leq c\: \Big[ \E\Big( \sup_{s\leq t} \lVert M^L(s)\rVert^p \Big) + \lVert a^L\rVert^p \,t^p \Big],
    \end{align*}
    where the constant $c$ may vary from term to term.
    Since $M^L(t) = L(t) - a^Lt$, we have that the term $\E\big( \sup_{s\leq t} \lVert M^L(s)\rVert^p \big)$ is finite iff $\E\big( \sup_{s\leq t} \lVert L(s)\rVert^p \big) < \infty$. Applying Theorem 25.18 from \cite{Sato1999} to the submultiplicative function $x \mapsto (\lVert x \rVert \lor 1)^p$, this is equivalent to $\E\big( \lVert L(1) \rVert^p\big) < \infty$, which proves $\E \big( \big\lVert \int_0^t \e^{Qs} \dd L(s) \big\rVert^p\big) < \infty$. Together with Proposition \ref{prop: OU_strongsolution} this implies the first claim. For the second suppose that $\alpha(Q) > 0$ and define $M_j \coloneq \int_{j}^{j + 1} \e^{Q(s - j)} \dd L(s)$. By the just proven fact and since $L$ is a Lévy process, $M_j$,  $j \in \N$ are independent and identically distributed random variables with $\E(\lVert M_0 \rVert^p) < \infty$. This yields (assuming w.l.o.g. that $t \in \N$)
    \begin{align*}
        \Big\lVert \int_0^t \e^{-Q(t-s)} \dd L(s) \Big\rVert_{L^p} &= \Big\lVert \sum_{j=0}^{t-1} \e^{-Q(t-j)} M_j \Big\rVert_{L^p} \leq \lVert M_0 \rVert_{L^p} \sum_{j=0}^{t-1} \big\lVert \e^{-Q(t-j)}\big\rVert \\
        &= \lVert M_0 \rVert_{L^p} \sum_{j=1}^{t} \big\lVert \e^{-Qj}\big\rVert.
    \end{align*}
    Elementary properties of the matrix exponential show that $\rho(\e^{-Q}) < 1$ if $\alpha(Q) > 0$. By Corollary \ref{linalg} we conclude that $\lVert \e^{-Qj}\rVert \to 0$ at a geometric rate. A fortiori, the sum on the right hand side above is bounded in $t$, which proves that $X$ is bounded in $L^p$ by using the explicit representation of $X$ from Proposition \ref{prop: OU_strongsolution}.
\end{proof}

\begin{proof}[Proof of Proposition \ref{prop: OU_expressions}]
    Since $\E_\theta(X(t)) = \E_\theta(X(0)) + a^L_\theta t - \int_0^t Q^\theta \E_\theta(X(s)) \dd s$, we obtain that $\mu^\theta(t) \coloneq \E_\theta(X(t))$ satisfies the ordinary differential equation $\frac{\mathrm{d}}{\mathrm{d}t} \mu^\theta(t) = a^L_\theta - Q^\theta \mu^\theta(t)$ together with $\mu^\theta(0) = \E_\theta(X(0))$. It follows easily from the properties of the matrix exponential that the unique solution to this initial value problem is given by $\mu^\theta(t) = \e^{-Q^\theta t} \mu^\theta(0) + (\mathrm{I}_d - \e^{-Q^\theta t}) \mu^\theta(\infty)$, from which the expressions for $A^\theta$ and $a^\theta$ follow. To derive \eqref{eq: OU_BB}, we first search for an expression for $S^\theta(t) \coloneq \E_\theta(X(t) X(t)^\top)$. Integration by parts for matrix-valued semimartingales yields
    \begin{fitequation}\label{eq: dXX_OU}
        \mathrm{d}(X(t) X(t)^\top) = -\mathrm{Sym}(Q^\theta X(t) X(t)^\top) \dd t + X(t) \dd L(t)^\top + \mathrm{d} L(t) X(t)^\top + \mathrm{d}[L, L](t).
    \end{fitequation}
    As in \cite{Eberlein2019}, Example 3.3, one has $\E_\theta([L, L](t)) = \mathrm{Cov}_\theta(L(t)) = c^L_\theta t$. Moreover, the stochastic integrals in \eqref{eq: dXX_OU} are of the form $\int_0^t \dd L(s) \, X(s)^\top = a^{L}_\theta \, \int_0^t X(s)^\top \dd s + \int_0^t \dd \widetilde L^\theta(s) \, X(s)^\top$, where $\widetilde L^\theta(t) = L(t) - a^L_\theta t$ is a centred Lévy process under $\PP_\theta$ and hence a martingale with finite second moments by the integrability assumption and by \cite{Eberlein2019}, Theorem 2.21. It follows from \cite{Jacod2003}, Theorem III.6.4(d), that the process $\int \dd \widetilde L^\theta \, X^\top$ is a square-integrable martingale on any interval $[0, t]$ with $t \geq 0$ because $X(t)$ has finite second moments for any $t \in \R_+$ by Proposition \ref{prop: OU_lp}, and so $\E_\theta(\int_0^t \dd L(s) \, X(s)^\top) = a_\theta^L \, \int_0^t \E_\theta[X(s)]^\top \dd s$. It follows that $S^\theta$ satisfies the matrix-valued initial value problem
    \begin{equation*}
        \frac{\mathrm{d}}{\mathrm{d}t} S^\theta(t) = \mathrm{Sym}\big[a^L_\theta \mu^\theta(t)^\top -  Q^\theta S^\theta(t)\big] + c^L_\theta
    \end{equation*}
    with $S^\theta(0) = \E_\theta(X(0) X(0)^\top)$. This matrix-valued Sylvester equation can be solved using the results in \cite{Fausett2009} to yield 
    \[S^\theta(t) = \e^{-Q^\theta t} S^\theta(0) \e^{-Q^{\theta^\top} t} + \int_0^t \e^{-Q^\theta (t - s)} \big[\mathrm{Sym}(a^L_\theta \mu^\theta(s)^\top) + c^L_\theta\big] \e^{-Q^{\theta^\top} (t - s)} \dd s.\]
    Using the fact that $M(t) \coloneq \int_0^t \e^{As} C \e^{B s} \dd s$ satisfies the equation $\e^{At} C \e^{B t} = A M(t) + M(t) B + C$ for any square matrices $A$, $B$ and $C$, it is then possible to write down a closed-form expression for $S^\theta$ involving the Kronecker sum and product. To obtain \eqref{eq: OU_BB}, one just needs to simplify the expression $C^\theta(t) = S^\theta(t) - A^\theta S^\theta(t) A^{\theta^\top} - \mathrm{Sym}[A^\theta \mu^\theta(t) a^{\theta^\top}] - a^\theta a^{\theta^\top}$, which follows easily from expanding the definition of $C^\theta(t) = \mathrm{Cov}_\theta\big[X(t) - \E_\theta(X(t) \mid X(t-1)) \big]$.
\end{proof}

\begin{appendix}
\section{Tools}\label{appn}

\subsection{Ergodic theorems for Markov chains}\label{appA: Markov}

This section provides a quick overview over some notions from the ergodic theory of discrete-time Markov chains that are needed for the proof of the ergodicity propositions \ref{prop: ergod} and \ref{ergod2} as well as for the proof of the geometric ergodicity property \ref{prop: geom_ergod_y}. It closely follows the standard work of \cite{Meyn2009}. For the course of this section, $X = (X(t))_{t \in \N}$ denotes a Markov chain with values in $E$, where $(E, \mathscr{E})$ is a measurable space with countably generated $\sigma$-algebra, and $P(x, \cdot) = P_1(x, \cdot)$ stands for the one-step transition measures of $X$, $x \in E$, that fully determine the transition semigroup $(P_t)_{t \in \N}$, which is connected to $X$ via $\E(f(X(s+t))|\F_s)=(P_tf)(X(s))$ for all $s,t\in\N$ and bounded, measurable $f:E\to\R$. As usual, we use the slightly ambiguous notation $P_t(x, A)\coloneq (P_t1_A)(x)=\PP_x(X(t)\in A)$ for the $t$-step transition function, where $\PP_\nu$ denotes the measure under which $X(0) \sim \nu$ and $\PP_x \coloneq \PP_{\varepsilon(x)}$ for $x \in E$.
We first define the necessary notions for the following ergodicity theorems.

$X$ is called \textbf{\textit{irreducible}} with respect to a non-trivial measure $\varphi$ on $\mathscr{E}$ if, for all $x \in E$ and all $A \in \mathscr{E}$ with $\varphi(A) > 0$, there exists some $t \in \N^*$ such that $P_t(x, A) > 0$. If $X$ is irreducible with respect to some measure $\varphi$, \cite{Meyn2009}, Proposition 4.2.2, yields that there exists some probability measure $\psi$ on $\mathscr{E}$ such that, for any other measure $\varphi'$, $X$ is irreducible with respect to $\varphi'$ if and only if $\varphi'$ is absolutely continuous with respect to $\psi$. In this case $\psi$ is termed a \textbf{\textit{maximal irreducibility measure}} for $X$ and $X$ is simply called \textbf{\textit{$\bm{\psi}$-irreducible}}.

If $X$ is a $\psi$-irreducible Markov chain, then there exists some $k \in \N^*$ and disjoint sets $D_1, \dots, D_k \in \mathscr{E}$ such that $\psi\big[ \big(\bigcup_{i=1}^k D_i\big)^c\big] = 0$ and $P(x, D_{i + 1}) = 1$ for any $x \in D_i$ and $i = 0, \dots, k - 1$ (mod $k$), see \cite{Meyn2009}, Theorem 5.4.4. $D_1, \dots, D_k$ is termed a \textbf{\textit{$\bm{k}$-cycle}} for $X$. If $k = 1$ is the only integer for which there exists a $k$-cycle for $X$, then $X$ is termed \textbf{\textit{aperiodic}}.

In line with the Feller property for continuous-time Markov processes we define a weak Feller property for discrete-time Markov chains. From now on, suppose that $E$ is equipped with a locally compact, separable and metrisable topology and that $\mathscr{E}$ is the corresponding Borel $\sigma$-algebra. Then $X$ is called a \textbf{\textit{weak Feller chain}} if $\smash{P(x, \cdot) \xrightarrow{x \to x_0} P(x_0, \cdot)}$ weakly for all $x_0 \in E$, see also Proposition 6.1.1 in \cite{Meyn2009}. A set $A \in \mathscr{E}$ is called \textbf{\textit{petite}} if there exists some non-trivial measure $\nu$ on $\mathscr{E}$ and $(\lambda_n)_{n \in \N} \in [0, 1]^{\N}$ with $\sum_{n =0}^\infty \lambda_n = 1$ such that $\sum_{n =0}^\infty \lambda_n P_n(x, B) \geq \nu(B)$ holds for all $x \in A$ and all $B \in \mathscr{E}$. The Markov chain $X$ is called \textbf{\textit{bounded in probability on average}} if for each $x \in E$ the sequence $(\overline{P}_n(x, \cdot))_{n \in \N^*}$ is tight, where $\overline{P}_n(x, \cdot) \coloneq \frac{1}{n}\sum_{k=1}^n P_k(x, \cdot)$. If $X$ is $\psi$-irreducible, we call $X$ a \textbf{\textit{T-chain}} if any compact subset of $E$ is petite, see Theorem 6.2.5 in \cite{Meyn2009}.

The notion of $\psi$-irreducibility of a Markov chain entails the behaviour that any set in the support of $\psi$ is visited with positive probability at some future point in time. The property of \textbf{\textit{Harris recurrence}} strengthens this condition by defining that a $\psi$-irreducible Markov chain $X$ is called Harris recurrent if $\PP_x(X \in A \text{ infinitely often}) = 1$ for any $A \in \mathscr{E}$ with $\psi(A) > 0$ and any $x \in A$. Additionally, we call a Harris recurrent Markov chain $X$ \textbf{\textit{positive Harris}} if it admits an invariant probability measure.

Using the above concepts we are now ready to cite the proof of an ergodicity theorem for $L^p$-bounded Markov chains, which forms the basis of any ergodicity proof in Section \ref{sec4: Estimation}. In order to do so, we further restrict the choice of state spaces $E$ to closed subsets of $\R^d$.

\begin{theorem}\label{theo: f-ergod}
    Let $X$ be an $E$-valued $\psi$-irreducible, aperiodic weak Feller chain such that $\mathrm{supp}(\psi)$ has non-empty interior. If $X$ is bounded in $L^p(\PP_x)$ for some $p \in [1,\infty)$ and any $x \in E$, then there exists a unique stationary distribution $\mu$ for $X$. Moreover, $X$ is strongly $f$-ergodic with respect to $\mu$ under any $\PP_\nu$ and for any $\mu$-integrable $f$ we have
    \[\lim_{t\to\infty}\E_x\big[f(X(t))\big] = \int f \dd \mu\quad\mbox{for $\psi$-almost every $x \in E$}.\]
    If $\big(\E_x[f(X(t))]\big)_{t \in \N}$ is bounded for some $x$, then $f$ is $\mu$-integrable.
\end{theorem}
\begin{proof}
    To show that there exists an invariant probability measure $\mu$ on $\mathscr{E}$, it suffices to verify that the distributions of $X$ under any $\PP_x$ are tight, i.e.\ for any $\varepsilon > 0$ and any $x \in E$ there exists a compact $C \subseteq E$ such that $\limsup_{t \in \N^*} \PP_x(X(t) \not \in C) \leq \varepsilon$. Indeed, then $X$ is bounded in probability on average, and any weak Feller chain bounded in probability on average admits an invariant probability measure, see \cite{Meyn2009}, Theorem 12.1.2(ii). Since $E \subseteq \R^d$ is closed and any compact set in the relative topology on $E$ is also compact in $\R^d$, it suffices to show that $\limsup_{t \in \N^*} \PP_x(\lVert X(t)\rVert \geq M)$ converges to 0 as $M \to \infty$. But
    \[\limsup_{t \to \infty}\PP_x(\lVert X(t)\rVert \geq M) \leq \limsup_{t \to \infty}\frac{\E_x(\lVert X(t) \rVert^p)}{M^p} \xrightarrow{M \to \infty} 0\]
     by Markov's inequality because $X$ is bounded in $L^p$. Moreover, $X$ is Harris-recurrent. To wit, note that $X$ is a $\psi$-irreducible T-chain, which follows from Theorem 6.0.1(iii) of \cite{Meyn2009} because $\mathrm{supp}(\psi)$ has non-empty interior. Moreover, $\smash{\PP_x(\lVert X(t) \rVert \xrightarrow{t \to \infty} \infty) = 0}$ for any $x \in E$ because otherwise $X$ could not be $L^p$-bounded by Fatou's lemma. Theorem 9.0.2 in \cite{Meyn2009} then yields that $X$ is Harris-recurrent.

     Together with the already shown properties, this yields that $X$ is an aperiodic positive Harris chain. In particular, Theorems 13.3.3, 14.0.1, and 17.0.1 in \cite{Meyn2009} yield, for any $\mu$-integrable function\footnote{Theorem 14.0.1 in \cite{Meyn2009} yields convergence of expectations only for $\mu$-integrable functions $f \geq 1$. This however poses no problem by considering positive and negative parts and shifting.} $f$, that $\smash{\PP_\nu^{X(t)} \xrightarrow{w} \mu}$ and $\frac{1}{t}\sum_{s=1}^t f(X(s)) \to \int f \dd \mu$ almost surely under any measure $\PP_\nu$, and $\E_{x}\big[f(X(t))\big] \to \int f \dd \mu$ for $\psi$-almost all $x \in E$. If $f$ is such that the above sequence of expectations is bounded for some $x \in E$, Theorem 14.3.3(i) in \cite{Meyn2009} yields moreover that the function $f$ is $\mu$-integrable.
\end{proof}

\begin{remark}
    Theorem 13.3.3 in \cite{Meyn2009}, which yields the weak convergence $\PP_\nu^{X(t)} \xrightarrow{w} \mu$ in the above proof, also yields the stronger result that $\PP_\nu^{X(t)}$ converges to $\mu$ with respect to the total variation norm. Since we do not need this stronger concept of convergence with respect to the total variation norm, we confine ourselves to weak convergence as in Theorem \ref{theo: f-ergod}.
\end{remark}

We end this section with the well-known Foster--Lyapunov sufficient condition for the concept of \textit{geometric ergodicity} of a Markov chain $X$ in the following sense:
\begin{definition}\label{defi: geom_ergod}
    Define $\lVert \nu \rVert_f \coloneq \sup\{|\int g \dd \nu|: g \text{ measurable with } |g| \leq f\}$ for any signed measure $\nu$ on $\mathscr{E}$, where $f: E \to [1,\infty)$ is measurable. A Markov chain $X$ with semigroup $(P_t)_{t\in\N}$ and invariant measure $\mu$ is \textbf{$\bm{f}$-geometrically ergodic} if there are $c \in \R_+$, $\gamma \in [0, 1)$ with
    \[\big\lVert P_t(x, \cdot)  - \mu\big\rVert_f \leq c f(x) \gamma^{t}.\]
\end{definition}
We need the following result in the context of an identifiability condition in Section \ref{sec4.2.3: Ident}.

\begin{theorem}\label{theo: geom-ergod}
    Let $X$ be an $E$-valued $\psi$-irreducible, aperiodic weak Feller chain with  transition operator $P$ and unique invariant probability measure $\mu$ such that $\mathrm{supp}(\psi)$ has non-empty interior. Then $X$ is $f$-geometrically ergodic for any measurable $f: E \to [1,\infty)$ fulfilling the Foster--Lyapunov condition $Pf - f \leq \beta f + b \mathbf{1}_C$ with some compact $C \subseteq E$ and $b \in \R$, $\beta < 0$.
\end{theorem}
\begin{proof}
    Since $X$ is a $\psi$-irreducible weak Feller chain such that $\mathrm{supp}(\psi)$ has non-empty interior, Proposition 6.2.8 in \cite{Meyn2009} implies that $X$ is a $\psi$-irreducible T-chain, i.e.\ every compact set $C\subseteq E$ is petite. The remaining part follows from \cite{Meyn2009}, Theorem 15.0.1.
\end{proof}

\subsection{Results from matrix analysis and linear systems theory}\label{appB: Matrix}

In this section we provide several useful results from linear algebra, matrix analysis, and dynamic linear systems theory. The most important result proved here is the stability lemma \ref{lem: lin_systems}, which provides (uniform) convergence conditions for linear and Lyapunov-type functional discrete dynamical systems. We start with a simple bound for inverse matrices. Recall that $\lVert \cdot\rVert$ denotes the spectral norm for matrices.

\begin{lemma}\label{lem: strange_bound}
    Let $A \in \R^{d\times d}$ be a non-singular matrix. Then $\lVert A^{-1}\rVert \leq \frac{\lVert A\rVert^{d-1}}{|\det(A)|}$.
\end{lemma}
\begin{proof}
    Let $\sigma_1, \dots, \sigma_d$ denote the singular values of $A$ in decreasing order and recall that this implies $|\det(A)| = \prod_{k=1}^d \sigma_d$, that $\norm A = \sigma_1$, and that $\lVert A^{-1}\rVert^{-1} = \sigma_d$. Then
    \begin{equation*}
        \lVert A^{-1}\rVert = \sigma_d^{-1} \leq \sigma_d^{-1} \bigg( \prod_{k=1}^{d-1} \frac{\sigma_1}{\sigma_k} \bigg) = \frac{\sigma_1^{d-1}}{\prod_{k=1}^{d} \sigma_k} = \frac{\lVert A\rVert^{d-1}}{|\det(A)|}.
    \end{equation*}
\end{proof}

The main content of the following lemmata consists in the fact that powers of a matrix $A$ converge geometrically fast to 0 whenever its spectral radius $\rho(A)$ is strictly smaller than 1. This basic result extends even to uniform convergence of matrix-valued functions.

\begin{lemma}\label{lem: linalg_unif_bounded}
    Let $E$ be a set and $A: E \to \R^{d\times d}$ bounded with $\rho(A(x)) \leq \alpha < 1$ for all $x \in E$. Then $A^m \to 0$ uniformly at a geometric rate, i.e.\ there are $c \in \R_+$ and $\gamma \in [0, 1)$ with $$\sup_{x\in E} \lVert A(x)^m\rVert \leq c \gamma^m, \qquad m \in \N^*.$$
\end{lemma}
\begin{proof}
    Fix some constant $C$ such that $\lVert A(x)\rVert \leq C$ for all $x\in E$ and some $0 < \varepsilon < \frac{1 - \alpha}{C}$. In \cite{Green1996}, Corollary 3.3.2, it is shown that $\smash{\big[\lVert B^m\rVert^{\frac{1}{m}} - \rho(B)\big] / \norm{B}}$ converges to 0 uniformly in $0 \neq B \in \R^{d \times d}$. Hence we can find $N \in \N^*$ such that for all $m \geq N$ and all $x \in E$ we have
    \[\lVert A(x)^m \rVert^{\frac{1}{m}} < \rho(A(x)) + \varepsilon\lVert A(x)\rVert \leq \alpha + C \varepsilon \eqqcolon \gamma < 1.\]
    for all $x \in E$. In particular, $\lVert A(x)^m\rVert < \gamma^m$ for all $x \in E$ and $m \geq N$. By defining the constant $c \coloneq \max\big\lbrace 1, \gamma^{-N}\max_{1 \leq n \leq N}\sup_{x \in E} \lVert A(x)^n\rVert\big\rbrace$, we obtain the desired result for all $m \in \N^*$.
\end{proof}

We can extend the preceding result to uniform convergence on compacta of products of matrix-valued functions which converge to a matrix-valued function fulfilling the assumptions of Lemma \ref{lem: linalg_unif_bounded} if it is additionally assumed that this limiting function is continuous.

\begin{lemma}\label{lem: linalg_unif}
    Let $E$ be a compact space and suppose that $A_t: E \to \R^{d\times d}$ is bounded for any $t\in \N^*$. Suppose also that $A: E \to \R^{d\times d}$ is a continuous function with $\rho(A(x)) < 1$ for all $x \in E$ and such that $\smash{A_t \xrightarrow{t\to \infty} A}$ uniformly on $E$. Then all chains in $(A_t)_{t \in \N^*}$ converge uniformly at a geometric rate to 0, i.e.\ there exist some $c \in \R_+$ and $\gamma \in [0, 1)$ such that $$\sup_{x \in E} \lVert A_{s+m}(x) A_{s+m-1}(x) \dots A_{s+1}(x) \rVert \leq c \gamma^m$$ for any $s \in \N$ and $m \in \N^*$. In particular, $A^m \to 0$ uniformly on $E$ at a geometric rate.
\end{lemma}
\begin{proof}
    Since $A$ is continuous, since the spectral radius is a continuous function on $\R^{d\times d}$, and since $E$ is compact, there exists some $\alpha \in [0, 1)$ such that $\rho(A(x)) < \alpha$ for all $x \in E$. It follows from \cite{Horn2012}, Lemma 5.6.10, that, for any $x \in E$, there exists a matrix norm $\lVert \cdot \rVert_x$ such that $\lVert A(x) \rVert_x < \alpha$. Then $D_x \coloneq \{y \in E: \: \lVert A(y)\rVert_x < \alpha\}$ is open and, since $x \in D_x$, $\bigcup_{x \in E} D_x$ is an open cover of $E$. By compactness of $E$ there exists a finite subcover $E = \bigcup_{k=1}^K D_{x_k}$ for some $x_1, \dots, x_K \in E$ and some $K \in \N^*$. In other words, for any $x \in E$ there exists a $k(x) \in \{1, \dots, K\}$ such that $\lVert A(x)\rVert_{k(x)} \coloneq \lVert A(x)\rVert_{x_{k(x)}} < \alpha$. Now fix some $0 < \varepsilon < 1 - \alpha$. By uniform convergence we deduce that for any $k \in \{1, \dots, K\}$ there is $N(k) \in \N^*$ such that $\sup_{x \in E}\lVert A_n(x) - A(x)\rVert_k < \varepsilon$ for all $n \geq N(k)$. Let $N\coloneq \max_{k \in \{1, \dots, K\}} N(k)$. Then \linebreak
    \[\lVert A_n(x) \rVert_{k(x)} \leq \lVert A(x)\rVert_{k(x)} + \lVert A_n(x) - A(x)\rVert_{k(x)} < \alpha + \varepsilon \eqqcolon \gamma < 1\]
    for all $n \geq N$ and all $x \in E$. Now let $\norm{\cdot}$ denote the spectral norm. By the equivalence of norms on $\R^{d\times d}$, there exist $C_1, \dots, C_K \in \R_+$ such that $\norm{\cdot} \leq C_k \norm{\cdot}_k$. Then, for any $n \geq N$,
    \[\norm{A_n(x) A_{n-1}(x) \dots A_N(x)} \leq C_{k(x)} \norm{A_n(x)}_{k(x)} \dots \norm{A_N(x)}_{k(x)} \leq \underbrace{\max_{k \in \{1, \dots, K\}} C_k}_{\eqqcolon C} \gamma^{n - N + 1},\]
    where $x \in E$. This proves the claim for all chains starting at $N$ or later. Finally, consider a chain $A_n(x) \dots A_{N - l}(x)$ of length $m$ starting at some index $N - l$ for $l < N$. Then
    \begin{align*}
        \norm{A_n(x) \dots A_{N - l}(x)} \leq C \sup_{x \in E} \prod_{k=1}^{N - 1} \norm{A_k(x)} \gamma^{m - l} \leq \underbrace{C \sup_{x \in E} \prod_{k=1}^{N - 1} \norm{A_k(x)} \gamma^{-N}}_{\eqqcolon C'} \gamma^m,
    \end{align*}
    where the supremum is finite by boundedness of the functions $A_k$. Setting $c \coloneq \max\{C, C'\}$ finishes the proof. The final comment follows by taking the constant sequence $\smash{A_t \equiv A}$.
\end{proof} 

If the compact space $E$ in the preceding lemma is chosen to be a singleton, the result translates verbatim to ordinary sequences of matrices instead of matrix-valued functions:

\begin{corollary}\label{linalg}
    Suppose that $(A_t)_{t \in \N^*}$ is a sequence of matrices that converges to some $A \in \R^{d \times d}$ with $\rho(A) < 1$. Then all chains in $(A_t)_{t \in \N^*}$ converge at a geometric rate to 0, i.e.\ there exist some constants $c \in \R_+$ and $\gamma \in [0, 1)$ such that $\norm{A_{s+m} A_{s+m-1} \dots A_{s+1}} \leq c \gamma^m$ for any $s \in \N$ and $m \in \N^*$. In particular, $A^m \to 0$ at a geometric rate as $m \to \infty$.
\end{corollary}

Alternatively, we can also uniformly bound partial derivatives of chains of differentiable matrix functions by a geometric rate, as the following corollary of Lemma \ref{lem: linalg_unif} shows:

\begin{corollary}\label{coro: linalg_diff}
    Let $E$ be a compact subset of $\R^k$ and suppose that functions $A_t: E \to \R^{d\times d}$ and $A: E \to \R^{d\times d}$ as in Lemma \ref{lem: linalg_unif} are given such that $A_t \in \mathrm{C}^1(E, \R^{d\times d})$ for any $t \in \N^*$ and $\partial_{x_j} A_t(x)$ is bounded uniformly in $t\in\N^*$ and $x\in E$ for any $j \in \{1, \dots, k\}$. Then there exist constants $c \in \R_+$ and $\gamma \in [0, 1)$ such that for any $j \in \{1, \dots, k\}$ $$\norm{\partial_{x_j} \big[A_{s+m}(x) A_{s+m-1}(x) \dots A_{s+1}(x)\big]} \leq c m \gamma^m$$ for any $s \in \N$ and $m \in \N^*$. In particular $\partial_{x_j} A^m \to 0$ uniformly on $E$ at a geometric rate.
\end{corollary}
\begin{proof}
    By the product rule $\partial_{j} \big[A_{s+m}(x) \dots A_{s+1}(x)\big] =: \partial_{j} A_{s, m}(x)$ has the form
    \begin{fitequation*}
    [\partial_{j} A_{s+m}(x)] A_{s, m-1}(x) + A_{s+m}(x)[\partial_{j} A_{s+m-1}(x)] A_{s, m-2,}(x) + \dots + A_{s + 1, m - 1}(x) [\partial_{j} A_{s+1}(x)].
    \end{fitequation*}
    If $M \geq 0$ is chosen such that $\lVert \partial_{x_j}A_t(x)\rVert \leq M$ for all $x \in E$, $t \in \N^*$ and $j \in \{1, \dots, k\}$, and if $c \in [1,\infty)$ and $\gamma \in [0, 1)$ are chosen as in Lemma \ref{lem: linalg_unif}, then $\|\partial_{j} \big[A_{s, m}(x)\big]\| \leq c^2 M m \gamma^{m-1} = (c^2 M \gamma^{-1}) m \gamma^m$. This finishes the proof.
\end{proof}

As already noted we will now state the most important lemma of this section concerning (uniform) convergence of linear and Lyapunov type functional discrete dynamical systems:
\begin{lemma}\label{lem: lin_systems}
    Let $E$ be a compact space. Suppose that $A_t: E \to \R^{d\times d}$ is bounded for any $t \in \N^*$ and uniformly convergent to some continuous function $A: E \to \R^{d\times d}$ with $\rho(A(x)) < 1$ for all $x \in E$. Moreover, let $B_t: E \to \R^{d\times d}$ and $B: E \to \R^{d\times d}$ as well as $b_t: E \to \R^d$ and $b: E \to \R^d$ be bounded functions. Consider the following linear systems:
    \begin{align}
        X_{t+1} &= A_t X_t A_t^\top + B_t, \label{eq: System1}\\
        X_{t+1} &= A X_t A^\top + B, \label{eq: System2}\\
        x_{t+1} &= A_t x_t + b_t, \label{eq: System3}\\
        x_{t+1} &= A x_t + b \label{eq: System4}
    \end{align}
    for $\R^{d\times d}$-valued functions $(X_t)_{t \in \N^*}$ and $\R^d$-valued functions $(x_t)_{t \in \N^*}$ on $E$ such that the initial values $X_1: E \to \R^{d\times d}$ and $x_1: E \to \R^{d}$ are bounded functions. \mbox{Then the following statements hold:}
    \begin{itemize}
        \item[1.] If the functions $(B_t)_{t \in \N^*}$, $B$, $(b_t)_{t \in \N^*}$, and $b$ are uniformly bounded on $E$, then the respective solutions $(X_t)_{t\in \N^*}$ and $(x_t)_{t\in \N^*}$ to (\ref{eq: System1}--\ref{eq: System4}) are uniformly bounded on $E$.
        \item[2.] The systems \eqref{eq: System2} and \eqref{eq: System4} possess unique fixed points $X^*: E \to \R^{d\times d}$ and $x^*: E \to \R^d$, respectively, and we have $X_t \to X^*$ as well as $x_t \to x^*$ uniformly on $E$ at a geometric rate. If $B$ is pointwise positive (semi-)definite, then $X^*$ is also pointwise positive \mbox{(semi-)}definite.
        \item[3.] If additionally $B_t \to B$ (respectively $b_t \to b$) and $A_t \to A$ uniformly at a geometric rate, then the system \eqref{eq: System1} (respectively \eqref{eq: System3}) converges uniformly at a geometric rate to the fixed point of \eqref{eq: System2} (respectively \eqref{eq: System4}).
    \end{itemize}
\end{lemma}

Similar to Corollary \ref{linalg}, Lemma \ref{lem: lin_systems} has an immediate counterpart for ordinary matrices instead of matrix-valued functions by taking $E$ to be a singleton, which we omit at this point.

\begin{proof}
    The systems \eqref{eq: System1} and $\eqref{eq: System2}$ can be reduced to \eqref{eq: System3} and \eqref{eq: System4} by vectorisation because e.g.\ $\mathrm{vec}(X_{t+1}) = (A \otimes A) \mathrm{vec}(X_t) + \mathrm{vec}(B)$ and because $\rho(A \otimes A) < 1$ whenever $\rho(A) < 1$. Hence, except for the definiteness in 2., it suffices to prove the claims 1.--3. for the vector-valued systems \eqref{eq: System3} and \eqref{eq: System4}. Moreover, by taking $A_t \equiv A$ and $b_t \equiv b$, it is sufficient to focus on \eqref{eq: System3}. By iterating this system we obtain the closed form solution
    \begin{equation}\label{eq: B_closedform}
        x_{t+1} = \Big(\prod_{s=1}^t A_s\Big) x_1 + \sum_{s=1}^t \Big(\prod_{r=s+1}^t A_r\Big) b_s.
    \end{equation}
    By Lemma \ref{lem: linalg_unif} we can bound the norms of the products in \eqref{eq: B_closedform} such that
    \begin{equation}\label{eq: B_closedform2}
        \norm{x_{t+1}(x)} \leq \widetilde c \gamma^t \norm{x_1(x)} + \widetilde c\sum_{s=1}^t \gamma^{t-s} \norm{b_s(x)},\quad x\in E
    \end{equation}
    holds for some constants $\widetilde c \in \R_+$ and $\gamma \in [0, 1)$. Since $b$ is uniformly bounded and since $x_1$ is also bounded, we obtain $\norm{x_{t+1}(x)} \leq c \gamma^t + c \sum_{s=0}^{t-1} \gamma^s$, $x \in E$ for some $c \in \R_+$, where the right-hand side is independent of $x$ and bounded in $t$. This proves 1.
    
    To show 2. for \eqref{eq: System4}, it is easy to see that $x^* \coloneq (\mathrm{I}_d - A)^{-1} b$ is the unique fixed point for \eqref{eq: System4}. In particular, since $A$ is continuous and $\rho(A(x)) < 1$ for all $x \in E$, there is $0 \leq \alpha < 1$ with $\rho(A(x)) \leq \alpha < 1$ for all $x \in E$, i.e.\ $\rho(A)$ is uniformly bounded below 1. It follows that the eigenvalues of $\mathrm{I}_d - A$ are uniformly bounded away from 0, which explains why the inverse $(\mathrm{I}_d - A(x))^{-1}$ is well defined for all $x \in E$. Moreover, by Lemma \ref{lem: strange_bound} we obtain the bound $\lVert (\mathrm{I}_d - A(x))^{-1}\rVert \leq \frac{\lVert \mathrm{I}_d - A(x)\rVert^{d-1}}{|\det(\mathrm{I}_d - A(x))|}$. Since the eigenvalues of $\mathrm{I}_d - A$ are uniformly bounded away from 0, also $|\det(\mathrm{I}_d - A)|$ is uniformly bounded away from 0. Since $\mathrm{I}_d - A$ is continuous on the compact set $E$, it is bounded and so $(\mathrm{I}_d - A)^{-1}$ is also bounded on $E$. Similar to \eqref{eq: B_closedform} we have
    \begin{equation}\label{eq: B_closedform3}
    x_{t+1} = A^t x_1 + \Big(\sum_{s=1}^{t} A^{t- s}\Big) b = A^t x_1 + (\mathrm{I}_d - A)^{-1} (\mathrm{I}_d - A^t) b,
    \end{equation}
    hence $\norm{x^*(x) - x_{t+1}(x)} \leq \lVert{A(x)^t}\rVert \norm{x_1} + \lVert (\mathrm{I}_d - A(x))^{-1}\rVert \lVert A(x)^t\rVert \norm{b} \leq c \lVert A(x)^t\rVert$ for some constant $c \in \R_+$. By Lemma \ref{lem: linalg_unif_bounded} it follows readily that $x_t \to x^*$ uniformly at a geometric rate. By \eqref{eq: B_closedform3} we can write $x^* = \big(\sum_{t=0}^\infty A^t\big) b$ for the fixed point of the system \eqref{eq: System4} and $X^* = \sum_{t=0}^\infty A^t B A^{t^\top}$ for the fixed point of the system \eqref{eq: System2}. The last identity proves that $X^*(x)$ is positive {(semi-)} definite whenever $B(x)$ is positive (semi-)definite, finishing 2.
    
    For 3. find $c \in \R_+$ as well as $\gamma \in [0, 1)$ such that $\lVert b_t(x) - b(x) \rVert \leq \gamma^t$ and $\lVert A_t(x) - A(x)\rVert \leq \gamma^t$. Suppose moreover that $\gamma$ and $c$ have been chosen large enough for chains in $(A_t)_{t \in N}$ and powers of $A$ to converge uniformly to 0 at the geometric rate $\gamma$ with constant $c$, which is guaranteed by Lemma \ref{lem: linalg_unif}.
    In the upcoming inequalities we use estimates of the form $ct\gamma^t \leq \widetilde c \widetilde \gamma^t$. Here, $\gamma$ and $\widetilde \gamma$ always denotes constants in $(0, 1)$. Let $(x_t)_{t \in \N^*}$ denote the solution to \eqref{eq: System3} and $(\overline{x}_t)_{t\in \N^*}$ the solution to \eqref{eq: System4}. Combining \eqref{eq: B_closedform} and \eqref{eq: B_closedform3} we obtain
    \begin{equation*}
        x(t+1) - \overline{x}(t+1) = \underbrace{\sum_{s=1}^t \Big[\Big(\prod_{r=s+1}^t A_r\Big) b_s - A^{t-s} b\Big]}_{(1)} + \underbrace{\Big(\prod_{s=1}^t A_s\Big) x_1  - A^t\overline{x}_1}_{(2)}.
    \end{equation*}
    We first focus on (1), which can be written as
    $\sum_{s=1}^t \big( \prod_{r=s+1}^t A_r - A^{t-s}\big) b_s + \sum_{s=1}^t A^{t-s} (b_s - b)$. The second sum in this decomposition is bounded in norm by $c \sum_{s=1}^t \gamma^{t-s} \gamma^s \leq \widetilde c \widetilde \gamma^t$ uniformly on $E$, which shows that it converges uniformly to 0 at a geometric rate. To deal with the first sum in the decomposition, a telescoping sum argument can be used to deduce
    \begin{equation*}
    \prod_{r=s+1}^t A_r - A^{t-s} = \sum_{r = s + 1}^t \Big[A^{t - r} (A_r - A) \prod_{u = s + 1}^{r - 1} A_u \Big].        
    \end{equation*}
    The last product is a chain of length $r - s - 1$ of the matrices $A_u$. The above arguments yield
    \begin{align}
        \Big\lVert \prod_{r=s+1}^t A_r(x) - A^{t-s}(x) \Big\rVert &\leq c \gamma^{t-s-1}\sum_{r = s+ 1}^t \lVert A_r(x) - A(x)\rVert \nonumber\\
        &\leq c \gamma^{t-s-1}(t-s) \gamma^{s+1} \leq c\gamma^t \label{eq: bound_for_product}
    \end{align}
    for all $x \in E$. Since $(b_t)_{t\in \N^*}$ is uniformly bounded on $E$ as a uniformly convergent sequence of bounded functions, it follows for the first summand in the decomposition of (1) that
    \begin{align*}
        \Big\lVert \sum_{s=1}^t \big( \prod_{r=s+1}^t A_r(x) - A^{t-s}(x)\big) b_s(x) \Big\rVert \leq c\gamma^t
    \end{align*}
    for all $x \in E$, which proves that (1) goes to 0 uniformly on $E$ at a geometric rate. Finally, it follows immediately from Lemma \ref{lem: linalg_unif} that (2) converges also to 0 uniformly at a geometric rate because the initial functions $x_1$ and $\overline{x}_1$ are bounded on $E$. All in all, we have uniform convergence to 0 of $x(t) - \overline{x}(t)$ at a geometric rate, which proves claim 3.
\end{proof}

The following lemma uses the Jordan canonical form to deduce that the convergence $A^n x \to 0$ entails that the map that $A$ induces on $E \subseteq \R^d$ has spectral radius less than 1:

\begin{lemma}\label{lem: power_convergence}
    Let $A \in \R^{d\times d}$ and let $E \subseteq \R^d$ such that $A^n x \to 0$ as $n \to \infty$ for any $x \in E$. Then there exists a matrix $\widetilde A \in \R^{d\times d}$ such that $\rho(\widetilde A) < 1$ and $\widetilde A x = Ax$ for any $x \in E$.
\end{lemma}
\begin{proof}
    Fix a basis $V = \{v_1, \dots, v_d\}$ of $\C^d$ such that $A$ is given in Jordan canonical form $J$ with respect to that basis, with $J$ being a complex block diagonal matrix of the form
    \[J = \begin{pmatrix}
        J_1 & 0 & \dots & 0\\
        0 & J_2 & \dots & 0 \\
        \vdots & \vdots & \ddots & \vdots \\
        0 & 0 & \dots & J_k
    \end{pmatrix} \in \C^{d\times d} \qquad \text{ with } \qquad J_i = \begin{pmatrix}
        \lambda_i & 1 & 0 & \dots & 0\\
        0 & \lambda_i & 1 & \dots & 0 \\
        \vdots & \vdots & \ddots & \ddots & \vdots \\
        0 & 0 & \dots & \lambda_i & 1 \\ 
        0 & 0 & \dots & 0 & \lambda_i \\ 
    \end{pmatrix} \in \C^{d_i \times d_i}\]
    for $i \in \{1, \dots, k\}$ and $k \in \{1, \dots, d\}$, where $\lambda_i$ are the eigenvalues of $A$. Since $d_i$ denotes the dimension of each Jordan block $J_i$, we have $m_k = d$ for $m_i = d_1 + \dots + d_i$. For any $x \in E$, fix the representation $x = k_1(x) v_1 + \dots + k_d(x) v_d$. Assume now that $A$ has some eigenvalue $\lambda_i$ with $|\lambda_i| \geq 1$. Since the map $A^n$ is given by $J^n$ with respect to the basis $V$, it follows that the components $m_{i-1} + 1, \dots, m_i$ of $A^n x$ are given by $J_i^n (k_{m_{i-1}}(x), \dots, k_{m_i}(x))^\top$ with respect to the basis $V$. Since $|\lambda_i| \geq 1$, this can only converge to 0 for all $x \in E$ if $k_{m_{i-1}}(x) = \dots = k_{m_i}(x) = 0$ for all $x \in E$. In this case the action of $A$ on any $x \in E$ is not altered if one replaces the Jordan block $J_i$ by $\widetilde J_i = 0 \in \C^{d_i \times d_i}$. For any eigenvalue $\lambda_j$ with $|\lambda_j| < 1$, define $\widetilde J_j \coloneq J_j$. Let
    \[\widetilde J \coloneq \begin{pmatrix}
        \widetilde J_1 & 0 & \dots & 0\\
        0 & \widetilde J_2 & \dots & 0 \\
        \vdots & \vdots & \ddots & \vdots \\
        0 & 0 & \dots & \widetilde J_k
    \end{pmatrix} \in \C^{d\times d} \qquad \text{ and } \qquad \widetilde A = V\widetilde J V^{-1}.\]
    By construction, $\rho(\widetilde A) < 1$ and $\widetilde Ax = Ax$ for all $x \in E$, which finishes the proof.
\end{proof}

The last theorem we state in this section provides an extension of the well-known mean-value theorem from elementary calculus to the case of differentiable matrix-valued functions:

\begin{theorem}\label{theo: matrix_mvt}
    Let $E \subseteq \R^k$ be a convex set and suppose that $F: E \to \R^{m\times n}$ is a continuously differentiable matrix-valued function. For arbitrary $y_1, y_2 \in E$ let $U(y_1, y_2)$ denote the line segment $U(y_1, y_2) \coloneq \{ty_1 + (1 - t)y_2: t \in [0, 1]\} \subseteq E$. Then
    \begin{equation*}
        \lVert F(y_1) - F(y_2) \rVert \leq \sqrt{\min\{m, n\}} \sup_{x \in U(y_1, y_2)} \sqrt{\lVert \partial_{x_1} F(x) \rVert^2 + \dots + \lVert \partial_{x_k} F(x) \rVert^2} \lVert y_1 - y_2\rVert
    \end{equation*}
    for any $y_1, y_2 \in E$.
\end{theorem}
\begin{proof}
    By definition we have $\lVert \mathrm{vec}(A)\rVert = \norm{A}_{\mathrm{Frob}}$, where $\norm{A}_{\mathrm{Frob}}$ denotes the Frobenius norm of a matrix. It follows that $\lVert F(y_1) - F(y_2) \rVert \leq \lVert \mathrm{vec}\,F(y_1) - \mathrm{vec}\,F(y_2) \rVert$ because $\lVert \cdot \rVert \leq \norm{\cdot}_{\mathrm{Frob}}$. By the mean value theorem for vector-valued functions we obtain
    \begin{equation}\label{eq: mvt1}
        \lVert F(y_1) - F(y_2) \rVert \leq \lVert \mathrm{vec}\,F(y_1) - \mathrm{vec}\,F(y_2) \rVert \leq \sup_{x \in U(y_1, y_2)} \lVert \nabla_x\mathrm{vec}\,F(x) \rVert \cdot \lVert y_1 - y_2 \rVert,        
    \end{equation} 
    where $\nabla_x\mathrm{vec}\,F(x)$ denotes the Jacobian matrix of $\mathrm{vec}\,F$ at $x$. Since the columns of $\nabla_x\mathrm{vec}\,F(x)$ are given by $\partial_{x_j} \mathrm{vec}\, F(x)$ for $j \in \{1, \dots, k\}$, it follows for any $x \in E$ that
    \begin{fitequation}\label{eq: mvt2}
        \lVert \nabla_x\mathrm{vec}\,F(x) \rVert \leq \lVert \nabla_x\mathrm{vec}\,F(x) \rVert_{\mathrm{Frob}} = \sqrt{\lVert\partial_{x_1} \mathrm{vec}\, F(x)\rVert^2 + \dots + \lVert\partial_{x_k} \mathrm{vec}\, F(x)\rVert^2}.        
    \end{fitequation}
    As $\norm{A}_{\mathrm{Frob}} \leq \sqrt{\min\{m, n\}} \norm{A}$ for $A \in \R^{m \times n}$, we have $\lVert\partial_{j} \mathrm{vec}\, F(x)\rVert = \lVert \mathrm{vec}(\partial_{j}F(x))\rVert = \lVert \partial_{j}F(x)\rVert_{\mathrm{Frob}} \linebreak\leq \sqrt{\min\{m, n\}} \lVert \partial_{j}F(x)\rVert$. Together with \eqref{eq: mvt1}, \eqref{eq: mvt2} yields the result.
\end{proof}

\subsection{Results from probability theory and analysis}\label{appC: Prob}

The purpose of this section is to collect some technical results from probability theory and from analysis. For this section, we fix a probability space $(\Omega, \F,\PP)$.

\begin{lemma}\label{lem: argmax}
    Let $E$ be a compact metric space and $f: E \to \R$ a continuous function such that $f(x^*) = \sup_{x \in E} f(x)$ holds for a unique maximiser $x^* \in E$. Let $\mathrm{arg}\max: \mathrm{C}(E, \R) \to \mathscr{P}(E)$ denote the set-valued function which maps a continuous function to the set of its maximisers. If $(f_n)_{n \in \N^*}$ is a sequence in $\mathrm{C}(E, \R)$ with $\lVert f_n - f\rVert_\infty \to 0$ and $(x_n)_{n \in \N^*}$ is a sequence in $E$ such that $x_n \in \mathrm{arg}\max(f_n)$ for any $n \in \N^*$, then $x_n \to x^*$ as $n \to \infty$.
\end{lemma}
\begin{proof}
    Since $E$ is a compact metric space, any sequence in $E$ has a subsequence converging in $E$, and since convergence to some element $x^*$ in a topological space is equivalent to the fact that any subsequence contains a further subsequence converging to $x^*$, it suffices to prove that any convergent subsequence of $(x_n)_{n \in \N^*}$ converges to $x^*$. Let $(x_{n_k})_{k \in \N^*}$ be a subsequence converging to some $x \in E$. We need to prove $x = x^*$. We can decompose
    \begin{align*}
        |f(x) - f(x^*)| &\leq |f(x) - f_{n_k}(x_{n_k})| + |f_{n_k}(x_{n_k}) - f(x^*)| \\
        &\leq |f(x) - f(x_{n_k})| + \lVert f_{n_k} - f \rVert_\infty + |f_{n_k}(x_{n_k}) - f(x^*)|.
    \end{align*}
    The first two summands on the right converge to 0 as $k \to \infty$ by continuity of $f$ and by uniform convergence. But also the third summand on the right above converges to 0 because
    \begin{align*}
        f_{n_k}(x_{n_k}) - f(x^*) \leq \hspace{0.3cm}f(x_{n_k}) - \hspace{0.3cm} f(x^*) + \lVert f_{n_k} - f\rVert_\infty \leq \hphantom{-}\lVert f_{n_k} - f\rVert_\infty \to 0, \\
        f_{n_k}(x_{n_k}) - f(x^*) \geq f_{n_k}(x_{n_k}) - f_{n_k}(x^*) - \lVert f_{n_k} - f\rVert_\infty \geq -\lVert f_{n_k} - f\rVert_\infty \to 0,
    \end{align*}
    and so we obtain $f(x) = f(x^*)$, which implies $x = x^*$ by uniqueness of the maximiser.
\end{proof}

The following shows that moving averages of $L^p$-bounded processes are also $L^p$-bounded:
\begin{lemma}\label{lem: Hamilton}
    Let $\smash{X = \big(X(t)\big)_{t \in \N}}$ be an $\R^d$-valued stochastic process bounded in $L^p$ for some $p \in [1,\infty)$. Let $E$ be a set and let $F^{(t)}_s: E \to  \R^{d \times d}$ be functions for each $t \in \N$ and $s \in \{0, \dots, t\}$ with $\smash{\sup_{t \in \N} \sum_{s=0}^t \sup_{x \in E} \lVert F^{(t)}_s(x) \rVert < \infty}$. Define processes $\smash{\big(Y_x(t) \big)_{t \in \N}}$ by $\smash{Y_x(t) \coloneq \sum_{s=0}^t F^{(t)}_s(x) X(t-s)}$ for any $x \in E$. Then $\sup_{t \in \N} \E\big[ \sup_{x \in E} \lVert Y_x(t)\rVert^p\big] < \infty$.
\end{lemma}
\begin{proof}
    $\smash{\vertiii{M} \coloneq \E\big[\sup_{x \in E} \lVert M_x\rVert^p\big]^{\frac{1}{p}}}$ defines a norm on the space of all collections $M = (M_x)_{x \in E}$ of random variables. To wit, we can generalise the Minkowski inequality by the following calculation for $\vertiii{M + N} \neq 0$:
    \begin{align*}
        \vertiii{M + N}^p &= \E\big[\sup_{x \in E} \lVert M_x + N_x\rVert^p\big] = \E\big[\sup_{x \in E} \big(\lVert M_x + N_x\rVert \cdot \lVert M_x + N_x\rVert^{p-1}\big)\big] \\
        &\leq \E\big[\sup_{x \in E} \big(\lVert M_x\rVert \sup_{x \in E}\lVert M_x + N_x\rVert^{p-1}\big)\big] + \E\big[\sup_{x \in E} \big(\lVert N_x\rVert \sup_{x \in E}\lVert M_x + N_x\rVert^{p-1}\big)\big] \\
        &\leq \Big(\E\big[ \sup_{x \in E} \lVert M_x\rVert^p\big]^{\frac{1}{p}} + \E\big[ \sup_{x \in E} \lVert N_x\rVert^p\big]^{\frac{1}{p}}\Big) \E\Big[\big(\sup_{x\in E} \lVert M_x + N_x \rVert^{p - 1} \big)^{\frac{p}{p - 1}}\Big]^{\frac{p - 1}{p}} \\
        &= \big( \vertiii{M} + \vertiii{N}\big) \vertiii{M + N}^{p - 1},
    \end{align*}
    where we used Hölder's inequality for the second inequality. This implies that $\vertiii{M + N} \leq \vertiii{M} + \vertiii{N}$ for $\vertiii{M + N} \neq 0$, as desired. Now, we can use this to obtain
    
    \begin{fitalign}
        \E\big[\sup_{x \in E} &\lVert Y_x(t)\rVert^p \big]^{\frac{1}{p}} = \E\Big[\sup_{x \in E} \Big\lVert \sum_{s=0} ^t F^{(t)}_s(x) X(t-s) \Big\rVert^p\Big]^{\frac{1}{p}} \leq \sum_{s=0}^t \E\Big[\sup_{x \in E} \big\lVert F^{(t)}_s(x) X(t-s) \big\rVert^p\Big]^{\frac{1}{p}} \\
        &\leq \sum_{s=0}^t \sup_{x \in E}\lVert F^{(t)}_s(x)\rVert \E\big[\lVert X(t-s)\rVert^p\big]^{\frac{1}{p}} \leq \sup_{t\in \N} \E\big[\lVert X(t)\rVert^p\big]^{\frac{1}{p}} \sup_{t \in \N}\sum_{s=0}^t \sup_{x \in E} \lVert F^{(t)}_s(x) \rVert
    \end{fitalign}
    
    \noindent for all $t \in \N$, where the right side is finite and independent of $t$, which finishes the proof.
\end{proof}
We now prove some elementary auxiliary results about convergence at a geometric rate:
\begin{lemma}\label{lem: geometric}
    Let $E$ be a unital normed algebra with respect to an inner multiplication operation. Let $a, b \in E$ and $(a_t)_{t\in \N}$ and $(b_t)_{t \in \N}$ be sequences in $E$ such that $a_t \to a$  and $b_t \to b$ at a geometric rate. Finally, let $F$ denote an arbitrary normed space and let $f: E \to F$. Then the following statements hold true:
    \begin{itemize}
        \item[1.] We have $a_t + b_t \to a + b$ at a geometric rate.
        \item[2.] We have $a_tb_t \to ab$ at a geometric rate. In particular, we have that $a_t^n \to a^n$ at a geometric rate for $n$-th powers with $n \in \N$.
        \item[3.] If $a$ and $a_t$ are invertible for all $t \in \N$, then $a_t^{-1} \to a^{-1}$ at a geometric rate.
    \end{itemize}
\end{lemma}
\begin{proof}
   \begin{itemize}
        \item[1.] This is obvious.
        \item[2.] This follows from part 1. by decomposing $a_tb_t - ab = a_t(b_t - b) + (a_t - a)b$ and by noting that the sequence $(a_t)_{t \in \N^*}$ is bounded in $E$ because it is a convergent sequence.
        \item[3.] We can decompose $a_t^{-1} - a^{-1} = a_t^{-1} (a - a_t) a^{-1}$. Since the map $x \mapsto x^{-1}$ is continuous on the group of invertible elements of a unital normed algebra, the sequence $(a_t^{-1})_{t \in \N}$ is bounded in $E$ because it is convergent. The result then follows upon taking norms.\qedhere
    \end{itemize}
\end{proof}

In the following two proofs we write $\norm{\cdot}_p$ for $\norm{\cdot}_{L^p}$ for better readability.
\begin{lemma}\label{lem: Boundedness_of_Polynomials}
    Let $\smash{X = \big(X(t)\big)_{t \in \N}}$ be an $\R^d$-valued stochastic process that is bounded in $L^p$ for some $p \in [1,\infty)$. If $g: \R^d \to \R^k$ is a $k$-dimensional polynomial of order $q$ for some $q \in \N^*$ with $q \leq p$ and $k \in \N^*$, then the process $g(X)$ is bounded in $L^{\frac{p}{q}}$.
\end{lemma}
\begin{proof}
    Since $\lVert Y \rVert_p \leq \lVert Y_1\rVert_p + \dots + \rVert Y_k\rVert_p$ holds for any $k$-dimensional random variable $Y$ and any $p \in [1,\infty)$, it suffices to prove that each component of the $k$-dimensional process $g(X)$ is bounded in $\smash{L^{\frac{p}{q}}}$. Hence assume that $k=1$ and let the function $g: \R^d \to \R$ be of the form $g(x) = \sum_{|\lambda| \leq q} \alpha_\lambda x^\lambda$. Then we obtain that $\smash{\lVert g(X(t))\rVert_{\frac{p}{q}} \leq \sum_{|\lambda| \leq q} |\alpha_\lambda|  \lVert X(t)^\lambda\rVert_{\frac{p}{q}}}$. Thus it suffices to show boundedness in $\smash{L^{\frac{p}{q}}}$ of monomials $X(t)^\lambda$ for $|\lambda| \geq 1$, which can be written as $X(t)^\lambda = \prod_{j=1}^d X_j(t)^{\lambda_j}$. By Hölder's inequality, we can deduce the bound
    \vspace*{-0.4cm}
    \begin{align*}
        \E\big(|X(t)^\lambda|^{\frac{p}{q}}\big) = \E\Big( \prod_{j=1}^d \big\lvert X_j(t)^{\frac{p\lambda_j}{q}}\big\rvert\Big) \leq \prod_{j=1}^d \big\lVert X_j(t)^{\frac{p\lambda_j}{q}}\big\rVert_{\frac{|\lambda|}{\lambda_j}} = \prod_{\lambda_j \neq 0} \big\lVert X_j(t)^{\frac{p\lambda_j}{q}}\big\rVert_{\frac{|\lambda|}{\lambda_j}}
    \end{align*}
    because $\sum_{j=1}^d \frac{\lambda_j}{|\lambda|} = 1$, where we set $\frac{|\lambda|}{\lambda_j} \coloneq \infty$ whenever $\lambda_j = 0$. The last identity follows because all factors with $\lambda_j = 0$ are equal to 1. For $\lambda_j \neq 0$ we can bound the factors by
    \begin{align}\label{eq: another_lengthy_term}
        \big\lVert X_j(t)^{\frac{p\lambda_j}{q}}\big\rVert_{\frac{|\lambda|}{\lambda_j}} = \E\big( |X_j(t)|^{\frac{p |\lambda|}{q}}\big)^{\frac{\lambda_j}{|\lambda|}} = \lVert X_j(t)\rVert_{\frac{p|\lambda|}{q}}^{\frac{p\lambda_j}{q}} \leq \lVert X_j(t)\rVert_p^{\frac{p\lambda_j}{q}},
    \end{align}
    where the last inequality follows from $|\lambda| \leq q$. Since $X$ is $L^p$-bounded, this ends the proof.
\end{proof}
\begin{lemma}\label{lem: Conv_of_Polynomials}
    Let $\smash{X = \big(X(t)\big)_{t \in \N}}$ and $\smash{Y = \big(Y(t)\big)_{t \in \N}}$ be $\R^d$-valued stochastic processes such that $X$ or $Y$ are bounded in $L^p$ for some $p \in [1,\infty)$. Suppose moreover that $\smash{X(t) - Y(t) \xrightarrow{L^p} 0}$ and let $g: \R^d \to \R$ be any polynomial of order $q \in \N^*$ with $q \leq p$. Then $g(X(t)) - g(Y(t)) \to 0$ in $\smash{L^{\frac{p}{q}}}$. If the convergence $\smash{X(t) - Y(t) \xrightarrow{L^p} 0}$ occurs at a geometric rate of order $r \in \R_+$, then also $\smash{g(X(t)) - g(Y(t)) \to 0}$ in $\smash{L^{\frac{p}{q}}}$ at a geometric rate of order $r$. 
\end{lemma}
\begin{proof}
    First, if one of $X$ and $Y$ is bounded in $L^p$, then the other one is clearly also bounded in $L^p$ because their difference converges to 0 in $L^p$. By the same reasoning as in Lemma \ref{lem: Boundedness_of_Polynomials}, we can assume that $k=1$. Let the function $g: \R^d \to \R$ be again of the form $g(x) = \sum_{|\lambda| \leq q} \alpha_\lambda x^\lambda$ As in the proof of Lemma \ref{lem: Boundedness_of_Polynomials}, we have $\lVert g(X(t)) - g(Y(t))\big\rVert_{\frac{p}{q}} \leq \smash{\sum_{|\lambda| \leq q} |\alpha_\lambda| \lVert X(t)^\lambda - Y(t)^\lambda\rVert_{\frac{p}{q}}}$; so it suffices the show the claim for monomials $g(x) = x^\lambda$ with $|\lambda| \geq 1$. Now
    \begin{fitequation*}
        X(t)^\lambda - Y(t)^\lambda = \prod_{j=1}^d X_j(t)^{\lambda_j} - \prod_{j=1}^d Y_j(t)^{\lambda_j} = \sum_{j=1}^d \Big[ \prod_{i=1}^{j-1} X_i(t)^{\lambda_i} (X_j(t)^{\lambda_j} - Y_j(t)^{\lambda_j}) \prod_{i=j+1}^d Y_j(t)^{\lambda_j} \Big]
    \end{fitequation*}
    \noindent for any $t \in \N$. Consequently, the triangle inequality for any $t \in \N$ yields that
    \begin{align*}
        \lVert X(t)^\lambda - Y(t)^\lambda\rVert_{\frac{p}{q}} \leq \sum_{j=1}^d \bigg\lVert \prod_{i=1}^{j-1} X_i(t)^{\lambda_i} (X_j(t)^{\lambda_j} - Y_j(t)^{\lambda_j}) \prod_{i=j+1}^d Y_j(t)^{\lambda_j}\bigg\rVert_{\frac{p}{q}}.
    \end{align*}
    By Hölder's inequality, since $\sum_{j=1}^d \frac{\lambda_j}{|\lambda|} = 1$, and using the convention $\frac{1}{0} \coloneq \infty$, it follows that we can bound the $j$-th summand on the right-hand side of the preceding inequality by
    \begin{align}\label{eq: lengthy_term}
        \prod_{i=1}^{j-1} \big\lVert X_i(t)^{\frac{p\lambda_i}{q}}\big\rVert_{\frac{|\lambda|}{\lambda_i}}^{\frac{q}{p}} \big\lVert (X_j(t)^{\lambda_j} - Y_j(t)^{\lambda_j})^{\frac{p}{q}}\big\rVert_{\frac{|\lambda|}{\lambda_j}}^{\frac{q}{p}} \prod_{i=j+1}^{d} \big\lVert Y_i(t)^{\frac{p\lambda_i}{q}}\big\rVert_{\frac{|\lambda|}{\lambda_i}}^{\frac{q}{p}}.
    \end{align}
    Focus on the products on the left of \eqref{eq: lengthy_term}. Each factor $\smash{\big\lVert X_i(t)^{\frac{p\lambda_i}{q}}\big\rVert_{\frac{|\lambda|}{\lambda_i}}}$ is uniformly bounded in $t$ by the same reasoning as in \eqref{eq: another_lengthy_term}. Likewise, the whole product on the right-hand side of \eqref{eq: lengthy_term} is uniformly bounded in $t$. \mbox{Since $|\lambda| \leq q$, the middle term in \eqref{eq: lengthy_term} is bounded by}
    \begin{align*}
        \big\lVert (X_j(t)^{\lambda_j} - Y_j(t)^{\lambda_j})^{\frac{p}{q}}\big\rVert_{\frac{|\lambda|}{\lambda_j}}^{\frac{q}{p}} \leq \big\lVert (X_j(t)^{\lambda_j} - Y_j(t)^{\lambda_j})^{\frac{p}{q}}\big\rVert_{\frac{q}{\lambda_j}}^{\frac{q}{p}} = \big\lVert X_j(t)^{\lambda_j} - Y_j(t)^{\lambda_j}\big\rVert_{\frac{p}{\lambda_j}}.
    \end{align*}
    If $\lambda_j = 0$, this term is equal to 1. Otherwise, we can write $x^n - y^n = (x - y) \sum_{j=1}^n x^{n - j} y^{j - 1}$ for any $x, y \in \R$ and $n \in \N^*$. We now apply Hölder's inequality $\lVert fg\rVert_r \leq \norm{f}_{rk} \norm{g}_{rl}$, which holds for any $r \in [1,\infty)$ and $\frac{1}{k} + \frac{1}{l} = 1$, to $r \coloneq \frac{p}{\lambda_j}$, $k \coloneq \lambda_j$ and $l \coloneq \frac{\lambda_j}{\lambda_j - 1}$. This gives
    \begin{align}\label{eq: third_lenghty_term}
        \big\lVert X_j(t)^{\lambda_j} - Y_j(t)^{\lambda_j}\big\rVert_{\frac{p}{\lambda_j}} &\leq \lVert X_j(t) - Y_j(t)\rVert_{p}\Big\lVert \sum_{i=1}^{\lambda_j} X_j(t)^{\lambda_j - i} Y_j(t)^{i - 1}\Big\rVert_{\frac{p}{\lambda_j - 1}}.
    \end{align}
     The sum on the right-hand side is a polynomial of order $\lambda_j - 1$ of the $2d$-dimensional process $(X, Y)$, which is bounded in $L^p$. Lemma \ref{lem: Boundedness_of_Polynomials} hence yields that this polynomial is bounded in $L^{\frac{p}{\lambda_j - 1}}$. Since moreover $\lVert X_j(t) - Y_j(t) \rVert_p \leq \lVert X(t) - Y(t) \rVert_p$, we obtain that the left-hand side of \eqref{eq: third_lenghty_term} converges to 0 (at a geometric rate of order $r$) whenever $X(t) - Y(t) \to 0$ in $L^p$ (at a geometric rate of order $r$). This finishes the proof.
\end{proof}

In the course of Section \ref{s:proofs}, we need a suitable version of a central limit theorem for discrete-time martingales. There is plenty of literature concerning martingale central limit theorems, see for example \cite{Hall1980}. Since however most of these results are stated in one dimension, we include a multivariate version of the martingale central limit theorem below:

\begin{theorem}\label{theo: MCLT}
Suppose that, for any $n \in \N^*$, a positive integer $k_n \in \N^*$ as well as some filtration $\smash{\F^{(n)} \coloneq (\F^{(n)}_k)_{k \in \{0, \dots, k_n\}}}$ on $(\Omega, \F)$ are given  and $\smash{(U^{(n)}_k)_{k \in \{1, \dots, k_n\}}}$ is some $d$-dimensional square-integrable martingale difference sequence with respect to $\smash{\F^{(n)}}$, i.e.\ $\smash{\E(U^{(n)}_{k} \mid \F^{(n)}_{k - 1}) = 0}$ for $k \in \{1, \dots, k_n\}$ and $n \in \N^*$. Assume that there exists some positive semidefinite matrix $\Sigma \in \R^{d\times d}$ such that the following holds:
\begin{itemize}
    \item[1.] $\sum_{k=1}^{k_n} U^{(n)}_k {U^{(n)\top}_k} \xrightarrow{\PP} \Sigma$ as $n \to \infty$,
    \item[2.] $\smash{\sum_{k=1}^{k_n} \E\big[\lVert U^{(n)}_k\rVert^2 \mathbf{1}_{\{\lVert U^{(n)}_k\rVert \geq \varepsilon\}} \mid \mathscr{F}_{k-1}^{(n)}\big] \xrightarrow{\PP} 0}$ for every $\varepsilon > 0$.
\end{itemize}
Then $M_n \coloneq \sum_{k=1}^{k_n} U_k^{(n)}$ satisfies $M_n \xrightarrow{d} Z$, where $Z \sim \mathscr{N}(0, \Sigma)$. 
\end{theorem}
\begin{proof}
This follows by \cite[Theorem VIII.3.33]{Jacod2003}.
\end{proof}

The main technical purpose of the following short lemma is to establish an elementary condition under which the uncountable supremum of real-valued random variables is again a well-defined random variable. We make occasional use of this result in Section \ref{sec4.2.2: ULLN}.

\begin{lemma}\label{lem: measurability}
    Let $E$ be a separable topological space. Suppose that for each $x \in E$ a random variable $Y(x)$ taking values in a normed space $F$ is given such that for each fixed $\omega$ the function $x \mapsto Y(x)(\omega)$ is continuous. Then the function $\sup_{x \in E} \lVert Y(x) \rVert$ is measurable.
\end{lemma}
\begin{proof}
    Since $E$ is separable, there exists a countable dense subset $E_0 \subseteq E$. By $\omega$-wise continuity in $x$ we then have $\sup_{x \in E} \lVert Y(x) \rVert = \sup_{x \in E_0} \lVert Y(x)\rVert$, where $\lVert \cdot \rVert$ denotes the norm in $F$. This is measurable as a countable supremum of measurable functions.
\end{proof}

The next result stated below contains a standard condition for \textit{uniform convergence in probability}, which is often termed a stochastic equicontinuity condition, see for example \cite{Andrews1992}. In the context of parameter estimation from Section \ref{sec4.2: Proof}, we need this condition to prove uniform weak laws of large numbers. The stochastic equicontinuity condition is quite old and dates back to \cite{Prokhorov1956}.

\begin{theorem}\label{theo: Pötscher_Prucha}
    Let $(E, d)$ be a compact metric space and suppose that for each $x \in E$ a sequence of $d$-dimensional random variables $(Y_x(t))_{t\in \N^*}$ as well as a deterministic $Y(x) \in \R^d$ are given such that $Y_x(t)$ and $Y(x)$ are continuous in $x$. Suppose that $\smash{Y_x(t) \xrightarrow{\PP} Y(x)}$ for each $x \in E$ and that for each $x \in E$ the uniform stochastic equicontinuity condition
    \begin{equation*}
        \lim_{\alpha \to 0} \limsup_{t \to \infty} \PP\Big(\sup_{d(x,y) < \alpha} \lVert Y_x(t) - Y_y(t) \rVert > \varepsilon\Big) = 0
    \end{equation*}
    holds for each $\varepsilon > 0$. Then $\sup_{x \in E} \lVert Y_x(t) - Y(x)\rVert \xrightarrow{t \to \infty} 0$ in probability.
\end{theorem}
\begin{proof}
    All suprema over uncountable sets occurring in the theorem are well-defined random variables by Lemma \ref{lem: measurability} because $E$ is separable as a compact space. The proof now closely follows the proof of Theorem 1 in \cite{Andrews1992}. Since $x \mapsto Y(x)$ is continuous on the compact space $E$, it is uniformly continuous. Let $\alpha > 0$ be so small that $\limsup_{t \to \infty} \PP\big(\sup_{d(x, y) < \alpha} \lVert Y_x(t) - Y_y(t) \rVert > \frac{\varepsilon}{3}\big) < \eta$ and that $\sup_{d(x, y) < \alpha}\lVert Y(x) - Y(y) \rVert < \frac{\varepsilon}{3}$ for fixed $\varepsilon, \eta > 0$. Using the total boundedness of $E$ we can find a finite set $\{x_1, \dots, x_K\} \subseteq E$ such that $\smash{E = \bigcup_{j=1}^K \mathrm{B}(x_j, \alpha)}$, where $\mathrm{B}(x_j, \alpha)$ denotes the open ball of radius $\alpha$ around $x_j$, as usual. We then obtain for $k(\varepsilon) \coloneq \limsup_{t \to \infty} \PP\big(\sup_{x \in E} \lVert Y_x(t) - Y(x) \rVert > \varepsilon\big)$ that
    
    \begin{fitalign}
        k(\varepsilon) &\leq \limsup_{t\to \infty} \PP\Big( \max_{j \in \{1, \dots, K\}} \sup_{x \in \mathrm{B}(x_j, \alpha)} \big[ \lVert Y_x(t) - Y_{x_j}(t)\rVert + \lVert Y_{x_j}(t) - Y(x_j)\rVert + \lVert Y(x_j) - Y(x)\big]\rVert > \varepsilon\Big) \\
        &\leq \limsup_{t \to \infty} \PP\Big(\sup_{d(x,y) < \alpha} \lVert Y_x(t) - Y_{y}(t)\rVert > \frac{\varepsilon}{3}\Big) + \limsup_{t \to \infty}\PP\Big(\max_{j \in \{1, \dots, K\}}\lVert Y_{x_j}(t) - Y(x_j)\rVert > \frac{\varepsilon}{3}\Big) \\
        &< \eta + \sum_{j=1}^K \limsup_{t \to \infty} \PP\Big( \lVert Y_{x_j}(t) - Y(x_j) \rVert > \frac{\varepsilon}{3K}\Big) = \eta.
    \end{fitalign}

\noindent Since $\eta > 0$ was arbitrary, it follows that $k(\varepsilon) = 0$, which finishes the proof.
\end{proof}

Another result that is needed for the proofs of Section \ref{sec4.2.3: Ident} is given by the following lemma, which states that a sequence of random functions converges uniformly in probability whenever it converges pointwise in probability and its partial derivatives converge uniformly.

\begin{lemma}\label{lem: tao}
    Let $E\subseteq \R^n$ be bounded and convex. Suppose that for each $x \in E$, a sequence of $m$-dimensional random variables $(Y_x(t))_{t \in \N^*}$ is given such that $Y_x(t)$ is continuously differentiable in $x$, and that there exists a deterministic function $Y: E \to \R^m$ such that $\smash{Y_{x}(t) \xrightarrow{\PP} Y(x)}$ for any $x \in E$. Suppose further that there exists some continuous $W: E \to \R^{m\times n}$ such that $\smash{\sup_{x \in E} \lVert \nabla_x Y_x(t) - W(x) \rVert \xrightarrow{\PP} 0}$, where $\nabla_x$ denotes the $(m \times n)$-Jacobian matrix with respect to $x$. Then $Y$ is continuously differentiable, $\smash{\sup_{x \in E} \lVert Y_x(t) - Y(x)\rVert \xrightarrow{\PP} 0}$, and $\nabla_x Y = W$.
\end{lemma}
\begin{proof}
    We first show that $Y_x(t)$ is a uniform Cauchy sequence in probability, i.e., for every $\varepsilon > 0$ there is a $T \in \N^*$ such that $\PP(\sup_{x \in E} \lVert Y_x(s) - Y_x(t)\rVert > \varepsilon) < \varepsilon$ for all $s, t \geq T$. Let $x_0 \in E$. By applying the multivariate mean value theorem to $x \mapsto Y_x(t) - Y_x(s)$ we have
    \begin{align*}
        \sup_{x \in E} \lVert Y_x(t) - Y_x(s) \rVert &\leq \sup_{x \in E}\big\lVert \big(Y_x(t) - Y_x(s)\big) - \big(Y_{x_0}(t) - Y_{x_0}(s)\big) \big\rVert + \big\lVert Y_{x_0}(t) - Y_{x_0}(s)\big\rVert \\
        &\leq \sup_{x \in E}\big\lVert \nabla_x Y_x(t) - \nabla_x Y_x(s) \big\rVert \sup_{x \in E} \lVert x - x_0 \rVert + \big\lVert Y_{x_0}(t) - Y_{x_0}(s)\big\rVert \\
        &\leq \mathrm{diam}(E) \sup_{x \in E}\big\lVert \nabla_x Y_x(t) - \nabla_x Y_x(s) \big\rVert + \big\lVert Y_{x_0}(t) - Y_{x_0}(s)\big\rVert.
    \end{align*}
    Since $\nabla_x Y_x(t) \to W(x)$ uniformly in probability, $(\nabla_x Y_x(t))_{t \in \N^*}$ is a uniform Cauchy sequence in probability, i.e.\ for $\varepsilon > 0$ we can find $T \in \N^*$ such that $\PP\big(\sup_{x \in E}\big\lVert \nabla_x Y_x(t) - \nabla_x Y_x(s) \big\rVert > \frac{\varepsilon}{2 \mathrm{diam}(E)}\big) < \frac{\varepsilon}{2}$ for all $s, t \geq T$. Likewise, since $Y_{x_0}(t)$ converges in probability for any $x_0 \in E$, it is a Cauchy sequence in probability, so we can choose $T$ large enough such that also $\PP\big(\big\lVert Y_{x_0}(t) - Y_{x_0}(s)\big\rVert > \frac{\varepsilon}{2}\big) < \frac{\varepsilon}{2}$ for all $s, t \geq T$. By the mean value theorem, it follows that
    \begin{fitequation*}
    \PP\Big( \sup_{x \in E} \lVert Y_x(t) - Y_x(s) \rVert > \varepsilon\Big) \leq \PP\Big( \sup_{x \in E} \lVert \nabla_x Y_x(t) - \nabla_x Y_x(s) \rVert \mathrm{diam}(E) + \lVert Y_{x_0}(t) - Y_{x_0}(s) \rVert > \varepsilon\Big) \leq \varepsilon
    \end{fitequation*}
    \noindent for $s, t \geq T$, so $(Y_x(t))_{t\in \N^*}$ is uniformly Cauchy in probability. Consequently, there exists some (possibly random) function $\widetilde Y: E \to \R^m$ such that $\smash{\sup_{x \in E} \lVert Y_x(t) - \widetilde Y(x) \rVert \xrightarrow{\PP} 0}$. But by the uniqueness of limits in probability, it follows that $\widetilde Y = Y$ almost surely.

    Fix now $i \in \{1, \dots, m\}$ and $j \in \{1, \dots, n\}$. For any $x \in E$ let $H_j(x)$ denote the set of all $h \in \R$ such that $x' \coloneq x + he_j^{(n)} \in E$, where $e_j^{(n)}$ is the $j$-th unit vector. Then, for $x \in E$ and $h \in H_j(x)$,
    \[\Big\lvert\frac{1}{h}\big[Y_{x'}(t)_i - Y_x(t)_i\big] - \frac{1}{h}\big[Y_{x'}(s)_i - Y_x(s)_i\big]\Big\rvert \leq \sup_{x \in E} |\nabla_x Y_x(t)_{ij} -  \nabla_x Y_x(s)_{ij}|,\]
    so the quotient sequence is Cauchy in probability uniformly in $x \in E$ and $h \in H_j(x)$. Therefore $\sup_{x \in E, h \in H_j(x)} \big\lvert \frac{1}{h}\big[Y_{x'}(t)_i - Y_x(t)_i\big] - \frac{1}{h}\big[Y(x')_i - Y(x)_i\big] \big\rvert \to 0$ in probability and in particular almost surely along a subsequence $(t_k)_{k \in \N^*}$. Take a further subsequence $(t_{k_l})_{l \in \N^*}$ such that also $\smash{\sup_{x \in E} \lVert \nabla_x Y_x(t_{k_l}) - W(x) \rVert \to 0}$ almost surely. With a slight abuse of notation, we will suppress the indices $k_l$ of these subsequences in the following. We will now show that $Y$ is continuously differentiable and $\nabla_x Y = W$. For $x \in E$ and $x' = x + h e_j$ with $h \in H_j(x)$ we get
    \begin{equation}
        \Big\lvert \frac{1}{h} \big[Y(x')_i - Y(x)_i\big] - W(x)_{ij}\Big\rvert \leq A + B + C,
    \end{equation}
    where $A = \big\lvert \frac{1}{h} \big[Y(x')_i - Y(x)_i\big] - \frac{1}{h}\big[Y_{x'}(t)_i - Y_x(t)_i\big]\big\rvert$, $B = \big\lvert \frac{1}{h}\big[Y_{x'}(t)_i - Y_x(t)_i\big] - \nabla_x Y_x(t)_{ij}\big\rvert$, and $C = \big\lvert \nabla_x Y_x(t)_{ij} - W(x)_{ij}\big\rvert$. Let $\varepsilon > 0$. By the arguments above, the terms $A$ and $C$ become almost surely less than $\frac{\varepsilon}{3}$ for some large enough $t$ independent of $x$ and $h$, while the term $B$ becomes less than $\frac{\varepsilon}{3}$ for any fixed $t$ if $|h|$ is small enough. It follows that $\partial_{x_j} Y(x)_i = W(x)_{ij}$.
\end{proof}

In Section \ref{sec5.3: OU} we need multivariate extensions of the well-known Burkholder--Davis--Gun\-dy inequality (see \cite{Dellacherie1982}, Theorem VII.92) and of Émery's inequality (see \cite{Protter2004}, Theorem V.3) for semimartingales. These multivariate extensions follow easily from their univariate counterparts, and we include the proofs here for the sake of completeness. In the following, let $\lVert M \rVert_{\max} \coloneq \max\big\lbrace |M_{ij}|: \: i \in \{1, \dots, m\}, \,j \in \{1, \dots, n\}\big\rbrace$ denote the maximum-element norm of a matrix $M \in \R^{m \times n}$ and recall that $\lVert M \rVert_{\max} \leq \lVert M \rVert \leq \sqrt{m n} \lVert M \rVert_{\max}$, where the latter can be derived from simple relations between the spectral norm and the Frobenius norm of a matrix. We start with the Burkholder--Davis--Gundy inequality and, as always, fix some right-continuous filtration $(\F(t))_{t \in \R_+}$ on the space $(\Omega, \F)$.

\begin{theorem}\label{theo: BDG}
    Let $X = (X(t))_{t \in \R_+}$ be an $\R^d$-valued local martingale with $X(0) = 0$ and let $p \in [1,\infty)$. Then there are $c_p > 0$ and $C_p > 0$ such that for all stopping times\footnote{Here $\lVert[X, X](\infty) \rVert \coloneq \lim_{t \to \infty} \lVert [X, X](t) \rVert$, which is well-defined because $[X, X]$ is almost surely increasing in the Loewner order because $v^\top [X, X] v = [v^\top X, v^\top X]$ is almost surely increasing for any $v \in \R^d$.} $\tau$
    \begin{equation*}
        c_p \E\big( \lVert [X, X](\tau)\rVert^{\frac{p}{2}}\big) \leq \E\Big( \sup_{s \leq \tau} \lVert X(s)\rVert^p\Big) \leq C_p\E\big( \lVert [X, X](\tau)\rVert^{\frac{p}{2}}\big).
    \end{equation*}
    In particular, the constants $c_p$ and $C_p$ can be chosen as $c_p = \big[6p  d^{\frac{p}{2} + 1}\big]^{-1}$ and $C_p = 4p d^{\frac{p}{2} + 1}$.
\end{theorem}
\begin{proof}
    First, note that $\max_{1 \leq i \leq n} \E(X_i) \leq \E(\max_{1 \leq i \leq n} X_i) \leq n \max_{1 \leq i \leq n} \E(X_i)$ for any nonnegative random variables $X_1, \dots, X_n$. Moreover, we use the univariate Burkholder--Davis--Gundy inequality $\frac{1}{6p} \E([X, X](\tau)^{\frac{p}{2}}) \leq \E(\sup_{s \leq \tau} |X(s)|^p) \leq 4p \E([X, X](\tau)^{\frac{p}{2}})$ in the case $d=1$, found e.g. in \cite{Dellacherie1982}, Theorem VII.92. To show the first inequality, note that
    \begin{align*}
        \E\Big( \lVert [X, X](\tau)\rVert^{\frac{p}{2}}\Big) &\leq d^{\frac{p}{2}} \E\Big( \lVert [X, X](\tau)\rVert_{\max}^{\frac{p}{2}}\Big) = d^{\frac{p}{2}} \E\Big(\max_{1 \leq i, j \leq d} |[X_i, X_j](\tau)|^{\frac{p}{2}}\Big) \\
        &\leq d^{\frac{p}{2}} \E\Big( \max_{1 \leq i, j \leq d} [X_i, X_i](\tau)^{\frac{p}{4}} [X_j, X_j](\tau)^{\frac{p}{4}} \Big) \leq d^{\frac{p}{2}} \E\Big(\max_{1 \leq i \leq d} [X_i, X_i](\tau)^{\frac{p}{2}} \Big) \\
        &\leq 6pd^{\frac{p}{2} + 1} \max_{1 \leq i \leq d} \E\Big( \sup_{s \leq \tau} |X_i(s)|^p\Big) \leq 6pd^{\frac{p}{2} + 1}\E\Big( \sup_{s \leq \tau} \lVert X(s)\rVert_{\max}^p\Big) \\
        &\leq 6pd^{\frac{p}{2} + 1}\E\Big( \sup_{s \leq \tau} \lVert X(s)\rVert^p\Big),
    \end{align*}
    where we used the Kunita--Watanabe (or the Cauchy--Schwarz) inequality in the third step, see for example \cite{Protter2004}, Theorem II.25. For the second inequality, we obtain the upper bound
    \begin{align*}
        \E\Big(\sup_{s \leq \tau} \lVert X(s)\rVert^p \Big) &\leq d^{\frac{p}{2}} \E\Big(\sup_{s \leq \tau} \lVert X(s)\rVert_{\max}^p \Big) \leq d^{\frac{p}{2} + 1} \max_{1 \leq i \leq d} \E\Big(\sup_{s \leq \tau} |X_i(s)|^p\Big) \\
        &\leq 4pd^{\frac{p}{2} + 1} \E\Big(\lVert [X, X](\tau)\rVert_{\max}^{\frac{p}{2}}\Big) \leq 4pd^{\frac{p}{2} + 1} \E\Big( \lVert [X, X](\tau)\rVert^{\frac{p}{2}}\Big). \qedhere
    \end{align*}
\end{proof}

On the space of $\R^{m \times n}$-valued càdlàg (or càglàd) adapted processes, we define the (possibly infinite) $S^p$ norm by $\lVert X \rVert_{S^p} \coloneq \big\lVert \sup_{t \in \R_+} \lVert X(t)\rVert \big\rVert_{L^p}$ for $1 \leq p \leq \infty$. For an $\R^d$-valued semimartingale $X$, we can then define the (possibly infinite) $H^p$ norm by $\lVert X \rVert_{H^p} \coloneq \inf_{M, A} \big\lbrace \big\lVert \lVert X(0)\rVert + \lVert [M, M](\infty)\rVert^{\frac{1}{2}} + \lVert A\rVert(\infty) \big\rVert_{L^p}\big\rbrace$, where the infimum is taken over all decompositions $X = X(0) + M + A$ with $M$ being a local martingale and $A$ an adapted, càdlàg process of finite variation with $M(0) = A(0) = 0$ almost surely. Here $\lVert A \rVert$ denotes the total variation process of $A$. The Hardy spaces $S^p$ and $H^p$ of processes with finite $S^p$ and $H^p$ norm can be shown to be Banach spaces, see \cite{Dellacherie1982}, VII.64 and VII.98(e), in the case $d=1$ and $p < \infty$. The $S^p$ and $H^p$ norms satisfy the following weak inequality:

\begin{theorem}\label{theo: Sp-Hp-domination}
    Let $p \in [1,\infty)$ and $d \in \N^*$. Then there exists $c_p > 0$ depending only on $p$ and $d$ such that $\lVert X \rVert_{S^p} \leq c_p \lVert X \rVert_{H^p}$ for all $d$-dimensional semimartingales $X$.
\end{theorem}
\begin{proof}
    Applying the Burkholder--Davis--Gundy inequality \ref{theo: BDG}, the fact that the total variation process is increasing, and $|a+b+c|^p \leq 3^{p-1}(|a|^p + |b|^p + |c|^p)$, we obtain
    \begin{align*}
        \lVert X\rVert^p_{S^p} &\leq 3^{p-1} \E\Big( \lVert X(0) \rVert^p + \sup_{t \in \R_+} \lVert M(t)\rVert^p + \lVert A \rVert(\infty)^p\Big) \\
        &\leq C_p \E\Big[\Big(\lVert X(0)\rVert + \lVert [M, M](\infty)\rVert^{\frac{1}{2}} + \lVert A\rVert(\infty) \Big)^p\Big]
    \end{align*}
    for some $C_p$ depending only on $p$ and $d$. Taking $p$-th roots and infima finishes the proof.
\end{proof}

The following theorem provides the multivariate version of Émery's inequality, attributed to \cite{Emery1978}, which gives a handy upper bound for $H^p$ norms of arbitrary stochastic integrals.

\begin{theorem}\label{theo: Emery}
    Let $H$ be a càglàd, adapted $\R^{m \times n}$-valued process and $X$ an $\R^n$-valued semimartingale. Moreover, let $\frac{1}{p} + \frac{1}{q} = \frac{1}{r}$ for $1 \leq p \leq \infty$ and $1 \leq q \leq \infty$. Then
    \begin{equation*}
        \Big\lVert \int H \dd X \Big\rVert_{H^r} \leq n\sqrt{m} \lVert H \rVert_{S^p} \lVert X - X(0) \rVert_{H^q}.
    \end{equation*}
\end{theorem}
\begin{proof}
    Fix a semimartingale decomposition $X = X(0) + M + A$. Then $\int H \dd X = \int H \dd M + \int H \dd A$ is a decomposition for $\int H \dd X$ Recall that $[\int H \dd M, \, \int H \dd M] = \int H \dd[M, M] H^\top$. We will now show $\lVert \int H \dd A \rVert(t) \leq n\sqrt{m} \sup_{s \leq t} \lVert H(s) \rVert \: \lVert A \rVert(t)$. Since each component $A_j$ of $A$ is a univariate process of finite variation, there are unique increasing $A^+_j$ and $A^-_j$ such that $A_j = A^+_j - A^-_j$ and $$\Big\lVert \int_a^b H(s) \dd A(s)\Big\rVert \leq \Big\lVert \int_a^b H(s) \dd A^+(s)\Big\rVert + \Big\lVert \int_a^b H(s) \dd A^-(s)\Big\rVert$$ for all $a \leq b$, where $A^+$ and $A^-$ denote the $n$-dimensional processes with components $A^+_j$ and $A^-_j$. Then, as $\lVert x \rVert \leq \sqrt{m} \lVert x\rVert_{\max}$,
    \begin{align*}
        \Big\lVert \int_a^b H(s) \dd A(s)\Big\rVert &\leq \sqrt{m} \max_{1 \leq i \leq m} \sum_{1 \leq j \leq n} \bigg( \bigg| \int_a^b H_{ij}(s) \dd A_j^+(s)\bigg| + \bigg| \int_a^b H_{ij}(s) \dd A_j^-(s)\bigg|\bigg)  \\
        &\leq \sqrt{m} \sup_{a \leq t \leq b} \lVert H(t)\rVert_{\max} \sum_{1 \leq j \leq n} \Big( |A_j|(b) - |A_j|(a)\Big).
    \end{align*}
     Hence $$\Big\lVert \int H \dd A\Big\rVert(t) = \sup_{\pi \in \mathscr{P}_t}\sum_{t_i \in \pi} \Big\lVert \int_{t_i}^{t_{i+1}} H(s) \dd A(s) \Big\rVert \leq \sqrt{m} \sup_{s \leq t} \Big\lVert H(s) \Big\rVert_{\max} \sum_{j=1}^n |A_j|(t),$$ which is itself bounded by $n \sqrt{m} \sup_{s \leq t} \lVert H(s) \rVert \: \lVert A \rVert(t)$. We have the further bound

    \begin{align*}
        \Big\lVert &\int_0^t H(s) \dd[M, M](s) \, H(s)^\top \Big\rVert \leq m \Big\lVert \int_0^t H(s) \dd[M, M](s) \, H(s)^\top \Big\rVert_{\max} \\
        &\leq m \max_{1 \leq i,j \leq m} \sum_{1 \leq l, k \leq n} \Big\lvert \int_0^t H_{il}(s) H_{jk}(s) \dd [M_l, M_k](s) \Big\rvert \\
        &\leq m \max_{1 \leq i,j \leq m} \sum_{1 \leq l, k \leq n} \Big( \int_0^t H_{il}(s)^2 \dd [M_l, M_l](s) \int_0^t H_{jk}(s)^2 \dd [M_k, M_k](s) \Big)^{\frac{1}{2}} \\
        &\leq m n^2 \max_{1 \leq i, j \leq m} \Big( \sup_{s \leq t} |H_{ij}(s)| \Big)^2 \max_{1 \leq k \leq n}[M_k, M_k](t) \\
        &= m n^2 \Big( \sup_{s \leq t} \lVert H(s) \rVert_{\max}\Big)^2 \: \lVert [M, M](t) \rVert_{\max} \\
        &\leq m n^2 \Big( \sup_{s \leq t} \lVert H(s) \rVert \Big)^2 \: \lVert [M, M](t) \rVert,
    \end{align*}
    where we used the Kunita--Watanabe inequality in the third step, see for example \cite{Protter2004}, Theorem II.25. This yields 
    \[\left\lVert \int H \dd X\right\rVert_{H^r} \leq n \sqrt{m} \left\lVert \Bigl(\sup_{t \in \R_+} \lVert H(t)\rVert\Bigr) \lVert [M, M](\infty)\rVert^{\frac{1}{2}} +  \Bigl(\sup_{t \in \R_+} \lVert H(t)\rVert\Bigr) \lVert A \rVert(\infty) \right\rVert_{L^r}.\]
    An application of Hölder's inequality yields the result.
\end{proof}
\end{appendix}

\begin{acks}
The authors would like to thank two anonymous referees for their valuable comments.
\end{acks}
\pagebreak

\printbibliography[heading=bibintoc]

\end{document}